\documentclass[a4paper, 10pt]{article}
\usepackage[applemac]{inputenc} 
\usepackage[T1]{fontenc}      
\usepackage[hmargin=1.5cm, vmargin=2cm]{geometry}         
\usepackage[english]{babel}  
\usepackage{bm,amsmath,amssymb}
\def\cat#1{{\mathfrak{#1}}}

\usepackage{floatflt}

\usepackage{multicol}
\usepackage{cwpuzzle}
\renewcommand{\PuzzleBlackBox}{\rule{.0\PuzzleUnitlength}%
{.0\PuzzleUnitlength}} 
\usepackage{pst-all}
\usepackage{pgf}
\usepackage{enumerate}
\usepackage{qsymbols}
\usepackage{multirow}
\usepackage[all]{xy}   
\usepackage{floatflt}
\usepackage{ntheorem}

\theoremstyle{break}
\theorembodyfont{\upshape \itshape}
\theoremheaderfont{\bfseries\itshape}
\theoremseparator{.} 
\newtheorem{theo}{Theorem}[subsection]
\newtheorem{prop}[theo]{Proposition}
\newtheorem{lemme}[theo]{Lemma}
\newtheorem{cor}[theo]{Corollary}

\theoremstyle{break}
\theorembodyfont{\upshape}
\theoremheaderfont{\bfseries\itshape}
\theoremseparator{.} 
\newtheorem*{nota}{Notation}
\newtheorem{dfn}[theo]{Definition}

\newtheorem*{convsanss}{Convention}

\theoremstyle{plain}
\theorembodyfont{\upshape}
\theoremheaderfont{\bfseries\itshape}
\theoremseparator{.} 
\newtheorem*{obs}{Observation}
\newtheorem{obsn}{Observation}
\newtheorem{ex}[theo]{Example}
\newtheorem{rmq}[theo]{Remark}

\usepackage[dvips,ps2pdf]{hyperref}

\newenvironment{itemizedot}{\begin{itemize} }{\end{itemize}}

\newenvironment{changemargin}[2]{\begin{list}{}{%
\setlength{\topsep}{0pt}%
\setlength{\leftmargin}{0pt}%
\setlength{\rightmargin}{0pt}%
\setlength{\listparindent}{0pt}
\setlength{\itemindent}{\parindent}%
\setlength{\parsep}{0pt plus 1pt}%
\addtolength{\leftmargin}{#1}%
\addtolength{\rightmargin}{#2}%
}\item }{\end{list}}

\newenvironment{preuve}{ \begin{changemargin}{0cm}{0cm}{\flushleft\bf Proof :} 
}{\flushright $\square$\par\end{changemargin}}
  
\newenvironment{preuve3}{ \begin{changemargin}{0cm}{0cm}{\flushleft\bf Proof of Theorem \ref{1.4} :}
}{\flushright $\square$\par\end{changemargin}}  
\newenvironment{preuve4}{ \begin{changemargin}{0cm}{0cm}{\flushleft\bf Proof of Theorem \ref{reciproque contrainte} :} 
}{\flushright $\square$\par\end{changemargin}} 
\newenvironment{preuve5}{ \begin{changemargin}{0cm}{0cm}{\flushleft\bf Proof of the observation :} 
}{\flushright $\square$\par\end{changemargin}}
      
\title{Cauchon diagrams for quantized enveloping algebras}
\author{Antoine M\'eriaux \footnote{\textsl{Acknowledgements :} The author wishes to express gratitude to S. Launois for useful discussions and valuable advices.}\\
\\
\it Laboratoire d'\'equations aux d\'eriv\'ees partielles et physique
math\'ematique,\\ \it U.F.R.
Sciences, B.P. 1039, 51687 Reims Cedex 2, France.\\
antoine.meriaux@univ-reims.fr}
\date{\ }

\begin{document}

\maketitle
$ $ \\
$ $ \\
\textbf{Abstract}\\

Let $\mathfrak{g}$ be a finite dimensional complex simple Lie algebra, $\mathbb{K}$ a commutative field and $q$ a nonzero element of $\mathbb{K}$ which is not a root of unity. 
To each reduced decomposition of the longest element $w_0$ of the Weyl group $W$ corresponds a PBW basis of the quantised enveloping algebra $\mathcal{U}_q^+(\mathfrak{g})$, and one can apply the theory of deleting-derivation to this iterated Ore extension. In particular, for each decomposition of $w_0$, this theory constructs a bijection between the set of prime ideals in $\mathcal{U}_q^+(\mathfrak{g})$ that are invariant under a natural torus action and certain combinatorial objects called Cauchon diagrams. In this paper, we give an algorithmic description of these Cauchon diagrams when the chosen reduced decomposition of $w_0$ corresponds to a good ordering (in the sense of Lusztig \cite{Lu2}) of the set of positive roots. This algorithmic description is based on the constraints that are coming from Lusztig's admissible planes \cite{Lu2}: each admissible plane leads to a set of constraints that a diagram has to satisfy to be Cauchon. Moreover, we explicitely describe the set of Cauchon diagrams for explicit reduced decomposition of $w_0$ in each possible type. In any case, we check that the number of Cauchon diagrams is always equal to the cardinal of $W$. In \cite{CM}, we use these results to prove that Cauchon diagrams correspond canonically to the positive subexpressions of $w_0$. So the results of this paper also give an algorithmic description of the positive subexpressions of any reduced decomposition of $w_0$ corresponding to a good ordering.  \\

$ $ \\

\section{Introduction}
Let $\mathfrak{g}$ be a finite dimensional complex simple Lie algebra of rank $n$, $\mathbb{K}$ a commutative field and $q$ a nonzero element of $\mathbb{K}$ which is not a root of unity. 
We follow the notation and convention of  \cite{J} for the quantum group $\mathcal{U}_q(\mathfrak{g})$. In particular, to each choice of a reduced decomposition of the longest Weyl word  $w_0$ of the Weyl group $W$ corresponds a generating system $(X_\beta)_{ \beta \in \Phi^+}$) of the positive part $\mathcal{U}_q^+(\mathfrak{g})$ of $\mathcal{U}_q(\mathfrak{g})$ (see Section 3), where $\Phi^+$ denotes the set of positive roots associated to $\mathfrak{g}$. 

The natural action of an $n$-dimensional torus on $\mathcal{U}_q^+(\mathfrak{g})$ induces a stratification of the prime spectrum Spec$(\mathcal{U}_q^+(\mathfrak{g}))$ of $\mathcal{U}_q^+(\mathfrak{g})$ via the so-called Stratification Theorem (see \cite{MR1615971}). 
In this stratification, the primitive ideals are easily identified: they are the primes of $\mathcal{U}_q^+(\mathfrak{g})$ that are maximal in their strata. This stratification 
was recently used in \cite{dumas}, \cite{Launois1} and \cite{Launois2} in order to describe the automorphism group of 
$\mathcal{U}_q^+(\mathfrak{g})$ in the case where $\mathfrak{g}$ is of type 
A$_2$ and B$_2$.    \\

As $\mathcal{U}_q^+(\mathfrak{g})$ can be presented as a skew-polynomial algebra, this stratification can also be described via the deleting-derivations theory of Cauchon \cite{MR1967309}. In particular, in this theory, the strata are in a natural bijection with certain combinatorial objects, called Cauchon diagrams, and their geometry is completely described by the associated diagram. In fact, in the above situation, Cauchon diagrams are 
distinguished subsets of the set of positive roots $\Phi^+$. (For this reason, we often refer to subsets of $\Phi^+$ as diagrams.) Note that to each reduced decomposition of $w_0$ corresponds a PBW basis of $\mathcal{U}_q^+(\mathfrak{g})$ and so a notion of Cauchon diagrams.\\

The main aim of this paper is to give an algorithmic description of Cauchon diagrams in the case where the reduced decomposition of $w_0$ corresponds to a good order of $\Phi^+$ (in the sense of \cite{Lu2}). Moreover, in each type, we exhibit a reduced decomposition of $w_0$ for which we are able to describe explicitly the corresponding Cauchon diagrams. \\

Our first ingredient in order to obtain an algorithmic description of Cauchon diagrams is the commutation relation between two generators $X_\beta$ and $X_{\beta'}$ given by Levendorskii and Soibelmann \cite{MR1116413}. These formulas are not explicitly known, so that one cannot easily use them in order to perform the deleting-derivations algorithm. As a consequence, the description of Cauchon diagrams does not seem accessible in the general case. For this reason, we will limit ourselves to the case where the reduced decomposition of $w_0$ corresponds to a good order on $\Phi^+$ (see \cite{Lu2}). We recall this notion in Section 2. Although we still do not know explicitly all the commutation relations between the generators of $\mathcal{U}_q^+(\mathfrak{g})$, the situation is better as we control enough commutation relations. More precisely, in this case, the commutation relation between two variables $X_\beta$, $X_{\beta'}$ is known when $\beta$ and $\beta'$ span a so-called admissible plane \cite{Lu2} (see Section 3.4). Those relations allow the (algorithmic) construction of a set of necessary conditions, called \textit{implications}, for a diagram $\Delta$ to be a Cauchon diagram (see Section 5.1). In Section 5.2, we prove that these conditions are necessary and sufficient (see Theorem \ref{reciproque contrainte}), so that we get an algorithmic description of Cauchon diagrams.\\
$ $\\
In Section 6, we use this theorem to give an explicit description of these implications and these diagrams for special choices of the reduced decomposition of $w_0$ . More precisely, in each type, we exhibit a reduced decomposition of $w_0$ for which we explicitly describe the corresponding Cauchon diagrams. As a corollary, we prove that in each type the number of diagrams is equal to the size $| W |$ of the Weyl group. As the strata do not depend on the choice of the reduced decomposition of $w_{0}$, this implies that the number of strata is always equal to $| W |$. This result was first proved by Gorelik \cite{MR1754232} by using different methods, but with the additional assumption that $q$ is transcendental. \\

In \cite{CM}, we use the results of this paper in order to show that Cauchon diagrams $\Delta$ are in one-to-one correspondence with positive sub-expressions \textbf{w}$^\Delta$  of $w_0$ as defined by Marsh and Rietsch \cite{MR2058727}. More precisely, assume $w_0$ has a reduced expression of the form $ w_{0} = s_{\alpha_1} \circ ... \circ s_{\alpha_N}$ and 
that this decomposition corresponds to a good order on $\Phi^+$. For all $i \in \{1, \cdots , N\}$, we set $\beta_i=s_{\alpha_1} \circ ... \circ s_{\alpha_{i-1}}(\alpha_i)$, so that 
$\Phi^+= \{ \beta_1 < \cdots < \beta_N\}$. For each diagram $\Delta = \{\beta_{i_1} < ... < \beta_{i_t} \}\subseteq \Phi^+$, we set $w^{\Delta} := s_{\alpha_{i_1}} \circ ... \circ s_{\alpha_{i_t}}$. Then we have the following results (see \cite{CM}).
$ $ \\
\begin{itemize}
\item[$\bullet$] If $\Delta$ is a Cauchon diagram, the above decomposition of $w^\Delta$ is reduced.
\item[$\bullet$] The map $\Delta \rightarrow w^{\Delta}$ is a bijection from the set $\mathcal{D}$ of Cauchon diagrams to $W$.
\end{itemize}
\section{Root systems}
\subsection{Classical results on root systems}
Let $\mathfrak{g}$ be a simple complex Lie algebra. Let follow the notations of \cite[Chap 4]{J}.\\
\par
\begin{nota} 
\begin{itemizedot}
\item We denote by $\Phi$ a root system and E = Vect($\Phi$) (dim E = n). When $\Pi := \{\alpha_{1}, ..., \alpha_{n}\}$ is a basis of $\Phi$, one has a decomposition $\Phi = \Phi^+   \sqcup \Phi^-$, where $\Phi^+$ (resp. $\Phi^-$) is the set of positive (resp. negative) roots .
\item Denote by W the Weyl group associated to the root system $\Phi$ ; it is generated by the reflections $s_{\alpha_{i}} (:= s_{i}) , 1 \leq i \leq n$. The longest Weyl word in W is written $w_{0}$.
\end{itemizedot}
\end{nota} 
\begin{dfn}
A root system $\Phi$ is \emph{reducible} if $\Phi = \Phi_1   \sqcup \Phi_2$ where $\Phi_1$ and $\Phi_2$ are two orthogonal root systems. Otherwise $\Phi$ is called \emph{irreducible}.\\  
\end{dfn}
Let us recall that there is a one to one correspondence between the irreducible root systems and the simple complex Lie algebras of finite dimension. We say that $\mathfrak{g}$  is of a given type if the associated root system is of the same type. The following definitions and results are taken from \cite{Lu2}.
\begin{dfn}
Let $\Pi =\{\alpha_1, \alpha_2, ... , \alpha_n\}$ be a basis of $\Phi$ and $j \in \llbracket 1, n \rrbracket \ (:= \{1, 2, ..., n \} )$.
\begin{enumerate}
\item \emph{The column j} is the set $C_j := \{\beta \in \Phi^+ | \ \beta = k_1 \alpha_1 + \ldots + k_j \alpha_j, k_i \in \mathbb{N}, k_j \neq 0\}$;
\item A root $\beta = k_1 \alpha_1 + \ldots + k_j \alpha_j \in C_{j}$ is called \emph{ordinary} if $k_j =1$; it is called \emph{exceptional} if $k_j = 2$;
\item A column  $C_j$ is called \emph{ordinary} if each root $\beta$ of $C_j$ is ordinary; this column is called \emph{exceptional}
if every root $\beta$ of $C_j$ is ordinary except a unique one ($\beta_{ex}$) which is exceptional.
\end{enumerate}
\end{dfn}
\begin{dfn}
A numbering $\Pi =\{\alpha_1, \alpha_2, ... , \alpha_n\}$ is \emph{good} if all columns $C_j$ are ordinary or exceptional.
\end{dfn}
\begin{ex}[The $G_2$ case]\label{exG2}
The root system pf type $G_2$ has rank 2, there are two simple roots $\alpha_1$ and $\alpha_2$ such that $\|\alpha_{2}\| = \sqrt{3} \|\alpha_{1}\|$. $\Pi =\{\alpha_1, \alpha_2 \}$ is a base for this roots system. The numbering $\Pi =\{\alpha_1, \alpha_2 \}$ is good in this case because $C_1 = \{ \alpha_1 \}$ is ordinary and $C_2 = \{\alpha_2, \alpha_1 + \alpha_2, 2\alpha_1 + \alpha_2, 3\alpha_1 + \alpha_2, 3\alpha_1 + 2\alpha_2  \}$ is exceptional.\\
On the contrary, the numbering $\Pi =\{\alpha_2 , \alpha_1 \}$ is not good.  For this numbering, $C_1 = \{ \alpha_2 \}$ is ordinary but $C_2 = \{\alpha_1, \alpha_2 + \alpha_1,\alpha_2 + 2\alpha_1, \alpha_2 +\mathbf{ 3}\alpha_1, 2\alpha_2 + \mathbf{3}\alpha_1  \}$ is neither ordinary nor exceptional.
\end{ex}
\begin{prop}
Les $\mathfrak{g}$ be a simple Lie algebra of finite dimension. The following numberings of the associated root system $\Pi$ are examples of good numberings.
\begin{itemizedot}
\item If $\mathfrak{g}$ is of type $A_n$, with Dynkin diagram : $ \alpha_1 - \alpha_2 - \cdots - \alpha_{n-1} - \alpha_n $, then\\
\[\Pi = \{\alpha_1 , \alpha_2 , \cdots , \alpha_{n-1} , \alpha_n\} \ \text{is a good numbering}.\] 
\item If $\mathfrak{g}$ is of type $B_n$, with Dynkin diagram : $\alpha_1 \Leftarrow \alpha_2 - \cdots - \alpha_{n-1} - \alpha_n $, then\\
\[\Pi = \{\alpha_1 , \alpha_2 , \cdots , \alpha_{n-1} , \alpha_n\} \ \text{is a good numbering}. \]
\item If $\mathfrak{g}$ is of type $C_n$, with Dynkin diagram : $\alpha_1 \Rightarrow \alpha_2 - \cdots - \alpha_{n-1} - \alpha_n $, then \\
\[\Pi = \{\alpha_1 , \alpha_2 , \cdots , \alpha_{n-1} , \alpha_n\} \ \text{is a good numbering}.\]
\item If $\mathfrak{g}$ is of type $D_n$, with Dynkin diagram : 
\scalebox{.8}{$\begin{array}{ccccccccccc}
\alpha_1 && & &   & &&&  &     &   \\ 
& \diagdown & & &   & &&&  &     & \\
&& \alpha_3 &\text{---~}& \alpha_{4} &\text{---~}& \cdots &\text{---~}& \alpha_{n-1}&\text{---~} & \alpha_n \\
& \diagup & & &   & &&&  &     & \\
\alpha_2&& & &  &  &&&&  &
\end{array}$}, then
\[ \Pi = \{\alpha_1 , \alpha_2 , \cdots , \alpha_{n-1} , \alpha_n\} \ \text{is a good numbering}.\]
\item If $\mathfrak{g}$ is of type $G_2$, with Dynkin diagram : $ \alpha_1 \Lleftarrow \alpha_2 $, then \[\Pi = \{\alpha_1 , \alpha_2 \} \ \text{is a good numbering}.\]
\item If $\mathfrak{g}$ is of type $F_4$, with Dynkin diagram : $\alpha_1 \ \text{---~} \alpha_2 \Rightarrow \alpha_3 \ \text{---~} \alpha_4$, then \\
\[\Pi = \{\alpha_4 , \alpha_3, \alpha_2, \alpha_1 \} \ \text{is a good numbering}.\]
\item If $\mathfrak{g}$ is of type $E_6$, with Dynkin diagram : \scalebox{.8}{$\begin{array}{ccccccccc}
 					& 					&   				 & 					 & \alpha_2	&&  &     &   \\ 
 					& 					&   				 & 					 & |					&&  &     & \\
 \alpha_1 &\text{---~}& \alpha_{3} &\text{---~}& \alpha_4 &\text{---~}& \alpha_{5}&\text{---~} & \alpha_6 \\
\end{array}$} then, 
\[\Pi = \{\alpha_2 , \alpha_5, \alpha_{4}, \alpha_3, \alpha_1, \alpha_6 \} \ \text{is a good numbering}.\]
\item If $\mathfrak{g}$ is of type $E_7$, with Dynkin diagram :\scalebox{.8}{$\begin{array}{ccccccccccc}
 					& 					&   				 & 					 & \alpha_2	&&  &     &&&  \\ 
 					& 					&   				 & 					 & |					&&  &     & &&\\
 \alpha_1 &\text{---~}& \alpha_{3} &\text{---~}& \alpha_4 &\text{---~}& \alpha_{5}&\text{---~} & \alpha_6 &\text{---~}& \alpha_7 \\
\end{array}$}, then
\[\Pi = \{\alpha_2 , \alpha_5, \alpha_{4}, \alpha_3, \alpha_1, \alpha_6, \alpha_7 \} \ \text{is a good numbering}.\]
\item If $\mathfrak{g}$ is of type $E_8$, with Dynkin diagram : \scalebox{.8}{$\begin{array}{ccccccccccccc}
 					& 					&   				 & 					 & \alpha_2	&&  &     &&& &&  \\ 
 					& 					&   				 & 					 & |					&&  &     &&& &&\\
 \alpha_1 &\text{---~}& \alpha_{3} &\text{---~}& \alpha_4 &\text{---~}& \alpha_{5}&\text{---~} & \alpha_6 &\text{---~}& \alpha_7& \text{---~}& \alpha_8\\
\end{array}$}, then
\[\Pi = \{\alpha_2 , \alpha_5, \alpha_{4}, \alpha_3, \alpha_1, \alpha_6, \alpha_7, \alpha_8 \} \ \text{is a good numbering}.\]
\end{itemizedot}
\end{prop}
The corresponding columns with these numberings are given explicitly in the section 3 and one could verify that each column is ordinary or exceptional. \textbf{In the following, the chosen numbering on $\Pi$ is always a good one. From Section 6
, we use the numbering from the previous proposition.}

\subsection{Lusztig order}
The following definition is from \cite[section 4.3]{Lu2}.
\begin{dfn}
\begin{itemizedot}
\item For a root $\beta = k_1 \alpha_1 + \ldots + k_j \alpha_j \in C_j$, the \emph{height of} $\beta$ is the positive integer $h(\beta) := k_1 + \cdots + k_j$ ; the \emph{Lusztig height} of $\beta$ is the rational number $h'(\beta) := \frac{1}{k_j} h(\beta)$. 
\item If $t \in h'(C_j)$, then the set $B^{j,t}:=\{\beta \in C_j | h'(\beta)=t\}$ is called \emph{the box of height t in the column} $C_j$. 
\end{itemizedot}
\end{dfn}
\begin{rmq}
 \[\displaystyle C_{j} = \bigsqcup_{t \in \mathbb{N}^{*}} B^{j,t}.\]
\end{rmq}
\begin{dfn}[Lusztig order on $\Phi^+$]
We define a partial order on $\Phi^+$ as follows. Let $\beta_1$ and $\beta_2$ be two positive roots,
\begin{itemizedot}
\item If $\beta_1 \in C_{j_1}$ and $\beta_2 \in C_{j_2}$ with $j_1 < j_2$, then $\beta_1 < \beta_2$.
\item If $\beta_1$ and $\beta_2$ are in the same column $C_j$ and if $h'(\beta_2) < h'(\beta_1)$, then $\beta_1 < \beta_2$.
\end{itemizedot}
One can refine the previous partial order in a total one by choosing arbitrarily an order inside the boxes. Such a total order on $\Phi^+$ is called \emph{"a" Lusztig order}.
\end{dfn}
\begin{obs}
\begin{itemizedot}
\item $\alpha_{j}$ is the greatest root in $C_{j}$ for any Lusztig order.
\item The positive roots of a box are consecutive for any Lusztig order, that is, $B^{j,t} = \{\beta_{p}, \beta_{p+1}, ..., \beta_{p+l}\}$.
\item Any Lusztig order induces an order on boxes. For example, the box before $B^{j,t}$ in the column $C_j$ is $B^{j,t+1}$.
\end{itemizedot}
\end{obs}

\begin{prop}\label{prop2.8}
Let $j \in \{2, ..., n\}$. Assume $C_j$ is an exceptional column, we denote by $\beta_{ex}$ its exceptional root. Then :
\begin{enumerate}
\item $\beta_{ex} \bot (C_1 \sqcup ... \sqcup C_{j-1})$
\item If $D = <\beta_{ex}>$ and if $s_D $ is the orthogonal against D, then :
\begin{itemizedot}
\item $s_D(C_j) = C_j$ and for any $ \beta \in C_{j}\setminus\{\beta_{ex}\}$, we have $\beta + s_{D}(\beta) = \beta_{ex}$.
\item Let $B^{j,t}$ be a box different from the box which contain $\beta_{ex}$. Then $s_D$ transforms $B^{j,t}$ into $B^{j,h(\beta_{ex})-t}$.
\end{itemizedot}
\end{enumerate}
\end{prop}
\begin{preuve}
\begin{enumerate}
\item Let $\beta \in C_1 \cup ... \cup C_{j-1}$. If $\beta$ is not orthogonal to $\beta_{ex}$, then $s_{\beta}(\beta_{ex}) = \beta_{ex} + k\beta \ (k \in \mathbb{Z}\setminus\{0\})$ is a root of $C_j$ whose  coordinate on $\alpha_{j}$ is equal to 2. This is a contradiction with the unicity of the exceptional root.
\item We observe that $s_D = -s_{\beta_{ex}}$, so that $s_D(\Phi) = \Phi$.
\begin{itemizedot}
\item Let $\beta$ be an ordinary root of $C_j$. We can decompose this root $$\beta = a_1 \alpha_1 + ... + a_{j-1} \alpha_{j-1} + \frac{1}{2}  \beta_{ex} \ (a_{i} \in \mathbb{Q}).$$ From 1. we deduce that $s_D(\beta) =   - a_1 \alpha_1  ...  - a_{j-1} \alpha_{j-1} + \frac{1}{2}  \beta_{ex} = \beta_{ex} - \beta$. This is a root from the previous observation. This root is clearly in $C_{j}$ since $\beta$ is in $C_{j}\setminus\{\beta_{ex}\}$.\\
\item By the previous assertion, $s_{D}$ transforms two element from $B^{j,t}$ into two roots of the same height. We deduce from this fact (and from the fact that $s_{D}$ is an involution) that $s_{D}(B^{j,t})$ is a box denoted by $B^{j,s}$. The formula $t + s = h(\beta_{ex})$ is a consequence of the previous assertion.
\end{itemizedot}
\end{enumerate}
\end{preuve}
\begin{dfn}
The support of a root $\beta = a_1 \alpha_1 + ... + a_{j} \alpha_{j} \in C_{j}$ is the set \emph{Supp} $(\beta) := \{\alpha_{i} \in \Pi \ | \ a_{i} \neq 0 \}$. In particular, for $\beta \in C_{j}$, we have $\textrm{\emph{Supp}}(\beta) \subset \{1, ..., j\}$. 
\end{dfn}
We are now ready to prove that the box containing the exceptional root of an exceptional column is reduced to the exceptional root.
\begin{prop}\label{boîteex}
Let $j \in \{2, ..., n\}$ Assume $C_j$ is an exceptional column and denote by $\beta_{ex}$ its exceptional root. Then $h'(\beta_{ex}) \notin \mathbb{N}$, so that $\beta_{ex}$ is alone in its box.
\end{prop}
\begin{preuve}
Denote $\Pi_{j} = \{\alpha_{1},..., \alpha_{j}\}$ and $\Phi_{j} = \Phi \cap \text{Vect}(\Pi_{j})$. Then $\Phi_{j}$ is a root system with basis $\Pi_{j}$ and $\Phi_{j}^{+} = \Phi^{+} \cap \text{Vect}(\Pi_{j})$. \\
\par
Let us consider the case where $\Phi_{j}$ irreducible. Then we have \\
\par
\textbf{Observation 1.} If $\beta$ is a root of $\Phi_{j}^{+}$ of maximal height then $\beta \in C_{j}$.\\
\textbf{Proof of Observation 1. :}\\
Assume that $\beta \in C_{i}$ with $i < j$. In the Dynkin diagram of $\Phi_{j}$ which is convex as $\Phi_{j}$ is irreducible, we can construct a path from $\alpha_{i}$ to $\alpha_{j}$. Denote this path by $P = (\alpha_{i_{1}}, ..., \alpha_{i_{t}})$, where $i_{1} = i$ and $i_{t} = j$. We know that $\alpha_{i} \in$ Supp$(\beta)$ and that $\alpha_{j} \notin$ Supp$(\beta)$. So there is a smallest index $l$ such that $\alpha_{i_{l}} \in Supp(\beta)$ and $\alpha_{i_{l+1}} \notin Supp(\beta)$. Thus, for all $ \alpha \in Supp(\beta)$, we have $\langle \alpha , \alpha_{i_{l+1}} \rangle  \leq 0$ and, since $\alpha_{i_{l}}$ and $\alpha_{i_{l+1}}$ are two consecutive elements from P, we have $\langle \alpha_{i_{l}} , \alpha_{i_{l+1}} \rangle < 0$. We deduce that, $\langle \beta, \alpha_{i_{l+1}} \rangle < 0$ thus $\beta + \alpha_{i_{l+1}} \in \Phi_{j}^+$ which contradicts the maximality of the height of $\beta$.\\
\par
\textbf{Observation 2.} $\beta_{ex}$ is the largest root of $\Phi_{j}$.\\
\textbf{Proof of Observation 2. :}
Let $\beta$ be a largest root in $\Phi_{j}$. We assume that $\beta \neq \beta_{ex}$. By the previous observation, we know that $\beta \in C_{j}$ and, by Proposition \ref{prop2.8}, $\beta_{ex} = \beta + s_{D}(\beta)$ is a sum of two positive roots, thus its height is greater than the height of $\beta$. So $\beta$ is equal to $\beta_{ex}$.\\
\par
We note that yhe existence of an exceptional root implies that $\Phi_{j}$ is not of type $A_{j}$. So $\Phi_{j}$ is of type $B_{j}, C_{j}, D_{j}, E_{6}, E_{7}, E_{8}, F_{4}$ or $G_{2}$ and, we deduce from  \cite{MR0240238} that the height of the largest root is odd. Hence it follows from Observation 2. that the height of $\beta_{ex}$ is odd, so that$h'(\beta_{ex}) \notin \mathbb{N}$.\\
\par
Let us now assume that $\Phi_{j}$ is reducible. Denote by $\Gamma_{j}$ the Dynkin diagram whose vertices are $\alpha_{1}, ..., \alpha_{j}$, and whose edges come from the Dynkin diagram of $\Phi$. We note $\Pi'$ the connected component of $\alpha_{j}$ in $\Gamma_{j}$; 
\[\Pi':=\{\alpha_{i} \in \Pi_{j} \ | \ \textrm{there exists a path contained in} \ \Gamma_{j} \ \textrm{connecting} \ \alpha_{i} \ \textrm{to} \ \alpha_{j}\}.\]
We note $\Phi' = \Phi \cap \text{Vect}(\Pi')$. It is a root system, whith basis $\Pi'$, and we have $\Phi'^+ = \Phi^+ \cap \text{Vect}(\Pi')$.\\
\par
\textbf{Observation 3.} $C_{j} \subseteq \Phi'^+$.\\
\textbf{Proof of Observation 3. :}\\
Otherwise, there is a root in $ C_{j} \setminus \Phi'^{+}$. If $\beta$ is such a root, there is a simple roots in its support which is also in the set $\Pi_{j}\setminus \Pi'$. As the support of $\beta$ contains $\alpha_{j} \in \Pi'$ too, we can write $\beta = u + v$ with $ u = \alpha_{i_{1}} + ... + \alpha_{i_{l}}$ whose support is a subset of $\Pi_{j}\setminus \Pi'$ and $v = \alpha_{i_{l+1}} + ... + \alpha_{i_{p}}$ whose support is a subset of $\Pi'$. Let us choose $\beta$ such that the integer $l$ is minimal.
\begin{itemizedot}
\item If $l = 1,$ then  $\beta = \alpha_{i_{1}} + v$. As $\alpha_{i_{1}} \notin \Pi'$, there is no link between $\alpha_{i_{1}}$ and the element of Supp$(v)$. Then $s_{i_{1}}(\beta) = -\alpha_{i_{1}} + v \in \Phi$, which is impossible because the coordinates of this root in the basis $\Pi$ do not have the same sign. 
\item So $l \geq 2$. As $\langle u, u \rangle > 0$, there exists a root in Supp$(u)$, for example $\alpha_{i_{l}}$, such that $\langle u, \alpha_{i_{l}} \rangle > 0$. As above, we have :
$$\langle v, \alpha_{i_{l}} \rangle = 0 \Rightarrow \langle \beta, \alpha_{i_{l}} \rangle > 0 \Rightarrow \beta' =\beta-\alpha_{i_{l}} \in C_{j}\setminus \Phi'^{+},$$ which is a contradiction with the minimality of $l$.
\end{itemizedot}
So we can conclude that $C_{j} \subseteq \Phi'^+$.\\
\par
Hence $C_{j}$ is an exceptional column of $\Phi'$ which is irreducible by construction. The proof above also shows that the exceptional root $\beta_{ex}$ satisfies $h'(\beta_{ex}) \notin \mathbb{N}$.\\
\end{preuve}

We can now proove that any Lusztig order is a convex order. 
\begin{prop}
"<" is a convex order over $\Phi^+$
\end{prop}
\begin{preuve}
Let $\beta_1 < \beta_2$ be two positives roots such that $\beta_1 + \beta_2 \in \Phi^+$. 
\begin{itemizedot}
\item If the two roots $\beta_{1}$ and $\beta_{2}$ do not belong to the same column, then $\beta_{1} + \beta_{2}$ is in the same column as $\beta_{2}$. In this case, neither $\beta_{2}$, nor $\beta_{1} + \beta_{2}$ are exceptional and 
\[h'(\beta_{1} + \beta_{2}) = h(\beta_{1} + \beta_{2}) = h(\beta_{1}) + h(\beta_{2}) > h(\beta_2) = h'(\beta_{2}).\]
Hence we have : $\beta_1 <\beta_1 + \beta_2 < \beta_2$.\\
\item If the two roots $\beta_1$ and $\beta_2$ belong to the same column, then $\beta_{1} + \beta_{2}$ is an exceptional root. We deduce from Proposition \ref{prop2.8} that $h'(\beta_{1} + \beta_{2}) = \frac{h'(\beta_{1})+h'(\beta_{2})}{2}$. Proposition \ref{boîteex} excludes the case where  $h'(\beta_{1} + \beta_{2}) = h'(\beta_{1}) = h'(\beta_{2})$ because the exceptional root is alone in its box. So we get $h'(\beta_{1}) > h'(\beta_{1} + \beta_{2}) > h'(\beta_{2})$, so that $\beta_1 <\beta_1 + \beta_2 < \beta_2$.
\end{itemizedot}
\end{preuve}
Consider a reduced decomposition of $w_0 =  s_{i_{1}} \circ s_{i_2} \circ ... \circ s_{i_N}$ of the longest Weyl word $w_0$. For all $j \in \llbracket 1 , N \rrbracket$, we set $\beta_j :=s_{i_1} \circ s_{i_2} \circ ... \circ s_{i_{j-1}}(\alpha_{i_{j}})$. Then it is well known (cf, for example, \cite[I.5.1]{MR1898492}) that $\{\beta_1,...,\beta_N\} = \Phi^+$. For each integer $j \in \llbracket 1, N \rrbracket$, we say that $\alpha_{i_{j}}$ is \emph{the simple root associated to the positive root} $\beta_{j}$.\\
We define an order on $\Phi^+$ by setting $\beta_i \prec \beta_j$ when $i<j$. We say that "$\prec$" is \emph{the order associated to the reduced decomposition} $w_0 = s_{i_{1}} \circ s_{i_2} \circ ... \circ s_{i_N}$ of $w_0$.\\
In \cite[Theorem and remark page 662]{Papi}, it is shown that this is a convex order and that this leads to a one to one correspondence between reduced decompositions of $w_0$ and convex orders on $\Phi^+$.\\
\textbf{Hence, as the Lusztig order "<" is convex, there is a unique reduced decomposition} \\
$\mathbf{w_0 = s_{i'_{1}} \circ s_{i'_2} \circ ... \circ s_{i'_N}}$ of $w_0$ \textbf{whose associated order is "<". In this article, we always choose this decomposition for} $\mathbf{w_0}$.\\
\par
The following proposition of Lusztig \cite[section 4.3]{Lu2} explains the behaviour of the positive roots inside (non-exceptional) boxes.
\begin{prop}\label{boîte}
 Inside each ordinary box (box which does not contain the exceptional root), roots are pairwise orthogonal. Moreover, simple roots associated to the positive roots of a given box are pairwise orthogonal.
\end{prop}
\begin{preuve}
For the type $G_{2}$, explicit computations leads to the result. We now assume that $\mathfrak{g}$ is a finite dimensional simple Lie algebra which is not of type $G_2$.  Recall that the positive roots of a box are consectutives. Let $\beta_{1}$ and $\beta_{2}$ be two consecutive roots of a box B in the column $C_{j}$. We note $\alpha_{i_{1}}$ and $\alpha_{i_{2}}$ the associated simple roots.\\
Suppose that $\alpha_{i_{1}}$ is not orthogonal to  $\alpha_{i_{2}}$, hence $\lambda = -<\alpha_{i_{1}}^{\vee},\alpha_{i_{2}}> = 1$ or $2$ (recall that $\mathfrak{g}$ is not of type $G_2$). So we can write $\beta_{2}  = w \circ s_{i_{1}} (\alpha_{i_{2}}) = w(\lambda \alpha_{i_{1}} + \alpha_{i_{2}}) = \lambda \beta_{1} + w(\alpha_{i_{2}})$.As  $w(\alpha_{i_{2}}) \in \Phi$, we must have $\lambda = 2$, otherwise $h(w(\alpha_{i_{2}})) = h(\beta_{2}) - h(\beta_{1}) = 0 $, which is absurd.\\
In this case, $\gamma = - w(\alpha_{i_{2}}) = 2\beta_{1} - \beta_{2} \in C_{j}$ and $h(\gamma) =  2h(\beta_{1})- h(\beta_{2}) = h(\beta_1)$. As $\beta_1$ and $\beta_2$ are distinct roots, so they are not collinear. So the set $\Phi' = \Phi \cap$  Vect$(\beta_1, \beta_2)$ is a root system of rank 2 which contains $\beta_1$, $\beta_2$, $\gamma$ and their opposites. The equality $2\beta_{1}  = \gamma  + \beta_{2}$ allows to state that $\Phi'$ is of type $B_2$ and that the situation is the following one :
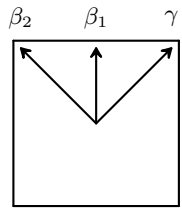
\begin{center}
\begin{pgfpicture}{-0.5cm}{-1cm}{2.9cm}{2.2cm}%
\pgfsetroundjoin \pgfsetroundcap%
\pgfsetlinewidth{0.8pt} 
\pgfsetfillcolor{black}
\pgfxyline(0.1,1.3)(2.3,1.3)
\pgfxyline(2.3,1.3)(2.3,-0.9)
\pgfxyline(2.3,-0.9)(0.1,-0.9)
\pgfxyline(0.1,-0.9)(0.1,1.3)
\pgfxyline(1.2,0.2)(2.2,1.2)
\pgfmoveto{\pgfxy(2.0712,1.1655)}\pgflineto{\pgfxy(2.2,1.2)}\pgflineto{\pgfxy(2.1655,1.0712)}\pgflineto{\pgfxy(2.1592,1.1592)}\pgfclosepath\pgffillstroke
\pgfxyline(1.2,0.2)(1.2,1.2)
\pgfmoveto{\pgfxy(1.1333,1.0845)}\pgflineto{\pgfxy(1.2,1.2)}\pgflineto{\pgfxy(1.2667,1.0845)}\pgflineto{\pgfxy(1.2,1.1423)}\pgfclosepath\pgffillstroke
\pgfxyline(1.2,0.2)(0.2,1.2)
\pgfmoveto{\pgfxy(0.2345,1.0712)}\pgflineto{\pgfxy(0.2,1.2)}\pgflineto{\pgfxy(0.3288,1.1655)}\pgflineto{\pgfxy(0.2408,1.1592)}\pgfclosepath\pgffillstroke
\pgfputat{\pgfxy(2.2,1.6)}{\pgfnode{rectangle}{center}{\small $\gamma$}{}{\pgfusepath{}}}
\pgfputat{\pgfxy(1.2,1.6)}{\pgfnode{rectangle}{center}{\small $\beta_{1}$}{}{\pgfusepath{}}}
\pgfputat{\pgfxy(0.2,1.6)}{\pgfnode{rectangle}{center}{\small $\beta_{2}$}{}{\pgfusepath{}}}
\end{pgfpicture} 
\end{center}
So $\gamma - \beta_{1} \in \Phi$, with $h(\gamma - \beta_1) = h(\gamma) - h(\beta_1) = 0 $. This is impossible, and so $\alpha_{i_1} \ \bot \ \alpha_{i_2}$.\\
Then we get $<\beta_{1},\beta_{2}> \ =  \ <w(\alpha_{i_{1}}),w(s_{i_{1}}(\alpha_{i_{2}}))> = <\alpha_{i_{1}},s_{i_{1}}(\alpha_{i_{2}})> \ = \ <\alpha_{i_{1}},\alpha_{i_{2}}> \ = 0$, as desired. This finishes the case where the two roots are consecutive. One concludes using an induction on the "distance" between the two roots $\beta_1$ and $\beta_2$.
\end{preuve}
\begin{convsanss}
For $j \in \llbracket 1, n \rrbracket$, denote $\delta_j$ the smallest root of $C_j$. Let us recall that $\alpha_j$ is the largest root of $C_j$.
\end{convsanss}
\begin{prop}\label{firstroot}
$\delta_j$ and $\alpha_j$ are alone in their boxes.
\end{prop}
\begin{preuve}
The root $\alpha_{j}$ is alone in its box because it is the only roots of $C_{j}$ whose height is equal to 1.\\
To prove that $\delta_j$ is alone in its box, we need the following result. 
\begin{lemme}
Let $1 \leq l \leq N$ and $1 \leq m \leq n$. Set $\Pi_m := \{\alpha_1, ..., \alpha_m \}$. If $\beta_l = s_{i_1}...s_{i_{l-1}}(\alpha_{i_l})$ is in the column $C_m$, then $\alpha_{i_j} \in \Pi_m$ for $j \in \llbracket 1, l \rrbracket$. 
\end{lemme}
\begin{preuve}
The proof is by induction on $l$.\\
If $l = 1$, then $\beta_1 = \alpha_1 \in \Pi_1$.\\
If $l \geq 2$, then $\beta_{l-1}$ is in the column $C_m$ or $C_{m-1}$. By the induction hypothesis, it implies that $\alpha_{i_t} \in \Pi_m$ (or $\alpha_{i_t} \in \Pi_{m-1} \subset \Pi_m$) for $t \in \llbracket 1, l-1 \rrbracket$. We observe that $\beta_l = s_{i_1}...s_{i_{l-1}}(\alpha_{i_l}) = \alpha_{i_l} + n_{l-1}\alpha_{i_{l-1}} + ... + n_1 \alpha_{i_{1}}$ where each $n_t$ is in $\mathbb{Z}$. As $\beta_l \in C_m$, necessarily $\alpha_{i_l} \in \Pi_m$.  
\end{preuve}
\textbf{Back to the proof of Proposition \ref{firstroot} :}\\
There is an integer $1 \leq l \leq N$ such that $\delta_{j} = \beta_{l} = s_{i_{1}} \circ s_{i_2} \circ ... \circ s_{i_{l-1}}(\alpha_{i_{l}})$. As above, $\beta_l =  \alpha_{i_l} + n_{l-1}\alpha_{i_{l-1}} + ... + n_1 \alpha_{i_{1}} \ (n_t \in \mathbb{Z})$ with $\alpha_{i_1}, ..., \alpha_{i_{l-1}}$ in $ \Pi_{j-1}$ since $\beta_{l-1} \in C_{j-1}$. As $\beta_l \in C_j$, it implies that $\alpha_{i_{l}} = \alpha_{j}$.\\
If $\delta_{j} (=\beta_{l}) $ is not alone in its box, then $\beta_{l+1}$ is also in this box and one has (Proposition \ref{boîte}) $\alpha_{i_{l}} \ \bot \ \alpha_{i_{l+1}}$. By the previous lemma, it implies that $\alpha_{i_{l+1}} \in \Pi_{j} \setminus \{\alpha_j\} = \Pi_{j-1}$ and $\beta_{l+1} = s_{i_{1}} \circ s_{i_2} \circ ...\circ s_{i_{l-1}} \circ s_{i_{l}}(\alpha_{i_{l+1}}) = s_{i_{1}} \circ s_{i_2} \circ ...\circ s_{i_{l-1}} (\alpha_{i_{l+1}}) = \alpha_{i_{l+1}} + n'_{l-1}\alpha_{i_{l-1}} + ... + n'_1 \alpha_{i_{1}} \ (n'_t \in \mathbb{Z})$, which contradicts the hypothesis $\beta_{l+1} \in C_{j}$.
\end{preuve}

Let us recall the following result (see, for example, \cite[Lemma 9.4]{Humphreys}).
\begin{lemme}\label{lemme2.9}
Let $\beta$ and $\delta$ be two distinct roots of $\Phi^+$ such that $\langle \beta,\delta \rangle \neq 0$
\begin{itemizedot}
\item If $\langle \beta,\delta \rangle > 0$, then $\beta -\delta \in \Phi$.
\item If $\langle \beta,\delta \rangle < 0$, then $\beta + \delta \in \Phi$.
\end{itemizedot}
\end{lemme}

\begin{prop}\label{prop2.10}
Let $\beta$ be an ordinary root of a column $C_j$. Denote (as in the proof of Proposition \ref{boîteex}) by $\Gamma_{j}$ the diagram whose vertices are $\alpha_{1}, ..., \alpha_{j}$, and whose edges are the edges from the Dynkin diagram of $\Phi$. Denote by $\Omega_j$ the connected component of $\alpha_j$ in $\Gamma_j$ .
\begin{enumerate}
\item If $\beta \neq \alpha_{j}$ then there exists $\epsilon \in \{\alpha_{1},...,\alpha_{j-1}\}$ such that $\beta - \epsilon \in C_j$.
\item Supp $\beta \subset \Omega_j$.
\item If $\beta \neq \delta_{j}$ then there exists $\epsilon \in \{\alpha_{1},...,\alpha_{j-1}\}$ such that $\beta + \epsilon \in C_j$.
\end{enumerate} 
\end{prop}
\begin{preuve}
One verifies that this proposition is true in the case $G_2$ thanks to the description of the columns given in Example \ref{exG2}. Assume from now on that the root system is not of type $G_2$. Set $\Pi_{j} = \{\alpha_{1},..., \alpha_{j}\}$.\\
\begin{enumerate}
\item If it exists $\epsilon \in \Pi_{j} \setminus \{\alpha_j\}$ such that $\langle \beta , \epsilon \rangle > 0$, one concludes with Lemma \ref{lemme2.9}. So let us assume that for all $\epsilon \in \Pi_{j} \setminus \{\alpha_j\}$, we have $\langle \beta , \epsilon \rangle \leq 0$.\\
As $\langle \beta, \beta \rangle > 0$ then $\langle \beta, \alpha_{j} \rangle > 0$ and $\beta - \alpha_{j} = \gamma_{1} \in \Phi^{+}$ with $\alpha_{j} \notin$ Supp$(\gamma_{1})$. One can write $\beta = \beta_{1} + \gamma_{1}$ with $\beta_{1} = \alpha_j \in C_{j}$ and, since $\beta$ is ordinary, we get $\alpha_{j} \notin Supp(\gamma_{1})$. Since $\beta \neq \alpha_j$, one has $h(\beta) \geq 2$ and we will prove by induction that, for each integer $i \in \llbracket 1, h(\beta)-1 \rrbracket$, one has $\beta = \beta_i + \gamma_i$ with $\beta_i \in C_j$, $\gamma_i \in \Phi^+$ and $h(\beta_i) = i$.\\
 \begin{itemizedot}
\item As $h(\beta_1) = h(\alpha_j) = 1$, the result is proved for $i = 1$.
\item Let us assume that the result is proved for $i$ with $1 \leq i < h(\beta)-1$ and observe that, since $\beta$ is ordinary, $\gamma_i \notin C_j$ so that Supp$(\gamma_i) \ \subset \Pi_{j-1}$. As $\langle \gamma_i, \gamma_i \rangle > 0$, there is $\epsilon \in \Pi_{j-1}$ such that $\langle \gamma_{i}, \epsilon \rangle > 0$. As $i < h(\beta)-1$, one has $h(\gamma_i) > 1$. Hence $\gamma_i \neq \epsilon$, and so we deduce from Lemma \ref{lemme2.9} that $\gamma_{i+1}:= \gamma_i - \epsilon$ belongs to $\Phi^+$.\\
So $\langle \beta , \epsilon \rangle = \langle \beta_{i} , \epsilon \rangle + \langle \gamma_{i} , \epsilon \rangle$ and, since $\langle \beta, \epsilon \rangle \leq 0$, one has $\langle \beta_{i} , \epsilon \rangle < 0$. So we deduce from Lemma \ref{lemme2.9} that $\beta_{i+1} := \beta_{i} + \epsilon \in C_{j}$.
Hence, $\beta = \beta_{i+1} + \gamma_{i+1}$ with $\beta_{i+1} \in C_j$ and $h(\beta_{i+1}) = h(\beta_{i})+1$.\\
\par
So we have proved the result and, for $i = h(\beta) - 1$, one has $h(\gamma_i) = 1$. So, since $\beta$ is ordinary, we have $\epsilon := \gamma_i \in \Pi_{j-1}$ and $\beta - \epsilon = \beta_i \in C_j$. 
\end{itemizedot}
\item The proof is by induction on $h(\beta)$.\\
If $h(\beta) = 1$, one has $\beta = \alpha_j \in \Omega_j$.\\
If $h(\beta) > 1$, we deduce from 1. that it exists $\epsilon \in \Pi_{j-1}$ such that $ \beta' := \beta - \epsilon \in C_j$. So $\beta'$ is ordinary and $h(\beta') = h(\beta) - 1$ so that, by the induction hypothesis, Supp $\beta' \subset \Omega_j$. It remains to prove that $\epsilon \in \Omega_j$.\\
If $\epsilon \notin \Omega_j$, then $\epsilon \ \bot \ \alpha$ for all $\alpha \in $ Supp $\beta'$.  As a result, $\epsilon \ \bot \ \beta'$. Let us consider the plane $P = <\beta, \epsilon>$ and observe that $\Phi_P = \Phi \cap P$ is a root system of rank 2. As $\Phi$ is not of type $G_2$, $\Phi_P$ is not of type $G_2$ as well. This forces $\Phi_P$ to be of type $B_2$. So one has the following configuration :
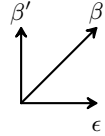
\begin{center}
\begin{pgfpicture}{-0.5cm}{-0.5cm}{2.9cm}{2.2cm}%
\pgfsetroundjoin \pgfsetroundcap%
\pgfsetlinewidth{0.8pt} 
\pgfxyline(1.2,0.2)(2.2,0.2)
\pgfsetfillcolor{black}
\pgfmoveto{\pgfxy(2.0845,0.2667)}\pgflineto{\pgfxy(2.2,0.2)}\pgflineto{\pgfxy(2.0845,0.1333)}\pgflineto{\pgfxy(2.1423,0.2)}\pgfclosepath\pgffillstroke
\pgfxyline(1.2,0.2)(2.2,1.2)
\pgfmoveto{\pgfxy(2.0712,1.1655)}\pgflineto{\pgfxy(2.2,1.2)}\pgflineto{\pgfxy(2.1655,1.0712)}\pgflineto{\pgfxy(2.1592,1.1592)}\pgfclosepath\pgffillstroke
\pgfxyline(1.2,0.2)(1.2,1.2)
\pgfmoveto{\pgfxy(1.1333,1.0845)}\pgflineto{\pgfxy(1.2,1.2)}\pgflineto{\pgfxy(1.2667,1.0845)}\pgflineto{\pgfxy(1.2,1.1423)}\pgfclosepath\pgffillstroke
\pgfputat{\pgfxy(2.2,-0.1)}{\pgfnode{rectangle}{center}{\small $\epsilon$}{}{\pgfusepath{}}}
\pgfputat{\pgfxy(2.2,1.4)}{\pgfnode{rectangle}{center}{\small $\beta$}{}{\pgfusepath{}}}
\pgfputat{\pgfxy(1.2,1.4)}{\pgfnode{rectangle}{center}{\small $\beta'$}{}{\pgfusepath{}}}
\end{pgfpicture} 
\end{center}
As a result, $\beta' - \epsilon$ is a root. As Supp $\beta' \subset \Omega_j$, one has $\epsilon \notin $ Supp $\beta'$ which is a contradiction with $\beta' - \epsilon$ being a root. So, one has $\epsilon \in \Omega_j$ and then, Supp $\beta \subset \Omega_j$.\\
\item Let $\beta$ be a root of $C_{j}$, not equal to $\delta_{j}$. If it exists $\epsilon \in \Pi_{j-1}$ such that $\langle \beta , \epsilon \rangle < 0$ then one concludes thanks to Lemma \ref{lemme2.9}.\\
Let us assume that $\langle \beta , \epsilon \rangle \geq 0$ for all $\epsilon \in \Pi_{j-1}$. If $\langle \beta, \alpha_j \rangle < 0$, then $C_{j}$ is exceptional and $\beta + \alpha_j = \beta_{ex}$. We deduce from Proposition \ref{prop2.8} that $s_D$ switches $\alpha_j$ and $\delta_{j}$, so that, $\beta_{ex} = \alpha_j + \delta_{j}$. It implies that $\beta = \delta_{j}$ which is a contradiction with the hypothesis. So one has $\langle \beta , \alpha_j \rangle \geq 0$ and then, $\langle \beta , \epsilon \rangle \geq 0$ for all $ \epsilon \in \Pi_{j}$.\\
Before finishing the proof, we need to establish two intermediaire results.
\begin{obsn}
$\textrm{\emph{Supp}} (\delta_{j}) = \Omega_{j}$.
\end{obsn}
\textbf{Proof of Observation 1.}\\
From 2. one has Supp $(\delta_{j}) \subset \Omega_{j}$. Let us suppose that this inclusion is a strict one and consider $\alpha \in \Omega_{j} \setminus \textrm{Supp} (\delta_{j})$. As $\Omega_j$ is connected, there exist $\alpha'_1, \alpha'_2, ..., \alpha'_s$ in $\Omega_j$, with $s \geq 2$, such that $\alpha'_1 = \alpha, \alpha'_s \in \ \text{Supp}(\delta_j)$ and, for $1 \leq i < s$, $\langle \alpha'_i , \alpha'_{i+1} \rangle < 0$.\\
Let $k$ be the greatest integer such that $\alpha'_k \notin$ Supp $(\delta_j)$. One has $k < s$ and $\alpha'_{k+1} \in$ Supp $(\delta_j)$. As $\delta_j$ is then a linear combination whose coefficients are positive integers of simple roots which are not equal to $\alpha'_k$ and, as the scalar product of two distinct simple roots is nonpositive, one has :
\[\langle \alpha'_k , \delta_j \rangle \leq \langle \alpha'_k, \alpha'_{k+1} \rangle < 0.\]
Then we deduce from Lemma \ref{lemme2.9} that $\delta_j + \alpha'_k \in C_j \setminus \{\beta_{ex}\}$, which is a contradiction with the fact that $\delta_j$ is the smallest root of $C_j$.
\begin{obsn}
$\langle \beta, \delta_{j} \rangle > 0$
\end{obsn}
\textbf{Proof of Observation 2.}\\
As $\langle \beta, \beta \rangle > 0$, there is $\alpha \in $ Supp$(\beta)$ such that $\langle \beta , \alpha \rangle > 0$. By  2. and Observation 1, one has $\alpha \in $ Supp $\delta_j$. As we have already proved that $\langle \beta , \epsilon \rangle \geq 0$ for all $ \epsilon \in \Pi_j$, one has $\langle \beta , \epsilon \rangle \geq 0$ for all $ \epsilon \in$ Supp $\delta_j = \Omega_j \subset \Pi_j \Rightarrow \langle \beta, \delta_{j} \rangle \geq \langle \beta, \alpha \rangle > 0$.\\
\par
By Observation 2 and Lemma \ref{lemme2.9}, we have $\gamma_{1} := \delta_{j} - \beta \in \Phi^{+}$ and, since $\beta$ and $\delta_j$ are ordinary roots of $C_j$, $\textrm{Supp} (\gamma_{1}) \subset \Pi_{j-1}$. We prove now, by decreasing induction that, for each integer $i \in \llbracket 1, h(\delta_j - \beta) \rrbracket$, there exists $\rho_i \in \Phi^+$ with Supp $(\rho_i) \subset \Pi_{j-1}$, such that $h(\rho_i) = i$ and $\beta + \rho_i \in C_j$.\\
\begin{itemizedot}
\item For $i = h(\delta_j- \beta)$, the result holds with $\rho_i = \gamma_1$ .
\item Assume $i \in \llbracket 2, h(\delta_j- \beta) \rrbracket$ and assume that there exists $\rho_i$ such that $\eta_i = \beta + \rho_i \in C_j$. As $\langle \rho_i , \rho_i \rangle > 0$, there is $\epsilon \in $ Supp $(\rho_i) \subset \Pi_{j-1}$ such that $\langle \rho_i , \epsilon \rangle > 0$. As we have already proved that $\langle \beta, \epsilon \rangle \geq 0$, this implies that $\langle \eta_i , \epsilon \rangle > 0$. So we deduce from Lemma \ref{lemme2.9} that $\eta_{i-1} := \eta_i - \epsilon \in C_j$ and $\rho_{i-1} := \rho_i - \epsilon \in \Phi^+$. Clearly one has Supp $(\rho_{i-1}) \subset \Pi_{j-1}$, $h(\rho_{i-1}) = i-1$ and $\beta + \rho_{i-1} = \eta_{i-1} \in C_j$.\\
So the result holds and, for $i = 1$, $\epsilon := \rho_1 \in \Pi_{j-1}$ and $\beta + \epsilon \in C_j$.
\end{itemizedot} 
\end{enumerate}
\end{preuve} 
\begin{prop}\label{h'}
Let $j \in \llbracket 1, n \rrbracket$.
\begin{enumerate}
\item If $C_j$ is ordinary, then $h'(C_j)$ is an interval of the form $\llbracket 1, t \rrbracket$;
\item If $C_j$ is exceptional, then $h'(C_j \setminus \{\beta_{ex}\})$ is an interval of the form $\llbracket 1, 2t \rrbracket \ (t \in \mathbb{N}^\star)$.\\
Moreover we have $h'(\beta_{ex}) = t + \frac{1}{2}$.
\end{enumerate}
\end{prop}
\begin{preuve}
The fact that $h'(C_j)$ in the ordinary case (resp. $h'(C_j \setminus \{\beta_{ex}\})$ in the exceptional case) is an interval of integers comes from Proposition \ref{prop2.10}. It contains $1 = h(\alpha_j)$, and so the first case is proved.\\
Let us assume that $C_j$ is exceptional. Denote by $B_{1},...,B_{t}$ the boxes which contain the roots smaller than $\beta_{ex}$ for the Lusztig order. For these boxes, we have $h'(B_{i}) > h'(\beta_{ex})$. But the relation $h(B_{i}) + h(B_{i}') = h(\beta_{ex})$, for the image $B'_{i}$ of $B_{i}$ by $s_{D}$, implies $h'(B_{i}) > h'(\beta_{ex}) > h'(B'_{i})$. So we have exactly t boxes appearing after $\beta_{ex}$ and the interval $h'(C_j \setminus \{\beta_{ex}\})$ is of the form $\llbracket 1, 2t \rrbracket \ (t \in \mathbb{N}^\star)$.\\
Moreover $h(\beta_{ex}) = h(\alpha_{j} + s_{D}(\alpha_{j})) = 1 + 2t$ and finally $h'(\beta_{ex}) = t + \frac{1}{2}$.
\end{preuve}
We now recall the notion of admissible planes introduced by Lusztig in \cite[Section 6.1]{Lu2}.
\begin{dfn}\label{planadm}
We call \emph{admissible plane} $P := <\beta, \beta'>$ a plane spanned by two positive roots $\beta$ and $\beta'$ such that :
$\beta$ belongs to an exceptional column $C_{j}$ and $\beta' = s_{D}(\beta)$ is such that $|h'(\beta') - h'(\beta)| = 1$. (In this case $\beta + \beta' = \beta_{ex}$ and  $h'(\beta_{ex}) = t \pm \frac{1}{2}$).\\
\textbf{or} $\beta$ is an ordinary root in any column  $C_{j}$ and $\beta' = \alpha_{i}$ with $i < j$.
We set $\Phi_P := \Phi \cap P$ and $\Phi_P^+ := \Phi^+ \cap P$.
\end{dfn}
\begin{rmq}\label{systracine}
If $\Phi_P$ is of type $G_2$ then $\Phi = G_2$ (due to the lengths of the roots).\\
If $\Phi$ is not of type $G_2$ then the first condition leads to two different type of admissible planes, $\Phi_{P}^{+}$ is of one of the following type :\\
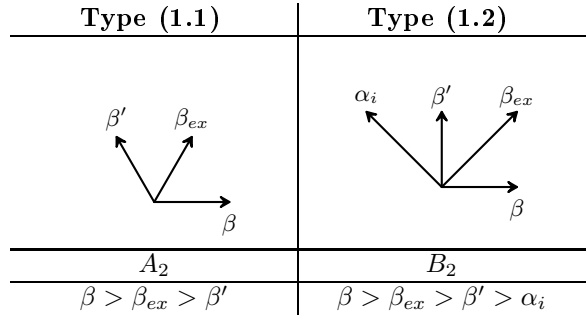
\begin{center}\begin{tabular}{c|c}
\textbf{Type (1.1)} \label{type 1.1} & \textbf{Type (1.2)} \label{type 1.2} \\
\hline
\begin{pgfpicture}{-0.5cm}{0cm}{2.9cm}{2.5cm}%
\pgfsetroundjoin \pgfsetroundcap%
\pgfsetlinewidth{0.8pt} 
\pgfxyline(1.2,0.5)(2.2,0.5)
\pgfsetfillcolor{black}
\pgfmoveto{\pgfxy(2.0845,0.5667)}\pgflineto{\pgfxy(2.2,0.5)}\pgflineto{\pgfxy(2.0845,0.4333)}\pgflineto{\pgfxy(2.1423,0.5)}\pgfclosepath\pgffillstroke
\pgfxyline(1.2,0.5)(1.7,1.366)
\pgfmoveto{\pgfxy(1.5845,1.2993)}\pgflineto{\pgfxy(1.7,1.366)}\pgflineto{\pgfxy(1.7,1.2327)}\pgflineto{\pgfxy(1.6711,1.316)}\pgfclosepath\pgffillstroke
\pgfxyline(1.2,0.5)(0.7,1.366)
\pgfmoveto{\pgfxy(0.7,1.2327)}\pgflineto{\pgfxy(0.7,1.366)}\pgflineto{\pgfxy(0.8155,1.2993)}\pgflineto{\pgfxy(0.7289,1.316)}\pgfclosepath\pgffillstroke
\pgfputat{\pgfxy(2.2,0.2)}{\pgfnode{rectangle}{center}{\small $\beta$}{}{\pgfusepath{}}}
\pgfputat{\pgfxy(1.7,1.6)}{\pgfnode{rectangle}{center}{\small $\beta_{ex}$}{}{\pgfusepath{}}}
\pgfputat{\pgfxy(0.7,1.6)}{\pgfnode{rectangle}{center}{\small $\beta'$}{}{\pgfusepath{}}}
\end{pgfpicture}
&
\begin{pgfpicture}{-0.5cm}{-0.5cm}{2.9cm}{2.2cm}%
\pgfsetroundjoin \pgfsetroundcap%
\pgfsetlinewidth{0.8pt} 
\pgfxyline(1.2,0.2)(2.2,0.2)
\pgfsetfillcolor{black}
\pgfmoveto{\pgfxy(2.0845,0.2667)}\pgflineto{\pgfxy(2.2,0.2)}\pgflineto{\pgfxy(2.0845,0.1333)}\pgflineto{\pgfxy(2.1423,0.2)}\pgfclosepath\pgffillstroke
\pgfxyline(1.2,0.2)(2.2,1.2)
\pgfmoveto{\pgfxy(2.0712,1.1655)}\pgflineto{\pgfxy(2.2,1.2)}\pgflineto{\pgfxy(2.1655,1.0712)}\pgflineto{\pgfxy(2.1592,1.1592)}\pgfclosepath\pgffillstroke
\pgfxyline(1.2,0.2)(1.2,1.2)
\pgfmoveto{\pgfxy(1.1333,1.0845)}\pgflineto{\pgfxy(1.2,1.2)}\pgflineto{\pgfxy(1.2667,1.0845)}\pgflineto{\pgfxy(1.2,1.1423)}\pgfclosepath\pgffillstroke
\pgfxyline(1.2,0.2)(0.2,1.2)
\pgfmoveto{\pgfxy(0.2345,1.0712)}\pgflineto{\pgfxy(0.2,1.2)}\pgflineto{\pgfxy(0.3288,1.1655)}\pgflineto{\pgfxy(0.2408,1.1592)}\pgfclosepath\pgffillstroke
\pgfputat{\pgfxy(2.2,-0.1)}{\pgfnode{rectangle}{center}{\small $\beta$}{}{\pgfusepath{}}}
\pgfputat{\pgfxy(2.2,1.4)}{\pgfnode{rectangle}{center}{\small $\beta_{ex}$}{}{\pgfusepath{}}}
\pgfputat{\pgfxy(1.2,1.4)}{\pgfnode{rectangle}{center}{\small $\beta'$}{}{\pgfusepath{}}}
\pgfputat{\pgfxy(0.2,1.4)}{\pgfnode{rectangle}{center}{\small $\alpha_i$}{}{\pgfusepath{}}}
\end{pgfpicture}
\\
\hline
$A_2$& $B_2$\\
\hline
  $\beta > \beta_{ex} > \beta'$ & $\beta > \beta_{ex} > \beta' > \alpha_i$
\end{tabular}
\end{center}
The second condition leads to four types of admissible planes, $\Phi_{P}^{+}$ is of one of the following types :\\
\begin{center}
\begin{tabular}{c|c|c|c}
\textbf{Type (2.1)} \label{type 2.1}&\textbf{Type (2.2)} \label{type 2.2}& \textbf{Type (2.3)} \label{type 2.3} & \textbf{Type (2.4)} \label{type 2.4}\\
\hline
\begin{pgfpicture}{-0.5cm}{0cm}{2.9cm}{2.5cm}%
\pgfsetroundjoin \pgfsetroundcap%
\pgfsetlinewidth{0.8pt} 
\pgfxyline(1.2,0.5)(2.2,0.5)
\pgfsetfillcolor{black}
\pgfmoveto{\pgfxy(2.0845,0.5667)}\pgflineto{\pgfxy(2.2,0.5)}\pgflineto{\pgfxy(2.0845,0.4333)}\pgflineto{\pgfxy(2.1423,0.5)}\pgfclosepath\pgffillstroke
\pgfxyline(1.2,0.5)(1.7,1.366)
\pgfmoveto{\pgfxy(1.5845,1.2993)}\pgflineto{\pgfxy(1.7,1.366)}\pgflineto{\pgfxy(1.7,1.2327)}\pgflineto{\pgfxy(1.6711,1.316)}\pgfclosepath\pgffillstroke
\pgfxyline(1.2,0.5)(0.7,1.366)
\pgfmoveto{\pgfxy(0.7,1.2327)}\pgflineto{\pgfxy(0.7,1.366)}\pgflineto{\pgfxy(0.8155,1.2993)}\pgflineto{\pgfxy(0.7289,1.316)}\pgfclosepath\pgffillstroke
\pgfputat{\pgfxy(2.2,0.2)}{\pgfnode{rectangle}{center}{\small $\beta_1$}{}{\pgfusepath{}}}
\pgfputat{\pgfxy(1.7,1.6)}{\pgfnode{rectangle}{center}{\small $\beta_{2}$}{}{\pgfusepath{}}}
\pgfputat{\pgfxy(0.7,1.6)}{\pgfnode{rectangle}{center}{\small $\alpha_i$}{}{\pgfusepath{}}}
\end{pgfpicture}  &
\begin{pgfpicture}{-0.5cm}{-0.5cm}{2.9cm}{2.2cm}%
\pgfsetroundjoin \pgfsetroundcap%
\pgfsetlinewidth{0.8pt} 
\pgfxyline(1.2,0.2)(2.2,0.2)
\pgfsetfillcolor{black}
\pgfmoveto{\pgfxy(2.0845,0.2667)}\pgflineto{\pgfxy(2.2,0.2)}\pgflineto{\pgfxy(2.0845,0.1333)}\pgflineto{\pgfxy(2.1423,0.2)}\pgfclosepath\pgffillstroke
\pgfxyline(1.2,0.2)(2.2,1.2)
\pgfmoveto{\pgfxy(2.0712,1.1655)}\pgflineto{\pgfxy(2.2,1.2)}\pgflineto{\pgfxy(2.1655,1.0712)}\pgflineto{\pgfxy(2.1592,1.1592)}\pgfclosepath\pgffillstroke
\pgfxyline(1.2,0.2)(1.2,1.2)
\pgfmoveto{\pgfxy(1.1333,1.0845)}\pgflineto{\pgfxy(1.2,1.2)}\pgflineto{\pgfxy(1.2667,1.0845)}\pgflineto{\pgfxy(1.2,1.1423)}\pgfclosepath\pgffillstroke
\pgfxyline(1.2,0.2)(0.2,1.2)
\pgfmoveto{\pgfxy(0.2345,1.0712)}\pgflineto{\pgfxy(0.2,1.2)}\pgflineto{\pgfxy(0.3288,1.1655)}\pgflineto{\pgfxy(0.2408,1.1592)}\pgfclosepath\pgffillstroke
\pgfputat{\pgfxy(2.2,-0.1)}{\pgfnode{rectangle}{center}{\small $\beta$}{}{\pgfusepath{}}}
\pgfputat{\pgfxy(2.2,1.4)}{\pgfnode{rectangle}{center}{\small $\beta_{ex}$}{}{\pgfusepath{}}}
\pgfputat{\pgfxy(1.2,1.4)}{\pgfnode{rectangle}{center}{\small $\beta'$}{}{\pgfusepath{}}}
\pgfputat{\pgfxy(0.2,1.4)}{\pgfnode{rectangle}{center}{\small $\alpha_i$}{}{\pgfusepath{}}}
\end{pgfpicture}
& \begin{pgfpicture}{-0.5cm}{-0.5cm}{2.9cm}{2.2cm}%
\pgfsetroundjoin \pgfsetroundcap%
\pgfsetlinewidth{0.8pt} 
\pgfxyline(1.2,0.2)(2.2,0.2)
\pgfsetfillcolor{black}
\pgfmoveto{\pgfxy(2.0845,0.2667)}\pgflineto{\pgfxy(2.2,0.2)}\pgflineto{\pgfxy(2.0845,0.1333)}\pgflineto{\pgfxy(2.1423,0.2)}\pgfclosepath\pgffillstroke
\pgfxyline(1.2,0.2)(2.2,1.2)
\pgfmoveto{\pgfxy(2.0712,1.1655)}\pgflineto{\pgfxy(2.2,1.2)}\pgflineto{\pgfxy(2.1655,1.0712)}\pgflineto{\pgfxy(2.1592,1.1592)}\pgfclosepath\pgffillstroke
\pgfxyline(1.2,0.2)(1.2,1.2)
\pgfmoveto{\pgfxy(1.1333,1.0845)}\pgflineto{\pgfxy(1.2,1.2)}\pgflineto{\pgfxy(1.2667,1.0845)}\pgflineto{\pgfxy(1.2,1.1423)}\pgfclosepath\pgffillstroke
\pgfxyline(1.2,0.2)(0.2,1.2)
\pgfmoveto{\pgfxy(0.2345,1.0712)}\pgflineto{\pgfxy(0.2,1.2)}\pgflineto{\pgfxy(0.3288,1.1655)}\pgflineto{\pgfxy(0.2408,1.1592)}\pgfclosepath\pgffillstroke
\pgfputat{\pgfxy(2.2,-0.1)}{\pgfnode{rectangle}{center}{\small $\alpha_i$}{}{\pgfusepath{}}}
\pgfputat{\pgfxy(2.2,1.4)}{\pgfnode{rectangle}{center}{\small $\beta_{3}$}{}{\pgfusepath{}}}
\pgfputat{\pgfxy(1.2,1.4)}{\pgfnode{rectangle}{center}{\small $\beta_2$}{}{\pgfusepath{}}}
\pgfputat{\pgfxy(0.2,1.4)}{\pgfnode{rectangle}{center}{\small $\beta_1$}{}{\pgfusepath{}}}
\end{pgfpicture}  & 
\begin{pgfpicture}{-0.7cm}{-0.7cm}{2cm}{2.5cm}%
\pgfsetroundjoin%
\pgfsetlinewidth{0.8pt} 
\pgfxyline(0,0)(1,0)
\pgfsetmiterjoin \pgfsetfillcolor{black}
\pgfmoveto{\pgfxy(0.8845,0.0667)}\pgflineto{\pgfxy(1,0)}\pgflineto{\pgfxy(0.8845,-0.0667)}\pgflineto{\pgfxy(0.9423,0)}\pgfclosepath\pgffillstroke
\pgfsetroundjoin%
\pgfxyline(0,0)(0,1.5)
\pgfsetmiterjoin \pgfmoveto{\pgfxy(-0.0667,1.3845)}\pgflineto{\pgfxy(0,1.5)}\pgflineto{\pgfxy(0.0667,1.3845)}\pgflineto{\pgfxy(0,1.4423)}\pgfclosepath\pgffillstroke
\pgfsetroundjoin%
\pgfsetlinewidth{0.2pt} 
\pgfputat{\pgfxy(1.2,0.3)}{\pgftext{\color{black}\small $\alpha_i$}}\pgfstroke
\pgfputat{\pgfxy(0.3,1.7)}{\pgftext{\color{black}\small $\beta$}}\pgfstroke
\end{pgfpicture}

 \\
\hline
 $A_2$ & $B_{2}$ with long $\alpha_i$ &$B_{2}$ with short $\alpha_{i}$ & $A_{1} \times A_{1}$  \\ 
 \hline
$\beta_1 > \beta_{2} > \alpha_i$ &$\beta > \beta_{ex} > \beta' > \alpha_i$ &$\beta_1 > \beta_{2} > \beta_3 > \alpha_i$ & $\beta > \alpha_i$ 
  \end{tabular}
  \end{center}
  We note that types (1.2) and (2.2) are the same.
\end{rmq}


\section{The quantized enveloping algebra $\mathcal{U}_q(\mathfrak{g})$}
Let $\mathbb{K}$ be a field of characteristic not equal to 2 and 3, and $q$ an element $\mathbb{K}^*$ which is not a root of unity. Firstly, we recall definitions about $\mathcal{U}_q(\mathfrak{g})$ and $\mathcal{U}_q^+(\mathfrak{g})$ using notations from \cite[chap4]{J}. We recall then the Poincare-Birkhoff-Witt bases of $\mathcal{U}_q(\mathfrak{g})$ construction using Lusztig automorphisms. There are several ways to construct the so called Lusztig automorphisms, we recall here three different methods. The Lusztig's one follows \cite[Section 3]{Lu2}, Jantzen's one, which is the same as De Concini, Kac et Procesi, is explained in \cite[Section 8.14]{J} and  \cite[Section 2.1]{MR1351503} and a third one is necessary to established a link between the two others. We will explain each method and then see the links between the obtained bases.
\subsection{Recalls on $\mathcal{U}_q(\cat{g})$}
\begin{nota}
For all $a$ and $n$ integers such that $ a \geq n \geq 0$, we set
\[ [n]_{q} = \frac{q^{n}-q^{-n}}{q-q^{-1}}, \ [n]^{!}_{q} = [n]_{q}[n-1]_{q}...[2]_{q}[1]_{q} , \ \left[ \begin{array}{c} a \\ n \end{array} \right]_{q} = \frac{[a]^{!}_{q}}{[n]^{!}_{q}[a-n]^{!}_{q}}.\]
Moreover for all $\alpha \in \Pi$, we set  $q_{\alpha} = q^{\frac{(\alpha,\alpha)}{2}}$.
\end{nota}
\begin{dfn}\label{Uqg}
\begin{itemizedot}
\item The quantized enveloping algebra $\mathcal{U}_q(\cat{g})$ is the $\mathbb{K}$-algebra with generators $E_{\alpha}$, $F_{\alpha}$, $K_{\alpha}$ and $K_{\alpha}^{-1}$ (for all $\alpha$ in $\Pi$) and relations (for all $\alpha , \beta \in \Pi$)
\begin{equation}
K_{\alpha}K_{\alpha}^{-1} = 1 = K_{\alpha}^{-1}K_{\alpha} \hspace{2cm}  K_{\alpha}K_{\beta}=K_{\beta}K_{\alpha}
\end{equation}  
\begin{equation}
K_{\alpha}E_{\beta}K_{\alpha}^{-1} = q^{(\alpha , \beta)}E_{\beta}
\end{equation} 
\begin{equation}
K_{\alpha}F_{\beta}K_{\alpha}^{-1} = q^{-(\alpha , \beta)}F_{\beta}
\end{equation} 
\begin{equation}
E_{\alpha}F_{\beta} - F_{\beta}E_{\alpha} = \delta_{\alpha \beta} \frac{K_{\alpha} - K_{\alpha}^{-1}}{q_{\alpha}-q_{\alpha}^{-1}}
\end{equation}  
where $\delta_{\alpha \beta}$ is the Kronecker symbol, and (for $\alpha \neq \beta$)
\begin{equation}
\sum_{s=0}^{1-a_{\alpha \beta}} (-1)^s \left[ \begin{array}{c} 1-a_{\alpha \beta} \\ i \end{array} \right]_{q_{\alpha}} E_{\alpha}^{1-a_{\alpha \beta}-s}E_{\beta}E_{\alpha}^s = 0
\end{equation} 
\begin{equation}
\sum_{s=0}^{1-a_{\alpha \beta}} (-1)^s \left[ \begin{array}{c} 1-a_{\alpha \beta} \\ i \end{array} \right]_{q_{\alpha}} F_{\alpha}^{1-a_{\alpha \beta}-s}F_{\beta}F_{\alpha}^s = 0
\end{equation} 
where $a_{\alpha \beta} = 2 (\alpha , \beta) / (\alpha , \alpha)$ for all $\alpha , \beta \in \Pi$.

\item $\mathcal{U}_q^+(\mathfrak{g})$ (resp. $\mathcal{U}_q^-(\mathfrak{g})$)  is the subalgebra of $\mathcal{U}_q(\mathfrak{g})$ generated by all $E_{\alpha}$ (resp. $F_{\alpha}$) with $\alpha \in \Pi$.
\end{itemizedot}
\end{dfn}
Let us recall two important results proven for example in \cite[Section I.6]{MR1898492}.
\begin{theo}
\begin{enumerate}
\item $\mathcal{U}_q^+(\mathfrak{g})$ is Noetherian domain.
\item $\mathcal{U}_q^+(\mathfrak{g})$ is graded by $\mathbb{Z}\Phi$ with $wt(E_\alpha)= \alpha, wt(F_\alpha)=- \alpha $ and $wt(K_\alpha^{\pm 1})= 0)$.
\end{enumerate}
\end{theo}
\subsection{Lusztig's construction}
\begin{dfn}[Lusztig's automorphisms]
For all $i \in \llbracket 1,n \rrbracket$ there is a unique automorphism $T_{\alpha_{i}}$ of the algebra $\mathcal{U}_q(\mathfrak{g})$ such that :
\[T_{\alpha_i}E_{\alpha_i} = -F_{\alpha_i}K_{\alpha_i}, \   T_{\alpha_i}F_{\alpha_i} = -K_{\alpha_i}^{-1} E_{\alpha_i} \ T_{\alpha_i}K_{\alpha_j} = K_{\alpha_j}K_{\alpha_i}^{-a_{ij}} \ (j \in \llbracket 1,n \rrbracket) \]
and for $j \neq i$ :
\[  T_{\alpha_i}E_{\alpha_j} = \sum_{r + s = - a_{ij}}(-1)^r q^{-d_i s} E_{\alpha_i}^{(r)}E_{\alpha_j}E_{\alpha_i}^{(s)} \]
\[ T_{\alpha_i}F_{\alpha_j} = \sum_{r + s = - a_{ij}}(-1)^r q^{d_i s} F_{\alpha_i}^{(s)}F_{\alpha_j}F_{\alpha_i}^{(r)}  \]
where $\displaystyle E_{\alpha_i}^{(n)} := \frac{E_{\alpha_i}^{n}}{[n]_{d_{i}}^{!}}$ and $d_i = \frac{(\alpha_i, \alpha_i)}{2}$. 
\end{dfn}
We now a fix a Lusztig order so that we can use the notations of columns and boxes as in the Section 2. The following result is given by Lusztig in \cite[Section 4.3]{Lu2} :
\begin{prop}
There is a unique map  from $\Phi^{+}$ to $\llbracket 1,n \rrbracket$, sending $\beta$ to $i_{\beta})$ such that the following properties are satisfied :
\begin{enumerate}
\item $s_{i_{\beta_{1}}}$ and $s_{i_{\beta_{2}}}$ commute in W whenever $\beta_{1}$ and $\beta_{2}$ are in the same box. Hence, for a box B, the product of all $s_{i_{\beta}}$ with $\beta \in B$ is a well defined element $s(B)$ in W, independent of the order of the factors.
\item $i_{\alpha_{j}} = j$ for $j \in \llbracket 1,n \rrbracket$.
\item If $\beta \in C_{j}$ and if $B_{1}, ..., B_{k}$ are the boxes in $C_{j}$ whose elements are strictly greater than $\beta$ for the Lusztig order then $s(B_{1})s(B_{2})...s(B_{k})(\alpha_{i_{\beta}}) = \beta$.
\end{enumerate}
We then set $w_{\beta} :=  s(B_{1})s(B_{2})...s(B_{k})$.
\end{prop}
We now recall the construction of a PBW basis of $\mathcal{U}_q^+(\mathfrak{g})$ due to Lusztig \cite[Theorem 3.2]{Lu2}.
\begin{theo}
Let $w \in W$ and $s_{i_{1}}...s_{i_{p}}$ be a reduced decomposition of $w$. Then the automorphism $T_{w} := T_{\alpha_{i_{1}}}...T_{\alpha_{i_{p}}}$ depends only on $w$ and not on the choice of reduced expression for it. Hence the $T_{\alpha_{i}}$ define a homomorphism 
of the braid group of $W$ in the group of automorphisms of the algebra $\mathcal{U}_q(\mathfrak{g})$.
\end{theo}
\begin{prop}
For all positive roots $\beta$, we define $E_{\beta} := T_{w_{\beta}}(E_{i_{\beta}}) \in \Phi^{+}$. These elements form a Poincaré-Birkhoff-Witt basis of $\mathcal{U}_q^+(\mathfrak{g})$ (\cite[proposition 4.2]{Lu2})
\end{prop}
To obtain commutation relation in this context, we can use admissible planes.
\begin{nota}
If $\beta > \beta'$, then we set $[E_{\beta},E_{\beta'}]_{q} = E_{\beta}E_{\beta'} - q^{(\beta,\beta')}E_{\beta'}E_{\beta}$.
\end{nota}
Our aim in the remaining of this paragraph is to exhibit the form of the commutation relation between two generators $E_\gamma$ and $E_{\gamma'}$, when $\gamma$ and $\gamma'$ belong to the same admissible plane P.\\
We first consider the case where $\Phi_P ( = \Phi \cap P)$ is of type $G_2$. In this case, $\Phi$ is also of type $G_2$ and the commutation relations have been computed in \cite[Section 5.2]{Lu2}. This leads us to the following result.
\begin{prop}\label{propG2}
Assume that $\Phi$ is of type $G_2$. Denote by $\alpha_1$ the short simple root and by $\alpha_2$ the long simple root. This is a good numbering of the set of simple roots (see Example \ref{exG2}). The corresponding reduced decomposition of $w_0$ is $s_1s_2s_1s_2s_1s_2$ $(s_i = s_{\alpha_i})$ and, describing the roots in the associated convex order, one has  
\[\Phi^+ = \{\beta_1 = \alpha_1 ,\beta_2 = 3\alpha_1 + \alpha_2, \beta_3 = 2\alpha_1 + \alpha_2, \beta_4 = 3\alpha_1 + 2\alpha_2, \beta_5 = \alpha_1 + \alpha_2, \beta_6 = \alpha_2\}.\]   
The first column $C_1$ is reduced to $\{\beta_1\}$, the second column $C_2 = \{\beta_2, \beta_3, \beta_4, \beta_5, \beta_6\}$ is exceptional with $\beta_{ex} = \beta_4$. One has : 
\begin{itemize}
\item $[E_{\beta_3},E_{\beta_1}]_{q} = \lambda E_{\beta_{2}}$ with $\lambda \neq 0$,
\item $[E_{\beta_4},E_{\beta_1}]_{q} = \lambda E_{\beta_{3}}^2$ with $\lambda \neq 0$,
\item $[E_{\beta_5},E_{\beta_1}]_{q} = \lambda E_{\beta_{3}}$ with $\lambda \neq 0$,
\item $[E_{\beta_6},E_{\beta_1}]_{q} = \lambda E_{\beta_{5}}$ with $\lambda \neq 0$,
\item $[E_{\beta_3},E_{\beta_1}]_{q} = \lambda E_{\beta_{2}}$ with $\lambda \neq 0$,
\item $[E_{\beta_5},E_{\beta_3}]_{q} = \lambda E_{\beta_{4}}$ with $\lambda \neq 0$.
\end{itemize}
\end{prop}
If $\Phi$ is not of type $G_2$, the commutation relations between the Lusztig's generators corresponding to two roots which are in the same admissible plane are known in several cases (\cite[Section 5.2]{Lu2}). In particular, we have the following relations.
\begin{prop}[$\Phi$ not of type $G_2$]\label{rellusztig}
\begin{itemizedot}
\item If $P = <\beta, \beta'>$ is an admissible plane of \hyperref[systracine]{type (1.1)}, then  $\Phi_{P}^{+} = \{\beta, \beta_{ex} = \beta + \beta', \beta'\}$ and:
\begin{itemize}
\item[$\circ$] $[E_{\beta},E_{\beta'}]_{q} = \lambda E_{\beta_{ex}}$ with $\lambda \neq 0$,
\item[$\circ$] $[E_{\beta},E_{\beta_{ex}}]_{q} = [E_{\beta_{ex}},E_{\beta'}]_{q} = 0$.
\end{itemize}

\item If $P = <\beta, \beta'>$ is an admissible plane of \hyperref[systracine]{type (1.2)}, then  $\Phi_{P}^{+} = \{\beta, \beta_{ex} = \beta + \beta', \beta', \alpha_{i}\}$ and the relations are:
\begin{itemize}
\item[$\circ$] $[E_{\beta},E_{\beta'}]_{q} = \lambda E_{\beta_{ex}}$ with $\lambda \neq 0$,
\item[$\circ$] $[E_{\beta_{ex}},E_{\alpha_i}]_{q} = \lambda' E_{\beta'}^{2}$ with $\lambda' \neq 0$,
\item[$\circ$] $[E_{\beta},E_{\alpha_i}]_{q} = \lambda'' E_{\beta'}$ with $\lambda'' \neq 0$,
\item[$\circ$] $[E_{\beta},E_{\beta_{ex}}]_{q} = [E_{\beta_{ex}},E_{\beta'}]_{q} = [E_{\beta'},E_{\alpha_{i}}]_{q} = 0$.
\end{itemize}
\item If $P = <\beta, \alpha_{i}>$ is an admissible plane of \hyperref[systracine]{type (2.1)}, then  $\Phi_{P}^{+} = \{\beta_{1}, \beta_{2} = \beta_{1} + \alpha_{i}, \alpha_{i}\}$ ($\beta = \beta_{1}$ or $\beta_{2}$) and:
\begin{itemize}
\item[$\circ$] $[E_{\beta_{1}},E_{\alpha_{i}}]_{q} = \lambda E_{\beta_{2}}$ with $\lambda \neq 0$.,
\item[$\circ$] $[E_{\beta_{1}},E_{\beta_{2}}]_{q} = [E_{\beta_{2}},E_{\alpha_{i}}]_{q} =  0. $
\end{itemize} 
\item If $P = <\beta,  \alpha_{i}>$ is an admissible plane of \hyperref[systracine]{type (2.2)}, then we have the same relations as in type (1.2).\\
\item If $P =$ \textrm{Vect}$(\beta, \alpha_{i})$ is an admissible plane of \hyperref[systracine]{type (2.3)}, then  $\Phi_{P}^{+} = \{\beta_{1}, \beta_{2} = \beta_{1} + \alpha_{i}, \beta_{3} = \beta_{1} + 2\alpha_{i}, \alpha_{i}\}$ ($\beta = \beta_{1}$,$\beta_{2}$ or $\beta_{3}$) and the relations are :
\begin{itemize}
\item[$\circ$] $[E_{\beta_{2}},E_{\alpha_{i}}]_{q} = \lambda E_{\beta_{3}}$ with $\lambda \neq 0$,
\item[$\circ$] $[E_{\beta_{1}},E_{\beta_{3}}]_{q} = \lambda' E_{\beta_{2}}^{2}$ with $\lambda' \neq 0$,
\item[$\circ$] $[E_{\beta_{1}},E_{\alpha_i}]_{q} = \lambda'' E_{\beta_{2}}$ with $\lambda'' \neq 0$,
\item[$\circ$] $[E_{\beta_{1}},E_{\beta_{2}}]_{q} = [E_{\beta_{2}},E_{\beta_{3}}]_{q} = [E_{\beta_{3}},E_{\alpha_{i}}]_{q} = 0$
\end{itemize}
\item If $P = <\beta, \alpha_{i}>$ is an admissible plane of \hyperref[systracine]{type (2.4)}, then $\Phi_{P}^{+} = \{\beta, \alpha_{i}\}$ with $\beta \ \bot \ \alpha_i$ and, if $\beta$ is ordinary, then $[E_{\beta},E_{\alpha_{i}}]_{q} = 0$.
\end{itemizedot}
\end{prop}
\begin{cor} \label{corplanadm}
Assume $\Phi$ is not of type $G_2$. Let $i, l$ be two integers such that $1 \leq i < l \leq n $ and  $\eta \in C_l$
\begin{enumerate}
\item If $(\eta, \alpha_i) > 0$, then $[E_{\eta},E_{\alpha_{i}}]_{q} = 0$.
\item If $\eta + \alpha_i = m \gamma$ with $\gamma \in \Phi^+$ and $m \in \mathbb{N}^{\star}$, then $[E_{\eta},E_{\alpha_{i}}]_{q} = \lambda E_{\gamma}^{m}, \ \text{with} \  \lambda \in \mathbb{K}^{\star}$.
\item If $\eta = \eta_1 + \eta_2$ with $\eta_1$ and $\eta_2$ in $C_l$ such that $h(\eta_1) + 1 = h(\eta_2)$ then  $[E_{\eta_1},E_{\eta_{2}}]_{q} = \lambda E_{\eta}, \ \text{with} \  \lambda \in \mathbb{K}^{\star}$.
\end{enumerate}
\end{cor}
\begin{preuve}
$P = \text{Vect}(\eta, \alpha_{i})$ is an admissible plane of type (2.1), (2.2), (2.3) or (2.4) by definition. 
\begin{enumerate}
\item $P$ is not of type (2.4) because $(\eta, \alpha_i) \neq 0$. We distinguish between three remaining cases.
\begin{itemizedot}
\item If $P$ is of type (2.1), then with the notations of Remark \ref{systracine}, we have $\eta = \beta_{2}$, and so the result follows from Proposition \ref{rellusztig}.
\item If $P$ is of type (2.2), then we have $\eta = \beta'$, and so the result follows from Proposition \ref{rellusztig}.
\item If $P$ is of type (2.3), then we have $\eta = \beta_{3}$, and so the result follows from proposition \ref{rellusztig}.
\end{itemizedot}

\item Since $m \neq 0$, we have $\gamma \in P \cap \Phi^{+} = \Phi_{p}^{+}$, so $P$ is not of type (2.4). We distinguish between three remaining cases.
\begin{itemizedot}
\item If $P$ is of type (2.1), then we deduce that $m=1$, $\eta = \beta_{1}$ and $\gamma = \beta_{2}$, and so the result follows from Proposition \ref{rellusztig}.
\item If $P$ is of type (2.2), then there are two possibilities:
\begin{itemize}
\item[$\circ$] $m=1$, $\eta = \beta$ and $\gamma = \beta'$. 
\item[$\circ$] $m=2$, $\eta = \beta_{ex}$ and $\gamma = \beta'$. 
\end{itemize}
In both cases, the result follows from Proposition \ref{rellusztig}.

\item If $P$ is of type (2.3), then we have $m=1$, $\eta = \beta_{1}$ (resp.  $\eta = \beta_{2}$) and $\gamma = \beta_{2}$  (resp. $\gamma = \beta_{3}$), and so the result follows from Proposition \ref{rellusztig}.

\end{itemizedot}
\item Let us consider the plane $P := <\eta_1, \eta_2>$. It is an admissible plane (see definition \ref{planadm}) and $\Phi_P^+ = \{\eta_1, \eta, \eta_2\}$ (if P is of type 1.1) or $\{\eta_1, \eta, \eta_2, \alpha_i\}$ with $i < l$ (if P is of type 1.2). The previous proposition implies that $[E_{\eta_1},E_{\eta_{2}}]_{q} = \lambda E_{\eta}, \ \text{with} \  \lambda \in \mathbb{K}^{\star}$ in both cases.
\end{enumerate}
\end{preuve}
\subsection{Jantzen's construction}
In \cite[section 8.14]{J}, a different construction of a PBW basis is explained which also uses the automorphisms $T_{\alpha}, (\alpha \in \Pi)$.\\
For a given reduced decomposition of $w_{0} = s_{i_{1}}...s_{i_{N}}$, we know that, for all $\beta \in \Phi^+$, there exists $i_{\beta} \in \llbracket 1,N \rrbracket$ such that $\beta = s_{i_{1}}...s_{i_{\beta}-1}(\alpha_{i_{\beta}})$.
\begin{dfn}
Let $\beta \in \Phi^+$, we set $w'_{\beta} := s_{i_{1}}...s_{i_{\beta}-1}$ and define $X_{\beta} := T_{w'_{\beta}}(E_{\alpha_{i_{\beta}}})$, $Y_{\beta} := T_{w'_{\beta}}(F_{\alpha_{i_{\beta}}})$.
\end{dfn}
The following result follows from \cite[theoreme 4.21 and 8.24]{J}.
\begin{theo}\label{PBW}
\begin{itemize}
\item If $\alpha \in \Pi$, then $X_{\alpha} = E_{\alpha}$ (\cite[Proposition 8.20]{J}).
\item The products $X_{\beta_{1}}^{k_{1}}...X_{\beta_{N}}^{k_{N}}$ ($k_i \in \mathbb{N}$)  form a basis of $\mathcal{U}_q^+(\mathfrak{g})$.
\item The products $X_{\beta_{1}}^{k_{1}}...X_{\beta_{N}}^{k_{N}}K_{\alpha_{1}}^{m_{1}}...K_{\alpha_{n}}^{m_{n}}Y_{\beta_{1}}^{l_{1}}...Y_{\beta_{N}}^{l_{N}}$ (resp. $K_{\alpha_{1}}^{m_{1}}...K_{\alpha_{n}}^{m_{n}}Y_{\beta_{1}}^{l_{1}}...Y_{\beta_{N}}^{l_{N}}X_{\beta_{1}}^{k_{1}}...X_{\beta_{N}}^{k_{N}}$,\\
resp. $Y_{\beta_{1}}^{l_{1}}...Y_{\beta_{N}}^{l_{N}}K_{\alpha_{1}}^{m_{1}}...K_{\alpha_{n}}^{m_{n}} X_{\beta_{1}}^{k_{1}}...X_{\beta_{N}}^{k_{N}}$), ($k_{i}, l_{i} \in \mathbb{N}, m_{i} \in \mathbb{Z})$ form a basis of $\mathcal{U}_q(\mathfrak{g})$.
\end{itemize}

\end{theo}
The following theorem was proved by Levendorski{\v\i} et Soibelman \cite[Proposition 5.5.2]{MR1116413} in a slightly different case. One can find other formulations in the literature (several containing small mistakes). That is why we give a proof of this result in \cite[Section 3.3]{moi}. We make this proof essentially by rewriting the one from \cite[Proposition 5.5.2]{MR1116413}.
\begin{theo}[of  Levendorski{\v\i} et Soibelman] \label{LS}
If $i$ and $j$ are two integers such that $1 \leq i < j \leq N$, then we have
\begin{equation*}
	X_{\beta_i}X_{\beta_j} - q^{(\beta_i , \beta_j)}X_{\beta_j}X_{\beta_i} = \sum_{\substack{ \beta_{i} < \gamma_{1} < ... < \gamma_{p} < \beta_{j}, \\ p \geq 1, \ k_i \in \mathbb{N}}}  c_{\mathbf{\overline{k},\overline{\gamma}}} X_{\gamma_{1}}^{k_{1}}...X_{\gamma_{p}}^{k_{p}},
\end{equation*}
 where  $c_{\mathbf{\overline{k},\overline{\gamma}}} \in \mathbb{K}$ and $c_{\mathbf{\overline{k},\overline{\gamma}}} \neq 0$ only if $wt(X_{\gamma_{1}}^{k_{1}}...X_{\gamma_{p}}^{k_{p}}) := k_{1} \times \gamma_{1}+...+k_{p} \times \gamma_{p} = \beta_i+\beta_j$.
\end{theo}

\subsection{Commutation relations between $X_{\gamma}$ in admissible planes}
The goal of this section is to show that the $X_{\gamma}$ satisfy analoguous relations to the $E_{\gamma}$ (see Section 3.2.). In order to achieve this aim, we start by introducing an intermediate generating system.
\subsubsection{Construction of a third generating system}
Let us recall the following well-known result :
\begin{lemme}[\cite{J}, Section 4.6] \label{tau}
\begin{enumerate}
\item There is a unique automorphism  $\omega$ of $\mathcal{U}_q(\mathfrak{g})$ such that $\omega(E_{\alpha}) = F_{\alpha}$, $\omega(F_{\alpha}) = E_{\alpha}$ and $\omega(K_{\alpha}) = K_{\alpha}^{-1}$. One has $\omega^{2} = 1$.
\item There is a unique anti-automorphism $\tau$ of $\mathcal{U}_q(\mathfrak{g})$ such that $\tau(E_{\alpha}) = E_{\alpha}$, $\tau(F_{\alpha}) = F_{\alpha}$ and $\tau(K_{\alpha}) = K_{\alpha}^{-1}$. One has $\tau^{2} = 1$.
\end{enumerate}
\end{lemme}
\begin{convsanss}
\begin{itemizedot}
\item Let $i$ be an integer of $\llbracket 1, n \rrbracket$. And set $T'_{\alpha_{i}} := \tau \circ T_{\alpha_{i}} \circ \tau$. This is an automorphism of $\mathcal{U}_q(\mathfrak{g})$ which satisfies the following conditions :
\[T'_{\alpha_i}E_{\alpha_i} = -K^{-1}_{\alpha_i}F_{\alpha_i}, \   T'_{\alpha_i}F_{\alpha_i} = -E_{\alpha_i} K_{\alpha_i}\ T'_{\alpha_i}K_{\alpha_j} = K_{\alpha_j}K_{\alpha_i}^{-a_{ij}} \ (j \in \llbracket 1,n \rrbracket) \]
and for $j \neq i$ :
\[  T'_{\alpha_i}E_{\alpha_j} = \sum_{r + s = - a_{ij}}(-1)^r q^{d_i s} E_{\alpha_i}^{(s)}E_{\alpha_j}E_{\alpha_i}^{(r)} \]
\[ T'_{\alpha_i}F_{\alpha_j} = \sum_{r + s = - a_{ij}}(-1)^r q^{-d_i s} F_{\alpha_i}^{(r)}F_{\alpha_j}F_{\alpha_i}^{(s)}  \]
\item If $w_{p}  \in W$ has a reduced decomposition given by $w_{p} = s_{i_{1}}...s_{i_{p}}$, then we set $T'_{w_{p}} :=  \tau \circ T_{w_{p}} \circ \tau$. We have $T'_{w_{p}} = T'_{\alpha_{i_{1}}}...T'_{\alpha_{i_{p}}} $.
\item If $\beta \in \Phi^+$, then we set $w_{\beta} := s_{i_{1}}...s_{i_{\beta}-1}$ and we define $X'_{\beta} := T'_{w_{\beta}}(E_{\alpha_{i_{\beta}}})$ and  $Y'_{\beta} := T_{w_{\beta}}(F_{\alpha_{i_{\beta}}})$. One has $X'_{\alpha} = E_{\alpha}$ and $Y'_{\alpha} = F_{\alpha}$ for $\alpha \in \Pi$.
\end{itemizedot}
\end{convsanss}
The Theorem of Levendorski{\v\i} and Soibelman can be rewritten as below. The proof can be found in \cite[Section 3.4]{moi} :
\begin{prop} \label{LS2}
If $i$ and $j$ are two integers such that $1 \leq i < j \leq N$ then we have
\begin{equation*}
	X'_{\beta_i}X'_{\beta_j} - q^{-(\beta_i , \beta_j)}X'_{\beta_j}X'_{\beta_i} = \sum_{\substack{ \beta_{i} < \gamma_{1} < ... < \gamma_{p} < \beta_{j}, \\ p \geq 1, \ k_i \in \mathbb{N}}}  c_{\mathbf{\overline{k},\overline{\gamma}}} X_{\gamma_{1}}^{\prime k_{1}}...X_{\gamma_{p}}^{\prime k_{p}}
\end{equation*}
 with  $c_{\mathbf{\overline{k},\overline{\gamma}}} \in \mathbb{K}$ and $c_{\mathbf{\overline{k},\overline{\gamma}}} \neq 0$ only if $wt(X_{\gamma_{1}}^{\prime k_{1}}...X_{\prime \gamma_{p}}^{k_{p}}) := k_{1} \times \gamma_{1}+...+k_{p} \times \gamma_{p} = \beta_i+\beta_j$
\end{prop}

\subsubsection{Relations between $E_{\beta}$ and $X'_{\beta}$.}
As in previous sections, $\Phi^+$ is provided with a given Lusztig order associated to a reduced decomposition of $w_{0} = s_{i_{1}}...s_{i_{N}}$. In this case, we can improve the Theorem of Levendorski{\v\i} and Soibelman.
\begin{theo} \label{1.4}
If $i$ and $j$ are two integers such that $1 \leq i < j \leq N$, then one has :
\begin{equation*}
	X'_{\beta_i}X'_{\beta_j} - q^{-(\beta_i , \beta_j)}X'_{\beta_j}X'_{\beta_i} = \sum_{\beta_{i} < \gamma_{1} < ... < \gamma_{p} < \beta_{j}}  C_{\mathbf{\overline{k},\overline{\gamma}}} X_{\gamma_{1}}^{\prime k_{1}}...X_{\gamma_{p}}^{\prime k_{p}}.
\end{equation*}
The monomials on the left-hand side whose coefficient $C_{\mathbf{\overline{k},\overline{\gamma}}}$ is not equal to zero satisfies:
\begin{itemizedot}
\item $wt(X_{\gamma_{1}}^{\prime k_{1}}...X_{\gamma_{p}}^{\prime k_{p}}) = \beta_i+\beta_j$;
\item $\gamma_{1}$ is not in the same box as $\beta_{i}$;
\item $\gamma_{p}$ is not in the same box as $\beta_{j}$.
\end{itemizedot}
\end{theo}
The proof of this theorem is essentially based on the following result:
\begin{lemme}\label{1.2}
Let $B=\{\beta_p, ..., \beta_{p+l}\}$ be a box and $\alpha_{i_{p}}, ..., \alpha_{i_{p+l}}$ be the corresponding simple roots. Then $\forall k \in \llbracket 0,l \rrbracket$, we have: 
\[ T'_{\alpha_{i_{p}}} ... T'_{\alpha_{i_{p+k-1}}}(E_{\alpha_{i_{p+k}}}) =  E_{\alpha_{i_{p+k}}} = T'_{\alpha_{i_{p}}} ...T'_{i_{\alpha_{p+k-1}}}T'_{i_{\alpha_{p+k+1}}}...T'_{\alpha_{i_{p+l}}}(E_{\alpha_{i_{p+k}}}). \]
\end{lemme}
\begin{preuve}
We already know that if $\alpha_{1}$ and $\alpha_{2}$ are two simple orthogonal roots, then $T_{\alpha_{1}}(E_{\alpha_{2}}) = E_{\alpha_{2}} = \tau(E_{\alpha_{2}})$, hence $T'_{\alpha_{1}}(E_{\alpha_{2}}) = E_{\alpha_{2}}$. As $\alpha_{i_{p}}, ..., \alpha_{i_{p+l}}$ are orthogonal to each others by Proposition \ref{boîte}, the formulas above are proved.
\end{preuve}
\begin{preuve3}
The first point is provided by Proposition \ref{LS2}. If in the reduced decomposition of $w_{0}$, we change the order of the reflexions associated to the simple roots coming from a single box B, we find a new reduced decomposition of $w_{0}$. The positive roots of B constructed with this new decomposition of $w_0$ are permuted as the simple roots are but the other roots are not moved. By the previous lemma, the $X'_{\beta}$, $\beta \in B$, are also permuted in the same way but are not modified, and the $X'_{\gamma}$, $\gamma \notin B$, are not modified. Thus, without lost of generality, we can assume that $\beta_i$ is maximal in its box and that $\beta_{j}$ is minimal in its box. As a result, if $\beta_i < \gamma_{1} < ... < \gamma_{p} < \beta_{j}$, then $\gamma_{1}$ is not in the same box as $\beta_{i}$ and  $\gamma_{p}$ is not in the same box as $\beta_{j}$.
\end{preuve3}
\begin{rmq}
The proof of the previous theorem can be rewritten with the elements $X_\beta (\beta \in \Phi_+)$ so that we also apply Theorem \ref{1.4} to those elements.
\end{rmq}
We can now establish a link between the $X'_\beta$'s and the $E_\beta$'s. 
\begin{theo}\label{theo3.6}
\[\forall \beta \in \Phi^+, \exists \lambda_\beta \in \mathbb{K}\setminus\{0\} \ \ \text{such that} \ \ X'_\beta = \lambda_\beta E_\beta\]
\end{theo}
\begin{preuve}
Let $\beta$ and $\beta'$ be two positive roots such that $\beta > \beta'$. Set $[X'_{\beta}, X'_{\beta'}]_q = X'_{\beta}X'_{\beta'} - q^{(\beta,\beta')} X'_{\beta'}X'_{\beta}$.\\
Let us deal first with the case where $\Phi$ is of type $G_2$. We keep the conventions of Proposition \ref{propG2}. It is known (Conventions 3.4.1.) that, since $\beta_1$ and $\beta_6$ are simple, one has $X'_{\beta_1} = E_{\beta_1}$ and $X'_{\beta_6} = E_{\beta_6}$.\\
Thus
\[[X'_{\beta_6},X'_{\beta_1}]_q = [E_{\beta_6},E_{\beta_1}]_q = \lambda E_{\beta_5} \ \text{avec} \ \lambda \in \mathbb{K}^\star. \]
By Theorem \ref{1.4},  one also has  $[X'_{\beta_6},X'_{\beta_1}]_q= \mu X'_{\beta_5}$ with $\mu \in \mathbb{K}$ and, then, $X'_{\beta_5} = \lambda_{\beta_5} E_{\beta_5}$ with $\lambda_{\beta_5} \in \mathbb{K}^\star$.\\
It implies that $[X'_{\beta_5},X'_{\beta_1}]_q = \lambda_{\beta_5} [E_{\beta_5},E_{\beta_1}]_q = \nu E_{\beta_3}$ with $\nu \in \mathbb{K}^\star$. We deduce as above that $X'_{\beta_3} = \lambda_{\beta_3} E_{\beta_3}$ with $\lambda_{\beta_3} \in \mathbb{K}^\star$.\\
Using the same method and considering $[X'_{\beta_3},X'_{\beta_1}]_q = \lambda_{\beta_3} [E_{\beta_3},E_{\beta_1}]_q$, one proves that $X'_{\beta_2} = \lambda_{\beta_2} E_{\beta_2}$ with $\lambda_{\beta_2} \in \mathbb{K}^\star$.\\
At last, one has $[X'_{\beta_5},X'_{\beta_3}]_q = \lambda_{\beta_5}\lambda_{\beta_3} [E_{\beta_5},E_{\beta_3}]_q = \nu E_{\beta_4}$ with $\nu \in \mathbb{K}^\star$, so it implies that $X'_{\beta_4} = \lambda_{\beta_4} E_{\beta_4}$ with $\lambda_{\beta_4} \in \mathbb{K}^\star$.\\
\par
Suppose now that $\Phi$ is of type $G_2$, and consider a column $C_{t}$ $(t \in \llbracket 1, n \rrbracket)$. We just prove the theorem for all the roots of $C_{t}$.\\
\par
We first study the case of ordinary roots.\\
Let $\beta \in C_{t}$ be an ordinary root. Let us prove the result by induction on $h(\beta)$.\\
If $h(\beta) = 1$, then $\beta = \alpha_{t}$ and as above $X'_{\alpha_{t}} = E_{\alpha_{t}}$.\\
Assume $h(\beta) > 1$ and the result proved for all $\delta \in C_{t}$ an ordinary root such that $h(\delta) < h(\beta)$.
By proposition \ref{prop2.10}, there is a simple root $\alpha_i$ ($i<t$) such that $\beta - \alpha_i = \gamma \in C_t$. Moreover, $\gamma$ is ordinary because, if it is not the case $\beta = \gamma + \alpha_{i}$ would be exceptional which contradicts the uniqueness of an exceptional root in a column. So $P := < \alpha_i, \beta>$  is an admissible plane of \hyperref[type 2]{type (2.1), (2.2) or (2.3)} and then $[E_\gamma,E_{\alpha_i}]_q = c E_\beta$ $(c \in \mathbb{K} \setminus\{ 0\})$ (see Section 3.2).\\
As $h(\gamma) = h(\beta)-1 < h(\beta)$,  one has $X'_\gamma = \lambda_{\gamma} E_{\gamma}$  $(\lambda_{\gamma } \in \mathbb{K}\setminus\{0\}$), and as $E_{\alpha_i} = X'_{\alpha_i}$, one has :
\[[X'_\gamma,X'_{\alpha_i}]_q = \lambda_{\gamma} [E_\gamma,E_{\alpha_i}]_q = \lambda_\gamma c E_\beta . \]
By Theorem \ref{1.4},  $E_\beta$ is a linear combination of monomials $X'_{\delta_1}...X'_{\delta_s}$ with \\
$\alpha_i < \delta_{1} \leq ... \leq \delta_{s} < \gamma$, $\delta_{s}$ not in the same box as $\gamma$, $\delta_{1}$ not in the same box as $\alpha_{i}$ and 
\[ \delta_1+...+\delta_s = \alpha_i + \gamma  = \beta  \ \ \ (\star) .\]
For all monomials, $\delta_s \in C_{t}$ and $\delta_{s}$ is ordinary (because $\beta \in C_{t}$ and $\beta$ is ordinary). As $\delta_s < \gamma $ and $\delta_{s}$ does not belong to the same box as $\gamma$, one has  $h(\delta_s) > h(\gamma)$. Hence $h(\delta_s) \geq h(\beta)$, so that $s = 1$ and $\delta_1 = \beta$. This implies $E_\beta =  aX'_\beta $ with $ a \in \mathbb{K}\setminus\{0\}$, and the result is proved.\\
\begin{tabular}{@{}lc@{}}
\begin{minipage}{.8\linewidth}
Let us now assume that $\beta$ is the exceptional root of $C_{t}$. Let $\gamma$ be the root of $C_{t}$ which precedes $\beta$ in the Lusztig order and let $\delta = s_{D}(\gamma)$, so that $\delta + \gamma = \beta$ (see Figure 1.). By Proposition \ref{h'}, one has $h'(\beta) = m+\frac{1}{2}$ with $m \in \mathbb{N}^{\star}$ and $h'(C_{t})= \llbracket 1, 2m \rrbracket$. If B is the box in $C_{t}$ which precedes $\beta$, then $h'(B) = h(B) = t+1$. As $\beta$ is alone in its box, we have $\gamma \in B$, so that $h(\gamma) = m+1$. Hence $h(\delta) = m$. Thus $P =$ Vect $(\gamma , \delta)$ is  an admissible plane of type (1.1) or (1.2), and $[E_{\delta},E_\gamma]_q = c E_\beta (c \neq 0)$ (see Section 3.2).\\
As $\gamma$ and $\delta$ are ordinary roots, we already know that $X'_\gamma =  \lambda_\gamma E_\gamma$ and $X'_{\delta} = \lambda_\delta E_{\delta}$ with $\lambda_{\gamma} \neq 0$ and  $\lambda_{\delta} \neq 0$. Thus, one has :
\[[X'_{\delta},X'_\gamma]_q =\lambda_\gamma \lambda_\delta  [E_{\delta},E_\gamma]_q = \lambda_\gamma \lambda_\delta c E_\beta \ (\lambda_{\gamma} \neq 0, \lambda_{\delta} \neq 0).\] 

\end{minipage} &
 \begin{tabular}{c|c|c}
\cline{2-2}
  & $\beta_k$ \\ 
 \cline{2-2}
     $C_t$ :  & \vdots \\ 
\cline{2-2}
     &$\gamma$ \\ 
 \cline{2-2}
    & $\beta$ \\ 
\cline{2-2}
     & $\delta $\\ 
   \cline{2-2}
     & \vdots \\ 
\cline{2-2}
     & $\beta_N$ \\ 
\cline{2-2}
\multicolumn{3}{c}{}\\
\multicolumn{3}{r}{Figure 1.}
  \end{tabular}\\
\end{tabular}\\
\par
As above,  $E_\beta$ is a linear combination of monomials $X'_{\delta_1}...X'_{\delta_s}$ with \\
$\gamma < \delta_{1} \leq ... \leq \delta_{s} < \delta$, $\delta_{s}$ not in the same box as $\delta$ and $\delta_{1}$ not in the same box as $\gamma$. As $\beta$  is the only root of $C_{t}$ which satisfies  $\gamma <\beta < \delta$, $\beta$ is not in the same box as $\delta$ and $\beta$ is not in the same box as $\gamma$. Hence $s=1$ and $\delta_{1} = \beta$. So that $E_\beta = a X'_\beta$ with $a \in \mathbb{K}\setminus\{0\}$.
\end{preuve}
From Theorems \ref{1.4} and \ref{theo3.6}, we deduce the following result.
\begin{cor} \label{LSE}
If $i$ and $j$ are two integers such that $1 \leq i < j \leq N$, one has:
\begin{equation*}
	E_{\beta_i}E_{\beta_j} - q^{-(\beta_i , \beta_j)}E_{\beta_j}E_{\beta_i} = \sum_{\substack{ \beta_{i} < \gamma_{1} < ... < \gamma_{p} < \beta_{j}, \\ p \geq 1, \ k_i \in \mathbb{N}}}  C'_{\mathbf{\overline{k},\overline{\gamma}}} E_{\gamma_{1}}^{k_{1}}...E_{\gamma_{p}}^{k_{p}}.
\end{equation*}
The monomials on the left-hand side whose coefficient $C'_{\mathbf{\overline{k},\overline{\gamma}}}$ is not equal to zero satisfies:
\begin{itemizedot}
\item $wt(X_{\gamma_{1}}^{\prime k_{1}}...X_{\gamma_{p}}^{\prime k_{p}}) = \beta_i+\beta_j$;
\item $\gamma_{1}$ is not in the same box as $\beta_{i}$;
\item $\gamma_{p}$ is not in the same box as $\beta_{j}$.
\end{itemizedot}\end{cor}

\subsubsection{Link with Jantzen's construction}
\begin{prop}\label{relcomX}
Let $\beta_{1} < \beta_{2}$ be two positive roots.
\begin{enumerate}
\item If $E_{\beta_{1}}E_{\beta_{2}}-q^{-(\beta_{1},\beta_{2})}E_{\beta_{2}}E_{\beta_{1}} = k E^{m}_{\gamma}$ ($k \neq 0, m \geq 1$ and  $\gamma \in \Phi^+$), then\\
$X_{\beta_{1}}X_{\beta_{2}}-q^{+(\beta_{1},\beta_{2})}X_{\beta_{2}}X_{\beta_{1}} = k' X^{m}_{\gamma}$ ($k' \neq 0$).
\item If $E_{\beta_{1}}E_{\beta_{2}}-q^{-(\beta_{1},\beta_{2})}E_{\beta_{2}}E_{\beta_{1}} = k E_{\gamma}E_{\delta}$ ($k \neq 0, \gamma, \delta \in \Phi^+, \gamma$ and $\delta$ belonging to the same box), then $X_{\beta_{1}}X_{\beta_{2}}-q^{+(\beta_{1},\beta_{2})}X_{\beta_{2}}X_{\beta_{1}} = k' X_{\gamma}X_{\delta}$ ($k' \neq 0$).
\end{enumerate}
\end{prop}
\begin{preuve}
Let $\beta \in \Phi^{+}$. Let us recall (see Section 3.4.1) that $X_{\beta} := T_{w'_{\beta}}(E_{\alpha_{i_{\beta}}})$, $X'_{\beta} := T'_{w'_{\beta}} (E_{\alpha_{i_{\beta}}})$, and that $T_{w'_{\beta}} = \tau \circ T'_{w'_{\beta}} \circ \tau $. So we have $X_{\beta} = \tau \circ T'_{w'_{\beta}} \circ \tau (E_{\alpha_{i_{\beta}}}) = \tau (X'_{\beta})$. Let us also recall (see Theoreme \ref{theo3.6}) that  $X'_{\beta} = \lambda_{\beta} E_{\beta}$ with $\lambda_{\beta} \in \mathbb{K}^{\star}$.\\
Let $\beta_{1} < \beta_{2}$ be two positive roots.
\begin{enumerate}
\item If $E_{\beta_{1}}E_{\beta_{2}}-q^{-(\beta_{1},\beta_{2})}E_{\beta_{2}}E_{\beta_{1}} = k E^m_{\gamma}$ ($k \neq 0, \gamma \in \Phi^+$), then :
\[\begin{array}{rl}
X_{\beta_{1}}X_{\beta_{2}}-q^{+(\beta_{1},\beta_{2})}X_{\beta_{2}}X_{\beta_{1}} &= \tau(X'_{\beta_{1}})\tau(X'_{\beta_{2}})-q^{+(\beta_{1},\beta_{2})}\tau(X_{\beta_{2}})\tau(X_{\beta_{1}})\\
\tiny\\
 &= \tau(X'_{\beta_{2}}X'_{\beta_{1}}-q^{+(\beta_{1},\beta_{2})}X'_{\beta_{1}}X'_{\beta_{2}})\\
 \tiny \\
 &=-q^{+(\beta_{1},\beta_{2})} \tau(X'_{\beta_{1}}X'_{\beta_{2}}-q^{-(\beta_{1},\beta_{2})}X'_{\beta_{2}}X'_{\beta_{1}})\\
 \tiny \\
 &=-q^{+(\beta_{1},\beta_{2})}\lambda_{\beta_{1}}\lambda_{\beta_{2}} \tau(E_{\beta_{1}}E_{\beta_{2}}-q^{-(\beta_{1},\beta_{2})}E_{\beta_{2}}E_{\beta_{1}})\\
 \tiny \\
 &=-q^{+(\beta_{1},\beta_{2})}\lambda_{\beta_{1}}\lambda_{\beta_{2}} \tau(k E^m_{\gamma})\\
 \tiny \\
 &=\displaystyle \frac{-q^{+(\beta_{1},\beta_{2})}\lambda_{\beta_{1}}\lambda_{\beta_{2}} k}{\lambda_{\gamma}}\tau((X'_{\gamma})^m)\\
 \tiny \\
 &= k' X^m_{\gamma} \ \text{with} \ k' \in \mathbb{K}^{\star}
 \end{array}\]
 \item If $E_{\beta_{1}}E_{\beta_{2}}-q^{-(\beta_{1},\beta_{2})}E_{\beta_{2}}E_{\beta_{1}} = k E_{\gamma}E_{\delta}$ ($k \neq 0, \gamma, \delta \in \Phi^+, \gamma$ and $\delta$ belonging two the same box) so, by doing the same computations as in 1., we obtain :
\[X_{\beta_{1}}X_{\beta_{2}}-q^{(\beta_{1},\beta_{2})}X_{\beta_{2}}X_{\beta_{1}} = k' \tau (X'_{\gamma}X'_{\delta}) = k' X_{\delta}X_{\gamma} \ (k' \neq 0)\]
As $\gamma$ and $\delta$ are in the same box, we know (see Proposition \ref{boîte}) that $(\delta, \gamma) = 0$, so that, by Theorem \ref{LS}, we get $X_{\gamma}X_{\delta} = X_{\gamma}X_{\delta} $, as desired.
\end{enumerate}
\end{preuve}

\section{Deleting derivations in $\mathcal{U}_q^+(\mathfrak{g})$}
\subsection{$\mathcal{U}_q^+(\mathfrak{g})$ is a CGL extension}
In this section, we set $A := \mathcal{U}_q^+(\mathfrak{g})$, $X_i := X_{\beta_{i}}$ for $1 \leq i \leq N$, and $\lambda_{i,j}:= q^{-(\beta_{j},\beta_{i})}$ for $1 \leq i,j \leq N$. We know from Proposition \ref{LS}) that, if $1 \leq i < j \leq N$, then one has:
\begin{equation}
X_{j}X_{i} - \lambda_{j,i}X_{i}X_{j} = P_{j,i}
\end{equation}
with
\begin{equation}\label{Pji}
P_{j,i} = \sum_{\overline{k} = (k_{i+1},...,k_{j-1})} c_{\bar{k}} X_{i+1}^{k_{i+1}}...X_{j-1}^{k_{j-1}},
\end{equation}
where $c_{\bar{k}} \in \mathbb{K}$.\\
Moreover, as $\mathcal{U}_q^+(\mathfrak{g})$  is $\Phi$-gradued, one has 
\begin{equation}\label{impli}
c_{\bar{k}} \neq 0 \Rightarrow \lambda_{l,i+1}^{k_{i+1}} ... \lambda_{l,j-1}^{k_{j-1}} = \lambda_{l,j}\lambda_{l,i} \ \text{for all} \ 1 \leq l \leq N
\end{equation}
Thus, $A$ satisfies  \cite[Hypothesis 6.1.1]{MR1967309}.\\
\par
From Theorem \ref{PBW}, ordered monomials in $X_{i}$ are a basis of $A$, so that we deduce from \cite[Proposition 6.1.1]{MR1967309}: 
\begin{prop}\label{ore}
\begin{enumerate}
\item $A$ is skew polynomial ring which could be expressed as:
\[A = \mathbb{K}[X_{1}][X_2;\sigma_2,\delta_2]...[X_N;\sigma_N,\delta_N],\]
where the $\sigma_j$'s are $\mathbb{K}$-linear automorphisms and the $\delta_j$'s are $\mathbb{K}$-linear $\sigma_j$-derivations such that, for $1 \leq i < j \leq N$, $\sigma_j(X_i) = \lambda_{j,i} X_i$ and $\delta_j(X_i) = P_{j,i}$.
\item If $1 \leq m \leq N$, then there is a (unique) automorphism $h_{m}$ of the algebra $A$ which satisfies $h_{m}(X_{i}) = \lambda_{m,i} X_{i}$ for  $1 \leq i \leq N$.
\end{enumerate}
\end{prop}
Moreover, we deduce from \cite[Proposition 6.1.2]{MR1967309} the following result.
\begin{prop}\label{conv}
\begin{enumerate}
\item $A$ satisfies conventions from \cite[Section 3.1]{MR1967309}, that is to say:
\begin{itemizedot}
\item For all $j \in \llbracket 2, N \rrbracket$, $\sigma_j$ is a $\mathbb{K}$-linear automorphism and $\delta_j$ is a $\mathbb{K}$-linear (left sided) $\sigma_{j}$-derivation and locally nilpotent. 
\item For all $j \in \llbracket 2, N \rrbracket$, one has $\sigma_{j} \circ \delta_{j} =q_{j} \delta_{j} \circ \sigma_{j}$ with $q_{j} = \lambda_{j,j} = q^{-||\beta_{j}||^2}$, and for all $i \in \llbracket 1, j-1 \rrbracket$, $\sigma_{j}(X_{i}) = \lambda_{j,i}X_{i}$.
\item None of the $q_{j} \ (2 \leq j \leq N )$ is a root of unity.
\end{itemizedot}
\item $A$ satisfies \cite[Hypothesis 4.1.2]{MR1967309}, that is to say:\\
The subgroup $H$ of the automorphisms group of A generated by the elements $h_{l}$ satisfies 
\begin{itemizedot}
\item For all h in H, the indeterminates $X_{1}, ..., X_{N}$ are h-eigenvectors.
\item The set $\{\lambda \in \mathbb{K}^{*} \ | \ (\exists h \in H)(h(X_{1}) = \lambda X_{1}\}$ is infinite.
\item If $m \in \llbracket 2, N \rrbracket$, there is $h_{m} \in H$ such that $h_{m}(X_{i}) = \lambda_{m,i}X_{i}$ if $1 \leq i < m$ and $h_{m}(X_{m}) = q_{m}X_{m}$. 
\end{itemizedot}
\end{enumerate}
\end{prop}
The previous proposition shows that $\mathcal{U}_q^+(\mathfrak{g})$ is a CGL extension in the sens of \cite{LLR} and so allows us to apply the deleting derivation theory (\cite{MR1967309}). We describe this theory in the following section.
\subsection{The deleting derivation algorithm}
It follows from Propositions \ref{ore} and \ref{conv}, that $A$ is an integral domain which is noetherian. Denote by $F$ its fields of fraction. We define, by induction, the families $X^{(l)} = (X_{i}^{(l)})_{1 \leq i \leq N}$ of elements of $F^{\star} := F \setminus \{0\}$, and the algebras $A^{(l)} := \mathbb{K}<X_{1}^{(l)},...,X_{N}^{(l)}>$ when $l$ decreases from $N + 1$ to 2 as in \cite[Section 3.2]{MR1967309}. So we have for all $l \in \llbracket 2, N+1 \rrbracket$ :

\begin{lemme}\label{4.1.1}
If $1 \leq i < j \leq N$, one has :
\begin{equation}\label{relationl}
X_{j}^{(l)}X^{(l)}_{i} - \lambda_{j,i}X^{(l)}_{i}X^{(l)}_{j} = P^{(l)}_{j,i}
\end{equation}
 with 
 \begin{equation}\label{Pjil}
P^{(l)}_{j,i} = \left\{ 
\begin{array}{lc}
0 & \text{si} \ j \geq l \\
\displaystyle \sum_{\overline{k} = (k_{i+1},...,k_{j-1})} c_{\bar{k}} (X_{i+1}^{(l)})^{k_{i+1}}...(X_{j-1}^{(l)})^{k_{j-1}} & \text{si} \ j < l,
\end{array}
\right.
\end{equation}
where $c_{\bar{k}}$ are \textbf{the same} as in the formula (\ref{Pji}), so that we also have the implication (\ref{impli}).
\end{lemme}
\begin{preuve}
See \cite[Théorème 3.2.1]{MR1967309}.
\end{preuve}
\begin{lemme}\label{4.1.2}
The ordered monomials on $X_{1}^{(l)},...,X_{N}^{(l)}$ form a basis $A^{(l)}$ as a $\mathbb{K}$-vectorial space.
\end{lemme}
\begin{preuve}
See \cite[Théorème 3.2.1]{MR1967309}.
\end{preuve}
From Lemmas \ref{4.1.1} and \ref{4.1.2} above and from \cite[Proposition 6.1.1]{MR1967309}, we deduce that:
\begin{lemme}\label{4.1.3}
\begin{enumerate}
\item $A^{(l)}$ is an iterated ore extension which can be written :
\[A^{(l)} = \mathbb{K}[X^{(l)}_{1}][X_2^{(l)};\sigma^{(l)}_2,\delta^{(l)}_2]...[X^{(l)}_N;\sigma^{(l)}_N,\delta^{(l)}_N]\]
where $\sigma_{j}^{(l)}$ are $\mathbb{K}$-linear automorphisms and $\delta_{j}^{(l)}$ are $\mathbb{K}$-linear (left sided) $\sigma_{j}^{(l)}$-derivations such that, for $1 \leq i < j \leq N$, $\sigma_{j}^{(l)}(X_{i}^{(l)}) = \lambda_{j,i} X_{i}^{(l)}$ and $\delta_{j}^{(l)}(X_{i}^{(l)}) = P_{j,i}^{(l)}$.
\item  $A^{(l)}$ is the $\mathbb{K}$ algebra generated by the elements $X_{1}^{(l)},...,X_{N}^{(l)}$ with relations (\ref{relationl}).
\end{enumerate}
\end{lemme}
Let us recall that the automorphisms $h_{m}$ ($1 \leq m \leq N$ ) of the algebra $A$ defined in Proposition \ref{ore} can be extended (uniquely) in automorphisms, also denoted by $h_{m}$, of the field $F$.  
\begin{lemme}\label{4.1.4}
If $1 \leq m,i \leq N$, one has  $h_{m}(X_{i}^{(l)}) = \lambda_{m,i} X_{i}^{(l)}$ so that $h_{m}$ induces (by restriction) an automorphism of the algebra $A^{(l)}$, denoted by $h_{m}^{(l)}$.
\end{lemme}
\begin{preuve}
See \cite[Lemme 4.2.1]{MR1967309}.
\end{preuve}
\begin{convsanss}
Denote by $H^{(l)}$ the subgroup of the automorphism group of $A^{(l)}$ generated by $h_{m}^{(l)}$ ($1 \leq m \leq N$) .
\end{convsanss}
By \cite[Proposition 6.1.2]{MR1967309}, one has :
\begin{lemme}\label{4.1.5}
The iterated Ore extension $A^{(l)} = \mathbb{K}[X^{(l)}_{1}][X_2^{(l)};\sigma^{(l)}_2,\delta^{(l)}_2]...[X^{(l)}_N;\sigma^{(l)}_N,\delta^{(l)}_N]$ satisfies the conventions of \cite[Section 3.1]{MR1967309} with, as above, $\lambda_{i,j} = q^{-(\beta_{i},\beta_{j})}$ and $q_{i} = \lambda_{i,i} = q^{-||\beta_{i}||^{2}}$ for $ 1\leq i,j \leq N$. It also satisfies the Hypothesis 4.1.2 of \cite{MR1967309} with $H^{(l)}$ replacing $H$.
\end{lemme}
\begin{cor}\label{Hprime}
If $J$ is a $H^{(l)}$-prime ideal of $A^{(l)}$ in the sense of \cite[II.1.9]{MR1898492}, then J is completely prime.
\end{cor}
\begin{preuve}
One has :
\begin{itemizedot}
\item $A^{(l)} = \mathbb{K}[X^{(l)}_{1}][X_2^{(l)};\sigma^{(l)}_2,\delta^{(l)}_2]...[X^{(l)}_N;\sigma^{(l)}_N,\delta^{(l)}_N]$ is an iterated Ore extension by Lemma \ref{4.1.3},
\item $X^{(l)}_{1},X_2^{(l)},...,X^{(l)}_N$ are $H^{(l)}$-eigenvectors by Lemma \ref{4.1.4}.
\item If $1 \leq i < j \leq N$, then one has $h_{i}^{(l)}(X_{j}^{(l)}) = \lambda_{j,i} X_{i}^{(l)} = \sigma_{j}^{(l)}(X_{i}^{(l)})$ and $h_{j}^{(l)}(X_{j}^{(l)}) = q_{j} X_{j}^{(l)}$ with $q_{j} = \lambda_{j,j} \in \mathbb{K}^{\star}$ is not a root of unity by Lemmas \ref{4.1.3} and \ref{4.1.4}. 
\end{itemizedot}
Hence we deduce from \cite[Theorem II.5.12]{MR1898492} that $J$ is completely prime.
\end{preuve}
From the construction of the deleting algorithm (\cite[Section 3.2]{MR1967309}), one has :
\begin{lemme}\label{4.1.6}
\begin{enumerate}
\item $X_{i}^{(N+1)} = X_{i}$ for all $i \in \llbracket 1, N \rrbracket$
\item If $2 \leq l \leq N$ and if $i \in \llbracket 1, N \rrbracket$, one has 
\begin{equation}\label{Xil}
X_{i}^{(l)} = 
\left\{
\begin{array}{lc}
 X_{i}^{(l+1)} & \text{si} \ i \geq l    \\
\displaystyle \sum_{n=0}^{+ \infty} \left[ \frac{(1-q_{l})^{-n}}{[n]!_{q_{l}}} \left(\delta_{l}^{(l+1)}\right)^{n}\circ \left(\sigma_{l}^{(l+1)}\right)^{-n}\left(X_{i}^{(l+1)} \right) \right]\left(X_{l}^{(l+1)} \right)^{-n} &   \text{si} \ i < l 
\end{array}
\right.
\end{equation}

\end{enumerate}

\end{lemme}
\begin{lemme}\label{4.1.7}
Let $J$ be an $H^{(l)}$-invariant (two sided) ideal of $A^{(l)}$. Let us consider an integer $j \in \llbracket 2, N \rrbracket$ and denote by\\ $B = \mathbb{K}[X^{(l)}_{1}][X_2^{(l)};\sigma^{(l)}_2,\delta^{(l)}_2]...[X^{(l)}_{j-1};\sigma^{(l)}_{j-1},\delta^{(l)}_{j-1}]$ the subalgebra of $A^{(l)}$ generated by $X^{(l)}_{1},...,X^{(l)}_{j-1}$. Then $\sigma_{j}^{(l)}(B \cap J) = B \cap J$ and $\delta_{j}^{(l)}(B \cap J) \subset B \cap J$. 
\end{lemme}
\begin{preuve}
By Lemmas \ref{4.1.3} and \ref{4.1.4}, one has for $1 \leq i < j$, 
\begin{equation}\label{HS}
\sigma_{j}^{(l)}(X_{i}^{(l)}) = \lambda_{j,i}X_{i}^{(l)} = h_{j}^{(l)}(X_{i}^{(l)})
\end{equation}
As a result, for all $b \in B$, $\sigma_{j}^{(l)}(b) = h_{j}^{(l)}(b)$. As $J$ is $H^{(l)}$-invariant, and as $B$ is  $\sigma_{j}^{(l)}$-invariant, we deduce that, for all $b \in B \cap J$, we have $\sigma_{j}^{(l)}(b) \in B \cap J$. So, $\sigma_{j}^{(l)}(B \cap J) \subset B \cap J$. From the equality ($\ref{HS}$), we get that:
\begin{equation}
\left( \sigma_j^{(l)} \right)^{-1}( X_{i}^{(l)} ) = \lambda_{j,i}^{-1} X_{i}^{(l)}  = \left( h_{i}^{(l)} \right)^{-1} ( X_{i}^{(l)} )
\end{equation}
As above, we deduce that $\left( \sigma_{j}^{(l)} \right)^{-1}(B \cap J) \subset B \cap J$, so that $\sigma_{j}^{(l)}(B \cap J) = B \cap J$. \\
Finally, if $b \in B \cap J$, then we have $\delta_{j}^{(l)}(b) = X_{j}^{(l)}b - \sigma_j^{(l)}(b)X_{j}^{(l)} \in B \cap J$.
\end{preuve}
If $l \in \llbracket 2, N \rrbracket$, then it follows from (\ref{Xil}) that $X_{l}^{(l)} = X_{l}^{(l+1)}$. This element is a nonzero element which belongs to the two algebras $A^{(l)}$ and $A^{(l+1)}$ (recall that none of the $X_{i}^{(l)}$ is null). So, the set $S_{l} := \{(X_{l}^{(l)})^{p} \ | \ p \in \mathbb{N} \}$ is a multiplicative system of regular elements of $A^{(l)}$ and $A^{(l+1)}$.
From \cite[Theorem 3.2.1]{MR1967309}, we deduce:
\begin{lemme}\label{4.1.8}
Let $l \in \llbracket 2, N \rrbracket$. Then $S_{l}$ is an Ore set in $A^{(l)}$ and also in $A^{(l+1)}$. Moreover, one has:
\[A^{(l)}S_{l}^{-1} = A^{(l+1)}S_{l}^{-1}\]
\end{lemme}

\subsection{Prime spectrum and diagrams}
Let us recall that the convention are $X_{i} = X_{\beta_{i}}$ for $1 \leq i \leq N$. Denote $\overline{A} := A^{(2)} = \mathbb{K}<T_{\beta_{1}},...,T_{\beta_{N}}>$ with $T_{\beta_{i}} = X_{i}^{(2)}$ for all $i$. By lemmas \ref{4.1.1} and \ref{4.1.3}, $\overline{A}$ is the quantum affine space generated by $T_{\beta_{i}}$ ($1 \leq i \leq N$) with relations $T_{\beta_{j}}T_{\beta_i} = \lambda_{j,i} T_{\beta_i}T_{\beta_j}$ for $1 \leq i < j \leq N$. \\
Let us consider an integer $l \in \llbracket 2, N \rrbracket$ and a prime ideal $P \in \text{Spec}(A^{(l+1)})$. \\
\par
\begin{itemizedot}
\item Assume $X_{l}^{(l+1)} \notin P$.\\
Then, by \cite[Lemmas 4.2.2 and 4.3.1]{MR1967309}, we have $S_{l} \cap P = \emptyset $ and $Q := A^{(l)} \cap P S_{l}^{-1} \in \text{Spec}(A^{(l)})$.
\item Assume $X_{l}^{(l+1)} \in P$.\\
Then, by \cite[Lemma 4.3.2]{MR1967309}, there is a (unique) surjective algebra homomorphism 
\[g:A^{(l)} \rightarrow \frac{A^{(l+1)}}{(X_{l}^{(l+1)})}\] which satisfies, for all i, $g(X_{i}^{(l)} )= \overline{X_{i}^{(l+1)}} (:= X_{i}^{(l+1)} + (X_{l}^{(l+1)}))$, so that $Q = g^{-1}(\frac{P}{(X_{l}^{(l+1)})}) \in \text{Spec}(A^{(l)})$.
\end{itemizedot}
\par
We define this way a map $\phi_{l} : \text{Spec}(A^{(l+1)}) \rightarrow \text{Spec}(A^{(l)})$ that maps  $P $ to $Q$ and, by composing these maps, we obtain a map $\phi = \phi_{2} \circ ... \circ \phi_{N} $ : Spec$(A) \rightarrow  \text{Spec}(\bar{A})$. By \cite[Proposition 4.3.1]{MR1967309}, one has :
\begin{lemme}
Each $\phi_{l}$ ($2 \leq l \leq N$) is injective, so that $\phi$ is injective.
\end{lemme}
We can now define the notion of diagrams and Cauchon diagrams
\begin{dfn}
\begin{enumerate}\label{diagcauchon}
\item We call \textit{diagram} a subset $\Delta $ of the set of positive roots $\Phi^+$, and we note:
 $$\text{Spec}_{\Delta}(\overline{A}) := \{Q \in \text{Spec}(\overline{A}) \ | \ Q \cap \{T_{\beta_1},...,T_{\beta_N}\} = \{T_{\beta} \ | \ \beta \in \Delta\}\}.$$
\item A diagram $\Delta$ is a \textit{Cauchon diagram} if there is $P \in \text{Spec}(A)$ such that $\phi(P) \in \text{Spec}_{\Delta}(\overline{A})$, that is to say, if $\phi(P) \cap  \{T_{\beta_1},...,T_{\beta_N}\} =  \{T_{\beta} \ |  \ \beta \in \Delta\}$. In this case, we set
$$\textrm{Spec}_{\Delta}(A) = \{P \in \textrm{Spec}(A) \ | \ \phi(P) \in \textrm{Spec}_{\Delta}(\overline{A})\}.$$
\end{enumerate}
\end{dfn}
By \cite[Theorems 5.1.1, 5.5.1 et 5.5.2]{MR1967309}, we have: 
\begin{prop}\label{Hstrat}
\begin{enumerate}
\item If $\Delta$ is a Cauchon diagram, then  $\phi (\textrm{Spec}_{\Delta}(A)) = \textrm{Spec}_{\Delta}(\overline{A})$ and $\phi$ induced a bi-increasing homeomorphism from $\textrm{Spec}_{\Delta}(A)$ onto $\textrm{Spec}_{\Delta}(\overline{A})$.
\item The family $\textrm{Spec}_{\Delta}(A)$ (with $\Delta$ Cauchon diagram) coincide  with the Goodearl-Letzter H-stratification of $\textrm{Spec}(A)$ (\cite{MR1898492}).
\end{enumerate}
\end{prop}
In the following section, we describe more precisely Cauchon Diagrams. In order to do this, the criteria in the next proposition will be needed. 
\begin{prop} \label{critereimage}
Let $P\it^{(m)}$ be an $H$-prime ideal of $A^{(m)}$.\\
$P\it^{(m)}$ $ $ $\in$ $ $ Im$(\phi_{m})$ $ $ if and only if one of two following conditions is satisfied \\
\begin{enumerate}
\item $ $ $X_{m}^{(m)}$ $ $ $\notin$ $ $ $P\it^{(m)}$.
\item $ $ $X_{m}^{(m)}$ $ $ $\in$ $ $ $P\it^{(m)}$ $ $ and $ $ $\Theta^{(m)}(\delta_{m}^{(m+1)}(X_{i}^{(m+1)})) $  $\in$ $ $ $P\it^{(m)}$ $ $ for $ $ $1$  $\leq$ $i$ $\leq$ $m-1$.\\
(where $ $ $\delta_{m}^{(m+1)}$($X_{i}^{(m+1)}$) =  $P_{m,i}^{(m+1)}$ ($X_{i+1}^{(m+1)}, \mbox{ } ... \mbox{ } , \mbox{ } X_{m-1}^{(m+1)}$) (lemma \ref{4.1.1}) and $\Theta^{(m)}: \mbox{ } \mathbb{K}<X_{1}^{(m+1)}, \mbox{ } ... \mbox{ } , \mbox{ } X_{m-1}^{(m+1)}> \mbox{ } \rightarrow \mbox{ } \mathbb{K}<X_{1}^{(m)}, \mbox{ } ... \mbox{ } , \mbox{ } X_{m-1}^{(m)}>$ $ $ is the homomorphism which send $ $ $X_{l}^{(m+1)}$ $ $ to $ $ $X_{l}^{(m)}$)
\end{enumerate}  
\end{prop}
\begin{preuve}
Assume that $ $ $P\it^{(m)}$ $ $ $\in$ $ $ Im$(\phi_{m})$, $ $ so that $ $ $P\it^{(m)}$ $ $ = $ $ $\phi_{m}$($P\it^{(m+1)}$) $ $ with  $ $ $P\it^{(m+1)}$ $ $ $\in$ $ $ Spec$(A^{(m+1)})$, $ $ and assume that condition 1. is not satisfied. This implies that $ $ $P\it^{(m)}$ $ $ = $ $ ker$(g)$ $ $ where $ $ $g: \mbox{ } A^{(m)} \mbox{ } \rightarrow \mbox{ } A^{(m+1)} / P\it^{(m+1)}$ $ $ is the homomorphism which sends $ $ $X_{i}^{(m)}$ $ $ to $ $ $x_{i}^{(m+1)}$ $ $ := $ $ $X_{i}^{(m+1)}$ $ $ + $ $ $P\it^{(m+1)}$.
Let $ $ $1$  $\leq$ $i$ $\leq$ $m-1$. \\ Recall that $ $ $\delta_{m}^{(m+1)}$($X_{i}^{(m+1)}$) $ $ = $ $ $P_{m,i}^{(m+1)}$($X_{i+1}^{(m+1)}, \mbox{ } ... \mbox{ } , \mbox{ } X_{m-1}^{(m+1)}$) and $ $ and that $ $ $\Theta^{(m)}: \mbox{ } k<X_{1}^{(m+1)}, \mbox{ } ... \mbox{ } , \mbox{ } X_{m-1}^{(m+1)}> \mbox{ } \rightarrow \mbox{ } k<X_{1}^{(m)}, \mbox{ } ... \mbox{ } , \mbox{ } X_{m-1}^{(m)}>$ $ $ is the homomorphism which transforms each $ $ $X_{l}^{(m+1)}$ $ $ in $ $ $X_{l}^{(m)}$. $ $ Since $ $ $X_{m}^{(m)}$ $ $ $\in$ $ $ $P\it^{(m)}$, $ $ we have $ $ $X_{m}^{(m+1)}$ $ $ $\in$ $ $ $P\it^{(m+1)}$ $ $ \cite{MR1967309}, Proposition 4.3.1. and so, $ $ $\delta_{m}^{(m+1)}$($X_{i}^{(m+1)}$) $ $ $\in$ $ $ $P\it^{(m+1)}$. $ $ Now, we have \\ $g$($\Theta^{(m)}$($\delta_{m}^{(m+1)}$($X_{i}^{(m+1)}$))) $ $ = $ $ $g$($\Theta^{(m)}$($P_{m,i}^{(m+1)}$($X_{i+1}^{(m+1)}, \mbox{ } ... \mbox{ } , \mbox{ } X_{m-1}^{(m+1)}$))) $ $ = $ $ $g$($P_{m,i}^{(m+1)}$($X_{i+1}^{(m)}, \mbox{ } ... \mbox{ } , \mbox{ } X_{m-1}^{(m)}$)) $ $ \\= $ $ $P_{m,i}^{(m+1)}$($x_{i+1}^{(m+1)}, \mbox{ } ... \mbox{ } , \mbox{ } x_{m-1}^{(m+1)})$ $ $ = $ $ $P_{m,i}^{(m+1)}$($X_{i+1}^{(m+1)}, \mbox{ } ... \mbox{ } , \mbox{ } X_{m-1}^{(m+1)})$ $ $ + $ $ $P\it^{(m+1)}$ $ $ = $ $ $0$. \\ This implies that $ $ $\Theta^{(m)}$($\delta_{m}^{(m+1)}$($X_{i}^{(m+1)}$))) $ $ $\in$ $ $ $ker(g)$ $ $ = $ $ $P\it^{(m)}$. \\ $ $ \\
If condition 1. is satisfied, then $ $ $P\it^{(m)}$ $ $ $\in$ $ $ $Im(\phi_{m})$ $ $ by \cite[Lemma 4.3.1.]{MR1967309}. \\ $ $ \\
Assume that condition 2. is satisfied. Let $ $ $1$  $\leq$ $i$ $\leq$ $m-1$, $ $. Then we have, as previously, $ $  $P_{m,i}^{(m+1)}$($X_{i+1}^{(m)}, \mbox{ } ... \mbox{ } , \mbox{ } X_{m-1}^{(m)}$) $ $ = $ $ $\Theta^{(m)}$($\delta_{m}^{(m+1)}$($X_{i}^{(m+1)}$)) $ $ $\in$ $ $ $P\it^{(m)}$. $ $ So, in $ $ $Q^{(m)}$ $ $ = $ $ $A^{(m)}/P\it^{(m)}$, $ $  we have  $ $  $P_{m,i}^{(m+1)}$($x_{i+1}^{(m)}, \mbox{ } ... \mbox{ } , \mbox{ } x_{m-1}^{(m)}$) $ $ = $ $ $0$. \\ Since $ $ $P_{m,i}^{(m)}$ $ $ = $ $ $0$ $ $ (see Lemma \ref{4.1.1}), we can write $ $ $ $ $x_{m}^{(m)}x_{i}^{(m)}$ $ $ $-$ $ $ $\lambda_{m,i}x_{i}^{(m)}x_{m}^{(m)}$ $ $ = $ $ $P_{m,i}^{(m)}$($x_{i+1}^{(m)}, \mbox{ } ... \mbox{ } , \mbox{ } x_{m-1}^{(m)}$) $ $ = $ $ $0$ $ $ = $ $ $P_{m,i}^{(m+1)}$($x_{i+1}^{(m)}, \mbox{ } ... \mbox{ } , \mbox{ } x_{m-1}^{(m)}$). \\  If  $ $ $1$  $\leq$ $i$ $\leq$ $j-1$ $ $ with $j$ $ $ $\neq$ $ $ $m$, $ $ it follows from Lemma \ref{4.1.1} that: \\ $ $ $x_{j}^{(m)}x_{i}^{(m)}$ $ $ $-$ $ $ $\lambda_{j,i}x_{i}^{(m)}x_{j}^{(m)}$ $ $ = $ $ $P_{j,i}^{(m)}$($x_{i+1}^{(m)}, \mbox{ } ... \mbox{ } , \mbox{ } x_{j-1}^{(m)}$) $ $ = $ $ $P_{j,i}^{(m+1)}$($x_{i+1}^{(m)}, \mbox{ } ... \mbox{ } , \mbox{ } x_{j-1}^{(m)}$). \\ So, by the universal property of algebras defined by generators and relations, there exists a (unique) homomorphism $ $ $\epsilon: \mbox{ } A^{(m+1)} \mbox{ } \rightarrow \mbox{ } Q^{(m)}$ $ $ which sends $ $ $X_{l}^{(m+1)}$ $ $ to $ $ $x_{l}^{(m)}$ for all $l$. $ $ This homomorphism is surjective,  and its kernel $ $ $ker(\epsilon)$ $ $ = $ $ $P\it^{(m+1)}$ $ $ is a prime ideal of $ $ $A^{(m+1)}$. $ $ We observe that, since $ $ $X_{m}^{(m)}$ $ $ $\in$ $ $ $P\it^{(m)}$, $ $ we have $ $ $X_{m}^{(m+1)}$ $ $ $\in$ $ $ $P\it^{(m+1)}$, $ $  and that $ $ $\epsilon$ $ $ induces an automorphism
$$ \overline{\epsilon}: \mbox{ } A^{(m+1)} / P\it^{(m+1)} \mbox{ } \rightarrow \mbox{ } Q^{(m)} \mbox{ } = \mbox{ } A^{(m)} / P\it^{(m)}$$ which sends  $ $ $x_{l}^{(m+1)}$ $ $ to $ $ $x_{l}^{(m)}$ for all $l$. $ $ Recall that $ $ $f_{m}: \mbox{ } A^{(m)} \mbox{ } \rightarrow \mbox{ } A^{(m)} / P\it^{(m)}$  denotes the canonical homomorphism. \\ So, $g \mbox{ } = \mbox{ } (\overline{\epsilon})^{-1} \circ f_{m}: \mbox{ } A^{(m)} \mbox{ } \rightarrow \mbox{ } A^{(m+1)} / P\it^{(m+1)}$ $ $ is the homomorphism which sends $ $ $X_{l}^{(m)}$ $ $ to $ $ $x_{l}^{(m+1)}$ for all $l$. $ $ As $ $ $ker(g)$ $ $ = $ $ $ker(f_{m})$ $ $ = $ $ $P\it^{(m)}$, $ $ we conclude that  $ $ $P\it^{(m)}$ $ $ = $ $ $\phi_{m}$($P\it^{(m+1)}$), as desired. 
\end{preuve}
\section{Cauchon diagrams in $\mathcal{U}_q^+(\mathfrak{g})$}
In \cite{MR1967310}, Cauchon uses a combinatorial tool to describe "admissible diagrams" (which are called "Cauchon diagrams" here) for the algebra $O\sb q(M\sb n(k))$ of quantum matrices. Thanks to Lusztig admissible planes theory (see Section 3.2), results from Section 3.3 and the deleting derivation theory, we describe those diagrams for $\mathcal{U}_q^+(\cat{g})$ (where $\cat{g}$ is a simple Lie algebra of finite dimension over $\mathbb{C}$). The goal of this section is to prove the following statement:\\
\par
\textbf{Theorem:}  A diagram $\Delta \subset \Phi^+$ satisfies all the implications from admissible planes (to be defined) if and only if $\Delta$ is a Cauchon diagram (in the sense of Definition \ref{diagcauchon}) .\\
\subsection{Implications in a diagram}
\begin{lemme}\label{lemme4.1}
Let $j \in \llbracket 1, N \rrbracket$, $l \in \llbracket 2, N \rrbracket$, $P^{(l+1)}$ be a prime ideal of $A^{(l+1)}$ and $P^{(l)} := \varphi_l(P^{(l+1)})$.
\begin{enumerate}
\item  If $X_{j}^{(l+1)} \in P^{(l+1)}$, then $X_{j}^{(l)} \in P^{(l)}$.
\item If $X_{j}^{(l+1)} = X_{j}^{(l)}$ (this is in particular the case if  $j \geq l$), then one has:  $X_{j}^{(l+1)} \in P^{(l+1)}$ if and only if $X_{j}^{(l)} \in P^{(l)}$.
\end{enumerate}
\end{lemme}
\begin{preuve}
The second point can be shown as in \cite[lemma 4.3.4]{MR1967309}. Let us show the first point when $j < l$.
\begin{description}
\item[$1^{st}$ case :] The pivot (in reference to Gaussian elimination) $\varpi:=X_{l}^{(l+1)}$ belongs to $P^{(l+1)}$. Recall (see Section 4.2) that there is a surjective homomorphism of algebra \[g:A^{(l)} \rightarrow \frac{A^{(l+1)}}{(X_{l}^{(l+1)})}\] which satisfies $g(X_{i}^{(l)} )= \overline{X_{i}^{(l+1)}} (:= X_{i}^{(l+1)} + (X_{l}^{(l+1)})$) for all $i \in \llbracket 1, N \rrbracket$. As $X_{j}^{(l+1)} \in P^{(l+1)}$, one has $g(X_{j}^{(l)}) \in \frac{P^{(l+1)}}{(X_{l}^{(l+1)})}$, so that $X_{j}^{(l)} \in g^{-1}(\frac{P^{(l+1)}}{(X_{l}^{(l+1)})}) =: P^{(l)}$.

\item[$2^{nd}$ case :] The pivot $\varpi = X_{l}^{(l+1)}$ does not belongs to $P^{(l+1)}$. Set $S_l := \{\varpi^n | n \in \mathbb{N}\}$. Recall (seeSection 4.2) that  we have $P^{(l)} = A^{(l)} \cap (P^{(l+1)}S_{l}^{-1})$.\\
Set  $\displaystyle J := \bigcap_{h \in H^{(l+1)}} h( P^{(l+1)})$ and observe that J is an $H^{(l+1)}$-invariant two-sided ideal by construction. As $A^{(l+1)}$ \cite[Hypothesis 4.1.2 ]{MR1967309} by Lemma \ref{4.1.5}, $X_{j}^{(l+1)}$ is an $H^{(l+1)}$-eigenvector. Thus, since $X_{j}^{(l+1)}$ belongs to $P^{(l+1)}$, it also belongs to J.\\
From Lemma \ref{4.1.7}, we deduce that $\left(\delta_{l}^{(l+1)}\right)^{n}\circ \left(\sigma_{l}^{(l+1)}\right)^{-n}\left(X_{j}^{(l+1)} \right) \in J \subset P^{(l+1)}$ for all $n \in \mathbb{N}$. As a result, we get :
\[X_{j}^{(l)} = \sum_{n=0}^{+ \infty} \left[ \frac{(1-q_{l})^{-n}}{[n]!_{q_{l}}}\left(\delta_{l}^{(l+1)}\right)^{n}\circ \left(\sigma_{l}^{(l+1)}\right)^{-n}\left(X_{j}^{(l+1)} \right) \right] \left(X_{l}^{(l+1)} \right)^{-n} \in  P^{(l+1)}S_{l}^{-1}.\] 
Thus, $X_{j}^{(l)} \in A^{(l)} \cap (P^{(l+1)}S_{l}^{-1}) = P^{(l)}$. 
\end{description}
\end{preuve}
\begin{lemme}\label{5.2}
Let $l \in \llbracket 2, N \rrbracket$ and $P^{(l+1)}$ be a prime ideal of $A^{(l+1)}$. Consider an integer $j$ with $2 \leq j < l$ and set $P^{(j)} = \varphi_j \circ ... \circ \varphi_l (P^{(l+1)})$.
\begin{enumerate}
\item Assume that $\beta_j$ is in the same box as $\beta_l$ or in the box before $\beta_l$'s one. Then
\begin{itemizedot}
\item $X_j^{(j+1)} = X_j^{(j+2)} = ... = X_j^{(l+1)}$,
\item $(X_j^{(j+1)} \in P^{(j+1)}) \Rightarrow (X_j^{(j+2)} \in P^{(j+2)}) \Rightarrow ... \Rightarrow (X_j^{(l+1)} \in P^{(l+1)})$.
\end{itemizedot}   
\item Assume that the boxes $B$ and $B'$ of $\beta_j$ and $\beta_l$ (respectively) are separated by a box $B''$ containing a unique element $\beta_e$ such that $X_e^{(e+1)} \in P^{(e+1)}$. Then $(X_j^{(j+1)} \in P^{(j+1)}) \Rightarrow (X_j^{(l+1)} \in P^{(l+1)})$.
\end{enumerate}
\end{lemme}
\begin{preuve}
\begin{enumerate}
\item Let $k \in  \llbracket j+1, l\rrbracket$ so that $\beta_k$ is, in the same box as as $\beta_j$, or in the same box as $\beta_l$. As these boxes are consecutive or equal, one has $X_kX_j = q^{-<\beta_k,\beta_j>}X_jX_k$, so that by Lemma \ref{4.1.1}, we have $X_k^{(k+1)}X_j^{(k+1)} = q^{-<\beta_k,\beta_j>}X_j^{(k+1)}X_k^{(k+1)}$. So one has  $\delta_k^{(k+1)}(X_j^{(k+1)}) = 0$ and, by \cite[Section 3.2]{MR1967309}, we get :
\[X_j^{(k)} = \sum_{s=0}^{+ \infty} \lambda_s \left( \delta_k^{(k+1)}\right)^s \circ \left( \sigma_k^{(k+1)} \right)^{-s}(X_j^{(k+1)})\left( X_k^{(k+1)}\right)^{-s} = \sum_{s=0}^{+ \infty} \lambda'_s \left( \delta_k^{(k+1)}\right)^s(X_j^{(k+1)})\left( X_k^{(k+1)}\right)^{-s} =  X_j^{(k+1)}\]
($\lambda_s, \lambda'_s \in \mathbb{K}$).\\
This shows the first point. The second point follows from Lemma \ref{lemme4.1}.
\item As $B$ and $B''$ are consecutive, 1. implies that $X_{j}^{(j+1)} = ... = X_{j}^{(e+1)}$ and that $(X_j^{(j+1)} \in P^{(j+1)}) \Rightarrow ... \Rightarrow (X_j^{(e+1)} \in P^{(e+1)})$. It just remains to show that:
\[X_j^{(k)} \in P^{(k)} \Rightarrow X_j^{(k+1)} \in P^{(k+1)} \ \text{pour} \ e+1 \leq k \leq l.\]
We do that by induction on $k$. As in the previous point, we have 
\[X_j^{(k)} = X_j^{(k+1)} + \sum_{s=1}^{+ \infty} \lambda_s \left( \delta_
k^{(k+1)}\right)^s \circ \left( \sigma_k^{(k+1)} \right)^{-s}(X_j^{(k+1)})\left( X_k^{(k+1)}\right)^{-s} \ \ (\lambda_s \in \mathbb{K})\]
\begin{itemizedot}
\item If $\delta_k^{(k+1)}(X_j^{(k+1)}) = 0$, then one has $X_j^{(k)} = X_j^{(k+1)}$ and we conclude thanks to Lemma \ref{lemme4.1}.
\item Otherwise, one has $\delta_k^{(k+1)}(X_j^{(k+1)}) = \lambda \left(X_e^{(k+1)}\right)^m$ $(m \in \mathbb{N}^\star, \lambda \in \mathbb{K}^\star$) by Lemma \ref{4.1.1} and, as $B'$ and $B''$ are consecutive,   
\[ \begin{array}{cl}
\delta_k^{(k+1)}\left( X_e^{(k+1)} \right) = 0 & \Rightarrow \left( \delta_k^{(k+1)}\right)^s(X_j^{(k+1)}) =  \lambda \left( \delta_k^{(k+1)}\right)^{s-1}\left( \left(X_e^{(k+1)}\right)^m \right) =0 \ \text{pour} \ s>1\\
 & \Rightarrow X_j^{(k)} = X_j^{(k+1)} + \lambda' \left(X_e^{(k+1)}\right)^m\left( X_k^{(k+1)}\right)^{-1} \ \text{avec} \ \lambda' \in \mathbb{K}^\star.
 \end{array}\]
\begin{itemizedot}
\item If $X_k^{k+1} \in P^{(k+1)}$, then consider the homomorphism $g : A^{(k)} \rightarrow \frac{A^{(k+1)}}{(X_k^{(k+1)})}$ which satisfies $g(X_i^{(k)}) = \overline{X_i^{(k+1)}} (:= X_i^{(k+1)} + (X_k^{(k+1)}) )$ for $i \in \llbracket 1, N \rrbracket$ (see Section 4.2). By definition of $\phi_k$ (\cite[Notation 4.3.1.]{MR1967309}), one has $P^{(k)} = g^{-1}\left( \frac{P^{(k+1)}}{(X_k^{(k+1)})} \right)$. So 
\[X_j^{(k)} \in P^{(k)} \Rightarrow g(X_j^{(k)}) = \overline{X_j^{(k+1)}} \in \frac{P^{(k+1)}}{(X_k^{(k+1)})} \Rightarrow  X_j^{(k+1)} \in P^{(k+1)}\]
\item  By 1., one has $X_e^{(e+1)}  = ... = X_e^{(k)}= X_e^{(k+1)}$ and $(X_e^{(e+1)} \in P^{(e+1)}) \Rightarrow ... \Rightarrow (X_e^{(k)} \in P^{(k)}) \Rightarrow (X_e^{(k+1)} \in P^{(k+1)})$. Set, as in \cite[Theorem 3.2.1]{MR1967309}, $S_k := \{\left(X_k^{(k+1)}\right)^n \ | \ n \in \mathbb{N} \}$ so that $P^{(k+1)} = A^{(k+1)}\cap (P^{(k)}S_k^{-1}) $ by definition of $\varphi_k$ \cite[Notation 4.3.1.]{MR1967309}. Then one has :\\
\[X_j^{(k+1)} = X_j^{(k)} - \lambda' \left(X_e^{(k+1)}\right)^m\left( X_k^{(k+1)}\right)^{-1} = X_j^{(k)} - \lambda' \left(X_e^{(k)}\right)^m\left( X_k^{(k+1)}\right)^{-1} \in P^{(k)}S_k^{-1}.\]
As $X_j^{(k+1)}$ is also in $A^{(k+1)}$, one has $X_j^{(k+1)} \in P^{(k+1)}$, as claimed.
\end{itemizedot} 
\end{itemizedot}
\end{enumerate}
\end{preuve}
We use \cite[Proposition 5.2.1]{MR1967309} to determine the shape of Cauchon diagrams. Let us rewrite this proposition in our notation: 
\begin{prop}\label{prop4.2}
Let $\Delta$ be a Cauchon diagram and let $P \in \text{Spec}(A)$. The ideal $P$ belongs to Spec$_{\Delta}(A)$ if and only if it satisfies the following criteria :
\[ (\forall \ l \in \llbracket 1, N \rrbracket ) \ \ (X_{l}^{(l+1)} \in P^{(l+1)} \Leftrightarrow \beta_{l} \in \Delta)  \]
\end{prop}
We can now prove the following proposition.
\begin{prop}\label{prop4.4.3}
Let $\Delta$ be a Cauchon diagram and $\beta_{l} \in \Delta$ ($1 \leq l \leq N$). Assume there is an integer $k \in \llbracket 1, l-1\rrbracket$ such that $X_{\beta_{l}}X_{\beta_{k}}-q^{-(\beta_{l},\beta_{k})}X_{\beta_{k}}X_{\beta_{l}} = c X_{\beta_{i_{1}}}...X_{\beta_{i_{s}}}$ with $c \in \mathbb{K}^{\star}, s \geq 1$ and $k < i_{1} \leq ... \leq i_{s} < l$. Then one of the $\beta_{i_{r}} \ (1 \leq r \leq s)$ belongs to $\Delta$.
\end{prop}
\begin{preuve}
Let $P \in \text{Spec}_{\Delta}(A)$. By Lemma \ref{4.1.1}, one has:
 $$X_{l}^{(l+1)}X^{(l+1)}_{k}-q^{-(\beta_{l},\beta_{k})}X^{(l+1)}_{k}X^{(l+1)}_{l} = c X^{(l+1)}_{\beta_{i_{1}}}...X^{(l+1)}_{\beta_{i_{s}}} := M.$$
By Proposition \ref{prop4.2}, one has $X_{l}^{(l+1)} \in P^{(l+1)}$ so that $M \in P^{(l+1)}$. As $P^{(l+1)}$ is a prime ideal of $A^{(l+1)}$, we know by \cite[II.6.9]{MR1898492} that $P^{(l+1)}$ is completely prime, so that there exists $r \in \llbracket 1, s \rrbracket$ such that $X_{\beta_{i_{r}}}^{(l+1)} \in P^{(l+1)}$. By Lemma \ref{lemme4.1}, we deduce that $X_{i_{r}}^{(i_{r}+1)} \in P^{(i_{r}+1)}$ and, by Proposition \ref{prop4.2}, we obtain $\beta_{i_{r}} \in \Delta$.
\end{preuve}
\begin{convsanss}
We say that a diagram $\Delta$ satisfies the implication
\begin{enumerate}
\item $\beta_{j_{0}} \rightarrow \beta_{j_{1}}$  if $(\beta_{j_{0}} \in \Delta) \Rightarrow (\beta_{j_{1}} \in \Delta)$.
\item \begin{tabular}{cc}
\scalebox{.8}{\begin{psmatrix}[rowsep=.2cm, colsep=1cm]
& $\beta_{j_{1}}$\\
$\beta_{j_{0}}$ & \vdots \\
& $\beta_{j_{s}}$
\end{psmatrix}
\psset{nodesep=2pt, linestyle=dashed, arrows=->}
\ncline{2,1}{1,2}
\ncline{2,1}{2,2}
\ncline{2,1}{3,2}}
&
if $(\beta_{j_{0}} \in \Delta) \Rightarrow ( \beta_{j_{1}} \in \Delta)$ or ... or $(\beta_{j_{1}} \in \Delta )$.
\end{tabular}
\end{enumerate}
\end{convsanss}
Proposition \ref{prop4.4.3} can be rewritten as follows:
\begin{prop} \label{prop4.4}
Let $\Delta$ be a Cauchon diagram and $\beta_{l} \in \Delta$ ($1 \leq l \leq N$). Assume that there exists an integer $k \in \llbracket 1, l-1\rrbracket$ such that $X_{\beta_{l}}X_{\beta_{k}}-q^{-(\beta_{l},\beta_{k})}X_{\beta_{k}}X_{\beta_{l}} = c X_{\beta_{i_{1}}}^{m_{1}}...X_{\beta_{i_{s}}}^{m_{s}}$ with $c \in \mathbb{K}^{\star}, s \geq 1, k < i_{1} < ... < i_{s} < l$ and  $m_{1},...,m_{s} \in \mathbb{N}^{\star}$.  
\begin{enumerate}
\item If $s = 1$, then the solid arrow $\beta_{l} \rightarrow \beta_{i_{1}}$ is an implication. 
\item If $s \geq 2$, then the system 
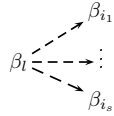
\begin{center}
\scalebox{.8}{\begin{psmatrix}[rowsep=.2cm, colsep=1cm]
& $\beta_{i_{1}}$\\
$\beta_{l}$ & \vdots \\
& $\beta_{i_{s}}$
\end{psmatrix}
\psset{nodesep=2pt, linestyle=dashed, arrows=->}
\ncline{2,1}{1,2}
\ncline{2,1}{2,2}
\ncline{2,1}{3,2}}
\end{center}
of dashed arrows is an implication.
\end{enumerate}
\end{prop}
In the three following propositions, denotes by $\Delta$ a Cauchon diagram.
\begin{prop}\label{contrainte elementaire}
Let $1 \leq l \leq n$ and $\beta \in C_{l}$. If there exists $i \in \llbracket 1, l-1 \rrbracket$ such that $\beta + \alpha_{i} = m \beta'$ with $m \in \mathbb{N}^{\star}$ and $\beta' \in \Phi^{+}$, then $\beta \rightarrow \beta'$ is an implication. 
\end{prop}
\begin{preuve}
We know (see Proposition \ref{propG2} when $\Phi$ is of type $G_2$, Corollary \ref{corplanadm} when $\Phi$ is not of type $G_2$) that we have in this case a commutation relation of the type $E_{\beta}E_{\alpha_{i}} - q^{(\beta,\alpha_{i})}E_{\alpha_{i}}E_{\beta} = k E^{m}_{\beta'}$ with $k \neq 0$ (where $E_{\gamma}$ are defined in Section 3.2).\\
Then, it follows from Proposition \ref{relcomX} that $X_{\beta}X_{\alpha_{i}} - q^{-(\beta,\alpha_{i})}X_{\alpha_{i}}X_{\beta} = k' X^{m}_{\beta'}$ with $k' \neq 0$. So we deduce from proposition \ref{prop4.4} that $\beta \rightarrow \beta'$ is an implication.
\end{preuve}
\begin{prop}\label{contrainte elementaire avant ex}
Let $C_{l}$ ($1 \leq l \leq n$) be an exceptional column. If $\beta \in C_{l}$ is in the box following the box of the exceptional root $\beta_{ex}$, then $\beta \rightarrow \beta_{ex}$ is an implication.  
\end{prop}
\begin{preuve}
Suppose first that $\Phi$ is of type $G_2$. With the notation of Proposition \ref{propG2}, one has $l=2$, $\beta_{ex} = \beta_4$, $\beta = \beta_5$ and one has a commutation formula of the type $E_{\beta_5}E_{\beta_3} - q^{(\beta_3,\beta_5)}E_{\beta_3}E_{\beta_5} = k E_{\beta_{4}}$ with $k \in \mathbb{K}^\star$. It implies, by Proposition \ref{prop4.4} that $\beta = \beta_5 \rightarrow \beta_{ex}= \beta_4$ is an implication.\\
\par 
Suppose now that $\Phi$ is not of type $G_2$. We know (see Proposition \ref{h'}) that $h'(\beta_{ex}) = t + \frac{1}{2}$ ($t \in \mathbb{N}^{\star}$), so that $h'(\beta) = h(\beta) = t$. We also know (see Proposition \ref{prop2.8}) that if $D = $Vect$(\beta_{ex})$, one has $\beta' = s_{D}(\beta) = \beta_{ex} - \beta \in C_{l}$, so that $h'(\beta') = h(\beta') = h(\beta_{ex}) - h(\beta) =  t+1$. As a result, $P = $Vect$(\beta, \beta')$ is an admissible plane of type (1.1) or (1.2). So, by Proposition \ref{rellusztig}, we have a commutation relation of the type $E_{\beta}E_{\beta'} - q^{(\beta,\beta')}E_{\beta'}E_{\beta} = k E_{\beta_{ex}}$ with $k \neq 0$. As in Proposition \ref{contrainte elementaire}, this implies that $\beta \rightarrow \beta_{ex}$ is an implication.  
\end{preuve}
\begin{prop}\label{contrainte elementaire apres ex}
Let $C_{l}$ ($1 \leq l \leq n$) be an exceptional column and $\beta_{ex}$ be its exceptional root. Assume that there exists $i \in \llbracket 1, l \llbracket$ such that $\beta_{ex} + \alpha_{i} = \beta'_{i_{1}} + \beta'_{i_{2}}$ with $\beta'_{i_{1}} \ \neq \ \beta'_{i_{2}}$ in the box which precedes $\beta_{ex}$. Then the system \begin{center}
\scalebox{.8}
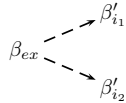
{\begin{psmatrix}[rowsep=.2cm, colsep=1cm]
& $\beta'_{i_{1}}$\\
$\beta_{ex}$ & \\
& $\beta'_{i_{2}}$
\end{psmatrix}
\psset{nodesep=2pt, linestyle=dashed, arrows=->}
\ncline{2,1}{1,2}
\ncline{2,1}{3,2}}
\end{center}
is an implication.
\end{prop}
\begin{preuve}
As, by hypothesis, $\beta'_{i_{1}} \neq \beta'_{i_{2}}$ are in the box preceding the box of $\beta_{ex}$, the root system is not of type $G_2$ (see Proposition \ref{propG2}).\\
As in the proof of Proposition \ref{contrainte elementaire}, it is enough to prove that:
$$[E_{\beta_{ex}},E_{\alpha_{i}}]_{q}:=E_{\beta_{ex}}E_{\alpha_{i}} - q^{(\beta_{ex},\alpha_{i})}E_{\alpha_{i}}E_{\beta_{ex}} = \lambda E_{\beta'_{i_{1}}}E_{\beta'_{i_{2}}} \ \text{avec} \ \lambda \in \mathbb{K}^{\star}$$ 
Recall from Proposition \ref{prop2.8} that $\beta_{ex} \ \bot \ \alpha_{i}$, so that:
$$(\alpha_{i} , \beta'_{i_{1}} + \beta'_{i_{2}} ) = (\alpha_{i}, \beta_{ex} + \alpha_{i}) = ||\alpha_{i}||^2 \Rightarrow (\alpha_{i} , \beta'_{i_{1}} ) > 0 \ \text{or} \ (\alpha_{i},\beta'_{i_{2}} ) > 0.$$ We can assume, without loss of generality, that $(\alpha_{i},\beta'_{i_{2}} ) > 0$, so that (Corollary \ref{corplanadm}) $[E_{\beta'_{i_{2}}},E_{\alpha_{i}}]_{q} = 0$.\\
\par
As in the proof of the previous proposition, one has:
\begin{itemizedot}
\item $h'(\beta_{ex}) = t + \frac{1}{2}$ ($t \in \mathbb{N}^{\star}$) and $h'(\beta'_{i_{1}}) = h'(\beta'_{i_{1}}) = t + 1$,
\item  $\beta_{i_{1}} = s_{D} (\beta'_{i_{1}})$ and $\beta_{i_{2}} = s_{D} (\beta'_{i_{2}})$ belong to $C_{l}$ and satisfy $h'(\beta_{i_{1}}) = h'(\beta_{i_{1}}) = t$,
\item $E_{\beta_{i_{2}}}E_{\beta'_{i_{2}}} - q^{(\beta_{i_{2}},\beta'_{i_{2}})}E_{\beta'_{i_{2}}}E_{\beta_{i_{2}}} = k E_{\beta_{ex}}$ with $k \neq 0$. \hspace{2cm} $(\star)$ 
\end{itemizedot}
$\mbox{ }$\\
By definition of $\beta_{i_{2}}$, one has $\beta_{ex} = \beta_{i_{2}} + \beta'_{i_{2}}$, so that \\
$$\beta'_{i_{1}} + \beta'_{i_{2}} =\beta_{ex} + \alpha_{i}  = \beta_{i_{2}} + \beta'_{i_{2}} + \alpha_{i} \Rightarrow \beta'_{i_{1}} = \beta_{i_{2}} + \alpha_{i}.$$
Thus, by Corollary \ref{corplanadm}, we have $\ [E_{\beta_{i_{2}}},E_{\alpha_{i}}]_{q}:= h E_{\beta'_{i_{1}}} \ (h \neq 0).$
\par
We know that $\mathcal{U}_q^+(\mathfrak{g})$ is $\mathbb{Z}\Phi$-graded. So there is a (unique) automorphism $\sigma$ of $\mathcal{U}_q^+(\mathfrak{g})$ 
such that for all $u \in \mathcal{U}_q^+(\mathfrak{g})$, homogeneous in degree $\beta$, $\sigma(u) = q^{(\beta,\alpha_{i})} u$.\\
Denote by $\delta$ the interior right-sided $\sigma$-derivation associated to $E_{\alpha_{i}}$, so that \\
 $$\delta(u) = uE_{\alpha_{i}} - E_{\alpha_{i}}\sigma(u) \ \ \ \ \ (\forall u \in \mathcal{U}_q^+(\mathfrak{g}))$$
If $\beta \in C_{l}$, one has $\delta(E_{\beta}) = E_{\beta}E_{\alpha_{i}} - q^{(\beta,\alpha_{i})}E_{\alpha_{i}}E_{\beta} = [E_{\beta},E_{\alpha_{i}}]_{q}$ and, this implies:
$\delta(E_{\beta'_{i_{2}}}) = 0$ and  $\delta(E_{\beta_{i_{2}}}) = h E_{\beta'_{i_{1}}}$.
We can show with $(\star)$ that : 
\[\begin{array}{rl}
k[E_{\beta_{ex}},E_{\alpha_{i}}]_{q} & = k \delta(\beta_{ex}) = \delta(E_{\beta_{i_{2}}}E_{\beta'_{i_{2}}})-q^{(\beta_{i_{2}},\beta'_{i_{2}})}\delta(E_{\beta'_{i_{2}}}E_{\beta_{i_{2}}}) \\
& = E_{\beta_{i_{2}}}\delta(E_{\beta'_{i_{2}}}) + \delta(E_{\beta_{i_{2}}})\sigma(E_{\beta'_{i_{2}}}) - q^{(\beta_{i_{2}},\beta'_{i_{2}})}(E_{\beta'_{i_{2}}}\delta(E_{\beta_{i_{2}}}) + \delta(E_{\beta'_{i_{2}}})\sigma(E_{\beta_{i_{2}}}))\\
& = h [q^{(\beta'_{i_{2}},\alpha_{i})}E_{\beta'_{i_{1}}}E_{\beta'_{i_{2}}} - q^{(\beta_{i_{2}},\beta'_{i_{2}})}E_{\beta'_{i_{2}}}E_{\beta'_{i_{1}}}]
\end{array}
\]
As $\beta'_{i_{2}}$ and $\beta'_{i_{1}}$ are in the same box, we know (Corollary \ref{LSE}) that $E_{\beta'_{i_{1}}}E_{\beta'_{i_{2}}} = E_{\beta'_{i_{2}}}E_{\beta'_{i_{1}}}$, so that:
\[k[E_{\beta_{ex}},E_{\alpha_{i}}]_{q} = h (q^{(\beta'_{i_{2}},\alpha_{i})} - q^{(\beta_{i_{2}},\beta'_{i_{2}})}) E_{\beta'_{i_{1}}}E_{\beta'_{i_{2}}}\]
Since $\beta_{i_{2}} + \beta'_{i_{2}} = \beta_{ex}$,  $P =$ Vect($\beta_{i_{2}},\beta'_{i_{2}}$) is an admissible plane of type (1.1) or (1.2) (see Remark \ref{systracine}) with $\{\beta_{i_{2}}, \beta'_{i_{2}}\} = \{\beta, \beta' \}$, so that $(\beta_{i_{2}},\beta'_{i_{2}}) \leq 0$.  As we have assumed that $(\alpha_{i},\beta'_{i_{2}} ) > 0$, this implies that:
$$[E_{\beta_{ex}},E_{\alpha_{i}}]_{q}:=E_{\beta_{ex}}E_{\alpha_{i}} - q^{(\beta_{ex},\alpha_{i})}E_{\alpha_{i}}E_{\beta_{ex}} = \lambda E_{\beta'_{i_{1}}}E_{\beta'_{i_{2}}} \ \text{with} \ \lambda \neq 0.$$ 
\end{preuve}

\subsection{Implications from an admissible plane}
We define the notion of implications coming from an admissible plane P, and we verify that all Cauchon diagrams satisfy all implications from admissible planes. Let us begin by showing some precise results on the exceptional root and near boxes behaviour. First, let us recall some notation introduced in Sections 2 and 3.
\begin{nota}
\begin{itemizedot}
\item  $C_1, ..., C_n$ denote the columns of $\Phi^+$ (relative to the chosen Lusztig order).
\item In the following, we consider a diagram $\Delta$, that is, $\Delta$ a subset of $\Phi^+$.
\item For any integer $j \in \llbracket 1, n \rrbracket$, we set $\Delta_{j} := \Delta \cap C_j = \{ \beta_{u_{1}}, ..., \beta_{u_{l}} \} \subset C_{j} = \{\beta_{k},...,\beta_{r}\}$. If the column $C_j$ is exceptional, $\beta_{ex}$ denotes the exceptional root and $B_{ex} := \{\beta_{ex}\}$ is its box. Then $B_1$ denote the box of $C_{j}$ which precedes $B_{ex}$ and $B'_1$ the one which follows $B_{ex}$ in the Lusztig order ; so that $s_{D}(B_1) = B_1'$.
\end{itemizedot}
\end{nota}
In Propositions \ref{contrainte  elementaire}, \ref{contrainte elementaire avant ex} and \ref{contrainte elementaire apres ex}, we proved the existence of implications thanks to admissible planes. We formalise this fact in the following definition of "implications coming from an admissible plane" :
\begin{dfn}\label{contrainteplanadm}
Let $\beta \in C_j$ with $h'(\beta) = l$ and, P be an admissible plane.
\begin{enumerate}
\item If $\Phi_P^+ = \{\beta, \beta + \alpha_i, \alpha_i\}$ with $i < j$ \hyperref[systracine]{type (2.1)}, then the implication coming from P is $\beta \rightarrow \beta + \alpha_i$.
\item $\Phi_P^+ = \{\beta, \beta + \alpha_i, \beta + 2\alpha_i, \alpha_i\}$ with $i < j$ \hyperref[systracine]{type (2.3)}, then the implications coming from P are $\beta \rightarrow \beta + \alpha_i$ and $\beta + \alpha_i \rightarrow \beta + 2\alpha_i$.
\item $\Phi_P^+ = \{\beta, \beta + \beta', \beta\}$ with $i < j $,$\beta' \in C_j$ and $h(\beta') = h(\beta) +1$ \hyperref[systracine]{type (1.1)}, then the implication coming from P is $\beta \rightarrow \beta + \beta'$.
\item $\Phi_P^+ = \{\alpha_i , \alpha_i + \beta, \alpha_i + 2\beta, \beta\}$ with $i < j$, $h'(\alpha_i + 2\beta) = \frac{2l +1}{2}$ and $h(\beta) = l$ \hyperref[systracine]{type (1.2)} or \hyperref[systracine]{type (2.2)}, , then the implications coming from P are $\beta \rightarrow \alpha_i + \beta $, $\beta \rightarrow \alpha_i + 2\beta$ and $\alpha_i + 2\beta \rightarrow \alpha_i +\beta $.
\item $\Phi_P^+ = \{\beta, \alpha_i\}$ with $i < j$, $\alpha_i \bot \beta$ and there are $\beta_1$ and $\beta_2$ in $C_j$ such that $\beta + \alpha_i = \beta_1 + \beta_2$ \hyperref[systracine]{type (2.4)}, then the implications coming from P are \scalebox{.8}
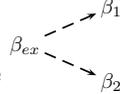
{\begin{psmatrix}[rowsep=.2cm, colsep=1cm]
& $\beta_{1}$\\
$\beta_{ex}$ & \\
& $\beta_{2}$
\end{psmatrix}
\psset{nodesep=2pt, linestyle=dashed, arrows=->}
\ncline{2,1}{1,2}
\ncline{2,1}{3,2}}\\ 
\item $\Phi_P^+ = \Phi^+ = \{\beta_1, ..., \beta_6\}$ is the positive part of a roots system of type $G_2$ (see Proposition \ref{propG2}), then the implications coming from P are $\beta_6 \rightarrow \beta_5, \ \beta_5 \rightarrow \beta_4, \ \beta_5 \rightarrow \beta_3, \ \beta_4 \rightarrow \beta_3, \ \beta_3 \rightarrow \beta_2$.
\end{enumerate}
\end{dfn}
\begin{lemme}\label{contraintedeP}
Let $\beta \in C_j$.
\begin{enumerate}
\item If $\beta$ belongs to a box which follows $\{\beta_{ex}\}$, then $\beta \rightarrow \beta_{ex}$ is an implication from an admissible plane.
\item If there is $i < j$ such that $\gamma = \beta + \alpha_i \in \Phi^+$ then $\beta \rightarrow \gamma$ is an implication from an admissible plane.
\end{enumerate}
\end{lemme}
\begin{preuve}
The results holds in the case where $\Phi$ is of type $G_2$. From now on, we assume that $\Phi$ is not of type $G_2$.
\begin{enumerate}
\item Let $P = <\beta, \beta_{ex}>$. It is an admissible plane of type 3 or 4 in the previous definition and in each case, $\beta \rightarrow \beta_{ex}$ is an implication coming from P.
\item Let $P = <\beta, \alpha_i>$. It is an admissible plane of type 1,2 or 4 in the previous definition and in each case, $\beta \rightarrow \gamma$ is an implication coming from P.
\end{enumerate}
\end{preuve}
\begin{prop}\label{direct}
Let $\Delta$ be a Cauchon diagram. Then $\Delta$ satisfies all the implication coming from admissible planes containing elements of $\Delta$.
\end{prop}
\begin{preuve}
Let $\beta \in \Delta$ and P be an admissible plane containing $\beta$. Recall (see Definition \ref{contrainteplanadm}) that $\Phi_P^+ = \Phi^+ \cap P$.
\begin{enumerate}
\item If $\Phi_P^+ = \{\beta, \beta + \alpha_i, \alpha_i\}$ with $i < j$, then it follows from Proposition \ref{contrainte elementaire} that $\Delta$ satisfies the implication $\beta \rightarrow \beta + \alpha_i$.
\item If $\Phi_P^+ = \{\beta, \beta + \alpha_i, \beta + 2\alpha_i, \alpha_i\}$ with $i < j$, then applying Proposition \ref{contrainte elementaire} to $\beta$ and $\beta + \alpha_i$, we get that $\Delta$ satisfies the implications $\beta \rightarrow \beta + \alpha_i$ and $\beta + \alpha_i \rightarrow \beta + 2\alpha_i$.
\item If $\Phi_P^+ = \{\beta, \beta + \beta', \beta\}$ with $i < j $, $\beta' \in C_j$ and $h(\beta') = h(\beta) +1$ then it follows from proposition \ref{contrainte elementaire avant ex} that $\Delta$ satisfies the implication $\beta \rightarrow \beta + \beta'$.
\item If $\Phi_P^+ = \{\alpha_i , \alpha_i + \beta, \alpha_i + 2\beta, \beta\}$ with $i < j$ and $h'(\alpha_i + 2\beta) = \frac{2l +1}{2}$, then it follows from Propositions \ref{contrainte elementaire}, \ref{contrainte elementaire avant ex} and \ref{contrainte elementaire apres ex} that $\Delta$ satisfies the implications $\beta \rightarrow \alpha_i + \beta $, $\beta \rightarrow \alpha_i + 2\beta$ and $\alpha_i + 2\beta \rightarrow \alpha_i +\beta $.
\item If $\Phi_P^+ = \{\beta, \alpha_i\}$ with $i < j$, $\alpha_i \bot \beta$ and there exist $\beta_1$ and $\beta_2$ in $C_j$ such that $\beta + \alpha_i = \beta_1 + \beta_2$, then it follows from proposition \ref{contrainte elementaire apres ex} that $\Delta$ satisfies the implication \scalebox{.8}
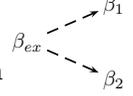
{\begin{psmatrix}[rowsep=.2cm, colsep=1cm]
& $\beta_{1}$\\
$\beta_{ex}$ & \\
& $\beta_{2}$
\end{psmatrix}
\psset{nodesep=2pt, linestyle=dashed, arrows=->}
\ncline{2,1}{1,2}
\ncline{2,1}{3,2}}.  \\ 
\item If $\Phi_P^+ = \Phi^+$ is of type $G_2$, Proposition \ref{contrainte elementaire} implies that $\Delta$ satisfies the implications  $\beta_6 \rightarrow \beta_5, \beta_5 \rightarrow \beta_3, \beta_4 \rightarrow \beta_3, \beta_3 \rightarrow \beta_2$. Moreover Proposition \ref{contrainte elementaire avant ex} implies that $\Delta$ satisfies the implication $\beta_5 \rightarrow \beta_4$.
\end{enumerate}
\end{preuve}

\subsection{The converse}
The goal of this section is to prove the converse of Proposition \ref{direct} , that is:
\begin{theo}\label{reciproque contrainte}
If $\Delta$ is a diagram which satisfies all the implications coming from admissible planes, then $\Delta$ is a Cauchon diagram.
\end{theo}
Let $\beta \in \Phi^+$ be a positive root of the column $C_j$. We denote by $B_{0}$ the box which contains $\beta$, by $B_{1}$ the box which precedes $B_{0}$ in the column $C_j$ (if it exists) and by $B_{2}$ the box which precedes $B_{1}$ in $C_j$ (if it exists).\\
Set $\Phi_\beta^+ = \{ \alpha_i \ | \ i < j\} \cup \{\gamma < \beta \ | \ \gamma$ is in the box of $\beta\} \cup B_1 \cup (B_2$ if $B_1 = \{\beta_{ex}\})$. If $\gamma \in \Phi^+$, then there exists $k \in \llbracket 1, N \rrbracket$ such that $\gamma = \beta_k$ and recall (see section 4.1) that $X_\gamma = X_k$.\\
Set $D_\beta := \mathbb{K}<X_\gamma \ | \ \gamma < \beta >$
\begin{lemme} \label{Dbeta}
\[D_\beta = \mathbb{K}<X_\gamma \ | \ \gamma \in \Phi_\beta^+>\]
\end{lemme}
\begin{preuve}
Set $D'_\beta := \mathbb{K}<X_\gamma \ | \ \gamma \in \Phi_\beta^+> \subset D_\beta $. Let us start by showing that, for $i < j$, we have $\{X_{\gamma}, \gamma \in C_{i}\} \subset D'_\beta$. If $\Phi$ is of type $G_2$, $\{X_{\gamma}, \gamma \in C_{i}\}$ is the empty set or it only contains $X_{\alpha_1} \in D'_\beta$. If $\Phi$ si not of type $G_2$, then we prove this result by induction on $h(\gamma)$. 
\begin{description}
\item[If $h(\gamma) = 1$ ,] then $\gamma = \alpha_{i}$ and $X_{\gamma} \in  D'_\beta$ by definition of $ \Phi_\beta^+$.
\item[If $h(\gamma) > 1$ and $\gamma$ ordinary,] then by Proposition \ref{prop2.10}, there exists $l < i$ such that $\gamma' = \gamma - \alpha_l \in \Phi^+$, so that, by Corollary \ref{corplanadm} and Proposition \ref{relcomX}, one has $X_\gamma \in \mathbb{K}<X_{\gamma'}, X_{\alpha_l}> \subset D'_\beta$ (by induction hypothesis).
\item[If $h(\gamma) > 1$ and $\gamma$ exceptional,] then we know (see Proposition \ref{prop2.8}) that in this case, there are two ordinary roots of $C_{i}$, denoted $\eta_1$ and $\eta_2$, such that $\eta_1+ \eta_2 = \gamma$ and $h(\eta_2) = h(\eta_1) + 1$. This implies by Corollary \ref{corplanadm} and Proposition \ref{relcomX} that $X_\gamma \in \mathbb{K}<X_{\eta_{1}}, X_{\eta_{2}}> \subset D'_\beta$ ($X_{\eta_{1}}$ and $X_{\eta_{2}}$ are in $D'_\beta$ because $\eta_1$ and $\eta_2$ are exceptional).
\end{description}
It just remains to show that $\{X_{\gamma} | \gamma \in C_{j}, \gamma < \beta \} \subset D'_\beta$.\\
If $h(\gamma) = h(\beta)$ with $\gamma < \beta$, then $\gamma \in \Phi_\beta^+$. So $X_{\gamma} \in D'_\beta$.\\
One uses again an induction to show that for each ordinary box $B$ of $C_{j}$ such that $B < B_0$ (i.e. all roots $\beta$ of $B$ are strictly less than all roots of $B_0$), one has $\{X_{\gamma} | \gamma \in B\} \subset D'_{\beta}$.
\begin{itemizedot}
\item \textbf{Assume that $B_1$ ordinary.}\\
The result is true for the box $B_{1}$ since $B_1 \subset \Phi_\beta^+$.\\
Let B be an ordinary root of $C_j$ such that $h(B) > h(B_{1})$ and $\gamma \in B$. By Proposition \ref{prop2.10}, there is $\alpha_{l} \in \Pi \ (l < j)$ such that $\gamma - \alpha_{l} \in \Phi^+$. Then $\gamma' := \gamma - \alpha_{l}$ is in an ordinary box $B'$ of $C_j$ such that $h(B) = h(B') + 1 > h(B') \geq h(B_{1}) > h(B_0)$ and one has $X_{\gamma'} \in D'_{\beta}$ by induction hypothesis.\\ 
If $\Phi$ is not of type $G_2$, then we deduce from orollary \ref{corplanadm} and Proposition \ref{relcomX} that $[X_{\gamma'},X_{\alpha_l}]_q = k X_\gamma$ with $k \in \mathbb{K}^\star$. As $X_{\alpha_{l}} \in D'_{\beta}$, this implies that $X_{\gamma} \in D'_{\beta}$.\\
If $\Phi$ is of type $G_2$, then we deduce from Propositions \ref{propG2} and \ref{relcomX} that $[X_{\gamma'},X_{\alpha_l}]_q = k X_\gamma$ with $k \in \mathbb{K}^\star$. As $X_{\alpha_{l}} \in D'_{\beta}$, this implies that $X_{\gamma} \in D'_{\beta}$.\\
\item \textbf{Assume that $B_1$ exceptional.}\\
The results is true for $B_{2}$ since, in this case, $B_2 \subset \Phi_\beta^+$.\\
This is the same proof as above with $B_{1}$ replaced by $B_{2}$.\\
\end{itemizedot}
It just remains to prove that if $B = \{\beta_{ex}\}$ is an exceptional box of $C_{j}$ such that  $B < B_{0}$, then one has $X_{\beta_{ex}} \in D'_{\beta}$.\\
If $B = B_1$, then one has $B \subset \Phi_\beta^+$, and the result is proved.\\
Assume that $B < B_1$. As above, one has $\beta_{ex} = \eta_{1} + \eta_{2}$ with $\eta_{1}$ and $\eta_{2}$ two exceptional roots of $C_j$ such that $h(\eta_{2}) = h(\eta_1) +1$. The boxes of $\eta_1$ and $\eta_2$ are ordinary, on each side of $B$, so less than or equal to $B_1$, so strictly less than $B_0$.  As the result holds for ordinary boxes, $X_{\eta_{1}} \in D'_{\beta}$ and $X_{\eta_{2}} \in D'_{\beta}$.\\
If $\Phi$ is not of type $G_2$, then we deduce (as above) from Corollary \ref{corplanadm} and Proposition \ref{relcomX} that $X_{\beta_{ex}} \in D'_{\beta}$.\\
If $\Phi$ is of type $G_2$, we deduce (as above) from Propositions \ref{propG2} and \ref{relcomX} that $X_{\beta_{ex}} \in D'_{\beta}$.\\
So we can conclude that $D_\beta = D'_\beta$.
\end{preuve}
Let us recall that $A = \mathcal{U}_q^+(\mathfrak{g}) = \mathbb{K}<X_{\beta_i} \ | \ i \in \llbracket 1 , N \rrbracket > := \mathbb{K}<X_i | i \in \llbracket 1 , N \rrbracket >$. Let $\beta_{r}$ and $\beta_{r+1}$ ($1 \leq r \leq N-1$) be two consecutive roots of $\Phi^+$ ($\beta_{r} < \beta_{r+1}$). Recall that $A^{(r+1)} = \mathbb{K}<X_{i}^{(r+1)}>$ and $A^{(r)} = \mathbb{K}<X_{i}^{(r)}>$ (1 < r < N) are the algebras deduced from $A$ by the deleting derivation algorithm of Section 4.
\begin{lemme} \label{relations}
Let  $\beta_r \in \Phi^+$ be a positive root of the column $C_j$ and $D_{\beta_r}^{(r+1)} := \mathbb{K}<X^{(r+1)}_\gamma \ | \ \gamma < \beta_r >$.
Then \[D^{(r+1)}_{\beta_r} = \mathbb{K}<X^{(r+1)}_\gamma \ | \ \gamma \in \Phi_{\beta_r}^+>.\]
\end{lemme}
\begin{preuve}
By Lemma \ref{4.1.1}, the commutation relations between the $X^{(r+1)}_\gamma$ with $\gamma \leq \beta_r$ are the same as the commutation relations between the $X_\gamma$ with $\gamma \leq \beta_r$. So the proof is the same as the proof of lemma \ref{Dbeta} but with $X_\gamma$ replaced by $X^{(r+1)}_\gamma$.
\end{preuve}
\textit{Denote, as in section 4, $\varphi : $ Spec $ A \hookrightarrow$ Spec $(\overline{A})$   $(\overline{A} = A^{(2)})$ the canonical injection, that is , the composition of canonical injections $\varphi_{r} : $ Spec $(A^{(r+1)}) \hookrightarrow$ Spec $(A^{(r)})$ for $r \in \llbracket 2, N \rrbracket$. Recall that a subset  $\Delta$ of $\Phi^+$ is a Cauchon diagram if and only if $(\exists P \in \text{Spec}(A))(\varphi (P) = <T_\gamma | \gamma \in \Delta>$).}\\ 
\begin{preuve4}
Let $\Delta \subset \Phi^+$ a diagram satisfying the implications coming from the admissible planes. Set $Q:= <T_\gamma | \gamma \in \Delta>$. By \cite[Section 5.5.]{MR1967309}, this is an $H^{(2)}$-prime ideal, so completely prime, of $A^{(2)} = \overline{A}$ and, if $\beta \in \Phi^+ \setminus \Delta$, then $T_\beta$ is regular modulo Q. So, $Q \cap \Phi^+ = \{ T_\gamma | \gamma \in \Delta \}$.\\
Let us show by induction, that for each $r \in \llbracket 2, N+1 \rrbracket $, there exists $P^{(r)} \in $ Spec $(A^{(r)})$ such that $Q = \varphi_{2} \circ ... \circ \varphi_{r-1} (P^{(r)})$ .\\
If $r= 2$, then in this case, one has $\varphi_{2} \circ ... \circ \varphi_{r-1} = Id_{\text{Spec}(\overline{A})}$ and $P^{(2)} = Q$.\\ 
\textbf{Consider an integer $r \in \llbracket 2, N \rrbracket $, assume that there exists $P^{(r)} \in \text{Spec} (A^{(r)})$ such that $\varphi_{2} \circ ... \circ \varphi_{r-1}(P^{(r)}) = Q$ and let us show there is $P^{(r+1)} \in \text{Spec} (A^{(r+1)})$ such that $\varphi_{r} (P^{(r+1)}) = P^{(r)}$( So that $\varphi_{2} \circ ... \circ \varphi_{r}(P^{(r+1)}) = Q$).}\\
$ $\\
\begin{itemizedot}
\item If $X_{r}^{(r)} \notin P^{(r)}$, then this follows from Proposition \ref{critereimage}.
\item Assume now that $X_{r}^{(r)} \in P^{(r)}$. From the second point of Proposition \ref{critereimage}, it is enough to show that 
$$\Theta^{(r)}\left(\delta_{r}^{(r+1)}(X_{i}^{(r+1)})\right) \in P\it^{(r)} \ \ \text{pour} \ 1  \leq i \leq  r-1.$$
\end{itemizedot}
\begin{obs}
It is enough to prove that $\Theta^{(r)}(\delta_{r}^{(r+1)}(X_{i}^{(r+1)})) $  $\in$ $ $ $P\it^{(r)}$ $ $ for $i \in \llbracket 1, r-1 \rrbracket$ such that $\beta_i \in \Phi_{\beta_r}^+$.
\end{obs}
\begin{preuve5}
Let $i \in \llbracket 1, r-1 \rrbracket$. It follows from Corollary \ref{relations} that:
\[ X_i^{(r+1)} = \sum_{j_1,...,j_s \in \Gamma} m_{i,r+1} X_{j_1}^{(r+1)}...X_{j_s}^{(r+1)} \ \ \text{où} \ \ \Gamma := \{ j \in \llbracket 1, r-1 \rrbracket | \beta_j \in \Phi_{\beta_r}^+ \}\]
Thus
\[\displaystyle \begin{array}{rl}
\delta_{r}^{(r+1)}(X_{i}^{(r+1)}) 	&= \sum  m_{i,r+1} \delta_{r}^{(r+1)}(X_{j_1}^{(r+1)} ... X_{j_s}^{(r+1)})\\
 & \\
 & \\
									&= \sum m_{i,r+1} \left[\delta_{r}^{(r+1)}\left(X_{j_1}^{(r+1)}\right)X_{j_2}^{(r+1)} ... \ X_{j_s}^{(r+1)} + \sigma_{r}^{(r+1)}\left(X_{j_1}^{(r+1)}\right)\delta_{r}^{(r+1)}\left(X_{j_2}^{(r+1)}\right) ...  \ X_{j_s}^{(r+1)}\right. \\
									& \\
									&\left. + \ ... + \sigma_{r}^{(r+1)}\left(X_{j_1}^{(r+1)} ... \ X_{j_{s-1}}^{(r+1)}\right)\delta_{r}^{(r+1)}\left(X_{j_s}^{(r+1)}\right)\right]\\
									& \\
									& \\
									&= \sum m_{i,r+1} \left[\delta_{r}^{(r+1)}\left(X_{j_1}^{(r+1)}\right)X_{j_2}^{(r+1)} ... \ X_{j_s}^{(r+1)} + \lambda_{r,j_1}X_{j_1}^{(r+1)}\delta_{r}^{(r+1)}\left(X_{j_2}^{(r+1)}\right) ... \ X_{j_s}^{(r+1)} \right.  \\
									& \\
									& \left. + \ ... + \lambda_{r,j_1} ... \lambda_{r,j_{s-1}}X_{j_1}^{(r+1)} ... \ X_{j_{s-1}}^{(r+1)}\delta_{r}^{(r+1)}\left(X_{j_s}^{(r+1)}\right) \right]
\end{array}\]
Then 
\[\begin{array}{rl}
\Theta^{(r)}\left(\delta_{r}^{(r+1)}(X_{i}^{(r+1)})\right) =
									& \sum m_{j_1, ..., j_s} \left[  \Theta^{(r)}\left(\delta_{r}^{(r+1)}(X_{j_1}^{(r+1)})\right)X_{j_2}^{(r)}... \ X_{j_s}^{(r)} \right.\\
									& \\
									& + \lambda_{r,j_1}X_{j_1}^{(r)}\Theta^{(r)}\left(\delta_{r}^{(r+1)}(X_{j_2}^{(r+1)})\right)... \ X_{j_s}^{(r)} +  ... \\
									& \\
									& \left. + \lambda_{r,j_1}...\lambda_{r,j_{s-1}}X_{j_1}^{(r)}... \ X_{j_{s-1}}^{(r)}\Theta^{(r)}\left(\delta_{r}^{(r+1)}(X_{j_s}^{(r+1)})\right)\right].
\end{array}\]
As each $\Theta^{(r)}\left(\delta_{r}^{(r+1)}(X_{j_l}^{(r+1)})\right) \in P^{(r)}$ by hypothesis, one has $\Theta^{(r)}\left(\delta_{r}^{(r+1)}(X_{i}^{(r+1)})\right) \in P^{(r)}$.
\end{preuve5}
\begin{tabular}{lc}
\begin{minipage}{.8\linewidth}
\bf{Back to the proof of the Theorem \ref{reciproque contrainte} :} \rm \\
For each $s \in \llbracket 2, r-1 \rrbracket$, set $P^{(s)} = \varphi_s \circ ... \varphi_{r-1}(P^{(r)})$.
\begin{obs}
$\beta_{r} \in \Delta$. 
\end{obs}
Indeed, as $X_{r}^{(r)} \in P^{(r)}$, Lemma \ref{lemme4.1} implies successively that $X_{r}^{(r-1)} \in P^{(r-1)}$, ..., $X_{r}^{(2)}\in P^{(2)} = Q$. Hence $T_{\beta_r} = X_{r}^{(2)} \in Q$ and so $\beta_{r} \in \Delta$.\\

Recall that, if $\beta_r \in C_j$, then $\Phi_{\beta_r}^+ = \{ \alpha_i \ | \ i < j\} \cup \{\gamma < \beta_r \ | \ \gamma \in B_0\} \cup B_1 \cup (B_2$ si $B_1 = \{\beta_{ex}\})$ ($B_0$ is the box containing $\beta_r$, $B_1$ is the box preceding $B_0$ if $C_j$ if exists and $B_2$ is the box preceding $B_1$ in $C_j$ if exists).
\end{minipage} 
&
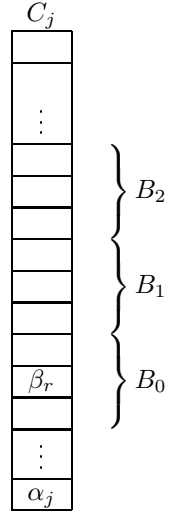
\begin{minipage}{.2\linewidth}
\hspace{.5cm} $C_j$\\
\begin{tabular}{l|c|c}
\cline{2-2}
&&\\
\cline{2-2}
&& \\
&$\vdots$ & \\
\cline{2-2}\cline{2-2}\cline{2-2}\cline{2-2}\cline{2-2}\cline{2-2}
 && \\
 \cline{2-2}
& & \\
 \cline{2-2}
 && \multirow{-3}{2cm}{ $\left. \begin{minipage}{0.1cm} $ $\\ $ $ \\ $ $ \\ \end{minipage} \right\} B_2$} \\
 \cline{2-2}\cline{2-2}\cline{2-2}\cline{2-2}\cline{2-2}\cline{2-2}\cline{2-2}
& & \\
 \cline{2-2}
 && \\
 \cline{2-2}
 && \multirow{-3}{2cm}{ $\left. \begin{minipage}{0.1cm} $ $\\ $ $ \\ $ $ \\ \end{minipage} \right\} B_1$} \\
 \cline{2-2}\cline{2-2}\cline{2-2}\cline{2-2}\cline{2-2}\cline{2-2}\cline{2-2}
 && \\
 \cline{2-2}
 &$\beta_{r}$& \\
 \cline{2-2}
 && \multirow{-3}{2cm}{ $\left. \begin{minipage}{0.1cm} $ $\\ $ $ \\ $ $ \\ \end{minipage} \right\} B_0$} \\
 \cline{2-2}\cline{2-2}\cline{2-2}\cline{2-2}\cline{2-2}\cline{2-2}\cline{2-2}
 &$\vdots $ \\
 \cline{2-2}
& $\alpha_j$ & \\
 \cline{2-2}
\end{tabular}
\end{minipage}
\end{tabular}\\
Let $i \in \llbracket 1, r-1 \rrbracket$ such that $\beta_i \in \Phi_{\beta_r}^+$.\\
\begin{itemizedot}
\item If $\beta_i \in B_0 \cup B_1$, then the Theorem \ref{1.4} implies that $\delta_{r}^{(r+1)}(X_{i}^{(r+1)}) = 0$. Hence $\Theta^{(r)}\left(\delta_{r}^{(r+1)}(X_{i}^{(r+1)})\right) = 0 \in P^{(r)}$.
\item Let us assume that $B_1 = \{\beta_{ex} \}$ with $\beta_{ex}= \beta_e$ ($e < r$), and that $\beta_i \in B_2$.\\ 
By Theorem \ref{1.4}, $\delta_{r}^{(r+1)}(X_{i}^{(r+1)}) = P_{r,i}^{(r+1)}$ is homogeneous of weight $\beta_r + \beta_i$ and the variables $X_l^{(r+1)}$ which appear in $P_{r,i}^{(r+1)}$ are such that $\beta_l \in B_1 = \{\beta_{e}\}$. So $P_{r,i}^{(r+1)}$ is equal to zero or is of the form $\lambda X_{e}^m$ with $\lambda \in \mathbb{K}^{\star}$ and $m\beta_{ex} = \beta_r + \beta_i$, so that (by comparing the coefficient on $\alpha_j$) one has $m = 1$. \\
If $P_{r,i}^{(r+1)} = 0$, then one has $\Theta^{(r)}\left(\delta_{r}^{(r+1)}(X_{i}^{(r+1)})\right) = 0 \in P^{(r)}$.\\
Otherwise, assume that $P_{r,i}^{(r+1)} = \lambda X_{e}^m$. As $\Delta$ satisfies the implications from admissible planes, Lemma \ref{contraintedeP} implies that $\Delta$ satisfies the implication $\beta_r \rightarrow \beta_{ex}$ and, as $\beta_r \in \Delta$, one has $\beta_{ex} \in \Delta$.\\
Then $X_e^{(2)} \in Q = P^{(2)}$ and by Lemma \ref{lemme4.1}, $X_e^{(e+1)} \in P^{(e+1)}$. As $\beta_{e}$ and $\beta_r$ are in consecutive boxes by construction, Lemma \ref{5.2} shows that 
\[X_e^{(e+1)} \in P^{(e+1)} \Rightarrow X_e^{(r)} \in P^{(r)}.\]
So, we deduce that $\Theta^{(r)}\left(\delta_{r}^{(r+1)}(X_{i}^{(r+1)})\right) = \Theta^{(r)}(\lambda  X_e^{(r+1)}) = \lambda X_e^{(r)} \in P^{(r)}$.
\item Consider now the case where $\beta_i = \alpha_k$ with $k<j$.\\ 
If $\delta_{r}^{(r+1)}(X_{i}^{(r+1)}) = 0$, then one has $\Theta^{(r)}\left(\delta_{r}^{(r+1)}(X_{i}^{(r+1)})\right) = 0 \in P^{(r)}$.\\
Assume that $\delta_{r}^{(r+1)}(X_{i}^{(r+1)}) \neq 0$. From Theorem \ref{1.4}, we get that
\[\delta_{r}^{(r+1)}(X_{i}^{(r+1)}) = \sum_{i < j_1 \leq ... \leq j_s < r} c_{j_1,...,j_s} X_{j_1}^{(r+1)}...X_{j_s}^{(r+1)} \ \ (c_{j_1,...,j_s} \in \mathbb{K})\]
\[c_{j_1,...,j_s} \in \mathbb{K}^* \Rightarrow \left( \beta_{j_1} + ... + \beta_{j_s} = \beta_r + \alpha_k \ \text{et} \ \beta_{j_1},...,\beta_{j_s} \notin B_0 \right) \]
This implies that
\[\Theta^{(r)}\left(\delta_{r}^{(r+1)}(X_{i}^{(r+1)})\right) = \sum_{i < j_1 \leq ... \leq j_s < r} c_{j_1,...,j_s} X_{j_1}^{(r)}...X_{j_s}^{(r)}\]
and that is enough to show that, if $c_{j_1,...,j_s} \in \mathbb{K}^*$, then one has $X_{j_1}^{(r)}...X_{j_s}^{(r)} \in P^{(r)}$.\\
\par 
So, take $(j_1,...j_s)$ such that $i < j_1 \leq ... \leq j_s < r$ and let us assume that $c_{j_1,...,j_s} \neq 0 $. Considering the coefficient of $\alpha_j$ in the following equality  
\begin{equation}\label{equaracine}
\beta_{j_1} + ... + \beta_{j_s} = \beta_r + \alpha_k,
\end{equation}
we deduce that $\beta_{j_s} \in C_j$. As $\beta_{j_s} \notin B_0$ and $j_s < r$, the box $B_1$ exists.\\
The proof splits into three cases.
\begin{itemizedot}
\item $B_0$ and $B_1$ are ordinaries. As $j_s < r$ and $\beta_{j_s} \notin B_0$, one has $h(\beta_r) < h(\beta_{j_s})$. By (\ref{equaracine}), $h(\beta_{j_s}) \leq h(\beta_r + \alpha_k) = h(\beta_r) + 1$. As a result, $s=1$ and $\beta_{j_s} \in B_1$. That is why, on has (Lemma \ref{contraintedeP}) the implication $\beta_r \rightarrow \beta_{j_s}$. Since $\beta_r \in \Delta$, one has $\beta_{j_s} \in \Delta$ and, as above, $X_{j_s}^{(j_s+1)} \in P^{(j_s+1)}$, so that $X_{j_s}^{(r)} \in P^{(r)}$. Hence the considered monomial whose coefficient $c_{j_1,...,j_s} \neq 0 $  is in $P^{(r)}$.\\
\item $B_0$ is ordinary and $B_1$ is exceptional so that $B_2$ exists. As in the previous case, one checks that $s =1$ and $\beta_{j_s} \in B_2$. So from Lemma \ref{contraintedeP}, there exists an implication $\beta_r \rightarrow \beta_{j_s}$. Also, From lemma \ref{contraintedeP}, one has the implication $\beta_r \rightarrow \beta_e$. Since $\beta_r \in \Delta$, one has $\beta_e, \beta_{j_s} \in \Delta$, so that $X_e^{(e+1)} \in P^{(e+1)}$ and $X_{j_s}^{(j_s+1)} \in P^{(j_s+1)}$. By the second point of Lemma \ref{5.2}, one deduces that $X_{j_s}^{(r)} \in P^{(r)}$. Thus the considered monomial is in $P^{(r)}$.\\
\item $B_0$ is exceptional. Since $\beta_{j_s} \notin B_0$, $\beta_{j_s}$ is ordinary in $C_j$. By the equality (\ref{equaracine}), one has $s \geq 2$ and $\beta_{j_{s-1}}$ is also ordinary in $C_j$. Set $h(\beta_r) := 2l +1 $ ($l \geq 1$). We knows that $h(\beta_{j_{s-1}}) \geq l+1$,  $h(\beta_{j_{s}}) \geq l+1$ and $h(\beta_r + \alpha_k) = 2l+2$. This implies that $s = 2$ and $\beta_{j_{s-1}}, \beta_{j_s} \in B_1$. The equality (\ref{equaracine}) can be then written as $\beta_r + \alpha_k = \beta_{j_{s-1}} + \beta_{j_s}$.
\begin{itemizedot}
\item Assume $\beta_{j_{s-1}} \neq \beta_{j_s}$, so that $\beta_{j_{s-1}}$ and $\beta_{j_s}$ are in the same box $B_1$, so they are orthogonal. As a result, $\Phi$ is not of the type $G_2$ (in the $G_2$ case, the boxes contain only one element). Set $P := <\beta_{j_s},\beta_{j_{s-1}}>$ the plane spanned by $\beta_{j_{s}}$, $\beta_{j_{s-1}}$, and assume $\Phi_P^+ \neq \{\beta_{j_{s-1}},\beta_{j_s}\}$. So, since $\Phi_P$ is not of type $G_2$, $\Phi_P$ is of type $A_2$ or $B_2$. As $\beta_{j_{s-1}}$ and $\beta_{j_s}$ are orthogonal, $\Phi_P$ is of type $B_2$ and there exists $\beta \in \Phi^+$ such that $\beta_r + \alpha_k = \beta_{j_{s-1}} + \beta_{j_s} = m \beta$ with $m = 1$ or $2$. \\
If $m = 1$, then $\beta$ and $\beta_r$ are two distinct exceptional roots of $C_j$, which is impossible.\\
Hence $m = 2$ and so $\beta_r + \alpha_k = \beta_{j_{s-1}} + \beta_{j_s} = 2 \beta$. This implies that $h(\beta) = l + 1$, so that $\beta$ is an element of $B_1$ too, different from $\beta_{j_{s-1}}$ and $\beta_{j_s} $. As a result, $\beta, \beta_{j_{s-1}}, \beta_{j_s}$ are pairwise orthogonal, which is a contradiction with the equality $\beta_{j_{s-1}} + \beta_{j_s} = 2 \beta$.\\
So one has $\Phi_P^+  =  \{\beta_{j_{s-1}},\beta_{j_s}\}$ and so we have the implication \scalebox{.8}{\begin{psmatrix}[rowsep=.2cm, colsep=1cm]
& $\beta_{j_{s-1}}$\\
$\beta_{r}$ & \\
& $\beta_{j_s}$
\end{psmatrix}
\psset{nodesep=2pt, linestyle=dashed, arrows=->}
\ncline{2,1}{1,2}
\ncline{2,1}{3,2}}. Hence one of the two roots $\beta_{j_s}$, $\beta_{j_{s-1}}$ is in $\Delta$. If, for example, $\beta_{j_s} \in \Delta$, one has, as in the first case, $X_{j_s}^{(j_s+1)} \in P^{(j_s+1)}$ and $X_{j_s}^{(r)} \in P^{(r)}$. The considered monomial is in $P^{(r)}$ as claimed.
\item If $\beta_{j_{s-1}} = \beta_{j_s}$, then the equality (\ref{equaracine}) becomes $\beta_r + \alpha_k = 2 \beta_{j_s}$. Set $\beta = s_D(\beta_{j_s}) = \beta_r - \beta_{j_s} \in \Phi^+$ and substract $\beta_{j_s}$ to each part of the previous equality, to obtain $\beta + \alpha_k = \beta_{j_s}$. Denote by $P$ the plane spanned by $\beta_r$ and $\beta_{j_s}$. \\
Assume that $\Phi$ is of type $G_2$. Then one has $\beta_r = \beta_4$, $\alpha_k = \beta_1$ and $\beta_{j_s} = \beta_3$. By Definition \ref{contrainteplanadm}, we have the implication $\beta_r \rightarrow \beta_{j_s}$.\\
Assume that $\Phi$ is not of type $G_2$. The equality $\beta_r + \alpha_k = 2 \beta_{j_s}$ implies that $\Phi_P$ is of type $B_2$, so that $\Phi_P^+ = \{\alpha_k, \alpha_k + \beta = \beta_{j_s}, \alpha_k + 2\beta = \beta_r, \beta \}$ with $h(\beta) = h(\beta_r) - h(\beta_{j_s}) = 2l+1 - (l+1) = l$. So $P$ is an admissible plane of type 4 in the sense of Definition \ref{contrainteplanadm}. So we have again the implication $\beta_r \rightarrow \beta_{j_s}$. \\
Thus, in all cases, one has $\beta_{j_s} \in \Delta$. So we have, as in the first case, $X_{j_s}^{(j_s+1)} \in P^{(j_s+1)}$ and $X_{j_s}^{(r)} \in P^{(r)}$. The considered monomial is again in $P^{(r)}$, as desired.
\end{itemizedot}
  
\end{itemizedot}
\end{itemizedot}
\end{preuve4}

\section{Cauchon diagrams for a particular decomposition of $w_{0}$.}
In this section, we give an explicit description of Cauchon diagrams for a chosen decomposition of $w_{0}$ in each type of simple Lie algebra of finite dimension. 
\subsection{Infinite series}
\subsubsection{Type $A_n$, $n \geq 1$}

\begin{convsanss}
The numbering of simple roots in the Dynkin diagram is as below : 
\begin{center}
	$ \alpha_1 - \alpha_2 - \cdots - \alpha_{n-1} - \alpha_n $.
\end{center}
We know (see for example \cite[Section 5]{MR1628449}) that $ s_{\alpha_1} \circ (s_{\alpha_2}\circ s_{\alpha_1})   \cdots \circ (s_{\alpha_n} \circ s_{\alpha_{n-1}} \circ \cdots \circ s_{\alpha_1} )$ is a reduced decomposition of  $w_0$ which induces the following order on positive roots. (We have arranged the roots in columns.)

\[\begin{array}{c|c|c|c|c|}
& C_{1} & C_{2} & & C_{n}\\
\cline{2-5}
 &   \beta_1 = \alpha_1 &  
 \beta_2 = \alpha_1 + \alpha_2  & \cdots & 
\beta_{N-n+1} = \alpha_1 + \cdots + \alpha_{n-1} + \alpha_n  \\
  \cline{2-5}
\multicolumn{2}{c|}{} & \multicolumn{1}{|c|}{
 \beta_3 = \alpha_2} & \multicolumn{1}{|c|}{\cdots} & \multicolumn{1}{|c|}{\vdots} \\
 \cline{3-5}
\multicolumn{3}{c|}{} & \multicolumn{1}{|c|}{\ddots} & \multicolumn{1}{|c|}{\vdots} \\
\cline{4-5}
\multicolumn{4}{c|}{} & \multicolumn{1}{|c|}{
\beta_N = \alpha_n}  \\
\cline{5-5}
\end{array}\]
\end{convsanss}
\par
This is a Lusztig order and none of the columns $C_{1}, ..., C_{n}$ is exceptional. Note that $|C_l| = l $. \\
\begin{lemme}\label{contraintesan}
If $1 \leq l \leq n$ and if $\beta_{j}, \beta_{j+1}$ are two consecutive roots of $C_{l}$, then $\beta_{j+1} \rightarrow \beta_{j} $ is an implication.
\end{lemme}
\begin{preuve}
If $\beta_{j}$ and $\beta_{j+1}$ are two consecutive roots of $C_{l}$, then there is a simple root $\alpha_i$ ($1 \leq i < l \leq n$) such that $\beta_{j} = \beta_{j+1} + \alpha_{i}$. From proposition \ref{contrainte elementaire}, $\beta \rightarrow \beta + \alpha_{i}$ is an implication.

\end{preuve}
It implies :
\begin{cor}\label{contraintecol}
If $C_{l} = \{\beta_{s}, \beta_{s+1},...,\beta_{r} = \alpha_{l}\}$ is the column $l$ with $1 \leq l \leq n$, then the implications coming from admissible planes are :
\[ \beta_{r} \rightarrow \beta_{r-1} \rightarrow ... \rightarrow \beta_{s+1} \rightarrow \beta_{s}\]
\end{cor}
\begin{convsanss}
\begin{itemize}
\item If $C_{l} = \{\beta_{s}, \beta_{s+1},...,\beta_{r} = \alpha_{l}\}$ is the column $l$ with $1 \leq l \leq n$, the truncated columns contained in $C_l$ are subset of the form $\{\beta_{s}, \beta_{s+1}, ..., \beta_{t}\}$, with $t \in \llbracket s, r \rrbracket$.
\item Denote by $\mathcal{D}$ the set of Cauchon diagrams, they satisfy the implications of Corollary \ref{contraintecol}.
\end{itemize}
\end{convsanss}
\begin{theo}
$\mathcal{D}$ is the set of all diagrams $\Delta$ which are unions of truncated columns.
\end{theo}
From now on, in order to represent a diagram $\Delta$, we fill the tabular used to describe the positive roots with black squares. The black squares represent the roots that belongs to the diagram. For example, the left-hand diagram below is $\Delta = \{\beta_1, \beta_4, \beta_5, \beta_6, \beta_7, \beta_8, \beta_{11}, \beta_{12}, \beta_{13}, \beta_{14} \}$.
\renewcommand{\PuzzleBlackBox}{\rule{.75\PuzzleUnitlength}%
{.75\PuzzleUnitlength}} 
\PuzzleUnsolved

\begin{center}
\begin{tabular}{cp{2cm}c}
\begin{Puzzle}{5}{5}%
|* |1 |*|* | * |.
|{} |1|*|*| * |.
|{}|{} |* | 1|* |.
|{} |{}|{} |1| * |.
|{} |{} |{}|{}|1 |.
\end{Puzzle}  
& &
\begin{Puzzle}{5}{5}%
|* |1 |*|* | * |.
|{} |1|* |*| * |.
|{}|{} |* | 1|* |.
|{} |{}|{} |*| * |.
|{} |{} |{}|{}|1 |.
\end{Puzzle}\\
$\Delta \in \mathcal{D} $& &
$\Delta \notin \mathcal{D}$
\end{tabular}
\end{center}
\begin{rmq}
The set of Cauchon diagrams $\mathcal{D}$ has the same cardinality as the Weyl group $W$.
\end{rmq}
\begin{preuve}
As $\mathcal{D}$ is the set of all diagrams $\Delta$ which are unions of truncated columns, one has $|\mathcal{D}| = (n+1)! = |W|$. 
\end{preuve}
\subsubsection{Type $B_n$, $n \geq 2$}
\begin{convsanss}
The numbering of simple roots in the Dynkin diagram is as below :
\[\alpha_1 \Leftarrow \alpha_2 - \cdots - \alpha_{n-1} - \alpha_n \]
We know (see for example \cite[Section 6]{MR1628449}) that \[s_{\alpha_1} \circ (s_{\alpha_2}\circ s_{\alpha_1}\circ s_{\alpha_2})   \cdots \circ (s_{\alpha_n} \circ s_{\alpha_{n-1}} \circ \cdots \circ s_{\alpha_2} \circ s_{\alpha_1}\circ s_{\alpha_2} \circ \cdots \circ s_{\alpha_n})\]
 is a reduced decomposition of $w_0$ which induces the following order on positive roots.
\end{convsanss}
\begin{tabular}{c|c|c|p{3cm}|c|}
\cline{5-5}
\multicolumn{4}{c|}{} & \multicolumn{1}{|c|}{$\beta_{(n-1)^2+1} = 2\alpha_1 + \cdots  + 2\alpha_{n-1} + \alpha_n$}  \\
\cline{4-5}
\multicolumn{3}{c|}{} & \multicolumn{1}{|c|}{} & \multicolumn{1}{|c|}{$\vdots$} \\
 \cline{3-5}
\multicolumn{2}{c|}{} & \multicolumn{1}{|c|}{$
 \beta_2 =2 \alpha_1 + \alpha_2$} & \multicolumn{1}{|c|}{} & \multicolumn{1}{|c|}{$\beta_{N-n} = 2\alpha_1 + \alpha_{2} + \cdots + \alpha_{n-1} + \alpha_n $} \\ 
\cline{2-5}
 &  $ \beta_1 = \alpha_1 $& $ \beta_3 = \alpha_1 + \alpha_2 $ &  & $
\beta_{N-n+1} = \alpha_1 + \cdots + \alpha_{n-1} + \alpha_n $ \\
  \cline{2-5}
\multicolumn{2}{c|}{} & \multicolumn{1}{|c|}{$
 \beta_4 = \alpha_2$} & \multicolumn{1}{|c|}{} & \multicolumn{1}{|c|}{$\beta_{N-n+2} = \alpha_2 + \cdots + \alpha_{n-1} + \alpha_n $} \\
 \cline{3-5}
\multicolumn{3}{c|}{} & \multicolumn{1}{|c|}{} & \multicolumn{1}{|c|}{$\vdots$} \\
\cline{4-5}
\multicolumn{4}{c|}{} & \multicolumn{1}{|c|}{$
\beta_N = \alpha_n$}  \\
\cline{5-5}
\end{tabular}\\
\par
This is a Lusztig order and none of the columns $C_{1}, ..., C_{n}$ is exceptional.\\
\begin{lemme}\label{contraintesbn}
If $1 \leq l \leq n$ and if $\beta_{j}, \beta_{j+1}$ are two consecutive roots of $C_{l}$, then $\beta_{j+1} \rightarrow \beta_{j} $ is an implication.
\end{lemme}
\begin{preuve}
If $\beta_{j}$ and $\beta_{j+1}$ are two consecutive roots of $C_{l}$, then there is a simple root $\alpha_i$ ($1 \leq i < l \leq n$) such that $\beta_{j} = \beta_{j+1} + \alpha_{i}$. From Proposition \ref{contrainte elementaire}, $\beta \rightarrow \beta + \alpha_{i}$ is an implication.
\end{preuve}
It implies :
\begin{cor}\label{contraintecolbn}
If $C_{l} = \{\beta_{s}, \beta_{s+1},...,\beta_{r} = \alpha_{l}\}$ is the column $l$ with $1 \leq l \leq n$, then the implications coming from admissible planes are :
\[ \beta_{r} \rightarrow \beta_{r-1} \rightarrow ... \rightarrow \beta_{s+1} \rightarrow \beta_{s}\]
\end{cor}
\begin{convsanss}
\begin{itemize}
\item If $C_{l} = \{\beta_{s}, \beta_{s+1},...,\beta_{r} = \alpha_{l}\}$ is the column $l$ with $1 \leq l \leq n$, the truncated columns contained in $C_l$ are subset of the form $\{\beta_{s}, \beta_{s+1}, ..., \beta_{t}\}$, with $t \in \llbracket s, r \rrbracket$.
\item Denote by $\mathcal{D}$ the set of Cauchon diagrams, they satisfy the implications of Corollary \ref{contraintecol}.
\end{itemize}
\end{convsanss}
\begin{theo}
$\mathcal{D}$ is the set of all diagrams $\Delta$ which are unions of truncated columns.
\end{theo}
\renewcommand{\PuzzleBlackBox}{\rule{.75\PuzzleUnitlength}%
{.75\PuzzleUnitlength}} 
\PuzzleUnsolved
\begin{center}
\begin{tabular}{cccc}
\begin{Puzzle}{7}{4}
|{}|{}|{}| *|.
|{}|{}| *| *| .
|{}| 1|* | 1| .
|*| 1| *|1 | . 
|{}| 1|* | 1| .
|{}|{}| 1| 1|.
|{}|{}|{}| 1| .
\end{Puzzle}  

& &

\begin{Puzzle}{7}{4}
|{}|{}|{}| 1| .
|{}|{}| *| 1| .
|{}| 1|* | 1|.
|*| 1| *|1 | . 
|{}| 1|* | *| .
|{}|{}| 1| *|.
|{}|{}|{}| 1|.
\end{Puzzle}  \\
&&&\\
&&&\\
&&&\\
&&&\\
&&&\\
 $\Delta \in \mathcal{D}$ & & $\Delta \notin \mathcal{D}$
\end{tabular}
\end{center}
\begin{rmq}
The set of Cauchon diagrams $\mathcal{D}$ has the same cardinality as the Weyl group $W$.
\end{rmq}
\begin{preuve}
As $\mathcal{D}$ is the set of all diagrams $\Delta$ which are unions of truncated columns, one has $|\mathcal{D}| = 2^{n+1}(n+1)! = |W|$. 
\end{preuve}

\subsubsection{Type $C_n$, $n \geq 3$}
\begin{convsanss}
The numbering of simple roots in the Dynkin diagram is as below :
\[\alpha_1 \Rightarrow \alpha_2 - \cdots - \alpha_{n-1} - \alpha_n  \]
We know (see for example \cite[Section 6]{MR1628449}) that
\[ s_{\alpha_1} \circ (s_{\alpha_2}\circ s_{\alpha_1}\circ s_{\alpha_2})   \cdots \circ (s_{\alpha_n} \circ s_{\alpha_{n-1}} \circ \cdots \circ s_{\alpha_2} \circ s_{\alpha_1}\circ s_{\alpha_2} \circ \cdots \circ s_{\alpha_n})\]
is a reduced decomposition of $w_0$ which induces the following order on positive roots.
\end{convsanss}
\begin{tabular}{c|c|c|p{3cm}|c|}
\cline{5-5}
\multicolumn{4}{c|}{} & \multicolumn{1}{|c|}{$
\beta_{(n-1)^2+1} = \alpha_1 + 2\alpha_2 \cdots  + 2\alpha_{n-1} + \alpha_n$}  \\
\cline{4-5}
\multicolumn{3}{c|}{} & \multicolumn{1}{|c|}{} & \multicolumn{1}{|c|}{$\vdots$} \\
 \cline{3-5}
\multicolumn{2}{c|}{} & \multicolumn{1}{|c|}{$ \beta_2 = \alpha_1 + \alpha_2 $} & \multicolumn{1}{|c|}{} & \multicolumn{1}{|c|}{$\beta_{N-n} = \alpha_1 + \alpha_{2} + \cdots + \alpha_{n-1} + \alpha_n $} \\ 
\cline{2-5}
 &  $ \beta_1 = \alpha_1$& $ \beta_3 = \alpha_1 + 2\alpha_2 $ &  & $
\beta_{N-n+1} = \alpha_1 +2\alpha_2 \cdots + 2\alpha_{n-1} + 2\alpha_n$ \\
  \cline{2-5}
\multicolumn{2}{c|}{} & \multicolumn{1}{|c|}{$
 \beta_4 = \alpha_2 $} & \multicolumn{1}{|c|}{} & \multicolumn{1}{|c|}{$\beta_{N-n+2} = \alpha_2 + \cdots + \alpha_{n-1} + \alpha_n$} \\
 \cline{3-5}
\multicolumn{3}{c|}{} & \multicolumn{1}{|c|}{} & \multicolumn{1}{|c|}{$\vdots$} \\
\cline{4-5}
\multicolumn{4}{c|}{} & \multicolumn{1}{|c|}{$\beta_N = \alpha_n$}  \\
\cline{5-5}
\end{tabular}\\

\par
This is a Lusztig order and all the columns  $C_{2}, ..., C_{n}$ are exceptional, the first one $C_1$ is ordinary.\\
\begin{lemme}\label{contraintescn}
If $1 \leq l \leq n$ and if $\beta_{j}, \beta_{j+1}$ are two consecutive roots of $C_{l}$, then $\beta_{j+1} \rightarrow \beta_{j} $ is an implication.
\end{lemme}
\begin{preuve}
Each column contains an odd number of elements, denoted by $C_{l} = \{\beta_{u_{1}},...,\beta_{u_{2k+1}}\}$, $\beta_{u_{k+1}}$ is an exceptional root and all boxes contain a single root.\\ 
\begin{itemizedot}
\item Let $j \in \llbracket 1, k-1 \rrbracket \cup \llbracket k+1, 2k+1 \rrbracket$, $\beta_{u_{j}}$ and $\beta_{u_{j+1}}$ be two consecutive roots of $C_{l}$. Then there is a simple root $\alpha_{i_{j}}$ ($1 \leq i_{j} < l \leq n$) such that $\beta_{u_{j}} = \beta_{u_{j+1}} + \alpha_{i_{j}}$. From Proposition \ref{contrainte elementaire}, $\beta_{u_{j+1}} \rightarrow \beta_{u_{j}}$ is an implication.
\item $\beta_{u_{k+2}}$ is in the box following the exceptional root, from Proposition \ref{contrainte elementaire avant ex}, $\beta_{u_{k+2}} \rightarrow \beta_{u_{k+1}}$ is an implication.
\item As $\beta_{u_{k+1}} + \alpha_{1} = 2 \beta_{u_{k-1}}$ so from Proposition \ref{contrainte elementaire apres ex}, $\beta_{u_{k+1}} \rightarrow \beta_{u_{k}}$ is an implication.
\end{itemizedot}
\end{preuve}
It implies :
\begin{cor}\label{contraintecolcn}
If $C_{l} = \{\beta_{s}, \beta_{s+1},...,\beta_{r} = \alpha_{l}\}$ is in the column $l$ with $1 \leq l \leq n$, the implications coming from admissible planes are :
\[ \beta_{r} \rightarrow \beta_{r-1} \rightarrow ... \rightarrow \beta_{s+1} \rightarrow \beta_{s}\]
\end{cor}
\begin{convsanss}
\begin{itemize}
\item If  $C_{l} = \{\beta_{s}, \beta_{s+1},...,\beta_{r} = \alpha_{l}\}$ is the column $l$ with $1 \leq l \leq n$, truncated columns contained in $C_l$ are subset of the following shape $\{\beta_{s}, \beta_{s+1}, ..., \beta_{t}\}$, with $t \in \llbracket s, r \rrbracket$.
\item Denote by $\mathcal{D}$ the set of Cauchon diagrams, they satisfy the implications of Corollary \ref{contraintecolcn}.
\end{itemize}
\end{convsanss}
\begin{theo}
$\mathcal{D}$ is the set of all diagrams $\Delta$ which are unions of truncated columns.
\end{theo}
\renewcommand{\PuzzleBlackBox}{\rule{.75\PuzzleUnitlength}%
{.75\PuzzleUnitlength}} 
\PuzzleUnsolved
\begin{center}
\begin{tabular}{cccc}
\begin{Puzzle}{7}{4}
|{}|{}|{}| *|.
|{}|{}| *| *| .
|{}| 1|* | 1| .
|*| 1| *|1 | . 
|{}| 1|* | 1| .
|{}|{}| 1| 1|.
|{}|{}|{}| 1| .
\end{Puzzle}  

& &

\begin{Puzzle}{7}{4}
|{}|{}|{}| 1| .
|{}|{}| *| 1| .
|{}| 1|* | 1|.
|*| 1| *|1 | . 
|{}| 1|* | *| .
|{}|{}| 1| *|.
|{}|{}|{}| 1|.
\end{Puzzle}  \\
&&&\\
&&&\\
&&&\\
&&&\\
&&&\\
 $\Delta \in \mathcal{D}$ & & $\Delta \notin \mathcal{D}$
\end{tabular}
\end{center}
\begin{rmq}
The set of Cauchon diagrams $\mathcal{D}$ has the same cardinality as the Weyl group $W$.
\end{rmq}
\begin{preuve}
As $\mathcal{D}$ is the set of all diagrams $\Delta$ which are unions of truncated columns, one has $|\mathcal{D}| = 2^{n+1}(n+1)! = |W|$.
\end{preuve}
\subsubsection{Type $D_n$, $n \geq 4$}
\begin{convsanss}
The numbering of simple roots in the Dynkin diagram is as below :
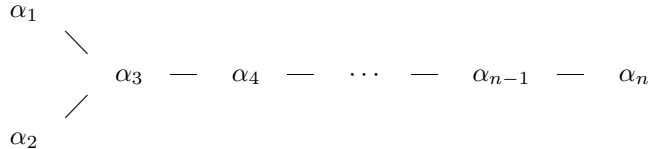
\[\begin{array}{ccccccccccc}
\alpha_1 && & &   & &&&  &     &   \\ 
& \diagdown & & &   & &&&  &     & \\
&& \alpha_3 &\text{---~}& \alpha_{4} &\text{---~}& \cdots &\text{---~}& \alpha_{n-1}&\text{---~} & \alpha_n \\
& \diagup & & &   & &&&  &     & \\
\alpha_2&& & &  &  &&&&  &
\end{array}\]
We know (see for example \cite[Section 6]{MR1628449}) that :
 \[s_{\alpha_1} \circ s_{\alpha_2} \circ (s_{\alpha_3}\circ s_{\alpha_1} \circ s_{\alpha_2} \circ s_{\alpha_3})   \cdots \circ (s_{\alpha_n} \circ s_{\alpha_{n-1}} \circ \cdots \circ s_{\alpha_3} \circ s_{\alpha_1} \circ s_{\alpha_2}\circ s_{\alpha_3} \circ \cdots \circ s_{\alpha_n}) \]
 is a reduced decomposition of $w_0$ which induces the following order on positive roots.\\
\begin{tabular}{c|c|c|c|c|}
\cline{5-5}
\multicolumn{4}{c|}{} & \multicolumn{1}{|c|}{$\beta_{N-2n+1} = \alpha_1 + \alpha_2 + 2 \alpha_3 \cdots  + 2\alpha_{n-1} + \alpha_n  $}  \\
\cline{4-5}
\multicolumn{3}{c|}{} & \multicolumn{1}{|c|}{} & \multicolumn{1}{|c|}{$\vdots$} \\
 \cline{3-5}
\multicolumn{2}{c|}{} & \multicolumn{1}{|c|}{$
 \beta_3 = \alpha_1 + \alpha_2 + \alpha_3   $} & \multicolumn{1}{|c|}{} & \multicolumn{1}{|c|}{$
\beta_{N-n-1} = \alpha_1 + \alpha_2 + \alpha_3 \cdots + \alpha_{n-1} + \alpha_n$} \\ 
\cline{2-5}
 &  $
 \beta_1 = \alpha_1 $& $  \beta_4 = \alpha_2 + \alpha_3 $ &  & $
\beta_{N-n} = \alpha_1 \ (\textrm{or} \ \alpha_2) \ + \alpha_3 \cdots + \alpha_{n-1} + \alpha_n$\\

\cline{2-5}
 &   $\beta_2 = \alpha_2$& $ \beta_5 = \alpha_1 + \alpha_3 $ &  & $
\beta_{N-n+1} = \alpha_2 \ (\textrm{or} \ \alpha_1) \ + \alpha_3 \cdots + \alpha_{n-1} + \alpha_n$ \\
  \cline{2-5}
\multicolumn{2}{c|}{} & \multicolumn{1}{|c|}{$
 \beta_6 = \alpha_3 $} & \multicolumn{1}{|c|}{} & \multicolumn{1}{|c|}{$
\beta_{N-n+2} = \alpha_3 \cdots + \alpha_{n-1} + \alpha_n$} \\
 \cline{3-5}
\multicolumn{3}{c|}{} & \multicolumn{1}{|c|}{} & \multicolumn{1}{|c|}{$\vdots$} \\
\cline{4-5}
\multicolumn{4}{c|}{} & \multicolumn{1}{|c|}{$
\beta_N = \alpha_n$}  \\
\cline{5-5}
\end{tabular}\\
The two first columns has been putted together to have a better view of roots.
\end{convsanss}
This is a Lusztig order and all the columns are ordinary.
\begin{lemme}\label{contraintesdn}
Let $C_l = \{\beta_{u_{1}},...,\beta_{u_{s}},\beta_{u_{s+1}}, ... ,\beta_{u_{2s}} \}$ (with $s = l-1$, $l \geq 3$), one has the implications :\\
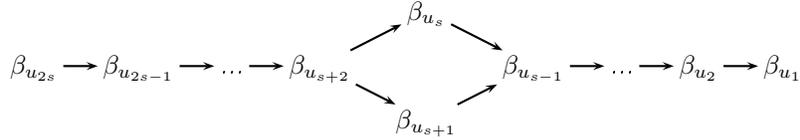
\begin{center}
\begin{psmatrix}[rowsep=.3cm, colsep=.6cm]
 				&				&		&			&$\beta_{u_{s}}$ \\ 
$\beta_{u_{2s}}$ &$ \beta_{u_{2s-1}} $& ... &$\beta_{u_{s+2}}$& 			&	$\beta_{u_{s-1}}$ &  ...& $\beta_{u_{2}}$  &$\beta_{u_{1}}$ \\
				& 						 &		&             &  $\beta_{u_{s+1}}$ 
\end{psmatrix}
\psset{arrows=->, nodesep=2pt}
\ncline{2,1}{2,2}
\ncline{2,2}{2,3}
\ncline{2,3}{2,4}
\ncline{2,4}{1,5}
\ncline{2,4}{3,5}
\ncline{1,5}{2,6}
\ncline{3,5}{2,6}
\ncline{2,6}{2,7}
\ncline{2,7}{2,8}
\ncline{2,8}{2,9}
\end{center}
\end{lemme}
\begin{preuve}
Observe first that the box $B = \{\beta_{u_{s}}, \beta_{u_{s+1}}\}$ is the only one which contains more than one element.
\begin{itemizedot}
\item For $j \in \llbracket 1, s-2 \rrbracket$, one has $\beta_{u_{j+1}} + \alpha_{s-j} = \beta_{u_{j}}$ and $\beta_{u_{2s-j+1}} + \alpha_{s-j+1} = \beta_{u_{2s-j}}$. Thus by Proposition \ref{contrainte elementaire}, one has the implications $\beta_{u_{j+1}} \rightarrow \beta_{u_{j}}$ et $\beta_{u_{2s-j+1}}\rightarrow \beta_{u_{2s-j}}$
\item One has $\beta_{u_{s+2}} + \alpha_{1} \in B$ and $\beta_{u_{s+2}} + \alpha_{2} \in B$ that is why from Proposition \ref{contrainte elementaire}, we obtain the implications $\beta_{u_{s+2}} \rightarrow \beta_{u_{s+1}}$ and $\beta_{u_{s+2}}\rightarrow \beta_{u_{s}}$.
\item One has $\beta_{u_{s-1}} - \alpha_{1} \in B$ and $\beta_{u_{s-1}} - \alpha_{2} \in B$ that is why from Proposition \ref{contrainte elementaire}, we obtain the implications $\beta_{u_{s+1}} \rightarrow \beta_{u_{s-1}}$ and $\beta_{u_{s}}\rightarrow \beta_{u_{s-1}}$.
\end{itemizedot}
\end{preuve}
\begin{convsanss}
\begin{itemize}
\item If  $C_{l} = \{\beta_{s}, \beta_{s+1},...,\beta_{r} = \alpha_{l}\}$ is the column $l$ with $1 \leq l \leq n$, the truncated columns with a hole in position $\beta_{m}$ ($s \leq m \leq r$) contained $C_l$ are subset of the shape $\{\beta_{s}, \beta_{s+1}, ...,\beta_{m-1},\beta_{m+1}, \beta_{t}\}$, $t \in \llbracket s, r \rrbracket$.
\item $\mathcal{D}$ is the set of Cauchon diagrams, they satisfy implications from Lemma \ref{contraintesdn} pour $1 \leq l \leq n$.
\end{itemize}
\end{convsanss}
\begin{theo}
$\mathcal{D}$ is the set of truncated columns and truncated columns with a hole in the line of $\alpha_1$ when this truncated column doesn't contains roots above the line of $\alpha_2$.  \\
\end{theo}
Thus :
\begin{center}
\begin{tabular}{cp{2cm}c}
\begin{Puzzle}{8}{4}
|{}|{}| *| .
|{}| *|* | .
|1| *| 1|. 
|*| *| *|. 
|{}| *|1 |.
|{}|{}| 1|.
\end{Puzzle}  
& &
\begin{Puzzle}{8}{4}
|{}|{}| *| .
|{}| *|* |.
|1| *| 1|.
|*| 1| *|.  
|{}| *|* | .
|{}|{}| *| .
\end{Puzzle}  \\
&&\\
&&\\
&& \\
 $\Delta \in \mathcal{D}$  & &$\Delta \notin \mathcal{D}$
\end{tabular}\\
\end{center}
\begin{rmq}
The set of Cauchon diagrams $\mathcal{D}$ has the same cardinality as the Weyl group $W$.
\end{rmq}
\begin{preuve}
As $\mathcal{D}$ is the set of truncated columns with or without a hole in a fixed position, we can compute its cardinality : 
\[|\mathcal{D}| = 4 \times 6 \times 8 \times ... \times 2n = 2^{n-1}n! = |W|.\]
\end{preuve}
\subsection{Exceptional cases}
\subsubsection{Type $G_2$}
\begin{convsanss}
The numbering of simple roots in the Dynkin diagram is as below :
\begin{center}
	$ \alpha_1 \Lleftarrow \alpha_2 $
\end{center}
We know that  $ s_{\alpha_1} \circ s_{\alpha_2} \circ s_{\alpha_1} \circ s_{\alpha_2} \circ s_{\alpha_1} \circ s_{\alpha_2} $ is a reduced decomposition of $w_0$ which induces the following order on positive roots.
\begin{center}
\begin{tabular}{c|c|c|}
\cline{3-3}
\multicolumn{2}{c|}{} & \multicolumn{1}{|c|}{$\beta_2 = 3\alpha_1 + \alpha_2$}\\
\cline{3-3}
\multicolumn{2}{c|}{} & \multicolumn{1}{|c|}{$\beta_3 = 2\alpha_1 + \alpha_2$} \\
\cline{2-3}
 &  $ \beta_1 = \alpha_1$ & $ \beta_4 = 3\alpha_1 + 2\alpha_2$ \\
 \cline{2-3}
\multicolumn{2}{c|}{} & \multicolumn{1}{|c|}{$\beta_5 = \alpha_1 + \alpha_2$}\\
 \cline{3-3}
\multicolumn{2}{c|}{} & \multicolumn{1}{|c|}{$\beta_6 = \alpha_2$}\\
\cline{3-3}
\end{tabular}
\end{center}
\end{convsanss}
\begin{lemme}
One has the following implications :
\scalebox{.8}{\begin{psmatrix}[rowsep=.2cm, colsep=1cm]
$\beta_{6}$ & $\beta_{5}$& $\beta_{4}$& $\beta_{3}$& $\beta_{2}$
\end{psmatrix}
\psset{nodesep=2pt, arrows=->}
\ncline{1,1}{1,2}
\ncline{1,2}{1,3}
\ncline{1,3}{1,4}
\ncline{1,4}{1,5}}
\end{lemme}
\begin{preuve}\label{contraintecolg2}
To prove this implications, we apply Propositions \ref{contrainte elementaire}, \ref{contrainte elementaire avant ex} and \ref{contrainte elementaire apres ex} with the following equalities ($\beta_{4}$ is an exceptional root):
\[\beta_{6} + \alpha_{1} = \beta_{5},\  h'(\beta_{5}) + 1 = \beta_{4},\  \beta_{4} + \alpha_{1} = 2 \beta_{3},\hspace{.2cm} \beta_{3} + \alpha_{1} = \beta_{2}.\]
\end{preuve}
\begin{convsanss}
$\mathcal{D}$ is the set of Cauchon diagrams, they satisfy implications from lemma \ref{contraintecolg2} .
\end{convsanss}
\begin{rmq}
The set of Cauchon diagrams $\mathcal{D}$ has the same cardinality as the Weyl group $W$.
\end{rmq}
\begin{preuve}
\[|\mathcal{D}| = 2 \times 6 = 12 = |W|.\]
\end{preuve}
\subsubsection{Type $F_4$}
\begin{convsanss}
The numbering of simple roots in the Dynkin diagram is as below : 
	\[\alpha_1 \ \text{---~} \alpha_2 \Rightarrow \alpha_3 \ \text{---~} \alpha_4\]
We choose the following reduced decomposition of $w_0$ :
 $$w_0 = s_4s_3s_4s_2s_3s_4s_2s_3s_2s_1s_2s_3s_4s_2s_3s_1s_2 s_1 s_3 s_4 s_2s_3s_2s_1 $$
This decomposition induces the following order on positive roots :\\
\begin{tabular}{cc}
\begin{tabular}{c|c|c|c|c|}
\cline{5-5}
\multicolumn{4}{c|}{} & \multicolumn{1}{|c|}{$
\beta_{10} (1, 3, 4, 2)$}\\
\cline{5-5}
\multicolumn{4}{c|}{} & \multicolumn{1}{|c|}{$
\beta_{11} (1, 2, 4, 2) $}\\
\cline{5-5}
\multicolumn{4}{c|}{} & \multicolumn{1}{|c|}{$
\beta_{12} (1, 2, 3, 2)$}\\
\cline{5-5}
\multicolumn{4}{c|}{} & \multicolumn{1}{|c|}{$
\beta_{13} (1, 2, 3, 1)$}\\
\cline{4-5}
\multicolumn{3}{c|}{} & \multicolumn{1}{|c|}{$\beta_4(0, 1, 2, 2 )$}&\multicolumn{1}{|c|}{$
\beta_{14} (1, 2, 2, 2)$}\\
\cline{4-5}
\multicolumn{3}{c|}{} & \multicolumn{1}{|c|}{$\beta_5( 0, 1, 2, 1)$} & \multicolumn{1}{|c|}{$\beta_{15} (1, 2, 2, 1)$} \\
 \cline{3-5}
\multicolumn{2}{c|}{} & \multicolumn{1}{|c|}{$\beta_2 (0, 0, 1, 1)$} & \multicolumn{1}{|c|}{$\beta_6 (0, 1, 1, 1)$} & \multicolumn{1}{|c|}{$\beta_{16} (1, 1, 2, 2)$} \\ 
\cline{2-5}
 &  $\beta_1 ( 0, 0, 0, 1) $& $ \beta_3 (0, 0, 1, 0)$  & $\beta_7 (0, 1, 2, 0) $& $
\beta_{17} (2, 3, 4, 2)$ \\
  \cline{2-5}
\multicolumn{3}{c|}{} & \multicolumn{1}{|c|}{$\beta_8 (0, 1, 1, 0)$} & \multicolumn{1}{|c|}{$\beta_{18} (1, 2, 2, 0)$} \\
 \cline{4-5}
\multicolumn{3}{c|}{} & \multicolumn{1}{|c|}{$\beta_9 (0, 1, 0, 0)$} & \multicolumn{1}{|c|}{$\beta_{19} (1, 1, 2, 1)$} \\
\cline{4-5}
\multicolumn{4}{c|}{} & \multicolumn{1}{|c|}{$\beta_{20} (1, 1, 1, 1)$}  \\
\cline{5-5}
\multicolumn{4}{c|}{} & \multicolumn{1}{|c|}{$\beta_{21} (1, 1, 2, 0)$}  \\
\cline{5-5}
\multicolumn{4}{c|}{} & \multicolumn{1}{|c|}{$\beta_{22} (1, 1, 1, 0)$}  \\
\cline{5-5}
\multicolumn{4}{c|}{} & \multicolumn{1}{|c|}{$\beta_{23} (1, 1, 0, 0)$}  \\
\cline{5-5}
\multicolumn{4}{c|}{} & \multicolumn{1}{|c|}{$\beta_{24} (1, 0, 0, 0)$}  \\
\cline{5-5}
\end{tabular}
&
\begin{tabular}{c|c|c|c|c|}
\cline{5-5}
\multicolumn{4}{c|}{} & \multicolumn{1}{|c|}{10}\\
\cline{5-5}
\multicolumn{4}{c|}{} & \multicolumn{1}{|c|}{9}\\
\cline{5-5}
\multicolumn{4}{c|}{} & \multicolumn{1}{|c|}{8}\\
\cline{5-5}
\multicolumn{4}{c|}{} & \multicolumn{1}{|c|}{7}\\
\cline{4-5}
\multicolumn{3}{c|}{} & \multicolumn{1}{|c|}{5}&\multicolumn{1}{|c|}{7}\\
\cline{4-5}
\multicolumn{3}{c|}{} & \multicolumn{1}{|c|}{4} & \multicolumn{1}{|c|}{6} \\
 \cline{3-5}
\multicolumn{2}{c|}{} & \multicolumn{1}{|c|}{2} & \multicolumn{1}{|c|}{3} & \multicolumn{1}{|c|}{6} \\ 
\cline{2-5}
 &  1&  1  & 3 & 11/2 \\
  \cline{2-5}
\multicolumn{3}{c|}{} & \multicolumn{1}{|c|}{2} & \multicolumn{1}{|c|}{5} \\
 \cline{4-5}
\multicolumn{3}{c|}{} & \multicolumn{1}{|c|}{1} & \multicolumn{1}{|c|}{5} \\
\cline{4-5}
\multicolumn{4}{c|}{} & \multicolumn{1}{|c|}{4}  \\
\cline{5-5}
\multicolumn{4}{c|}{} & \multicolumn{1}{|c|}{4}  \\
\cline{5-5}
\multicolumn{4}{c|}{} & \multicolumn{1}{|c|}{3}  \\
\cline{5-5}
\multicolumn{4}{c|}{} & \multicolumn{1}{|c|}{2}  \\
\cline{5-5}
\multicolumn{4}{c|}{} & \multicolumn{1}{|c|}{1}  \\
\cline{5-5}
\end{tabular}
\end{tabular}\\
\end{convsanss}
One checks that each column is ordinary or exceptional and then computes $h'(\beta_i)$ for all roots to verify that the order is a Lusztig one. We already know the form of diagrams for the two first columns. Thanks to commutation relations and Proposition \ref{contrainte elementaire}, we obtain the following implications for penultimate column :
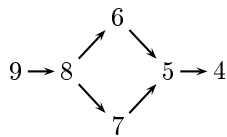
\begin{center}
\begin{psmatrix}[rowsep=.3cm, colsep=.5cm]
&  & 6 \\
 9 & 8 & & 5& 4 \\ 
 & & 7  \\  
\end{psmatrix}
\psset{arrows=->, nodesep=2pt}
\ncline{2,1}{2,2}
\ncline{2,2}{1,3}
\ncline{2,2}{3,3}
\ncline{1,3}{2,4}
\ncline{3,3}{2,4}
\ncline{2,4}{2,5}
\end{center}
So there are 8 possibilities to construct a diagrams satisfying the implications.\\

From Propositions \ref{contrainte elementaire}, \ref{contrainte elementaire avant ex} and \ref{contrainte elementaire apres ex}, one obtains the following implications for the last column: \\
\begin{tabular}{lc}
\begin{tabular}{c|c}
Diagrams beginning with case : & Counting diagrams :\\
\hline
none & 1\\
10					& 1	\\
11					& 1\\
12					& 1\\
13					& 1\\
14				& 2\\
15					& 1\\
16					& 3\\
17					& 2\\
18					& 2\\
19				& 2\\
20					& 2\\
21					& 2\\
22					& 1\\
23					& 1\\
24					& 1\\
\hline
 & \\
TOTAL 			& 24 
\end{tabular}
&
\begin{minipage}{3cm}
\begin{psmatrix}[rowsep=.3cm, colsep=.2cm]
 & 24 \\ 
 & 23 \\ 
 & 22 \\ 
20 &  & 21 \\ 
19 & &18 \\ 
&   17 \\ 
 16& & 15 \\ 
14 & & 13 \\ 
   & 12 \\ 
 & 11 &\\
 & 10&  
\end{psmatrix}
\psset{arrows=->, nodesep=2pt}
\ncline{10,2}{11,2}
\ncline{9,2}{10,2}
\ncline{8,3}{9,2}
\ncline{8,1}{9,2}
\ncline{7,3}{8,3}
\ncline{7,3}{8,1}
\ncline{7,1}{8,1}
\ncline{6,2}{7,3}
\ncline{5,1}{7,1}
\ncline{5,3}{6,2}
\ncline{5,1}{6,2}
\ncline{4,3}{5,3}
\ncline{4,1}{5,1}
\ncline{4,3}{5,1}
\ncline{3,2}{4,1}
\ncline{3,2}{4,3}
\ncline{2,2}{3,2}
\ncline{1,2}{2,2}
\end{minipage}
\end{tabular}\\
\par
So we count : $2 \times 3 \times 8 \times 24 = 2^7\times 3^2$ diagrams which is the cardinality of Weyl group.
\subsubsection{Type $E_6$}
\begin{convsanss}
The numbering of simple roots in the Dynkin diagram is as below :
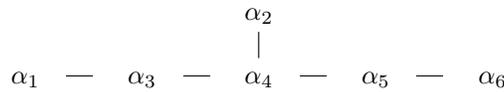
\[\begin{array}{ccccccccc}
 					& 					&   				 & 					 & \alpha_2	&&  &     &   \\ 
 					& 					&   				 & 					 & |					&&  &     & \\
 \alpha_1 &\text{---~}& \alpha_{3} &\text{---~}& \alpha_4 &\text{---~}& \alpha_{5}&\text{---~} & \alpha_6 \\
\end{array}\]
To describe the chosen reduced decomposition of $w_0$, we remark that the roots $\alpha_1$ to $\alpha_5$ span a roots system of $D_5$. Denote by $\tau$, the longest Weyl word used for $D_5$ then the decomposition  
 $$ \tau s_{6} s_5 s_4 s_2 s_3  s_1  s_4  s_3 s_5  s_4  s_6 s_2 s_5  s_4 s_3 s_1$$
is a reduced decomposition of $w_0$ which induces the following order on positive roots : \\
 \small
\begin{minipage}{20cm}
\begin{tabular}{p{0cm}@{}|c|c|c|c|c|}
\cline{6-6}
\multicolumn{5}{c|}{} & \multicolumn{1}{|c|}{ $ 
\beta_{21} = (1,2,2,3,2,1)$}\\
\cline{6-6}
\multicolumn{5}{c|}{} & \multicolumn{1}{|c|}{$ 
\beta_{22} = (1,1,2,3,2,1)$}\\
\cline{6-6}
\multicolumn{5}{c|}{} & \multicolumn{1}{|c|}{$
\beta_{23} = (1,1,2,2,2,1)$}\\
\cline{6-6}
\multicolumn{5}{c|}{} & \multicolumn{1}{|c|}{$
\beta_{24} = (1,1,2,2,1,1)$}\\
\cline{5-6}
\multicolumn{4}{c|}{} & \multicolumn{1}{|c|}{\tiny $\beta_{13}=\alpha_1+\alpha_2 + 2\alpha_3 + 2\alpha_4 + \alpha_5$} & \multicolumn{1}{|c|}{$\beta_{25} = (1,1,1,2,2,1)$} \\
 \cline{4-6}
\multicolumn{3}{c|}{} & \multicolumn{1}{|c|}{\tiny $\beta_7 =  \alpha_2 + \alpha_3 + 2\alpha_4 + \alpha_5$} & \multicolumn{1}{|c|}{\tiny $\beta_{14} = \alpha_1 + \alpha_2 + \alpha_3 + 2\alpha_4 + \alpha_5$} & \multicolumn{1}{|c|}{$\beta_{26} = (0,1,1,2,2,1)$} \\ 
\cline{3-6} 
\multicolumn{2}{c|}{} & \multicolumn{1}{|c|}{\tiny $\beta_3 =  \alpha_2 + \alpha_4+\alpha_5 $} & \multicolumn{1}{|c|}{\tiny $\beta_{8} =  \alpha_2 + \alpha_3 + \alpha_4+ \alpha_5 $} & \multicolumn{1}{|c|}{\tiny $\beta_{15} = \alpha_1 + \alpha_2 + \alpha_3 + \alpha_4 + \alpha_5$}& \multicolumn{1}{|c|}{$\beta_{27} = (1,1,1,2,1,1)$} \\ 
\cline{2-6}
 & \tiny $\beta_1 = \alpha_2 $& \tiny $ \beta_4 = \alpha_4 + \alpha_5$  & \tiny $\beta_9 = \alpha_2 + \alpha_3 + \alpha_4 $& \tiny $  \beta_{16} = \alpha_1 + \alpha_{3} + \alpha_4 + \alpha_5$ & $\beta_{28} =(0,1,1,2,1,1)$ \\
  \cline{2-6}
& \tiny $\beta_2 = \alpha_5 $& \tiny $ \beta_5 = \alpha_2 + \alpha_4$  & \tiny $\beta_{10} =  \alpha_3 + \alpha_4 +\alpha_5 $&\tiny  $ \beta_{17} = \alpha_1+ \alpha_2 + \alpha_{3} + \alpha_4 $ & $\beta_{29} = (1,1,1,1,1,1)$ \\
  \cline{2-6}
\multicolumn{2}{c|}{} & \multicolumn{1}{|c|}{\tiny $\beta_6 = \alpha_4$} & \multicolumn{1}{|c|}{\tiny $\beta_{11} = \alpha_3 + \alpha_4$} & \multicolumn{1}{|c|}{\tiny $\beta_{18} =\alpha_1 + \alpha_3 + \alpha_4  $} & \multicolumn{1}{|c|}{$\beta_{30} =  (0,1,1,1,1,1)$} \\
 \cline{3-6}
\multicolumn{3}{c|}{} & \multicolumn{1}{|c|}{\tiny $\beta_{12} = \alpha_3$} & \multicolumn{1}{|c|}{\tiny $\beta_{19} = \alpha_1 + \alpha_3$} & \multicolumn{1}{|c|}{ $\beta_{31} =  (1,0,1,1,1,1)$} \\
\cline{4-6}
\multicolumn{4}{c|}{} & \multicolumn{1}{|c|}{\tiny $\beta_{20} = \alpha_1$} & \multicolumn{1}{|c|}{$\beta_{32} =(0,1,0,1,1,1)$}  \\
\cline{5-6}
\multicolumn{5}{c|}{} & \multicolumn{1}{|c|}{$\beta_{33} = (0,0,1,1,1,1) $}  \\
\cline{6-6}
\multicolumn{5}{c|}{} & \multicolumn{1}{|c|}{$\beta_{34} =  (0,0,0,1,1,1)$}  \\
\cline{6-6}
\multicolumn{5}{c|}{} & \multicolumn{1}{|c|}{$\beta_{35} = (0,0,0,0,1,1)$}  \\
\cline{6-6}
\multicolumn{5}{c|}{} & \multicolumn{1}{|c|}{$\beta_{36} = (0,0,0,0,0,1)$}  \\
\cline{6-6}
\end{tabular}\\
\end{minipage}
\end{convsanss}
One obtains the following implications thanks to Proposition \ref{contrainte elementaire}, one counts the diagrams satisfying the implications :\\
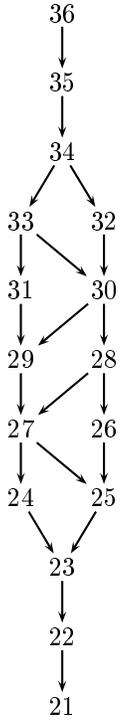
\begin{tabular}{lccr}
\begin{minipage}{3cm}
\begin{psmatrix}[rowsep=.5cm, colsep=.2cm]
 & 36 \\ 
 & 35 \\ 
 & 34 &  \\ 
33& &32\\
31 &  & 30 \\ 
29 &  & 28 \\ 
27 &  & 26 \\ 
24 &  & 25 \\
&23 \\
 & 22 \\ 
 & 21 
\end{psmatrix}
\psset{arrows=->, nodesep=2pt}
\ncline{1,2}{2,2}
\ncline{2,2}{3,2}
\ncline{3,2}{4,1}
\ncline{3,2}{4,3}
\ncline{4,1}{5,1}
\ncline{4,1}{5,3}
\ncline{4,3}{5,3}
\ncline{5,1}{6,1}
\ncline{5,3}{6,1}
\ncline{5,3}{6,3}
\ncline{6,1}{7,1}
\ncline{6,3}{7,1}
\ncline{6,3}{7,3}
\ncline{7,1}{8,3}
\ncline{7,1}{8,1}
\ncline{7,3}{8,3}
\ncline{8,3}{9,2}
\ncline{8,1}{9,2}
\ncline{9,2}{10,2}
\ncline{10,2}{11,2}
\end{minipage}
& &&
\begin{minipage}{6cm}
\begin{tabular}{c|c}
Diagrams beginning by : & Counting the diagrams :\\
\hline
none & 1\\
21					& 1	\\
22					& 1\\
23					& 1\\
24					& 1\\
25					& 2\\
26					& 2\\
27					& 2\\
28					& 1\\
29					& 3\\
30					& 1\\
31					& 4\\
32					& 2\\
33					& 2\\
34					& 1\\
35					& 1\\
36					& 1\\
\hline
 & \\
TOTAL 			& 27 = $3^3$
\end{tabular} \\
\end{minipage}
\end{tabular} \\
\par
For the first five columns, we have a system of type $D_5$, so $4 \times 6 \times 8 \times 10 = 2^7 \times 3 \times 5$ diagrams.\\
By adding the last one, one obtains $2^7 \times 3^4 \times 5$ diagrams which equal the cardinality of the Weyl group. 
So the set of Cauchon diagrams $\mathcal{D}$ has the same cardinality as the Weyl group $W$.
\subsubsection{Type $E_7$}
\begin{convsanss}
The numbering of simple roots in the Dynkin diagram is as below :
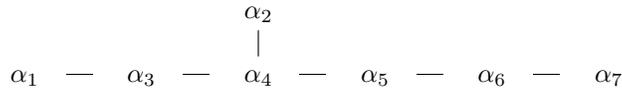
\[\begin{array}{ccccccccccc}
 					& 					&   				 & 					 & \alpha_2	&&  &     &&&  \\ 
 					& 					&   				 & 					 & |					&&  &     & &&\\
 \alpha_1 &\text{---~}& \alpha_{3} &\text{---~}& \alpha_4 &\text{---~}& \alpha_{5}&\text{---~} & \alpha_6 &\text{---~}& \alpha_7 \\
\end{array}\]
As the roots $\alpha_1$ to $\alpha_6$ span a roots system of type $E_6$, denote by $\sigma$ the longest Weyl word used for the type $E_6$. The decomposition
\[ \sigma s_{7} s_6 s_5 s_4 s_2 s_3  s_1  s_4  s_3 s_5  s_4  s_6 s_2 s_5  s_7 s_4 s_6 s_3 s_5 s_1 s_4 s_2 s_3 s_4 s_5 s_6 s_7\]
 is a reduced decomposition of $w_0$ which induces the following order on positive roots. To stock all the informations in a tabular, we replace vectors by their coordinates on the basis $\alpha_1, \cdots , \alpha_7$ omitting useless coordinates. \\
\begin{tabular}{p{0cm}@{}|c|c|c|c|c|c|}
\multicolumn{6}{c|}{} & \multicolumn{1}{|c|}{\tiny $ 
\mathbf{\alpha_1,\alpha_{2},\alpha_3,\alpha_4,\alpha_5,\alpha_6, \alpha_7}$}\\
\cline{7-7}
\multicolumn{6}{c|}{} & \multicolumn{1}{|c|}{\tiny $\beta_{37} (2, 2, 3, 4, 3, 2, 1)$}\\
\cline{7-7}
\multicolumn{6}{c|}{} & \multicolumn{1}{|c|}{\tiny $\beta_{38} (1, 2, 3, 4, 3, 2, 1)$}\\
\cline{7-7}
\multicolumn{6}{c|}{} & \multicolumn{1}{|c|}{\tiny $\beta_{39} (1, 2, 2, 4, 3, 2, 1)$}\\
\cline{7-7}
\multicolumn{6}{c|}{} & \multicolumn{1}{|c|}{\tiny $\beta_{40} (1, 2, 2, 3, 3, 2, 1)$}\\
\cline{7-7}
\multicolumn{6}{c|}{} & \multicolumn{1}{|c|}{\tiny $\beta_{41} ( 1, 1, 2, 3, 3, 2, 1 )$}\\
\cline{7-7}
\multicolumn{5}{c|}{} & \multicolumn{1}{|c|}{\tiny $ 
\mathbf{\alpha_1,\alpha_{2},\alpha_3,\alpha_4,\alpha_5,\alpha_6}$}& \multicolumn{1}{|c|}{\tiny $\beta_{42} (1, 2, 2, 3, 2, 2, 1)$}\\
\cline{6-7}
\multicolumn{5}{c|}{} & \multicolumn{1}{|c|}{\tiny $ 
\beta_{21} (1, 2, 2, 3, 2, 1)$}& \multicolumn{1}{|c|}{\tiny $\beta_{43} (1, 2, 2, 3, 2, 1, 1)$}\\
\cline{6-7}
\multicolumn{5}{c|}{} & \multicolumn{1}{|c|}{\tiny$ 
\beta_{22} (1, 1, 2, 3, 2, 1) $}& \multicolumn{1}{|c|}{\tiny $\beta_{44} (1, 1, 2, 3, 2, 2, 1)$} \\
\cline{6-7}
\multicolumn{5}{c|}{} & \multicolumn{1}{|c|}{\tiny$
\beta_{23} (1, 1, 2, 2, 2, 1)$}& \multicolumn{1}{|c|}{\tiny $\beta_{45} (1, 1, 2, 3, 2, 1, 1)$}\\
\cline{6-7}
\multicolumn{4}{c|}{} & \multicolumn{1}{|c|}{\tiny $ 
\mathbf{\alpha_1,\alpha_{2},\alpha_3,\alpha_4,\alpha_5,\alpha_6}$}& \multicolumn{1}{|c|}{\tiny$
\beta_{24}(1, 1, 2, 2, 1, 1)$}& \multicolumn{1}{|c|}{\tiny $\beta_{46} (1, 1, 2, 2, 2, 2, 1)$}\\
\cline{5-7}
\multicolumn{3}{c|}{} & \multicolumn{1}{|c|}{\tiny $ 
\mathbf{\alpha_2,\alpha_{3},\alpha_4,\alpha_5}$}& \multicolumn{1}{|c|}{\tiny $\beta_{13}(1,1,2,2,1)$} & \multicolumn{1}{|c|}{\tiny$\beta_{25} (1, 1, 1, 2, 2, 1)$}& \multicolumn{1}{|c|}{\tiny $\beta_{47} (1, 1, 2, 2, 2, 1, 1)$} \\
 \cline{4-7}
\multicolumn{2}{c|}{}& \multicolumn{1}{|c|}{\tiny $ 
\mathbf{\alpha_2,\alpha_{4},\alpha_5}$} & \multicolumn{1}{|c|}{\tiny $\beta_7 (1,1,2,1)$} & \multicolumn{1}{|c|}{\tiny $\beta_{14} (1,1,1,2,1)$} & \multicolumn{1}{|c|}{\tiny$\beta_{26} (0, 1, 1, 2, 2, 1)$}& \multicolumn{1}{|c|}{\tiny $\beta_{48}(1, 1, 1, 2, 2, 2, 1) $} \\ 
\cline{3-7} 
\multicolumn{1}{c|}{}& \multicolumn{1}{|c|}{\tiny $ 
\mathbf{\alpha_2,\alpha_{5}}$} & \multicolumn{1}{|c|}{\tiny $\beta_3 (1,1,1) $} & \multicolumn{1}{|c|}{\tiny $\beta_{8} (1,1,1,1) $} & \multicolumn{1}{|c|}{\tiny $\beta_{15} (1,1,1,1,1)$}& \multicolumn{1}{|c|}{\tiny$\beta_{27} (1, 1, 1, 2, 1, 1)$}& \multicolumn{1}{|c|}{\tiny $\beta_{49} (1, 1, 2, 2, 1, 1, 1)$} \\ 
\cline{2-7}
 & \tiny $\beta_1 (1,0) $& \tiny $ \beta_4 (0,1,1)$  & \tiny $\beta_9 (1,1,1,0) $& \tiny $  \beta_{16} (1,0,1,1,1)$ & \tiny $\beta_{28} (0, 1, 1, 2, 1, 1)$ & \tiny $\beta_{50} (1, 1, 1, 2, 2, 1, 1)$ \\
  \cline{2-7}
& \tiny $\beta_2 (0,1) $& \tiny $ \beta_5 (1,1,0)$  & \tiny $\beta_{10} (0,1,1,1) $&\tiny  $ \beta_{17} (1,1,1,1,0)$ & \tiny $\beta_{29} (1,1,1,1,1,1)$& \tiny $\beta_{51} (0, 1, 1, 2, 2, 2, 1)$ \\
  \cline{2-7}
\multicolumn{2}{c|}{} & \multicolumn{1}{|c|}{\tiny $\beta_6 (0,1,0)$} & \multicolumn{1}{|c|}{\tiny $\beta_{11} (0,1,1,0)$} & \multicolumn{1}{|c|}{\tiny $\beta_{18} (1,0,1,1,0)$} & \multicolumn{1}{|c|}{\tiny $\beta_{30} (0, 1, 1, 1, 1, 1)$}& \multicolumn{1}{|c|}{\tiny $\beta_{52} (1, 1, 1, 2, 1, 1, 1)$} \\
 \cline{3-7}
\multicolumn{3}{c|}{} & \multicolumn{1}{|c|}{\tiny $\beta_{12} (0,1,0,0)$} & \multicolumn{1}{|c|}{\tiny $\beta_{19}(1,0,1,0,0)$} & \multicolumn{1}{|c|}{\tiny $\beta_{31} (1, 0, 1, 1, 1, 1)$}& \multicolumn{1}{|c|}{\tiny $\beta_{53} (0, 1, 1, 2, 2, 1, 1)$} \\
\cline{4-7}
\multicolumn{4}{c|}{} & \multicolumn{1}{|c|}{\tiny $\beta_{20}(1,0,0,0,0)$} & \multicolumn{1}{|c|}{\tiny$\beta_{32}(0, 1, 0, 1, 1, 1) $}& \multicolumn{1}{|c|}{\tiny $\beta_{54}(1, 1, 1, 1, 1, 1, 1) $}  \\
\cline{5-7}
\multicolumn{5}{c|}{} & \multicolumn{1}{|c|}{\tiny$\beta_{33} (0, 0, 1, 1, 1, 1)$}& \multicolumn{1}{|c|}{\tiny $\beta_{55}(0, 1, 1, 2, 1, 1, 1) $}  \\
\cline{6-7}
\multicolumn{5}{c|}{} & \multicolumn{1}{|c|}{\tiny$\beta_{34} (0, 0, 0, 1, 1, 1)$}& \multicolumn{1}{|c|}{\tiny $\beta_{56}(1, 0, 1, 1, 1, 1, 1) $}  \\
\cline{6-7}
\multicolumn{5}{c|}{} & \multicolumn{1}{|c|}{\tiny$\beta_{35} ( 0, 0, 0, 0, 1, 1)$}& \multicolumn{1}{|c|}{\tiny $\beta_{57} (0, 1, 1, 1, 1, 1, 1)$}  \\
\cline{6-7}
\multicolumn{5}{c|}{} & \multicolumn{1}{|c|}{\tiny$\beta_{36} (0,0,0,0,0,1)$}& \multicolumn{1}{|c|}{\tiny $\beta_{58} (0, 1, 0, 1, 1, 1, 1)$}  \\
\cline{6-7}
\multicolumn{6}{c|}{} & \multicolumn{1}{|c|}{\tiny $\beta_{59} (0, 0, 1, 1, 1, 1, 1)$}\\
\cline{7-7}
\multicolumn{6}{c|}{} & \multicolumn{1}{|c|}{\tiny $\beta_{60} (0, 0, 0, 1, 1, 1, 1)$}\\
\cline{7-7}
\multicolumn{6}{c|}{} & \multicolumn{1}{|c|}{\tiny $\beta_{61} (0, 0, 0, 0, 1, 1, 1)$}\\
\cline{7-7}
\multicolumn{6}{c|}{} & \multicolumn{1}{|c|}{\tiny $\beta_{62} (0,0,0,0,0,1,1)$}\\
\cline{7-7}
\multicolumn{6}{c|}{} & \multicolumn{1}{|c|}{\tiny $\beta_{63} (0,0,0,0,0,0,1)$}\\
\cline{7-7}
\end{tabular}\\
\end{convsanss}
We already know the form of diagrams for the first six columns. We use the same methods as in type $E_6$ (Proposition \ref{contrainte elementaire}) to find the implications in the last columns of type $E_7$. One obtains : \\
\scalebox{.85}{
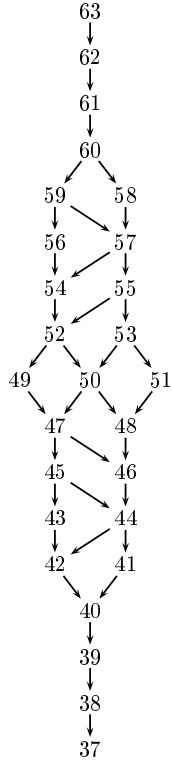
\begin{tabular}{cc}
\begin{minipage}{.5\linewidth}
\begin{center}
\begin{psmatrix}[rowsep=.3cm, colsep=.2cm]
 &  & 63 \\ 
 &  & 62 \\ 
 &  & 61 \\ 
 &  & 60 \\ 
 & 59 &  & 58 \\ 
&56 &  & 57 \\ 
 & 54 &  & 55 \\ 
 & 52 &  & 53 \\ 
  49 &  & 50 &  & 51 \\ 
  & 47 &  & 48 \\ 
 & 45 &  & 46 \\ 
&43 &  & 44 \\ 
 & 42 &  & 41 \\ 
 &  & 40 \\ 
 &  & 39 \\ 
 &  & 38 \\
 & & 37
\end{psmatrix}
\psset{arrows=->, nodesep=2pt}
\ncline{1,3}{2,3}
\ncline{2,3}{3,3}
\ncline{3,3}{4,3}
\ncline{4,3}{5,2}
\ncline{4,3}{5,4}
\ncline{5,2}{6,4}
\ncline{5,2}{6,2}
\ncline{5,4}{6,4}
\ncline{6,4}{7,4}
\ncline{6,4}{7,2}
\ncline{6,2}{7,2}
\ncline{7,4}{8,4}
\ncline{7,4}{8,2}
\ncline{7,2}{8,2}
\ncline{8,4}{9,5}
\ncline{8,4}{9,3}
\ncline{8,2}{9,3}
\ncline{8,2}{9,1}
\ncline{9,5}{10,4}
\ncline{9,3}{10,4}
\ncline{9,3}{10,2}
\ncline{9,1}{10,2}
\ncline{10,4}{11,4}
\ncline{10,2}{11,4}
\ncline{10,2}{11,2}
\ncline{11,4}{12,4}
\ncline{11,2}{12,4}
\ncline{11,2}{12,2}
\ncline{12,4}{13,2}
\ncline{12,4}{13,4}
\ncline{12,2}{13,2}
\ncline{13,2}{14,3}
\ncline{13,4}{14,3}
\ncline{14,3}{15,3}
\ncline{15,3}{16,3}
\ncline{16,3}{17,3}
\end{center}
\end{minipage}
&
\begin{minipage}{.5\linewidth}
\begin{center}
\begin{tabular}{c|c}
Diagrams beginning : & Counting the diagrams\\
\hline
none & 1\\
37					& 1	\\
38					& 1\\
39					& 1\\
40					& 1\\
41					& 1\\
42					& 2\\
43					& 2\\
44					& 2\\
45					& 1\\
46					& 3\\
47					& 1\\
48					& 4\\
49					& 2\\
50					& 2\\
51					& 6\\
52					& 2\\
53					& 4\\
54					&	3\\
55					& 2\\
56					& 4\\
57					&	2\\
58					& 2\\
59					& 2\\
60					& 1\\
61					& 1\\
62					& 1\\
63					&	1\\
\hline
 & \\
TOTAL 			& 56 = $2^3 \times 7$

\end{tabular} \\
\end{center}
\end{minipage}
\end{tabular}}\\
\par
Cauchon diagrams for the first six columns come from the type $E_{6}$, that is to say $2^7 \times 3^4 \times 5$ diagrams. With the last column, one has $2^{10} \times 3^4 \times 5 \times 7$ diagrams which equals the cardinality of the Weyl group. We have again $| \mathcal{D}| = |W|$.

\subsubsection{Type $E_8$}
\begin{convsanss}
The numbering of simple roots in the Dynkin diagram is as below :
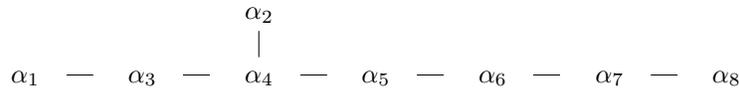
\[\begin{array}{ccccccccccccc}
 					& 					&   				 & 					 & \alpha_2	&&  &     &&& &&  \\ 
 					& 					&   				 & 					 & |					&&  &     &&& &&\\
 \alpha_1 &\text{---~}& \alpha_{3} &\text{---~}& \alpha_4 &\text{---~}& \alpha_{5}&\text{---~} & \alpha_6 &\text{---~}& \alpha_7& \text{---~}& \alpha_8\\
\end{array}\]
As the roots $\alpha_1$ to $\alpha_7$ span a roots system of type $E_7$, denote by $\sigma_{7}$, the longest Weyl word used for type $E_7$. The decomposition  
\par
$ \sigma_{7} s_8 s_{7} s_6 s_5 s_4 s_2 s_3  s_1  s_4  s_3 s_5  s_4  s_6 s_2 s_5  s_7 s_4 s_6
 s_8 s_3 s_5 s_7 s_1 s_4 s_6 s_3 s_2 s_5 s_4 s_{5} s_2 s_3 s_6 s_1 s_4  s_7  s_3  s_5 s_8s_4  s_6 s_2 s_5  s_7s_4 $\\$s_6 s_3 s_5 s_1 s_4 s_2 s_3 s_4 s_5 s_6 s_7 s_8$.\\
is a reduced decomposition of $w_0$ which induces the following order on positive roots. (Only the last column is given). In the second column of this tabular, we compute  $h'(\beta_i)$ : 
\end{convsanss}
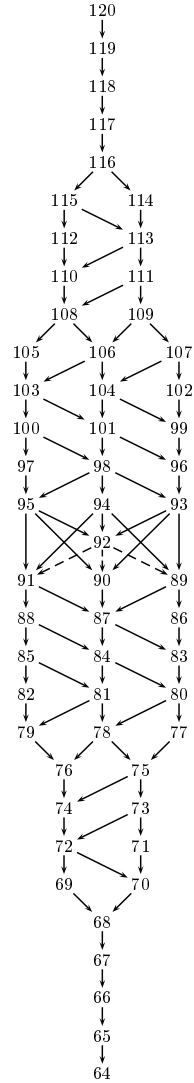
\begin{tabular}{lc}
\begin{minipage}{.6\linewidth}
\[
\begin{array}{ccc}
\begin{array}{c|c}
\beta_i & h'(\beta_i) \\
\hline
\beta_{64}( 2, 3, 4, 6, 5, 4, 3, 1 ) &	28 \\
\hline
\beta_{65}( 2, 3, 4, 6, 5, 4, 2, 1 ) &	27 \\
\hline
\beta_{66}( 2, 3, 4, 6, 5, 3, 2, 1 ) &	26 \\
\hline
\beta_{67}( 2, 3, 4, 6, 4, 3, 2, 1 ) &	25 \\
\hline
\beta_{68}( 2, 3, 4, 5, 4, 3, 2, 1 ) &	24 \\
\hline
\beta_{69}( 2, 2, 4, 5, 4, 3, 2, 1 ) &	23 \\
\beta_{70}( 2, 3, 3, 5, 4, 3, 2, 1 ) &	23 \\
\hline
\beta_{71}( 1, 3, 3, 5, 4, 3, 2, 1 ) & 	22 \\
\beta_{72}( 2, 2, 3, 5, 4, 3, 2, 1 )	& 	22 \\
\hline
\beta_{73}(1, 2, 3, 5, 4, 3, 2, 1 ) &	21 \\
\beta_{74}( 2, 2, 3, 4, 4, 3, 2, 1 ) 	& 	21 \\
\hline
\beta_{75}( 1, 2, 3, 4, 4, 3, 2, 1)&	20 \\
\beta_{76}(2, 2, 3, 4, 3, 3, 2, 1 )& 	20 \\
\hline
\beta_{77}( 1, 2, 2, 4, 4, 3, 2, 1 )&	19 \\
\beta_{78}( 1, 2, 3, 4, 3, 3, 2, 1 )& 	19 \\
\beta_{79}( 2, 2, 3, 4, 3, 2, 2, 1 )&	19 \\
\hline
\beta_{80}(1, 2, 2, 4, 3, 3, 2, 1 )&	18  \\
\beta_{81}( 1, 2, 3, 4, 3, 2, 2, 1 )&	18 \\
\beta_{82} (2, 2, 3, 4, 3, 2, 1, 1 )&18\\
\hline
\beta_{83}( 1, 2, 2, 3, 3, 3, 2, 1 )& 	17 \\
\beta_{84}( 1, 2, 2, 4, 3, 2, 2, 1 )&	17 \\
\beta_{85}( 1, 2, 3, 4, 3, 2, 1, 1 )& 	17 \\
\hline
\beta_{86}( 1, 1, 2, 3, 3, 3, 2, 1 )& 	16 \\
\beta_{87}( 1, 2, 2, 3, 3, 2, 2, 1 )&	16 \\
\beta_{88}( 1, 2, 2, 4, 3, 2, 1, 1 )&	16 \\
\hline
\beta_{89}( 1, 1, 2, 3, 3, 2, 2, 1 )&	15 \\
\beta_{90}( 1, 2, 2, 3, 2, 2, 2, 1 )& 	15 \\
\beta_{91}( 1, 2, 2, 3, 3, 2, 1, 1 )&	15 \\
\hline
\beta_{92}( 2, 3, 4, 6, 5, 4, 3, 2 )&29/2 \\
\hline
\end{array}
&&
\begin{array}{c|c}
\beta_i & h'(\beta_i) \\
\hline
\beta_{93}( 1, 1, 2, 3, 2, 2, 2, 1 )&	14 \\
\beta_{94}( 1, 1, 2, 3, 3, 2, 1, 1)&	14\\
\beta_{95}( 1, 2, 2, 3, 2, 2, 1, 1 )&	14\\
\hline 
\beta_{96}( 1, 1, 2, 2, 2, 2, 2, 1 )&	13\\
\beta_{97}( 1, 2, 2, 3, 2, 1, 1, 1 )&	13\\
\beta_{98}( 1, 1, 2, 3, 2, 2, 1, 1 )&	13 \\
\hline 
\beta_{99}( 1, 1, 1, 2, 2, 2, 2, 1 )& 	12\\
\beta_{100}	(1, 1, 2, 3, 2, 1, 1, 1)&	12\\
\beta_{101}	( 1, 1, 2, 2, 2, 2, 1, 1 )&	12\\
\hline
\beta_{102}( 0, 1, 1, 2, 2, 2, 2, 1 )&	11\\
\beta_{103}( 1, 1, 2, 2, 2, 1, 1, 1 )&	11\\
\beta_{104}	( 1, 1, 1, 2, 2, 2, 1, 1 )&	11\\
\hline 
\beta_{105}( 1, 1, 2, 2, 1, 1, 1, 1 )&	10 \\
\beta_{106}	(1, 1, 1, 2, 2, 1, 1, 1 )&10 \\
\beta_{107}	( 0, 1, 1, 2, 2, 2, 1, 1 )&	10 \\
\hline
\beta_{108}( 1, 1, 1, 2, 1, 1, 1, 1 )&	9 \\
\beta_{109}	(0, 1, 1, 2, 2, 1, 1, 1 )&	9 \\
\hline
\beta_{110}	(1, 1, 1, 1, 1, 1, 1, 1 )& 	8 \\
\beta_{111}	( 0, 1, 1, 2, 1, 1, 1, 1)&	8 \\
\hline
\beta_{112}	( 1, 0, 1, 1, 1, 1, 1, 1 )&	7  \\
\beta_{113}	( 0, 1, 1, 1, 1, 1, 1, 1 )&	7 \\
\hline
\beta_{114}	( 0, 1, 0, 1, 1, 1, 1, 1 )&	6 \\ 
\beta_{115}	( 0, 0, 1, 1, 1, 1, 1, 1 )&	6 \\
\hline
\beta_{116}	( 0, 0, 0, 1, 1, 1, 1, 1 )&	5 \\
\hline 
\beta_{117}	( 0, 0, 0, 0, 1, 1, 1, 1)&	4 \\
\hline
\beta_{118}	( 0, 0, 0, 0, 0, 1, 1, 1 )& 	3\\
\hline
\beta_{119}	( 0, 0, 0, 0, 0, 0, 1, 1 )&	2 \\
\hline
\beta_{120}	( 0, 0, 0, 0, 0, 0, 0, 1 )&	1 \\
\end{array}
\end{array}\]
We already know the shape of diagrams from the first seven columns. Thanks to Propositions \ref{contrainte elementaire}, \ref{contrainte elementaire avant ex} and \ref{contrainte elementaire apres ex}, one obtains the implications for the last column.\\
In particular, we obtain implications such as $\mathbf{(i \Rightarrow j \ \text{or} \ k)}$ :
\[(92 \Rightarrow 91 \ \text{or} \ 90) \ \text{and} \ (92 \Rightarrow 90 \ \text{or} \ 89) \ \text{and} \ (92 \Rightarrow 91 \ \text{or} \ 89).\]
Those implications appear with dashed arrows in the diagram : \\
\end{minipage}
& 
\begin{minipage}{.4\linewidth}
\begin{center}
\scalebox{.7}{
\begin{psmatrix}[rowsep=.3cm, colsep=.2cm]
 &  && 120 \\ 
 &  & &119 \\ 
 &  & &118 \\ 
 &  & &117 \\ 
 &  & &116 \\ 
 & &115 &  & 114 \\ 
& &112 & &113 \\ 
 & &110  & &111 \\ 
 &  & 108 &  & 109 \\
 & 105 &  & 106 &  & 107 \\ 
 & 103 &  & 104 &  & 102 \\ 
 & 100 &  & 101 &  & 99 \\ 
 & 97 &  & 98 &  & 96 \\ 
 & 95 &  & 94 &  & 93 \\ 
 &  &  & 92 \\ 
 & 91 &  & 90 &  & 89 \\ 
 & 88 &  & 87 &  & 86 \\ 
 & 85 &  & 84 &  & 83 \\ 
 & 82 &  & 81 &  & 80 \\
 & 79 &  & 78 &  & 77 \\ 
 && 76 &  & 75 \\ 
 && 74 &  & 73 \\ 
 && 72 &  & 71 \\ 
 && 69 &  & 70 \\ 
 &&  & 68 \\ 
 &&  & 67 \\ 
 &&  & 66 \\ 
 &&  & 65 \\ 
 &&  & 64 
\end{psmatrix}
\psset{arrows=->, nodesep=2pt}
\ncline{1,4}{2,4}
\ncline{2,4}{3,4}
\ncline{3,4}{4,4}
\ncline{4,4}{5,4}
\ncline{5,4}{6,3}
\ncline{5,4}{6,5}
\ncline{6,3}{7,5}
\ncline{6,3}{7,3}
\ncline{6,5}{7,5}
\ncline{7,5}{8,5}
\ncline{7,5}{8,3}
\ncline{7,3}{8,3}
\ncline{8,5}{9,5}
\ncline{8,5}{9,3}
\ncline{8,3}{9,3}
\ncline{9,5}{10,6}
\ncline{9,5}{10,4}
\ncline{9,3}{10,4}
\ncline{9,3}{10,2}
\ncline{10,6}{11,4}
\ncline{10,6}{11,6}
\ncline{10,4}{11,4}
\ncline{10,4}{11,2}
\ncline{10,2}{11,2}
\ncline{11,4}{12,4}
\ncline{11,4}{12,6}
\ncline{11,2}{12,4}
\ncline{11,2}{12,2}
\ncline{11,6}{12,6}
\ncline{12,4}{13,4}
\ncline{12,4}{13,6}
\ncline{12,2}{13,4}
\ncline{12,2}{13,2}
\ncline{12,6}{13,6}
\ncline{13,4}{14,2}
\ncline{13,4}{14,4}
\ncline{13,4}{14,6}
\ncline{13,2}{14,2}
\ncline{13,6}{14,6}
\ncline{14,2}{15,4}
\ncline{14,2}{16,2}
\ncline{14,2}{16,4}
\ncline{14,4}{15,4}
\ncline{14,4}{16,2}
\ncline{14,4}{16,6}
\ncline{14,6}{15,4}
\ncline{14,6}{16,4}
\ncline{14,6}{16,6}
\ncline{16,2}{17,2}
\ncline{16,2}{17,4}
\ncline{16,4}{17,4}
\ncline{16,6}{17,4}
\ncline{16,6}{17,6}
\ncline{17,2}{18,2}
\ncline{17,2}{18,4}
\ncline{17,4}{18,4}
\ncline{17,4}{18,6}
\ncline{17,6}{18,6}
\ncline{18,2}{19,2}
\ncline{18,2}{19,4}
\ncline{18,4}{19,4}
\ncline{18,4}{19,6}
\ncline{18,6}{19,6}
\ncline{19,2}{20,2}
\ncline{19,4}{20,2}
\ncline{19,4}{20,4}
\ncline{19,6}{20,4}
\ncline{19,6}{20,6}
\ncline{20,2}{21,3}
\ncline{20,4}{21,3}
\ncline{20,4}{21,5}
\ncline{20,6}{21,5}
\ncline{21,3}{22,3}
\ncline{21,5}{22,3}
\ncline{21,5}{22,5}
\ncline{22,3}{23,3}
\ncline{22,5}{23,3}
\ncline{22,5}{23,5}
\ncline{23,3}{24,5}
\ncline{23,3}{24,3}
\ncline{23,5}{24,5}
\ncline{24,5}{25,4}
\ncline{24,3}{25,4}
\ncline{25,4}{26,4}
\ncline{26,4}{27,4}
\ncline{27,4}{28,4}
\ncline{28,4}{29,4}
\psset{nodesep=1pt, linestyle=dashed}
\ncline{15,4}{16,2}
\ncline{15,4}{16,4}
\ncline{15,4}{16,6}}
\end{center}
\end{minipage}
\end{tabular}\\
With the previous implications, we can count the diagrams.\\
\begin{tabular}{cc}
\begin{tabular}{c|c}
Diagrams beginning :  & Counting diagrams \\
\hline
none & 1\\
64					& 1	\\
65					& 1\\
66					& 1\\
67					& 1\\
68					& 1\\
69					& 1\\
70					& 2\\
71					& 2\\
72					& 2\\
73					& 1\\
74					& 3\\
75					& 1\\
76					& 4\\
77					& 2\\
78					& 2\\
79					& 7\\
80					& 2\\
81					&	3\\
82					& 11\\
83					& 5\\
84					&	4\\
85					& 6\\
86					& 9\\
87					& 6\\
88					& 5\\
89					& 4\\
90					&	12\\
91					& 6\\
92					& 8
\end{tabular}
&
\begin{tabular}{c|c}
Diagrams beginning : & Counting diagrams \\
\hline
93					& 5\\
94					& 3\\
95					& 6\\
96					& 8\\
97					& 8\\
98					&4 \\
99					& 12\\
100					& 3\\
101					& 6\\
102					& 16\\
103					& 3\\
104					& 8\\
105					& 5\\
106					& 4\\
107					& 7\\
108					& 3\\
109					& 3\\
110					& 4\\
111					& 2\\
112					& 5\\
113					& 2\\
114					& 2\\
115					& 2\\
116					& 1 \\
117					& 1\\
118					& 1\\
119					& 1\\
120					& 1\\
\hline
 & \\
TOTAL 			& 240 = $2^4 \times 3 \times 5$
\end{tabular} 
\end{tabular} \\
\par
$ $ \\
$ $ \\
For the first seven columns, one have a system of type $E_7$, so $2^{10} \times 3^4 \times 5$ diagrams.\\
By adding the last column, we count $2^{10} \times 3^4 \times 5 \times 7$ diagrams which equals the cardinality of the Weyl group. We have as in the other cases $|\mathcal{D}|=|W|$.

\pagebreak
\tableofcontents
\pagebreak
\bibliographystyle{alpha}
\bibliography{biblio.bib}
\end{document}


\maketitle
$ $ \\
$ $ \\
$ $ \\
$ $ \\
$ $ \\
$ $ \\ 
\textbf{Abstract}\\

In this paper, we give an algorithmic description of Cauchon's diagrams in nilpotent quantum algebras $\mathcal{U}_q^+(\mathfrak{g})$ when the choosen reduced decomposition of the longest element $w_0$ of the Weyl group $W$ corresponds to a good ordering of the positive part $\Phi^+$ in the sense of G. Lusztig (\cite{Lu2}), where $\Phi$ denotes the root system. This algorithmic description is based on conditions, each one being defined by each Lusztig's admissible plane (\cite{Lu2}). Moreover, we give examples for explicit reduced decomposition of $w_0$ in each possible type of the finite dimensional complex simple Lie algebra $\mathfrak{g}$. We check that the number of Cauchon diagrams is always equal to the cardinal of $W$. In a future paper, we will prove that the Cauchon diagrams correspond canonically to the positive subexpressions of $w_0$. So the results of this paper also give an algorithmic description of those positive subexpressions when the reduced decomposition of $w_0$ satisfies the previous conditions. \\
$ $ \\

\section{Introduction}
Soit $\mathfrak{g}$ une algèbre de Lie complexe simple de dimension finie, $\mathbb{K}$ un corps commutatif et $q$ un élément de $\mathbb{K}^*$, $q$ non racine de l'unité. \\
$ $\\
On utilise les conventions de Jantzen (\cite{J}) pour la définition du groupe quantique $\mathcal{U}_q(\mathfrak{g})$ et, moyennant le choix d'une décomposition réduite de l'élément de plus grande longueur $w_0$ du groupe de Weyl $W$, pour la construction d'un système générateur ($X_\beta, \beta \in \Phi^+$) de la partie positive $\mathcal{U}_q^+(\mathfrak{g})$ (cf paragraphe 3). L'action naturelle du tore sur $\mathcal{U}_q^+(\mathfrak{g})$ induit une stratification de Spec$(\mathcal{U}_q^+(\mathfrak{g}))$ qui peut être décrite au moyen de l'algorithme d'effacement des dérivations (cf paragraphe 4.2), les strates étant en bijection naturelle avec les diagrammes de Cauchon, la géométrie d'une strate étant complètement décrite par le diagramme qui lui est associé. Il est donc naturel de chercher à décrire ces diagrammes. Cette description est étroitement liée aux formules de redressement de Levendorskii et Soibelmann \cite{MR1116413} entre les $X_\beta$ sur lesquelles nous revenons en détail dans le paragraphe 3.\\
$ $\\
Ces formules ne sont pas connues de manière explicite, de sorte que la description de ces diagrammes semble inaccessible dans le cas général. Nous nous limitons donc dans ce papier au cas d'une décomposition réduite de $w_0$ associée à un bon ordre sur l'ensemble $\Phi^+$ des racines positives, au sens de G. Lusztig \cite{Lu2} (cf paragraphe 2). Dans cette situation, nous disposons des formules de redressement entre deux variables $X_\beta$, $X_{\beta'}$ lorsque $\beta$ et $\beta'$ engendrent un plan admissible au sens de G. Lusztig (cf paragraphe 3.4). Cela nous permet de construire une famille de conditions nécessaires appelées \textit{contraintes} pour qu'un diagramme $\Delta$ soit un diagramme de Cauchon (cf paragraphe 5.1). Dans le paragraphe 5.2, nous démontrons que ces conditions sont en fait nécessaires et suffisantes (cf théorème \ref{reciproque contrainte}).\\
$ $\\
Dans le paragraphe 6, nous utilisons ce théorème pour donner une description explicite de ces contraintes et de ces diagrammes pour un choix particulier de la décomposition réduite de $w_{0}$, pour tous les types possibles du système de racines associé à $\mathfrak{g}$. On constate que dans tous les cas, le nombre de diagrammes est égal au cardinal $| W |$ du groupe de Weyl. Comme les strates ne dépendent pas du choix de la décomposition réduite de $w_{0}$, cela entraîne que le nombre de strates est toujours égal à $| W |$. Si on fait l'hypothèse supplémentaire que $q$ transcendant, ce résultat a aussi été démontré en utilisant d'autres méthodes par M. Gorelik (cf \cite{MR1754232} et \cite{dumas}).\\
$ $ \\
 Par ailleurs, dans un prochain article, nous démontrerons que les diagrammes de Cauchon $\Delta$ fournissent les sous-expressions positives \textbf{w}$^\Delta$  de $w_0$ au sens de R. Marsh et K. Rietsch (\cite{MR2058727}). En particulier, si $w_0$ a une écriture réduite $ w_{0} = s_{\alpha_1} \circ ... \circ s_{\alpha_N}$ et si, pour chaque diagramme $\Delta = \{\beta_{i_1} < ... < \beta_{i_t} \}$ on note $w^{\Delta} = s_{\alpha_{i_1}} \circ ... \circ s_{\alpha_{i_t}}$, on a :
 $ $ \\
\begin{itemize}
\item[$\bullet$] Si $\Delta$ est un diagramme de Cauchon, l'écriture de $w^\Delta$ ci-dessus est réduite.
\item[$\bullet$] L'application $\Delta \rightarrow w^{\Delta}$ est une bijection de l'ensemble $\mathcal{D}$ des diagrammes de Cauchon sur $W$.
\end{itemize}

\section{Systèmes de racines : rappels et compléments}
\subsection{Résultats classiques sur les systèmes de racines}
Soit $\mathfrak{g}$ une algèbre de Lie simple complexe. Nous commençons par introduire les notations utilisées pour les systèmes de racines ; nous suivrons celles de \cite[chap4]{J}.\\
\par
\begin{nota} 
\begin{itemizedot}
\item On note $\Phi$ un système de racines et E = Vect($\Phi$) (dim E = n). Quand on fixe une base $\Pi := \{\alpha_{1}, ..., \alpha_{n}\}$ de $\Phi$, il existe une décomposition $\Phi = \Phi^+   \sqcup \Phi^-$, où $\Phi^+$ (resp. $\Phi^-$) désigne, comme d'habitude, l'ensemble des racines positives (resp. négatives).
\item On note W le groupe de Weyl associé au système de racines $\Phi$, il est engendré par les réflexions $s_{\alpha_{i}} (:= s_{i}) , 1 \leq i \leq n$. L'élément de plus grande longueur de W est noté $w_{0}$.
\end{itemizedot}
\end{nota} \begin{dfn}
Un système de racines $\Phi$ est \emph{réductible} si $\Phi = \Phi_1   \sqcup \Phi_2$ où $\Phi_1$ et $\Phi_2$ sont deux systèmes de racines orthogonaux. Sinon $\Phi$ est dit \emph{irréductible}.\\  
\end{dfn}
Rappelons qu'il y a une correspondance bijective entre les systèmes de racines irréductibles et les algèbres de Lie simples complexes de dimension finie. On dira que $\mathfrak{g}$ est d'un type donné si le système de racines associé   $\mathfrak{g}$ est de ce type. Les définitions et résultats qui suivent proviennent de \cite{Lu2}.
\begin{dfn}
Soient $\Pi =\{\alpha_1, \alpha_2, ... , \alpha_n\}$ une base de $\Phi$ et $j$ un entier de $\llbracket 1, n \rrbracket$.
\begin{enumerate}
\item \emph{La colonne j} est l'ensemble $C_j := \{\beta \in \Phi^+ | \beta = k_1 \alpha_1 + \ldots + k_j \alpha_j, k_i \in \mathbb{N}, k_j \neq 0\}$;
\item On dit qu'une racine $\beta = k_1 \alpha_1 + \ldots + k_j \alpha_j \in C_{j}$ est \emph{ordinaire} si $k_j =1$; elle est dite \emph{exceptionnelle} si $k_j = 2$;
\item Une colonne  $C_j$ est appelée \emph{ordinaire} si toute racine  $\beta$ de $C_j$ est ordinaire; cette colonne est dite \emph{exceptionnelle}
si toute racine $\beta$ de $C_j$ est ordinaire sauf une seule racine note $\beta_{ex}$ qui est exceptionnelle.
\end{enumerate}
\end{dfn}
\begin{dfn}
Le numérotation $\Pi =\{\alpha_1, \alpha_2, ... , \alpha_n\}$ est \emph{bonne} si toute colonne $C_j$ est ordinaire ou exceptionnelle.
\end{dfn}
\begin{ex}[Le cas $G_2$]\label{exG2}
Le système de racines de type $G_{2}$ est de rang 2, il a 2 racines simples $\alpha_1$ et $\alpha_2$ avec $\|\alpha_{2}\| = \sqrt{3} \|\alpha_{1}\|$. $\Pi =\{\alpha_1, \alpha_2 \}$ est une base de ce système de racine. La numérotation $\Pi =\{\alpha_1, \alpha_2 \}$ est bonne car $C_1 = \{ \alpha_1 \}$ est ordinaire et $C_2 = \{\alpha_2, \alpha_1 + \alpha_2, 2\alpha_1 + \alpha_2, 3\alpha_1 + \alpha_2, 3\alpha_1 + 2\alpha_2  \}$ est exceptionnelle.\\
Par contre, la numérotation $\Pi =\{\alpha_2 , \alpha_1 \}$ n'est pas bonne. En effet, pour cette numérotation, $C_1 = \{ \alpha_2 \}$ est ordinaire mais $C_2 = \{\alpha_1, \alpha_2 + \alpha_1,\alpha_2 + 2\alpha_1, \alpha_2 +\mathbf{ 3}\alpha_1, 2\alpha_2 + \mathbf{3}\alpha_1  \}$ n'est ni ordinaire ni exceptionnelle.
\end{ex}
\begin{prop}
Soit $\mathfrak{g}$ une algèbre de Lie simple de dimension finie. Les numérotations ci-dessous du système des racines simples $\Pi$ sont des exemples de bonnes numérotations.
\begin{itemizedot}
\item Si $\mathfrak{g}$ est de type $A_n$, de diagramme de Dynkin : $ \alpha_1 - \alpha_2 - \cdots - \alpha_{n-1} - \alpha_n $,\\
\[\Pi = \{\alpha_1 , \alpha_2 , \cdots , \alpha_{n-1} , \alpha_n\}.\] 
\item Si $\mathfrak{g}$ est de type $B_n$, de diagramme de Dynkin : $\alpha_1 \Leftarrow \alpha_2 - \cdots - \alpha_{n-1} - \alpha_n $,\\
\[\Pi = \{\alpha_1 , \alpha_2 , \cdots , \alpha_{n-1} , \alpha_n\}.\]
\item Si $\mathfrak{g}$ est de type $C_n$, de diagramme de Dynkin : $\alpha_1 \Rightarrow \alpha_2 - \cdots - \alpha_{n-1} - \alpha_n $,\\
\[\Pi = \{\alpha_1 , \alpha_2 , \cdots , \alpha_{n-1} , \alpha_n\} .\]
\item Si $\mathfrak{g}$ est de type $D_n$, de diagramme de Dynkin : 
\scalebox{.8}{$\begin{array}{ccccccccccc}
\alpha_1 && & &   & &&&  &     &   \\ 
& \diagdown & & &   & &&&  &     & \\
&& \alpha_3 &\text{---~}& \alpha_{4} &\text{---~}& \cdots &\text{---~}& \alpha_{n-1}&\text{---~} & \alpha_n \\
& \diagup & & &   & &&&  &     & \\
\alpha_2&& & &  &  &&&&  &
\end{array}$},
\[ \Pi = \{\alpha_1 , \alpha_2 , \cdots , \alpha_{n-1} , \alpha_n\}.\]
\item Si $\mathfrak{g}$ est de type $G_2$, de diagramme de Dynkin : $ \alpha_1 \Lleftarrow \alpha_2 $, \[\Pi = \{\alpha_1 , \alpha_2 \}.\]
\item Si $\mathfrak{g}$ est de type $F_4$, de diagramme de Dynkin : $\alpha_1 \ \text{---~} \alpha_2 \Rightarrow \alpha_3 \ \text{---~} \alpha_4$, \\
\[\Pi = \{\alpha_4 , \alpha_3, \alpha_2, \alpha_1 \}.\]
\item Si $\mathfrak{g}$ est de type $E_6$, de diagramme de Dynkin : \scalebox{.8}{$\begin{array}{ccccccccc}
 					& 					&   				 & 					 & \alpha_2	&&  &     &   \\ 
 					& 					&   				 & 					 & |					&&  &     & \\
 \alpha_1 &\text{---~}& \alpha_{3} &\text{---~}& \alpha_4 &\text{---~}& \alpha_{5}&\text{---~} & \alpha_6 \\
\end{array}$}, 
\[\Pi = \{\alpha_2 , \alpha_5, \alpha_{4}, \alpha_3, \alpha_1, \alpha_6 \}.\]
\item Si $\mathfrak{g}$ est de type $E_7$, de diagramme de Dynkin :\scalebox{.8}{$\begin{array}{ccccccccccc}
 					& 					&   				 & 					 & \alpha_2	&&  &     &&&  \\ 
 					& 					&   				 & 					 & |					&&  &     & &&\\
 \alpha_1 &\text{---~}& \alpha_{3} &\text{---~}& \alpha_4 &\text{---~}& \alpha_{5}&\text{---~} & \alpha_6 &\text{---~}& \alpha_7 \\
\end{array}$},
\[\Pi = \{\alpha_2 , \alpha_5, \alpha_{4}, \alpha_3, \alpha_1, \alpha_6, \alpha_7 \}.\]
\item Si $\mathfrak{g}$ est de type $E_8$, de diagramme de Dynkin : \scalebox{.8}{$\begin{array}{ccccccccccccc}
 					& 					&   				 & 					 & \alpha_2	&&  &     &&& &&  \\ 
 					& 					&   				 & 					 & |					&&  &     &&& &&\\
 \alpha_1 &\text{---~}& \alpha_{3} &\text{---~}& \alpha_4 &\text{---~}& \alpha_{5}&\text{---~} & \alpha_6 &\text{---~}& \alpha_7& \text{---~}& \alpha_8\\
\end{array}$},
\[\Pi = \{\alpha_2 , \alpha_5, \alpha_{4}, \alpha_3, \alpha_1, \alpha_6, \alpha_7, \alpha_8 \}.\]
\end{itemizedot}
\end{prop}
Les colonnes correspondant à ces numérotations seront donnes explicitement dans la partie 3 et on constatera que chaque colonne est bien, soit ordinaire soit exceptionnelle. \textbf{Dans tout ce qui suit, la numérotation choisie sur $\Pi$ sera toujours bonne et  partir de la section 6, on prendra explicitement l'ordre décrit dans la proposition ci-dessus.}

\subsection{Ordre de Lusztig}
\begin{nota}
\begin{itemizedot}
\item Pour $\beta = k_1 \alpha_1 + \ldots + k_j \alpha_j \in C_j$, on appelle \emph{hauteur de }$\beta$ l'entier $h(\beta) := k_1 + \cdots + k_j$ et on appellera \emph{hauteur de Lusztig de} $\beta$, le nombre rationnel $h'(\beta) := \frac{1}{k_j} h(\beta)$. \cite[section 4.3]{Lu2}
\item Si $t \in h'(C_j)$, l'ensemble $B^{j,t}:=\{\beta \in C_j | h'(\beta)=t\}$ est appelé \emph{la boîte de hauteur t dans la colonne} $C_j$. 
\end{itemizedot}
\end{nota}
\begin{rmq}
 \[\displaystyle C_{j} = \bigsqcup_{t \in \mathbb{N}^{*}} B^{j,t}.\]
\end{rmq}
\begin{dfn}[Ordre de Lusztig sur $\Phi^+$]
On définit un ordre partiel sur $\Phi^+$ comme suit :\\
Soient $\beta_1$ et $\beta_2$ deux racines de $\Phi^+$,
\begin{itemizedot}
\item Si $\beta_1 \in C_{j_1}$ et $\beta_2 \in C_{j_2}$ avec $j_1 < j_2$,  alors $\beta_1 < \beta_2$.
\item Si $\beta_1$ et $\beta_2$ sont dans la même colonne $C_j$ et si $h'(\beta_2) < h'(\beta_1)$, alors $\beta_1 < \beta_2$.
\end{itemizedot}
On peut raffiner l'ordre partiel précèdent en un ordre total en choisissant un ordre arbitraire à l'intérieur des boîtes. Comme il y a plusieurs façons d'ordonner chaque boîte on dira que l'ordre obtenu est \emph{"un" ordre de Lusztig}.
\end{dfn}
\begin{obs}
\begin{itemizedot}
\item $\alpha_{j}$ est la plus grande racine de $C_{j}$ pour l'ordre ci-dessus.
\item Les racines positives d'une même boîte sont consécutives pour un tel ordre : $B^{j,t} = \{\beta_{p}, \beta_{p+1}, ..., \beta_{p+l}\}$.
\end{itemizedot}
\end{obs}

\begin{prop}\label{prop2.8}
Soient $C_j$ une colonne exceptionnelle et $\beta_{ex}$ sa racine exceptionnelle.
\begin{enumerate}

\item $\beta_{ex} \bot (C_1 \sqcup ... \sqcup C_{j-1})$
\item Si $D = <\beta_{ex}>$ et si $s_D $ est la symétrie orthogonale par rapport  D, on a :
\begin{itemizedot}
\item $s_D(C_j) = C_j$ et  $\forall \beta \in C_{j}\setminus\{\beta_{ex}\}$ on a  $\beta + s_{D}(\beta) = \beta_{ex}$.
\item Pour toute boîte $B^{j,t}$ différente de la boîte contenant $\beta_{ex}$, $s_D$ transforme $B^{j,t}$ en $B^{j,h(\beta_{ex})-t}$.
\end{itemizedot}
\end{enumerate}
\end{prop}
\begin{preuve}
\begin{enumerate}
\item Soit $\beta \in C_1 \cup ... \cup C_{j-1}$. Si $\beta$ n'est pas orthogonale  $\beta_{ex}$, alors $s_{\beta}(\beta_{ex}) = \beta_{ex} + k\beta \ (k \in \mathbb{Z}\setminus\{0\})$ est une racine de $C_j$ dont le coefficient en $\alpha_{j}$ est égal à 2. Ceci contredit l'unicité de la racine exceptionnelle. 
\item Observons que $s_D = -s_{\beta_{ex}}$, de sorte que $s_D(\Phi) = \Phi$.
\begin{itemizedot}
\item Soit $\beta$ une racine non exceptionnelle de $C_j$. On peut écrire $$\beta = a_1 \alpha_1 + ... + a_{j-1} \alpha_{j-1} + \frac{1}{2}  \beta_{ex} \ (a_{i} \in \mathbb{Q}).$$ Il résulte de 1. que $s_D(\beta) =   - a_1 \alpha_1  ...  - a_{j-1} \alpha_{j-1} + \frac{1}{2}  \beta_{ex} = \beta_{ex} - \beta$. Ceci est une racine d'après l'observation ci-dessus. Elle est dans $C_{j}$ puisque $\beta$ est dans $C_{j}\setminus\{\beta_{ex}\}$.\\
\item Par le point précèdent, $s_{D}$ transforme deux éléments de $B^{j,t}$ en deux racines de même hauteur. On en déduit (au moyen de l'involutivité de $s_{D}$) que $s_{D}(B^{j,t})$ est une boîte $B^{j,s}$. La formule $t + s = h(\beta_{ex})$ résulte immédiatement du premier point.
\end{itemizedot}
\end{enumerate}
\end{preuve}
\begin{dfn}
Le support d'une racine $\beta = a_1 \alpha_1 + ... + a_{j} \alpha_{j} \in C_{j}$ est l'ensemble \emph{Supp} $(\beta) := \{\alpha_{i} \in \Pi | a_{i} \neq 0 \}$. En particulier, pour $\beta \in C_{j}$, on a $\textrm{\emph{Supp}}(\beta) \subset \{1, ..., j\}$. 
\end{dfn}
On va à présent démontrer la
\begin{prop}\label{boîteex}
Supposons qu'il existe une colonne exceptionnelle $C_{j}$ ($1 < j \leq n$) et notons $\beta_{ex}$ sa racine exceptionnelle. Alors $h'(\beta_{ex}) \notin \mathbb{N}$ de sorte que $\beta_{ex}$ est seule dans sa boîte.
\end{prop}
\begin{preuve}
On note $\Pi_{j} = \{\alpha_{1},..., \alpha_{j}\}$ et $\Phi_{j} = \Phi \cap \text{Vect}(\Pi_{j})$. On vérifie facilement que $\Phi_{j}$ est un système de racines de base $\Pi_{j}$ et que $\Phi_{j}^{+} = \Phi^{+} \cap \text{Vect}(\Pi_{j})$. \\
\par
Considérons d'abord le cas où $\Phi_{j}$ est irréductible. On a alors \\
\par
\textbf{Observation 1.} Si $\beta$ est une racine de $\Phi_{j}^{+}$ de hauteur maximale alors $\beta \in C_{j}$.\\
On suppose $\beta \in C_{i}$ avec $i < j$. Dans le diagramme de Dynkin de $\Phi_{j}$ qui est connexe ($\Phi_{j}$ est irréductible), on peut construire un chemin de $\alpha_{i}$ à $\alpha_{j}$. On note ce chemin $P = (\alpha_{i_{1}}, ..., \alpha_{i_{t}})$, où $i_{1} = i$ et $i_{t} = j$. On sait que $\alpha_{i} \in Supp(\beta)$ et que $\alpha_{j} \notin Supp(\beta)$. Ainsi il existe un plus petit indice $l$ tel que $\alpha_{i_{l}} \in Supp(\beta)$ et $\alpha_{i_{l+1}} \notin Supp(\beta)$. Ainsi, pour tout $ \alpha \in Supp(\beta)$ on a $\langle \alpha , \alpha_{i_{l+1}} \rangle  \leq 0$ et, puisque $\alpha_{i_{l}}$ et $\alpha_{i_{l+1}}$ sont deux éléments consécutifs de P, $\langle \alpha_{i_{l}} , \alpha_{i_{l+1}} \rangle < 0$. De l, $\langle \beta, \alpha_{i_{l+1}} \rangle < 0$ donc $\beta + \alpha_{i_{l+1}} \in \Phi_{j}^+$ ce qui contredit la maximalité de la hauteur $\beta$.\\
\par
\textbf{Observation 2.} $\beta_{ex}$ est la racine de hauteur maximale dans $\Phi_{j}$.\\

Soit $\beta$ une racine de hauteur maximale dans $\Phi_{j}$. Supposons $\beta \neq \beta_{ex}$. Par l'observation précédente, $\beta \in C_{j}$ et, par la proposition \ref{prop2.8}, $\beta_{ex} = \beta + s_{D}(\beta)$ est une somme de deux racines positives donc de hauteur plus grande que $\beta$. Donc $\beta$ est nécessairement égale  $\beta_{ex}$.\\
\par
L'existence d'une racine exceptionnelle entraine que $\Phi_{j}$ n'est pas de type $A_{j}$. Par suite $\Phi_{j}$ est de type $B_{j}, C_{j}, D_{j}, E_{6}, E_{7}, E_{8}, F_{4}$ ou $G_{2}$ et, en se référant aux planches de \cite{MR0240238}, on voit que la hauteur de la racine de plus grande hauteur est toujours impair. Par l'observation 2. $\beta_{ex}$ est de hauteur impair, et donc $h'(\beta_{ex}) \notin \mathbb{N}$.\\
\par
Supposons à présent $\Phi_{j}$ réductible. Notons $\Gamma_{j}$ le diagramme dont les sommets sont $\alpha_{1}, ..., \alpha_{j}$, et dont les arêtes sont celles qui proviennent du diagramme de Dynkin de $\Phi$ et notons $\Pi'$ la composante connexe de $\alpha_{j}$ dans $\Gamma_{j}$; 
\[\Pi':=\{\alpha_{i} \in \Pi_{j} | \ \textrm{il existe un chemin dans} \ \Gamma_{j} \ \textrm{reliant} \ \alpha_{i} \ \textrm{et} \ \alpha_{j}\}.\]
Notons $\Phi' = \Phi \cap \text{Vect}(\Pi')$. C'est un système de racines de base $\Pi'$ et $\Phi'^+ = \Phi^+ \cap \text{Vect}(\Pi')$.\\
\par
\textbf{Observation 3.} $C_{j} \subset \Phi'^+$.\\
\par
Sinon, il existe des racines dans $ C_{j} \setminus \Phi'^{+}$. Si $\beta$ est une telle racine, son support contient des racines simples qui sont dans $\Pi_{j}\setminus \Pi'$. Comme le support de $\beta$ contient aussi $\alpha_{j} \in \Pi'$, on peut écrire $\beta = u + v$ avec $ u = \alpha_{i_{1}} + ... + \alpha_{i_{l}}$  support dans $\Pi_{j}\setminus \Pi'$ et $v = \alpha_{i_{l+1}} + ... + \alpha_{i_{p}}$  support dans $\Pi'$. Choisissons $\beta$ de manière que l'entier $l$ ainsi défini soit minimal.
\begin{itemizedot}
\item Si $l = 1, \beta = \alpha_{i_{1}} + v$. Comme $\alpha_{i_{1}} \notin \Pi'$, il n'existe aucune liaison entre $\alpha_{i_{1}}$ et les éléments du support de $v$. Donc $s_{i_{1}}(\beta) = -\alpha_{i_{1}} + v \in \Phi$, ce qui est impossible puisque les coordonnées dans $\Pi$ de cette racine ne sont pas toutes de même signe. 
\item Donc $l \geq 2$. Comme $\langle u, u \rangle > 0$, il existe une racine simple du support de $u$, par exemple $\alpha_{i_{l}}$, vérifiant $\langle u, \alpha_{i_{l}} \rangle > 0$. Comme ci dessus, on a :
$$\langle v, \alpha_{i_{l}} \rangle = 0 \Rightarrow \langle \beta, \alpha_{i_{l}} \rangle > 0 \Rightarrow \beta' =\beta-\alpha_{i_{l}} \in C_{j}\setminus \Phi'^{+},$$ ce qui contredit la minimalité de $l$.
\end{itemizedot}
On a donc bien $C_{j} \subset \Phi'^+$.\\
\par
Ainsi $C_{j}$ est une colonne exceptionnelle de $\Phi'$ qui est irréductible par construction. L'étude ci-dessus montre que sa racine exceptionnelle $\beta_{ex}$ vérifie $h'(\beta_{ex}) \notin \mathbb{N}$.\\
\end{preuve}

On peut maintenant montrer la
\begin{prop}
"<" est un ordre convexe sur  $\Phi^+$
\end{prop}
\begin{preuve}
Soient $\beta_1 < \beta_2$ deux racines positives telles que $\beta_1 + \beta_2 \in \Phi^+$. 
\begin{itemizedot}
\item Si les deux racines $\beta_{1}$ et $\beta_{2}$ ne sont pas dans la même colonne, alors $\beta_{1} + \beta_{2}$ est dans la même colonne que $\beta_{2}$. Dans ce cas, ni $\beta_{2}$, ni $\beta_{1} + \beta_{2}$ ne sont exceptionnelle et on a :
\[h'(\beta_{1} + \beta_{2}) = h(\beta_{1} + \beta_{2}) = h(\beta_{1}) + h'(\beta_{2}) > h'(\beta_{2}).\]
On en déduit que $\beta_1 <\beta_1 + \beta_2 < \beta_2$.\\
\item Si les deux racines sont dans la même colonne, alors $\beta_{1} + \beta_{2}$ est une racine exceptionnelle. D'après la proposition \ref{prop2.8}, on a $h'(\beta_{1} + \beta_{2}) = \frac{h'(\beta_{1})+h'(\beta_{2})}{2}$. La proposition \ref{boîteex} exclus le cas $h'(\beta_{1} + \beta_{2}) = h'(\beta_{1}) = h'(\beta_{2})$ car la racine exceptionnelle est seule dans sa boîte. On a donc $h'(\beta_{1}) > h'(\beta_{1} + \beta_{2}) > h'(\beta_{2})$ ce qui implique $\beta_1 <\beta_1 + \beta_2 < \beta_2$.
\end{itemizedot}
\end{preuve}
Si on considère une décomposition réduite de $w_0 =  s_{i_{1}} \circ s_{i_2} \circ ... \circ s_{i_N}$ (l'élément de plus grande longueur du groupe de Weyl), on sait (cf, par exemple, \cite[I.5.1]{MR1898492}) que $\beta_j :=s_{i_1} \circ s_{i_2} \circ ... \circ s_{i_{j-1}}(\alpha_{i_{j}})$ décrit $\Phi^+$ lorsque $j$ décrit $\llbracket 1 , N \rrbracket$. Pour chaque entier $j \in \llbracket 1, N \rrbracket$, on dit que $\alpha_{i_{j}}$ est \emph{la racine simple associée à la racine positive} $\beta_{j}$.\\
On ordonne alors $\Phi^+$ en posant $\beta_i \prec \beta_j$ lorsque $i<j$. On dit que "$\prec$" est \emph{l'ordre associé  la décomposition réduite de} $w_0 = s_{i_{1}} \circ s_{i_2} \circ ... \circ s_{i_N}$.\\
Dans \cite[Théorème et remarque page 662]{Papi}, il est démontré que cet ordre est convexe et que l'on définit, de cette manière, une correspondance biunivoque entre les décompositions réduites de $w_0$ et les ordres convexes sur $\Phi^+$.\\
\textbf{Ainsi, l'ordre "<"de Lusztig tant convexe, il existe une unique décomposition réduite de} \\
$\mathbf{w_0 = s_{i'_{1}} \circ s_{i'_2} \circ ... \circ s_{i'_N}}$ \textbf{dont l'ordre associé est "<". Dans ce texte, on choisira systématiquement cette décomposition pour} $\mathbf{w_0}$.\\
\par
La prochaine proposition provient de \cite[section 4.3]{Lu2} et précise comment sont disposes les racines  l'intérieur des boîtes.
\begin{prop}\label{boîte}
 A l'intérieur de chaque boîte ne contenant pas la racine exceptionnelle, les racines sont orthogonales 2 à 2. De plus, les racines simples associées aux racines d'une boîte donnée sont orthogonales 2 à 2.
\end{prop}
\begin{preuve}
Le cas d'un système de type $G_{2}$ se déduit simplement de l'étude de l'exemple \ref{exG2}. On suppose donc que $\mathfrak{g}$ est une algèbre de Lie simple de dimension finie et de type différent de $G_2$. Soient $\beta_{1}$ et $\beta_{2}$ deux racines consécutives d'une boîte B de la colonne $C_{j}$. On note $\alpha_{i_{1}}$ et $\alpha_{i_{2}}$ les racines simples respectivement associes  $\beta_{1}$ et $\beta_{2}$.\\
Supposons que $\alpha_{i_{1}}$ n'est pas orthogonale  $\alpha_{i_{2}}$, de sorte que $\lambda = -<\alpha_{i_{1}}^{\vee},\alpha_{i_{2}}> = 1$ ou $2$.\\
On peut donc écrire $\beta_{2}  = w \circ s_{i_{1}} (\alpha_{i_{2}})
 = w(\lambda \alpha_{i_{1}} + \alpha_{i_{2}})
 = \lambda \beta_{1} + w(\alpha_{i_{2}})$.
Comme  $w(\alpha_{i_{2}}) \in \Phi$, on a nécessairement $\lambda = 2$, sinon $h(w(\alpha_{i_{2}})) = h(\beta_{2}) - h(\beta_{1}) = 0 $, ce qui est absurde.\\
Mais alors $\gamma = - w(\alpha_{i_{2}}) = 2\beta_{1} - \beta_{2} \in C_{j}$ et $h(\gamma) =  2h(\beta_{1})- h(\beta_{2}) = h(\beta_1)$. $\beta_1$ et $\beta_2$ sont deux racines positives distinctes donc non colinéaires. La trace $\Phi'$ de $\Phi$ sur le plan Vect$(\beta_1, \beta_2)$ est donc un système de racines de rang 2 contenant $\beta_1$, $\beta_2$, $\gamma$ et leurs opposées. L'égalité $2\beta_{1}  = \gamma  + \beta_{2}$ permet d'affirmer que $\Phi'$ est de type $B_2$ et qu'on est dans la situation ci-dessous :
\begin{center}
\begin{pgfpicture}{-0.5cm}{-1cm}{2.9cm}{2.2cm}%
\pgfsetroundjoin \pgfsetroundcap%
\pgfsetlinewidth{0.8pt} 
\pgfsetfillcolor{black}
\pgfxyline(0.1,1.3)(2.3,1.3)
\pgfxyline(2.3,1.3)(2.3,-0.9)
\pgfxyline(2.3,-0.9)(0.1,-0.9)
\pgfxyline(0.1,-0.9)(0.1,1.3)
\pgfxyline(1.2,0.2)(2.2,1.2)
\pgfmoveto{\pgfxy(2.0712,1.1655)}\pgflineto{\pgfxy(2.2,1.2)}\pgflineto{\pgfxy(2.1655,1.0712)}\pgflineto{\pgfxy(2.1592,1.1592)}\pgfclosepath\pgffillstroke
\pgfxyline(1.2,0.2)(1.2,1.2)
\pgfmoveto{\pgfxy(1.1333,1.0845)}\pgflineto{\pgfxy(1.2,1.2)}\pgflineto{\pgfxy(1.2667,1.0845)}\pgflineto{\pgfxy(1.2,1.1423)}\pgfclosepath\pgffillstroke
\pgfxyline(1.2,0.2)(0.2,1.2)
\pgfmoveto{\pgfxy(0.2345,1.0712)}\pgflineto{\pgfxy(0.2,1.2)}\pgflineto{\pgfxy(0.3288,1.1655)}\pgflineto{\pgfxy(0.2408,1.1592)}\pgfclosepath\pgffillstroke
\pgfputat{\pgfxy(2.2,1.6)}{\pgfnode{rectangle}{center}{\small $\gamma$}{}{\pgfusepath{}}}
\pgfputat{\pgfxy(1.2,1.6)}{\pgfnode{rectangle}{center}{\small $\beta_{1}$}{}{\pgfusepath{}}}
\pgfputat{\pgfxy(0.2,1.6)}{\pgfnode{rectangle}{center}{\small $\beta_{2}$}{}{\pgfusepath{}}}
\end{pgfpicture} 
\end{center}
Il en résulte que $\gamma - \beta_{1} \in \Phi$, avec $h(\gamma - \beta_1) = h(\gamma) - h(\beta_1) = 0 $. Ce qui est impossible. On a donc $\alpha_1 \ \bot \ \alpha_2$.\\
On en déduit que $<\beta_{1},\beta_{2}> \ =  \ <w(\alpha_{i_{1}}),w(s_{i_{1}}(\alpha_{i_{2}}))> = <\alpha_{i_{1}},s_{i_{1}}(\alpha_{i_{2}})> \ = \ <\alpha_{i_{1}},\alpha_{i_{2}}> \ = 0$.
\end{preuve}
\par
\begin{convsanss}
Pour $j \in \llbracket 1, n \rrbracket$, on note $\delta_j$ la première racine de $C_j$. On rappelle que $\alpha_j$ est la dernière racine de $C_j$.
\end{convsanss}
\begin{prop}
$\delta_j$ et $\alpha_j$ sont seules dans leur boîtes.
\end{prop}
\begin{preuve}
La racine $\alpha_{j}$ est seule dans sa boîte car c'est la seule racine de $C_{j}$ de hauteur égale à 1.\\

Pour démontrer que $\delta_j$ est seule dans sa boîte, on va se servir du 
\begin{lemme}
Soit $1 \leq l \leq N$ et $1 \leq m \leq n$. On pose $\Pi_m := \{\alpha_1, ..., \alpha_m \}$. Si $\beta_l = s_{i_1}...s_{i_{l-1}}(\alpha_{i_l})$ est dans la colonne $C_m$, alors $\alpha_{i_j} \in \Pi_m$ pour $j \in \llbracket 1, l \rrbracket$. 
\end{lemme}
\begin{preuve}
On raisonne par récurrence sur $l$.\\
Si $l = 1$ : $\beta_1 = \alpha_1 \in \Pi_1$.\\
Si $l \geq 2$, $\beta_{l-1}$ est dans la colonne $C_m$ ou $C_{m-1}$. Par l'hypothèse de récurrence, on en déduit que $\alpha_{i_t} \in \Pi_m$ (ou $\alpha_{i_t} \in \Pi_{m-1} \subset \Pi_m$) pour $t \in \llbracket 1, l-1 \rrbracket$. On observe que $\beta_l = s_{i_1}...s_{i_{l-1}}(\alpha_{i_l}) = \alpha_{i_l} + n_{l-1}\alpha_{i_{l-1}} + ... + n_1 \alpha_{i_{1}}$ où chaque $n_t$ est dans $\mathbb{Z}$. Comme $\beta_l \in C_m$, nécessairement $\alpha_{i_l} \in \Pi_m$.  
\end{preuve}
\textbf{Retour à la démonstration de la proposition :}\\
Il existe un entier $1 \leq l \leq N$ tel que $\delta_{j} = \beta_{l} = s_{i_{1}} \circ s_{i_2} \circ ... \circ s_{i_{l-1}}(\alpha_{i_{l}})$. comme ci-dessus, $\beta_l =  \alpha_{i_l} + n_{l-1}\alpha_{i_{l-1}} + ... + n_1 \alpha_{i_{1}} \ (n_t \in \mathbb{Z})$ avec, puisque $\beta_{l-1} \in C_{j-1}$, $\alpha_{i_1}, ..., \alpha_{i_{l-1}}$ dans $ \Pi_{j-1}$. Comme $\beta_l \in C_j$, on en déduit que $\alpha_{i_{l}} = \alpha_{j}$.\\
Si $\delta_{j} (=\beta_{l}) $ n'est pas seule dans sa boîte, alors $\beta_{l+1}$ est aussi dans cette boîte et on a (proposition \ref{boîte}) $\alpha_{i_{l}} \bot \alpha_{i_{l+1}}$. Par le lemme précédent, on en déduit que $\alpha_{i_{l+1}} \in \Pi_{j} \setminus \{\alpha_j\} = \Pi_{j-1}$ et $\beta_{l+1} = s_{i_{1}} \circ s_{i_2} \circ ...\circ s_{i_{l-1}} \circ s_{i_{l}}(\alpha_{i_{l+1}}) = s_{i_{1}} \circ s_{i_2} \circ ...\circ s_{i_{l-1}} (\alpha_{i_{l+1}}) = \alpha_{i_{l+1}} + n'_{l-1}\alpha_{i_{l-1}} + ... + n'_1 \alpha_{i_{1}} \ (n'_t \in \mathbb{Z})$, ce qui contredit l'hypothèse  $\beta_{l+1} \in C_{j}$.
\end{preuve}

Rappelons le résultat classique suivant (cf par exemple, \cite[lemme 9.4]{Humphreys}).
\begin{lemme}\label{lemme2.9}
Soient $\beta$ et $\delta$ deux racines distinctes de $\Phi^+$ telles que $\langle \beta,\delta \rangle \neq 0$
\begin{itemizedot}
\item Si $\langle \beta,\delta \rangle > 0$, alors $\beta -\delta \in \Phi$.
\item Si $\langle \beta,\delta \rangle < 0$, alors $\beta + \delta \in \Phi$.
\end{itemizedot}
\end{lemme}
\begin{prop}\label{prop2.10}
Soit $\beta$ une racine ordinaire de la colonne $C_j$. Notons (comme dans la démonstration de la proposition \ref{boîteex}) $\Gamma_{j}$ le diagramme dont les sommets sont $\alpha_{1}, ..., \alpha_{j}$, et dont les arêtes sont celles qui proviennent du diagramme de Dynkin de $\Phi$. On note $\Omega_j$ la composante connexe de $\alpha_j$ dans $\Gamma_j$ .
\begin{enumerate}
\item Si $\beta \neq \alpha_{j}$, alors il existe $\epsilon \in \{\alpha_{1},...,\alpha_{j-1}\}$ telle que $\beta - \epsilon \in C_j$.
\item Supp $\beta \subset \Omega_j$.
\item Si $\beta \neq \delta_{j}$, alors il existe $\epsilon \in \{\alpha_{1},...,\alpha_{j-1}\}$ telle que $\beta + \epsilon \in C_j$.
\end{enumerate} 
\end{prop}
\begin{preuve}
On vérifie que cette proposition est vraie dans le cas $G_2$ à l'aide de la description des colonnes donnée dans l'exemple \ref{exG2}. On suppose maintenant que le système de racines considéré n'est pas de type $G_2$.\\
Dans cette démonstration, on note $\Pi_{j} = \{\alpha_{1},..., \alpha_{j}\}$.\\
\begin{enumerate}
\item S'il existe $\epsilon \in \Pi_{j} \setminus \{\alpha_j\}$ tel que $\langle \beta , \epsilon \rangle > 0$, on conclut avec le lemme \ref{lemme2.9}. On suppose donc que $ \forall \epsilon \in \Pi_{j} \setminus \{\alpha_j\}$, $\langle \beta , \epsilon \rangle \leq 0$.\\
 Comme $\langle \beta, \beta \rangle > 0$ alors, $\langle \beta, \alpha_{j} \rangle > 0$ et $\beta - \alpha_{j} = \gamma_{1} \in \Phi^{+}$ avec $\alpha_{j} \notin$ Supp$(\gamma_{1})$. On peut donc écrire $\beta = \beta_{1} + \gamma_{1}$ avec $\beta_{1} = \alpha_j \in C_{j}$ et, puisque $\beta$ est ordinaire, $\alpha_{j} \notin Supp(\gamma_{1})$. Puisque $\beta \neq \alpha_j$, on a $h(\beta) \geq 2$ et on se propose de montrer par récurrence que, pour chaque entier $i \in \llbracket 1, h(\beta)-1 \rrbracket$, on a $\beta = \beta_i + \gamma_i$ avec $\beta_i \in C_j$, $\gamma_i \in \Phi^+$ et $h(\beta_i) = i$.\\
 \begin{itemizedot}
\item Comme $h(\beta_1) = h(\alpha_j) = 1$, on a le résultat au rang $i = 1$.
\item Supposons le résultat démontré au rang $i$ avec $1 \leq i < h(\beta)-1$ et observons que, puisque $\beta$ est ordinaire, $\gamma_i \notin C_j$ de sorte que Supp$(\gamma_i) \ \subset \Pi_{j-1}$. Comme $\langle \gamma_i, \gamma_i \rangle > 0$, il existe $\epsilon \in \Pi_{j-1}$ telle que  $\langle \gamma_{i}, \epsilon \rangle > 0$. Comme $i < h(\beta)-1$, on a $h(\gamma_i) > 1 \Rightarrow \gamma_i \neq \epsilon \Rightarrow \gamma_{i+1}:= \gamma_i - \epsilon \in \Phi^+$ (lemme \ref{lemme2.9}).\\
Alors $\langle \beta , \epsilon \rangle = \langle \beta_{i} , \epsilon \rangle + \langle \gamma_{i} , \epsilon \rangle$ et, puisque $\langle \beta, \epsilon \rangle \leq 0$, on a $\langle \beta_{i} , \epsilon \rangle < 0$. Ceci entraine que (lemme \ref{lemme2.9}) $\beta_{i+1} := \beta_{i} + \epsilon \in C_{j}$.
Ainsi, $\beta = \beta_{i+1} + \gamma_{i+1}$ avec $\beta_{i+1} \in C_j$ et $h(\beta_{i+1}) = h(\beta_{i})+1$.\\
\par
On a donc bien le résultat annoncé et, pour $i = h(\beta) - 1$, on a $h(\gamma_i) = 1$. Donc, puisque $\beta$ est ordinaire, $\epsilon := \gamma_i \in \Pi_{j-1}$ et $\beta - \epsilon = \beta_i \in C_j$. 
\end{itemizedot}
\item On raisonne par récurrence sur $h(\beta)$.\\
Si $h(\beta) = 1$, on a $\beta = \alpha_j \in \Omega_j$.\\
Si $h(\beta) > 1$, il résulte de 1. qu'il existe $\epsilon \in \Pi_{j-1}$ telle que $ \beta' := \beta - \epsilon \in C_j$. $\beta'$ est alors ordinaire et de hauteur $h(\beta') = h(\beta) - 1$ de sorte que, par l'hypothèse de récurrence, Supp $\beta' \subset \Omega_j$. Il reste donc à montrer que $\epsilon \in \Omega_j$.\\
Si $\epsilon \notin \Omega_j$, alors $\epsilon \ \bot \ \alpha$ pour tout $\alpha \in $ Supp $\beta'$.  Il en résulte que $\epsilon \ \bot \ \beta'$. Considérons le plan $P = <\beta, \epsilon>$ et observons que $\Phi_P = \Phi \cap P$ est un système de racines de rang 2. Comme $\Phi \neq G_2$, on a $\Phi_P \neq G_2$ et, nécessairement, $\Phi_P$ est de type $B_2$. On a donc la configuration suivante :
\begin{center}
\begin{pgfpicture}{-0.5cm}{-0.5cm}{2.9cm}{2.2cm}%
\pgfsetroundjoin \pgfsetroundcap%
\pgfsetlinewidth{0.8pt} 
\pgfxyline(1.2,0.2)(2.2,0.2)
\pgfsetfillcolor{black}
\pgfmoveto{\pgfxy(2.0845,0.2667)}\pgflineto{\pgfxy(2.2,0.2)}\pgflineto{\pgfxy(2.0845,0.1333)}\pgflineto{\pgfxy(2.1423,0.2)}\pgfclosepath\pgffillstroke
\pgfxyline(1.2,0.2)(2.2,1.2)
\pgfmoveto{\pgfxy(2.0712,1.1655)}\pgflineto{\pgfxy(2.2,1.2)}\pgflineto{\pgfxy(2.1655,1.0712)}\pgflineto{\pgfxy(2.1592,1.1592)}\pgfclosepath\pgffillstroke
\pgfxyline(1.2,0.2)(1.2,1.2)
\pgfmoveto{\pgfxy(1.1333,1.0845)}\pgflineto{\pgfxy(1.2,1.2)}\pgflineto{\pgfxy(1.2667,1.0845)}\pgflineto{\pgfxy(1.2,1.1423)}\pgfclosepath\pgffillstroke
\pgfputat{\pgfxy(2.2,-0.1)}{\pgfnode{rectangle}{center}{\small $\epsilon$}{}{\pgfusepath{}}}
\pgfputat{\pgfxy(2.2,1.4)}{\pgfnode{rectangle}{center}{\small $\beta$}{}{\pgfusepath{}}}
\pgfputat{\pgfxy(1.2,1.4)}{\pgfnode{rectangle}{center}{\small $\beta'$}{}{\pgfusepath{}}}
\end{pgfpicture} 
\end{center}
Il en résulte que $\beta' - \epsilon$ est une racine. Comme Supp $\beta' \subset \Omega_j$, on a $\epsilon \notin $ Supp $\beta'$ ce qui est contradictoire avec le fait que $\beta' - \epsilon$ soit une racine. On a donc bien $\epsilon \in \Omega_j$ et par suite, Supp $\beta \subset \Omega_j$.\\
\item Soit $\beta$ une racine de $C_{j}$ différente de $\delta_{j}$. S'il existe $\epsilon \in \Pi_{j-1}$ telle que $\langle \beta , \epsilon \rangle < 0$ alors on conclut par le lemme \ref{lemme2.9}.\\
Supposons que $\langle \beta , \epsilon \rangle \geq 0$ pour tout $\epsilon \in \Pi_{j-1}$. Si $\langle \beta, \alpha_j \rangle < 0$, alors $C_{j}$ est exceptionnelle et $\beta + \alpha_j = \beta_{ex}$. Il résulte de la proposition \ref{prop2.8} que $s_D$ échange $\alpha_j$ et $\delta_{j}$, de sorte que, $\beta_{ex} = \alpha_j + \delta_{j}$. On en déduit que $\beta = \delta_{j}$ ce qui est contraire à l'hypothèse. On a donc $\langle \beta , \alpha_j \rangle \geq 0$ et, par suite, $\langle \beta , \epsilon \rangle \geq 0$ pour tout $ \epsilon \in \Pi_{j}$.\\
\begin{obsn}
$\textrm{\emph{Supp}} (\delta_{j}) = \Omega_{j}$.
\end{obsn}
\textbf{preuve :}\\
Par 2. on a Supp $(\delta_{j}) \subset \Omega_{j}$. Supposons cette inclusion stricte et considérons $\alpha \in \Omega_{j} \setminus \textrm{Supp} (\delta_{j})$. Comme $\Omega_j$ est connexe, il existe $\alpha'_1, \alpha'_2, ..., \alpha'_s$ dans $\Omega_j$, avec $s \geq 2$, $\alpha'_1 = \alpha, \alpha'_s \in \ \text{Supp}(\delta_j)$ et, pour $1 \leq i < s$, $\langle \alpha'_i , \alpha'_{i+1} \rangle < 0$.\\
Soit $k$ le plus grand entier tel que $\alpha'_k \notin$ Supp $(\delta_j)$. On a $k < s$ et $\alpha'_{k+1} \in$ Supp $(\delta_j)$. Comme $\delta_j$ est alors une combinaison linéaire à coefficients entiers positifs de racines simples différentes de $\alpha'_k$ et, comme le produit scalaire de deux racines simples distinctes est négatif ou nul, on a :
\[\langle \alpha'_k , \delta_j \rangle \leq \langle \alpha'_k, \alpha'_{k+1} \rangle < 0.\]
On en déduit par le lemme \ref{lemme2.9} que $\delta_j + \alpha'_k \in C_j \setminus \{\beta_{ex}\}$, ce qui contredit le fait que $\delta_j$ est la première racine de $C_j$.
\begin{obsn}
$\langle \beta, \delta_{j} \rangle > 0$
\end{obsn}
\textbf{preuve :}\\
Comme $\langle \beta, \beta \rangle > 0$, il existe $\alpha \in $ Supp$(\beta)$ telle que $\langle \beta , \alpha \rangle > 0$. Par le point 2. et l'observation 1, on a $\alpha \in $ Supp $\delta_j$. Comme (voir ci-dessus) $\langle \beta , \epsilon \rangle \geq 0$ pour tout $ \epsilon \in \Pi_j$, on a $\langle \beta , \epsilon \rangle \geq 0$ pour tout $ \epsilon \in$ Supp $\delta_j = \Omega_j \subset \Pi_j \Rightarrow \langle \beta, \delta_{j} \rangle \geq \langle \beta, \alpha \rangle > 0$.\\
\par
Par l'observation 2 et le lemme \ref{lemme2.9} $\gamma_{1} := \delta_{j} - \beta \in \Phi^{+}$ et, puisque $\beta$ et $\delta_j$ sont des racines ordinaires de $C_j$, $\textrm{Supp} (\gamma_{1}) \subset \Pi_{j-1}$.\\ 
On se propose de démontrer, par récurrence descendante, que, pour chaque entier $i \in \llbracket 1, h(\delta_j - \beta) \rrbracket$, il existe $\rho_i \in \Phi^+$ avec Supp $(\rho_i) \subset \Pi_{j-1}$, $h(\rho_i) = i$ et $\beta + \rho_i \in C_j$.\\
\begin{itemizedot}
\item Pour $i = h(\delta_j- \beta)$, on a le résultat avec $\rho_i = \gamma_1$ .
\item Supposons démontrée l'existence de $\rho_i$ avec $2 \leq i \leq h(\delta_j- \beta)$, de sorte $\eta_i = \beta + \rho_i \in C_j$. Comme $\langle \rho_i , \rho_i \rangle > 0$, il existe $\epsilon \in $ Supp $(\rho_i) \subset \Pi_{j-1}$ telle que $\langle \rho_i , \epsilon \rangle > 0$. Comme on a (voir ci-dessus), $\langle \beta, \epsilon \rangle \geq 0$, on en déduit que $\langle \eta_i , \epsilon \rangle > 0$. Donc (lemme \ref{lemme2.9}) $\eta_{i-1} := \eta_i - \epsilon \in C_j$ et $\rho_{i-1} := \rho_i - \epsilon \in \Phi^+$. On a alors clairement Supp $(\rho_{i-1}) \subset \Pi_{j-1}$, $h(\rho_{i-1}) = i-1$ et $\beta + \rho_{i-1} = \eta_{i-1} \in C_j$.\\
On a donc bien le résultat annoncé et, pour $i = 1$, $\epsilon := \rho_1 \in \Pi_{j-1}$ et $\beta + \epsilon \in C_j$.
\end{itemizedot} 
\end{enumerate}
\end{preuve} 
\begin{prop}\label{h'}
Soit un entier $j \in \llbracket 1, n \rrbracket$.
\begin{enumerate}
\item Si $C_j$ est ordinaire, alors $h'(C_j)$ est un intervalle de $\mathbb{N}^*$ de la forme  $\llbracket 1, t \rrbracket$;
\item Si $C_j$ est exceptionnelle, alors $h'(C_j \setminus \{\beta_{ex}\})$ est un intervalle de la forme $\llbracket 1, 2t \rrbracket \ (t \in \mathbb{N})$.\\
On a  $h'(\beta_{ex}) = t + \frac{1}{2}$.
\end{enumerate}
\end{prop}
\begin{preuve}
Le fait que $h'(C_j)$ dans le cas ordinaire (resp. $h'(C_j \setminus \{\beta_{ex}\})$ dans le cas exceptionnel) soit un intervalle de $\mathbb{N}$ résulte de la proposition \ref{prop2.10}. Il contient $1 = h(\alpha_j)$, ce qui démontre le premier point.\\
Supposons $C_j$ exceptionnelle. On note $B_{1},...,B_{t}$ les boîtes contenants les racines inférieurs à $\beta_{ex}$ dans l'ordre de Lusztig. Pour toutes ces boîtes, on a $h'(B_{i}) > h'(\beta_{ex})$. Mais la relation $h(B_{i}) + h(B_{i}') = h(\beta_{ex})$, pour l'image $B'_{i}$ de $B_{i}$ par $s_{D}$, implique $h'(B_{i}) > h'(\beta_{ex}) > h'(B'_{i})$. On a donc exactement t boîtes apparaissant après $\beta_{ex}$ et l'intervalle $h'(C_j \setminus \{\beta_{ex}\})$ est bien de la forme $\llbracket 1, 2t \rrbracket \ (t \in \mathbb{N})$.\\
De plus $h(\beta_{ex}) = h(\alpha_{n} + s_{D}(\alpha_{n})) = 1 + 2t$ et finalement $h'(\beta_{ex}) = t + \frac{1}{2}$.
\end{preuve}
On définit ci-dessous les plans admissibles introduit par G. Lusztig dans \cite[section 6.1]{Lu2}.
\begin{dfn}\label{planadm}
On appelle \emph{plan admissible} $P := <\beta, \beta'>$ tout plan engendré par deux racines positives $\beta$ et $\beta'$ telles que :
\begin{description}
\item[ou bien] $\beta$ appartient à une colonne exceptionnelle $C_{j}$ et $\beta' = s_{D}(\beta)$ est telle que $|h'(\beta') - h'(\beta)| = 1$ (dans ce cas on a $\beta + \beta' = \beta_{ex}$ et  $h'(\beta_{ex}) = t \pm \frac{1}{2}$).
\item[ou bien] $\beta$ est une racine non exceptionnelle qui appartient à une colonne quelconque $C_{j}$ et $\beta' = \alpha_{i}$ avec $i < j$.
\end{description}
On note $\Phi_P := \Phi \cap P$ et $\Phi_P^+ := \Phi^+ \cap P$.
\end{dfn}
\begin{rmq}\label{systracine}
Si $\Phi_P = G_2$ alors $\Phi = G_2$ (à cause de la longueur des racines).\\
Si $\Phi \neq G_2$ alors la première condition donne lieu  deux types de plans admissibles pour lesquels $\Phi_{P}^{+}$ est donné par :\\
\begin{center}\begin{tabular}{c|c}
\textbf{Type (1.1)} \label{type 1.1} & \textbf{Type (1.2)} \label{type 1.2} \\
\hline
\begin{pgfpicture}{-0.5cm}{0cm}{2.9cm}{2.5cm}%
\pgfsetroundjoin \pgfsetroundcap%
\pgfsetlinewidth{0.8pt} 
\pgfxyline(1.2,0.5)(2.2,0.5)
\pgfsetfillcolor{black}
\pgfmoveto{\pgfxy(2.0845,0.5667)}\pgflineto{\pgfxy(2.2,0.5)}\pgflineto{\pgfxy(2.0845,0.4333)}\pgflineto{\pgfxy(2.1423,0.5)}\pgfclosepath\pgffillstroke
\pgfxyline(1.2,0.5)(1.7,1.366)
\pgfmoveto{\pgfxy(1.5845,1.2993)}\pgflineto{\pgfxy(1.7,1.366)}\pgflineto{\pgfxy(1.7,1.2327)}\pgflineto{\pgfxy(1.6711,1.316)}\pgfclosepath\pgffillstroke
\pgfxyline(1.2,0.5)(0.7,1.366)
\pgfmoveto{\pgfxy(0.7,1.2327)}\pgflineto{\pgfxy(0.7,1.366)}\pgflineto{\pgfxy(0.8155,1.2993)}\pgflineto{\pgfxy(0.7289,1.316)}\pgfclosepath\pgffillstroke
\pgfputat{\pgfxy(2.2,0.2)}{\pgfnode{rectangle}{center}{\small $\beta$}{}{\pgfusepath{}}}
\pgfputat{\pgfxy(1.7,1.6)}{\pgfnode{rectangle}{center}{\small $\beta_{ex}$}{}{\pgfusepath{}}}
\pgfputat{\pgfxy(0.7,1.6)}{\pgfnode{rectangle}{center}{\small $\beta'$}{}{\pgfusepath{}}}
\end{pgfpicture}
&
\begin{pgfpicture}{-0.5cm}{-0.5cm}{2.9cm}{2.2cm}%
\pgfsetroundjoin \pgfsetroundcap%
\pgfsetlinewidth{0.8pt} 
\pgfxyline(1.2,0.2)(2.2,0.2)
\pgfsetfillcolor{black}
\pgfmoveto{\pgfxy(2.0845,0.2667)}\pgflineto{\pgfxy(2.2,0.2)}\pgflineto{\pgfxy(2.0845,0.1333)}\pgflineto{\pgfxy(2.1423,0.2)}\pgfclosepath\pgffillstroke
\pgfxyline(1.2,0.2)(2.2,1.2)
\pgfmoveto{\pgfxy(2.0712,1.1655)}\pgflineto{\pgfxy(2.2,1.2)}\pgflineto{\pgfxy(2.1655,1.0712)}\pgflineto{\pgfxy(2.1592,1.1592)}\pgfclosepath\pgffillstroke
\pgfxyline(1.2,0.2)(1.2,1.2)
\pgfmoveto{\pgfxy(1.1333,1.0845)}\pgflineto{\pgfxy(1.2,1.2)}\pgflineto{\pgfxy(1.2667,1.0845)}\pgflineto{\pgfxy(1.2,1.1423)}\pgfclosepath\pgffillstroke
\pgfxyline(1.2,0.2)(0.2,1.2)
\pgfmoveto{\pgfxy(0.2345,1.0712)}\pgflineto{\pgfxy(0.2,1.2)}\pgflineto{\pgfxy(0.3288,1.1655)}\pgflineto{\pgfxy(0.2408,1.1592)}\pgfclosepath\pgffillstroke
\pgfputat{\pgfxy(2.2,-0.1)}{\pgfnode{rectangle}{center}{\small $\beta$}{}{\pgfusepath{}}}
\pgfputat{\pgfxy(2.2,1.4)}{\pgfnode{rectangle}{center}{\small $\beta_{ex}$}{}{\pgfusepath{}}}
\pgfputat{\pgfxy(1.2,1.4)}{\pgfnode{rectangle}{center}{\small $\beta'$}{}{\pgfusepath{}}}
\pgfputat{\pgfxy(0.2,1.4)}{\pgfnode{rectangle}{center}{\small $\alpha_i$}{}{\pgfusepath{}}}
\end{pgfpicture}
\\
\hline
$A_2$& $B_2$\\
\hline
  $\beta > \beta_{ex} > \beta'$ & $\beta > \beta_{ex} > \beta' > \alpha_i$
\end{tabular}
\end{center}
La seconde condition donne lieu  quatre types de plans admissibles pour lesquels $\Phi_{P}^{+}$ est donné par :\\
\begin{center}
\begin{tabular}{c|c|c|c}
\textbf{Type (2.1)} \label{type 2.1}&\textbf{Type (2.2)} \label{type 2.2}& \textbf{Type (2.3)} \label{type 2.3} & \textbf{Type (2.4)} \label{type 2.4}\\
\hline
\begin{pgfpicture}{-0.5cm}{0cm}{2.9cm}{2.5cm}%
\pgfsetroundjoin \pgfsetroundcap%
\pgfsetlinewidth{0.8pt} 
\pgfxyline(1.2,0.5)(2.2,0.5)
\pgfsetfillcolor{black}
\pgfmoveto{\pgfxy(2.0845,0.5667)}\pgflineto{\pgfxy(2.2,0.5)}\pgflineto{\pgfxy(2.0845,0.4333)}\pgflineto{\pgfxy(2.1423,0.5)}\pgfclosepath\pgffillstroke
\pgfxyline(1.2,0.5)(1.7,1.366)
\pgfmoveto{\pgfxy(1.5845,1.2993)}\pgflineto{\pgfxy(1.7,1.366)}\pgflineto{\pgfxy(1.7,1.2327)}\pgflineto{\pgfxy(1.6711,1.316)}\pgfclosepath\pgffillstroke
\pgfxyline(1.2,0.5)(0.7,1.366)
\pgfmoveto{\pgfxy(0.7,1.2327)}\pgflineto{\pgfxy(0.7,1.366)}\pgflineto{\pgfxy(0.8155,1.2993)}\pgflineto{\pgfxy(0.7289,1.316)}\pgfclosepath\pgffillstroke
\pgfputat{\pgfxy(2.2,0.2)}{\pgfnode{rectangle}{center}{\small $\beta_1$}{}{\pgfusepath{}}}
\pgfputat{\pgfxy(1.7,1.6)}{\pgfnode{rectangle}{center}{\small $\beta_{2}$}{}{\pgfusepath{}}}
\pgfputat{\pgfxy(0.7,1.6)}{\pgfnode{rectangle}{center}{\small $\alpha_i$}{}{\pgfusepath{}}}
\end{pgfpicture}  &
\begin{pgfpicture}{-0.5cm}{-0.5cm}{2.9cm}{2.2cm}%
\pgfsetroundjoin \pgfsetroundcap%
\pgfsetlinewidth{0.8pt} 
\pgfxyline(1.2,0.2)(2.2,0.2)
\pgfsetfillcolor{black}
\pgfmoveto{\pgfxy(2.0845,0.2667)}\pgflineto{\pgfxy(2.2,0.2)}\pgflineto{\pgfxy(2.0845,0.1333)}\pgflineto{\pgfxy(2.1423,0.2)}\pgfclosepath\pgffillstroke
\pgfxyline(1.2,0.2)(2.2,1.2)
\pgfmoveto{\pgfxy(2.0712,1.1655)}\pgflineto{\pgfxy(2.2,1.2)}\pgflineto{\pgfxy(2.1655,1.0712)}\pgflineto{\pgfxy(2.1592,1.1592)}\pgfclosepath\pgffillstroke
\pgfxyline(1.2,0.2)(1.2,1.2)
\pgfmoveto{\pgfxy(1.1333,1.0845)}\pgflineto{\pgfxy(1.2,1.2)}\pgflineto{\pgfxy(1.2667,1.0845)}\pgflineto{\pgfxy(1.2,1.1423)}\pgfclosepath\pgffillstroke
\pgfxyline(1.2,0.2)(0.2,1.2)
\pgfmoveto{\pgfxy(0.2345,1.0712)}\pgflineto{\pgfxy(0.2,1.2)}\pgflineto{\pgfxy(0.3288,1.1655)}\pgflineto{\pgfxy(0.2408,1.1592)}\pgfclosepath\pgffillstroke
\pgfputat{\pgfxy(2.2,-0.1)}{\pgfnode{rectangle}{center}{\small $\beta$}{}{\pgfusepath{}}}
\pgfputat{\pgfxy(2.2,1.4)}{\pgfnode{rectangle}{center}{\small $\beta_{ex}$}{}{\pgfusepath{}}}
\pgfputat{\pgfxy(1.2,1.4)}{\pgfnode{rectangle}{center}{\small $\beta'$}{}{\pgfusepath{}}}
\pgfputat{\pgfxy(0.2,1.4)}{\pgfnode{rectangle}{center}{\small $\alpha_i$}{}{\pgfusepath{}}}
\end{pgfpicture}
& \begin{pgfpicture}{-0.5cm}{-0.5cm}{2.9cm}{2.2cm}%
\pgfsetroundjoin \pgfsetroundcap%
\pgfsetlinewidth{0.8pt} 
\pgfxyline(1.2,0.2)(2.2,0.2)
\pgfsetfillcolor{black}
\pgfmoveto{\pgfxy(2.0845,0.2667)}\pgflineto{\pgfxy(2.2,0.2)}\pgflineto{\pgfxy(2.0845,0.1333)}\pgflineto{\pgfxy(2.1423,0.2)}\pgfclosepath\pgffillstroke
\pgfxyline(1.2,0.2)(2.2,1.2)
\pgfmoveto{\pgfxy(2.0712,1.1655)}\pgflineto{\pgfxy(2.2,1.2)}\pgflineto{\pgfxy(2.1655,1.0712)}\pgflineto{\pgfxy(2.1592,1.1592)}\pgfclosepath\pgffillstroke
\pgfxyline(1.2,0.2)(1.2,1.2)
\pgfmoveto{\pgfxy(1.1333,1.0845)}\pgflineto{\pgfxy(1.2,1.2)}\pgflineto{\pgfxy(1.2667,1.0845)}\pgflineto{\pgfxy(1.2,1.1423)}\pgfclosepath\pgffillstroke
\pgfxyline(1.2,0.2)(0.2,1.2)
\pgfmoveto{\pgfxy(0.2345,1.0712)}\pgflineto{\pgfxy(0.2,1.2)}\pgflineto{\pgfxy(0.3288,1.1655)}\pgflineto{\pgfxy(0.2408,1.1592)}\pgfclosepath\pgffillstroke
\pgfputat{\pgfxy(2.2,-0.1)}{\pgfnode{rectangle}{center}{\small $\alpha_i$}{}{\pgfusepath{}}}
\pgfputat{\pgfxy(2.2,1.4)}{\pgfnode{rectangle}{center}{\small $\beta_{3}$}{}{\pgfusepath{}}}
\pgfputat{\pgfxy(1.2,1.4)}{\pgfnode{rectangle}{center}{\small $\beta_2$}{}{\pgfusepath{}}}
\pgfputat{\pgfxy(0.2,1.4)}{\pgfnode{rectangle}{center}{\small $\beta_1$}{}{\pgfusepath{}}}
\end{pgfpicture}  & 
\begin{pgfpicture}{-0.7cm}{-0.7cm}{2cm}{2.5cm}%
\pgfsetroundjoin%
\pgfsetlinewidth{0.8pt} 
\pgfxyline(0,0)(1,0)
\pgfsetmiterjoin \pgfsetfillcolor{black}
\pgfmoveto{\pgfxy(0.8845,0.0667)}\pgflineto{\pgfxy(1,0)}\pgflineto{\pgfxy(0.8845,-0.0667)}\pgflineto{\pgfxy(0.9423,0)}\pgfclosepath\pgffillstroke
\pgfsetroundjoin%
\pgfxyline(0,0)(0,1.5)
\pgfsetmiterjoin \pgfmoveto{\pgfxy(-0.0667,1.3845)}\pgflineto{\pgfxy(0,1.5)}\pgflineto{\pgfxy(0.0667,1.3845)}\pgflineto{\pgfxy(0,1.4423)}\pgfclosepath\pgffillstroke
\pgfsetroundjoin%
\pgfsetlinewidth{0.2pt} 
\pgfputat{\pgfxy(1.2,0.3)}{\pgftext{\color{black}\small $\alpha_i$}}\pgfstroke
\pgfputat{\pgfxy(0.3,1.7)}{\pgftext{\color{black}\small $\beta$}}\pgfstroke
\end{pgfpicture}

 \\
\hline
 $A_2$ & $B_{2}$ avec $\alpha_i$ longue &$B_{2}$ avec $\alpha_{i}$ courte & $A_{1} \times A_{1}$  \\ 
 \hline
$\beta_1 > \beta_{2} > \alpha_i$ &$\beta > \beta_{ex} > \beta' > \alpha_i$ &$\beta_1 > \beta_{2} > \beta_3 > \alpha_i$ & $\beta > \alpha_i$ 
  \end{tabular}
  \end{center}
  On remarque que les types (1.2) et (2.2) sont les mêmes.
\end{rmq}


\section{L'algèbre enveloppante quantique $\mathcal{U}_q(\mathfrak{g})$}
Soit $\mathbb{K}$ un corps de caractéristique différente de 2 et 3, et $q$ un élément de  $\mathbb{K}^*$ qui n'est pas une racine de l'unité. Tout d'abord, on va rappeler les définitions de $\mathcal{U}_q(\mathfrak{g})$ et $\mathcal{U}_q^+(\mathfrak{g})$ en utilisant les notations de \cite[chap4]{J}. On rappellera ensuite la construction des bases de Poincaré-Birkhoff-Witt de $\mathcal{U}_q(\mathfrak{g})$ au moyen des automorphismes de Lusztig. Cependant, il y a plusieurs fa\c cons de définir ces automorphismes dits ``de Lusztig'', on va ici décrire trois méthodes différentes. Celle de Lusztig provient de \cite[section 3]{Lu2}, celle de Jantzen, qui est la même que De Concini, Kac et Procesi, est développe dans \cite[section 8.14]{J} et  \cite[section 2.1]{MR1351503} et enfin une troisième qui permet de faire le lien entre les deux précédentes. Nous allons expliciter chaque méthode puis voir quelles sont les liens entre les bases de Poincaré-Birkhoff-Witt construites dans chacun des cas.
\subsection{Rappels sur $\mathcal{U}_q(\cat{g})$}
\begin{nota}[q-entiers et q-coefficients binomiaux]
on note pour $\alpha \in \Pi, a$ et $n$ entiers naturels avec $ a \geq n \geq 0$:
\[ [n]_{q} = \frac{q^{n}-q^{-n}}{q-q^{-1}}, \ [n]^{!}_{q} = [n]_{q}[n-1]_{q}...[2]_{q}[1]_{q} , \ \left[ \begin{array}{c} a \\ n \end{array} \right]_{q} = \frac{[a]^{!}_{q}}{[n]^{!}_{q}[a-n]^{!}_{q}}, \ \text{et} \ q_{\alpha} = q^{\frac{(\alpha,\alpha)}{2}}. \]
\end{nota}
\begin{dfn}[$\mathcal{U}_q(\cat{g})$ et $\mathcal{U}_q^+(\mathfrak{g})$]\label{Uqg}
\begin{itemizedot}
\item L'algèbre enveloppante quantique $\mathcal{U}_q(\cat{g})$ est la $\mathbb{K}$-algèbre engendrée par les générateurs $E_{\alpha}$, $F_{\alpha}$, $K_{\alpha}$ et $K_{\alpha}^{-1}$ pour tout $\alpha$ dans $\Pi$, soumis aux relations suivantes : ($\forall \alpha , \beta \in \Pi$)
\begin{equation}
K_{\alpha}K_{\alpha}^{-1} = 1 = K_{\alpha}^{-1}K_{\alpha} \hspace{2cm}  K_{\alpha}K_{\beta}=K_{\beta}K_{\alpha}
\end{equation}  
\begin{equation}
K_{\alpha}E_{\beta}K_{\alpha}^{-1} = q^{(\alpha , \beta)}E_{\beta}
\end{equation} 
\begin{equation}
K_{\alpha}F_{\beta}K_{\alpha}^{-1} = q^{-(\alpha , \beta)}F_{\beta}
\end{equation} 
\begin{equation}
E_{\alpha}F_{\beta} - F_{\beta}E_{\alpha} = \delta_{\alpha \beta} \frac{K_{\alpha} - K_{\alpha}^{-1}}{q_{\alpha}-q_{\alpha}^{-1}}
\end{equation}  
et aussi pour $\alpha \neq \beta$ :
\begin{equation}
\sum_{s=0}^{1-a_{\alpha \beta}} (-1)^s \left[ \begin{array}{c} 1-a_{\alpha \beta} \\ i \end{array} \right]_{q_{\alpha}} E_{\alpha}^{1-a_{\alpha \beta}-s}E_{\beta}E_{\alpha}^s = 0
\end{equation} 
\begin{equation}
\sum_{s=0}^{1-a_{\alpha \beta}} (-1)^s \left[ \begin{array}{c} 1-a_{\alpha \beta} \\ i \end{array} \right]_{q_{\alpha}} F_{\alpha}^{1-a_{\alpha \beta}-s}F_{\beta}F_{\alpha}^s = 0
\end{equation} 
où $a_{\alpha \beta} = 2 (\alpha , \beta) / (\alpha , \alpha)$ pour toutes racines simples $\alpha , \beta \in \Pi$.

\item On note $\mathcal{U}_q^+(\mathfrak{g})$ la sous-algèbre de $\mathcal{U}_q(\mathfrak{g})$ engendrée par les seuls $E_{\alpha}$ avec $\alpha \in \Pi$.
\end{itemizedot}
\end{dfn}
Rappelons deux résultats important démontrés par exemple dans \cite[section I.6]{MR1898492}.
\begin{theo}
\begin{enumerate}
\item $\mathcal{U}_q^+(\mathfrak{g})$ est noethérien.
\item $\mathcal{U}_q^+(\mathfrak{g})$ est gradué par $\mathbb{Z}\Phi$. ($wt(E_\alpha)= \alpha, wt(F_\alpha)=- \alpha $ et $wt(K_\alpha^{\pm 1})= 0)$
\end{enumerate}
\end{theo}
\subsection{La construction de Lusztig}
C'est à cette construction que les autres auteurs font référence malgré les modifications qu'ils effectuent dans leurs textes.
\begin{dfn}[Automorphismes de Lusztig]
Pour tout $i \in \llbracket 1,n \rrbracket$ il existe un unique automorphisme $T_{\alpha_{i}}$ de l'algèbre $\mathcal{U}_q(\mathfrak{g})$ tel que :
\[T_{\alpha_i}E_{\alpha_i} = -F_{\alpha_i}K_{\alpha_i}, \   T_{\alpha_i}F_{\alpha_i} = -K_{\alpha_i}^{-1} E_{\alpha_i} \ T_{\alpha_i}K_{\alpha_j} = K_{\alpha_j}K_{\alpha_i}^{-a_{ij}} \ (j \in \llbracket 1,n \rrbracket) \]
et pour $j \neq i$ :
\[  T_{\alpha_i}E_{\alpha_j} = \sum_{r + s = - a_{ij}}(-1)^r q^{-d_i s} E_{\alpha_i}^{(r)}E_{\alpha_j}E_{\alpha_i}^{(s)} \]
\[ T_{\alpha_i}F_{\alpha_j} = \sum_{r + s = - a_{ij}}(-1)^r q^{d_i s} F_{\alpha_i}^{(s)}F_{\alpha_j}F_{\alpha_i}^{(r)}  \]
où $\displaystyle E_{\alpha_i}^{(n)} := \frac{E_{\alpha_i}^{n}}{[n]_{d_{i}}^{!}}$. 
\end{dfn}
On utilise ensuite un ordre de Lusztig, les colonnes et les boîtes définies dans la partie 2. Le résultat suivant est énoncé par G. Lusztig dans \cite[section 4.3]{Lu2} :
\begin{prop}
Il existe une unique fonction $\Phi^{+} \rightarrow \llbracket 1,n \rrbracket, (\beta \rightarrow i_{\beta})$ telle que les propriétés suivantes sont vérifiées :
\begin{enumerate}
\item $s_{i_{\beta_{1}}}$ et $s_{i_{\beta_{2}}}$ commutent dans W si $\beta_{1}$ et $\beta_{2}$ sont dans la même boîte. Ainsi pour une boîte B, le produit des $s_{i_{\beta}}$ pour $\beta \in B$ est un élément bien définit $s(B)$ de W, indépendant de l'ordre de ses facteurs.
\item $i_{\alpha_{j}} = j$ pour $j \in \llbracket 1,n \rrbracket$.
\item Si $\beta \in C_{j}$ et si $B_{1}, ..., B_{k}$ sont les boîtes de $C_{j}$ dont les éléments sont strictement plus grands que $\beta$ pour l'ordre de Lusztig alors $s(B_{1})s(B_{2})...s(B_{k})(\alpha_{i_{\beta}}) = \beta$.
\end{enumerate}
On pose alors $w_{\beta} :=  s(B_{1})s(B_{2})...s(B_{k})$.
\end{prop}
On construit à présent $T_{w_{\beta}}$  l'aide du théorème suivant (\cite[théoréme 3.2]{Lu2})
\begin{theo}
Soit $w \in W$ et $s_{i_{1}}...s_{i_{p}}$ une décomposition réduite de $w$. Alors l'automorphisme $T_{w} := T_{\alpha_{i_{1}}}...T_{\alpha_{i_{p}}}$ ne dépends pas de la décomposition réduite choisie, il dépend uniquement de l'élément $w \in W$. Ainsi les $T_{\alpha_{i}}$ définissent un morphisme du groupe de tresse de W dans le groupe des automorphismes de $\mathcal{U}_q(\mathfrak{g})$.
\end{theo}
\begin{dfn}
Pour toute racine positive $\beta$, on définit $E_{\beta} := T_{w_{\beta}}(E_{i_{\beta}}) \in \Phi^{+}$. Ces éléments forment une base de PBW de $\mathcal{U}_q^+(\mathfrak{g})$ (\cite[proposition 4.2]{Lu2})
\end{dfn}
Pour obtenir des relations de commutation dans ce contexte, on utilise les plans admissibles.
\begin{nota}
Si $\beta > \beta'$, on note $[E_{\beta},E_{\beta'}]_{q} = E_{\beta}E_{\beta'} - q^{(\beta,\beta')}E_{\beta'}E_{\beta}$. Dans la suite, $\Phi_{P}$ désigne la trace du système de racines $\Phi$ sur un plan donné P.
\end{nota}
\textbf{Relations de commutation entre les $E_{\gamma}$ dans les plans admissibles :} \label{type 2}\\
Si $\Phi_P = G_2$ alors $\Phi = G_2$ (à cause de la longueur des racines) et les relations de commutation entre les générateurs de Lusztig sont connues (\cite[Section 5.2]{Lu2}). Elles permettent d'énoncer :
\begin{prop}\label{propG2}
Supposons $\Phi = G_2$, notons $\alpha_1$ la racine simple courte et $\alpha_2$ la racine simple longue. On obtient ainsi la bonne numérotation de l'ensemble des racines simples (cf exemple \ref{exG2}). $w_0$ a pour décomposition réduite $s_1s_2s_1s_2s_1s_2$ $(s_i = s_{\alpha_i})$ et, les racines étant écrites dans l'ordre convexe associée, on a 
\[\Phi^+ = \{\beta_1 = \alpha_1 ,\beta_2 = 3\alpha_1 + \alpha_2, \beta_3 = 2\alpha_1 + \alpha_2, \beta_4 = 3\alpha_1 + 2\alpha_2, \beta_5 = \alpha_1 + \alpha_2, \beta_6 = \alpha_2\}.\]   
La première colonne $C_1$ est réduite à $\{\beta_1\}$, la seconde colonne $C_2 = \{\beta_2, \beta_3, \beta_4, \beta_5, \beta_6\}$ est exceptionnelle ($\beta_{ex} = \beta_4$). On a alors : 
\begin{itemize}
\item $[E_{\beta_3},E_{\beta_1}]_{q} = \lambda E_{\beta_{2}}$ avec $\lambda \neq 0$,
\item $[E_{\beta_4},E_{\beta_1}]_{q} = \lambda E_{\beta_{3}}^2$ avec $\lambda \neq 0$,
\item $[E_{\beta_5},E_{\beta_1}]_{q} = \lambda E_{\beta_{3}}$ avec $\lambda \neq 0$,
\item $[E_{\beta_6},E_{\beta_1}]_{q} = \lambda E_{\beta_{5}}$ avec $\lambda \neq 0$,
\item $[E_{\beta_3},E_{\beta_1}]_{q} = \lambda E_{\beta_{2}}$ avec $\lambda \neq 0$,
\item $[E_{\beta_5},E_{\beta_3}]_{q} = \lambda E_{\beta_{4}}$ avec $\lambda \neq 0$.
\end{itemize}
\end{prop}
Si $\Phi \neq G_2$, les relations de commutation entre les générateurs de Lusztig correspondant aux racines d'un plan admissible sont partiellement connues (\cite[Section 5.2]{Lu2}). Elles permettent d'énoncer :
\begin{prop}[$\Phi \neq G_2$]\label{rellusztig}
\begin{itemizedot}
\item Si $P = <\beta, \beta'>$ est un plan admissible de \hyperref[systracine]{type (1.1)} avec $h'(\beta') = h'(\beta) + 1$, alors  $\Phi_{P}^{+} = \{\beta, \beta_{ex} = \beta + \beta', \beta'\}$ et :
\begin{itemize}
\item[$\circ$] $[E_{\beta},E_{\beta'}]_{q} = \lambda E_{\beta_{ex}}$ avec $\lambda \neq 0$,
\item[$\circ$] $[E_{\beta},E_{\beta_{ex}}]_{q} = [E_{\beta_{ex}},E_{\beta'}]_{q} = 0$.
\end{itemize}

\item Si $P = <\beta, \beta'>$ est un plan admissible de \hyperref[systracine]{type (1.2)} avec $h'(\beta') = h'(\beta) + 1$, alors  $\Phi_{P}^{+} = \{\beta, \beta_{ex} = \beta + \beta', \beta', \alpha_{i}\}$ et on a les relations :
\begin{itemize}
\item[$\circ$] $[E_{\beta},E_{\beta'}]_{q} = \lambda E_{\beta_{ex}}$ avec $\lambda \neq 0$,
\item[$\circ$] $[E_{\beta_{ex}},E_{\alpha_i}]_{q} = \lambda' E_{\beta'}^{2}$ avec $\lambda' \neq 0$,
\item[$\circ$] $[E_{\beta},E_{\alpha_i}]_{q} = \lambda'' E_{\beta'}$ avec $\lambda'' \neq 0$,
\item[$\circ$] $[E_{\beta},E_{\beta_{ex}}]_{q} = [E_{\beta_{ex}},E_{\beta'}]_{q} = [E_{\beta'},E_{\alpha_{i}}]_{q} = 0$.
\end{itemize}
\item Si $P = <\beta, \alpha_{i}>$ est un plan admissible de \hyperref[systracine]{type (2.1)}, alors  $\Phi_{P}^{+} = \{\beta_{1}, \beta_{2} = \beta_{1} + \alpha_{i}, \alpha_{i}\}$ ($\beta = \beta_{1}$ ou $\beta_{2}$) et :
\begin{itemize}
\item[$\circ$] $[E_{\beta_{1}},E_{\alpha_{i}}]_{q} = \lambda E_{\beta_{2}}$ avec $\lambda \neq 0$.,
\item[$\circ$] $[E_{\beta_{1}},E_{\beta_{2}}]_{q} = [E_{\beta_{2}},E_{\alpha_{i}}]_{q} =  0. $
\end{itemize} 
\item Si $P = <\beta,  \alpha_{i}>$ est un plan admissible de \hyperref[systracine]{type (2.2)}, alors on a les même relations que dans le type (1.2).\\
\item Si $P = \text{Vect}(\beta, \alpha_{i})$ est un plan admissible de \hyperref[systracine]{type (2.3)}, alors  $\Phi_{P}^{+} = \{\beta_{1}, \beta_{2} = \beta_{1} + \alpha_{i}, \beta_{3} = \beta_{1} + 2\alpha_{i}, \alpha_{i}\}$ ($\beta = \beta_{1}$,$\beta_{2}$ ou $\beta_{3}$) et on a les relations :
\begin{itemize}
\item[$\circ$] $[E_{\beta_{2}},E_{\alpha_{i}}]_{q} = \lambda E_{\beta_{3}}$ avec $\lambda \neq 0$,
\item[$\circ$] $[E_{\beta_{1}},E_{\beta_{3}}]_{q} = \lambda' E_{\beta_{2}}^{2}$ avec $\lambda' \neq 0$,
\item[$\circ$] $[E_{\beta_{1}},E_{\alpha_i}]_{q} = \lambda'' E_{\beta_{2}}$ avec $\lambda'' \neq 0$,
\item[$\circ$] $[E_{\beta_{1}},E_{\beta_{2}}]_{q} = [E_{\beta_{2}},E_{\beta_{3}}]_{q} = [E_{\beta_{3}},E_{\alpha_{i}}]_{q} = 0$
\end{itemize}
\item Si $P = <\beta, \alpha_{i}>$ est un plan admissible de \hyperref[systracine]{type (2.4)}, alors $\Phi_{P}^{+} = \{\beta, \alpha_{i}\}$ avec $\beta \ \bot \ \alpha_i$ et, si $\beta$ n'est pas exceptionnelle, alors $[E_{\beta},E_{\alpha_{i}}]_{q} = 0$.
\end{itemizedot}
\end{prop}
\begin{cor}[$\Phi \neq G_2$] \label{corplanadm}
Soient $i, l$ deux entiers tels que $1 \leq i < l \leq n $ et  $\eta \in C_l$
\begin{enumerate}
\item Si $(\eta, \alpha_i) > 0$, alors $[E_{\eta},E_{\alpha_{i}}]_{q} = 0$.
\item Si $\eta + \alpha_i = m \gamma$ avec $\gamma \in \Phi^+$ et $m \in \mathbb{N}^{\star}$, alors $[E_{\eta},E_{\alpha_{i}}]_{q} = \lambda E_{\gamma}^{m}, \ \text{avec} \  \lambda \in \mathbb{K}^{\star}$.
\item Si $\eta = \eta_1 + \eta_2$ avec $\eta_1$ et $\eta_2$ dans $C_l$ telles que $h(\eta_1) + 1 = h(\eta_2)$ alors  $[E_{\eta_1},E_{\eta_{2}}]_{q} = \lambda E_{\eta}, \ \text{avec} \  \lambda \in \mathbb{K}^{\star}$.
\end{enumerate}
\end{cor}
\begin{preuve}
$P = \text{Vect}(\eta, \alpha_{i})$ est un plan admissible de type (2.1), (2.2), (2.3) ou (2.4) par définition. 
\begin{enumerate}
\item $P$ n'est pas de type (2.4) puisque $(\eta, \alpha_i) \neq 0$. Il nous reste donc à envisager les cas suivants :
\begin{itemizedot}
\item $P$ est de type (2.1). On voit alors (en utilisant les notations \ref{systracine}) que, $\eta = \beta_{2}$. D'où le résultat par la proposition \ref{rellusztig}.
\item $P$ est de type (2.2). On voit alors que $\eta = \beta'$. D'où le résultat par la proposition \ref{rellusztig}.
\item $P$ est de type (2.3). On voit alors que $\eta = \beta_{3}$. D'où le résultat par la proposition \ref{rellusztig}.
\end{itemizedot}

\item Puisque $m \neq 0$, on a $\gamma \in P \cap \Phi^{+} = \Phi_{p}^{+}$ de sorte que $P$ n'est pas non plus de type (2.4). Il nous reste donc à envisager les cas suivants :
\begin{itemizedot}
\item $P$ est de type (2.1). On voit alors que $m=1$, $\eta = \beta_{1}$ et $\gamma = \beta_{2}$. D'où le résultat par la proposition \ref{rellusztig}.
\item $P$ est de type (2.2). On voit alors que l'on a deux possibilités :
\begin{itemize}
\item[$\circ$] $m=1$, $\eta = \beta$ et $\gamma = \beta'$. 
\item[$\circ$] $m=2$, $\eta = \beta_{ex}$ et $\gamma = \beta'$. 
\end{itemize}
D'où le résultat par la proposition \ref{rellusztig}.

\item $P$ est de type (2.3). On voit alors qu'on a $m=1$, $\eta = \beta_{1}$ (resp.  $\eta = \beta_{2}$) et $\gamma = \beta_{2}$  (resp. $\gamma = \beta_{3}$). D'où le résultat par la proposition \ref{rellusztig}.

\end{itemizedot}
\item Considérons le plan $P := <\eta_1, \eta_2>$. Il est admissible (cf définition \ref{planadm}) $\Phi_P^+$ est égal à $\{\eta_1, \eta, \eta_2\}$ (type 1.1) ou $\{\eta_1, \eta, \eta_2, \alpha_i\}$ avec $i < l$ (type 1.2). Il résulte alors de la proposition précèdente que $[E_{\eta_1},E_{\eta_{2}}]_{q} = \lambda E_{\eta}, \ \text{avec} \  \lambda \in \mathbb{K}^{\star}$.
\end{enumerate}
\end{preuve}
\subsection{La construction de Jantzen}
Cette construction est utilisée dans \cite[section 8.14]{J}, on utilise les même automorphismes $T_{\alpha}, \alpha \in \Pi$ mais d'une manière différente pour construire les éléments de la base de PBW.\\
Pour une décomposition donnée de $w_{0} = s_{i_{1}}...s_{i_{N}}$, on sait que pour tout $\beta \in \Phi^+$, il existe $i_{\beta} \in \llbracket 1,N \rrbracket$ tel que $\beta = s_{i_{1}}...s_{i_{\beta}-1}(\alpha_{i_{\beta}})$.
\begin{dfn}
Soit $\beta \in \Phi^+$, on pose $w'_{\beta} := s_{i_{1}}...s_{i_{\beta}-1}$ et on définit $X_{\beta} := T_{w'_{\beta}}(E_{\alpha_{i_{\beta}}})$, $Y_{\beta} := T_{w'_{\beta}}(F_{\alpha_{i_{\beta}}})$.
\end{dfn}
En utilisant les résultats \cite[théoréme 4.21]{J} et \cite[théoréme section 8.24]{J}, il vient :
\begin{theo}\label{PBW}
\begin{itemize}
\item Si $\alpha \in \Pi$, on a $X_{\alpha} = E_{\alpha}$ pour $\alpha \in \Pi$ (\cite[Proposition 8.20]{J}).
\item Les monômes $X_{\beta_{1}}^{k_{1}}...X_{\beta_{N}}^{k_{N}}$ ($k_i \in \mathbb{N}$)  forment une base de PBW de $\mathcal{U}_q^+(\mathfrak{g})$.
\item Les monômes $X_{\beta_{1}}^{k_{1}}...X_{\beta_{N}}^{k_{N}}K_{\alpha_{1}}^{m_{1}}...K_{\alpha_{n}}^{m_{n}}Y_{\beta_{1}}^{l_{1}}...Y_{\beta_{N}}^{l_{N}}$ (resp. $K_{\alpha_{1}}^{m_{1}}...K_{\alpha_{n}}^{m_{n}}Y_{\beta_{1}}^{l_{1}}...Y_{\beta_{N}}^{l_{N}}X_{\beta_{1}}^{k_{1}}...X_{\beta_{N}}^{k_{N}}$,\\
resp. $Y_{\beta_{1}}^{l_{1}}...Y_{\beta_{N}}^{l_{N}}K_{\alpha_{1}}^{m_{1}}...K_{\alpha_{n}}^{m_{n}} X_{\beta_{1}}^{k_{1}}...X_{\beta_{N}}^{k_{N}}$), ($k_{i}, l_{i} \in \mathbb{N}, m_{i} \in \mathbb{Z})$ forment une base de $\mathcal{U}_q(\mathfrak{g})$.
\end{itemize}

\end{theo}
La proposition ci-dessous a été démontré par Levendorski{\v\i} et Soibelman \cite[Proposition 5.5.2]{MR1116413} dans un cadre légèrement différent. On trouve d'autres formulations dans la littérature dont certaines comportent des inexactitudes. Pour cette raison, nous allons la redémontrer dans le cadre de ce travail. Les idées utilisées dans cette démonstration proviennent essentiellement de \cite[Proposition 5.5.2]{MR1116413}.
\begin{theo}[de Levendorski{\v\i} et Soibelman] \label{LS}
Avec les notations précédentes, si $i$ et $j$ sont deux entiers tels que $1 \leq i < j \leq N$
\begin{equation*}
	X_{\beta_i}X_{\beta_j} - q^{(\beta_i , \beta_j)}X_{\beta_j}X_{\beta_i} = \sum_{\beta_{i} < \gamma_{1} < ... < \gamma_{p} < \beta_{j}, \ p \geq 1, \ k_i \in \mathbb{N}}  c_{\mathbf{\overline{k},\overline{\gamma}}} X_{\gamma_{1}}^{k_{1}}...X_{\gamma_{p}}^{k_{p}}
\end{equation*}
 avec  $c_{\mathbf{\overline{k},\overline{\gamma}}} \in \mathbb{K}$ et $c_{\mathbf{\overline{k},\overline{\gamma}}} \neq 0 \Rightarrow wt(X_{\gamma_{1}}^{k_{1}}...X_{\gamma_{p}}^{k_{p}}) := k_{1} \times \gamma_{1}+...+k_{p} \times \gamma_{p} = \beta_i+\beta_j$
\end{theo}
Avant de démontrer ce théorème, rappelons quelques résultats classiques sur les automorphismes de Lusztig et les racines positives.
\begin{lemme}[\cite{J} section 4.6] \label{tau}
\begin{enumerate}
\item Il existe une unique automorphisme $\omega$ de $\mathcal{U}_q(\mathfrak{g})$ tel que $\omega(E_{\alpha}) = F_{\alpha}$, $\omega(F_{\alpha}) = E_{\alpha}$ et $\omega(K_{\alpha}) = K_{\alpha}^{-1}$. On a $\omega^{2} = 1$.
\item Il existe une unique anti-automorphisme $\tau$ de $\mathcal{U}_q(\mathfrak{g})$ tel que $\tau(E_{\alpha}) = E_{\alpha}$, $\tau(F_{\alpha}) = F_{\alpha}$ et $\tau(K_{\alpha}) = K_{\alpha}^{-1}$. On a $\tau^{2} = 1$.
\end{enumerate}
\end{lemme}
\begin{prop}[\cite{J} section 8.20] \label{auto lusztig}
Soient $w \in W$ et $\alpha \in \Pi$. Si $w(\alpha) > 0$, alors $T_{w}(E_{\alpha}) \in \mathcal{U}_q^{+}(\mathfrak{g})$,  $T_{w}(F_{\alpha}) \in \mathcal{U}_q^{-}(\mathfrak{g})$ et, de plus,  $T_{w}(E_{\alpha})$ (resp.  $T_{w}(F_{\alpha})$) est homogène de poids $w(\alpha)$ (resp. $-w(\alpha)$).
\end{prop}
Comme (toujours d'après \cite[section 8.18]{J}) on a $T_{w}^{-1} = \tau \circ T_{w^{-1}} \circ \tau$ et comme $\tau$ laisse stable $\mathcal{U}_q^{-}(\mathfrak{g})$ et $\mathcal{U}_q^{+}(\mathfrak{g})$, la proposition précédente peut se ré-écrire : 
\begin{prop}\label{auto lusztig2}
Soient $w \in W$ et $\alpha \in \Pi$. Si $w^{-1}(\alpha) > 0$, alors $T^{-1}_{w}(F_{\alpha}) \in \mathcal{U}_q^{-}(\mathfrak{g})$ et $T^{-1}_{w}(E_{\alpha}) \in \mathcal{U}_q^{+}(\mathfrak{g})$.
\end{prop}
\begin{lemme}\label{Uqb-}
Soient $\beta_{l}$ et $\beta_{p}$ deux racines positives telles que $\beta_{l} \leq \beta_{p}$ dans $\Phi^+$
\begin{enumerate}
\item $T_{w'_{\beta_{l}}}^{-1}(X_{\beta_{p}}) $ est dans $\mathcal{U}_q^{+}(\mathfrak{g})$  et homogène de poids $w^{\prime -1}_{\beta_{l}}(\beta_{p})$.
\item $T_{w'_{\beta_{p}}}^{-1}(X_{\beta_{l}}) $ est dans $\mathcal{U}_q(\mathfrak{b}^-)$.
\end{enumerate}

\end{lemme}
\begin{preuve}
\begin{enumerate}
\item Par la proposition \ref{auto lusztig}, $T_{w'_{\beta_l}}^{-1}(X_{\beta_{p}}) = T_{\alpha_{i_{l-1}}}^{-1}...T_{\alpha_{i_{1}}}^{-1}T_{\alpha_{i_{1}}}...T_{\alpha_{i_{p-1}}}(E_{\alpha_{i_{p}}}) = T_{\alpha_{i_{l}}}...T_{\alpha_{i_{p-1}}}(E_{\alpha_{i_{p}}})$ est un élément de $\mathcal{U}_q^{+}(\mathfrak{g})$, homogène de poids $s_{i_{l}} \circ ... \circ s_{i_{p-1}}(\alpha_{p}) = w^{\prime -1}_{\beta_{l}}(\beta_{p})$.
\item \[\begin{array}{ll}T_{w'_{\beta_p}}^{-1}(X_{\beta_{l}}) &= T_{\alpha_{i_{p}}}^{-1}...T_{\alpha_{i_{l}}}^{-1}T_{\alpha_{i_{l-1}}}^{-1}...T_{\alpha_{i_{1}}}^{-1}T_{\alpha_{i_{1}}}...T_{\alpha_{i_{l-1}}}(E_{\alpha_{i_{l}}})\\
&= T_{\alpha_{i_{p}}}^{-1}...T_{\alpha_{i_{l}}}^{-1}(E_{\alpha_{i_{l}}})\\
&= T_{\alpha_{i_{p}}}^{-1}...T_{\alpha_{i_{l-1}}}^{-1}(-K_{\alpha_{i_{l}}}^{-1}F_{\alpha_{i_{l}}})\\
&= T_{\alpha_{i_{p}}}^{-1}...T_{\alpha_{i_{l-1}}}^{-1}(-K_{\alpha_{i_{l}}}^{-1})T_{\alpha_{i_{p}}}^{-1}...T_{\alpha_{i_{l-1}}}^{-1}(F_{\alpha_{i_{l}}})
\end{array}\]
$T_{\alpha_{i_{p}}}^{-1}...T_{\alpha_{i_{l-1}}}^{-1}(-K_{\alpha_{i_{l}}}^{-1}) = - K_{s_{i_{p}}...s_{i_{l-1}}(\alpha_{i_{l}})}^{-1} $ est dans $\mathcal{U}_q^0(\mathfrak{g})$ car c'est un produit de $K_{\alpha}^{-1}$ (cf \cite[section 4.4]{J} pour la définition de $K_{\lambda}$ avec $\lambda \in \mathbb{Z}\Phi$).
\par
Comme $s_{i_{p}}...s_{i_{l-1}}(\alpha_{i_{l}})$ est une racine positive, on sait (cf proposition \ref{auto lusztig2}) que $T_{\alpha_{i_{p}}}^{-1}...T_{\alpha_{i_{l-1}}}^{-1}(F_{\alpha_{i_{l}}})$ est dans $\mathcal{U}_q^{-}(\mathfrak{g})$.

\end{enumerate}

\end{preuve}

\begin{lemme}\label{convexe}
On note ici $\leq$ l'ordre sur les racines positives associée à la décomposition de $w_{0}$ choisie ci-dessus.\\
Soient $\beta_{i} \leq \beta_{i+1} \leq ... \leq \beta_{i+p}$, p racines positives non nécessairement distinctes ($p \geq 2$). Si $\beta := \beta_{i} + ... + \beta_{i+p}$ est une racine positive, alors $\beta_{i} < \beta < \beta_{i+p}$.
\end{lemme}
\begin{preuve}
On démontre ce résultat par récurrence sur p :
\begin{description}
\item[Initialisation :] Rappelons que l'ordre sur les racines positives induit par une décomposition de $w_{0}$ est toujours convexe (voir par exemple \cite{Papi}). Ceci donne le résultat pour p = 1.
\item[p-1 $\Rightarrow$ p:] On suppose que $\beta := \beta_{i} + ... + \beta_{i+p} \in \Phi^+$, il existe $k \in \llbracket 1, p \rrbracket$ tel que $\langle \beta_{i+k}, \beta \rangle > 0$ (car $\langle \beta, \beta \rangle > 0$). D'après le lemme \ref{lemme2.9},  on  a $\gamma:= \beta - \beta_{i+k} \in \Phi$ et, puisque $\gamma = \beta_i + ...+ \widehat{\beta_{i+k}} + ... + \beta_{i+p}, \gamma \in \Phi^+$. Par le cas p = 1, on a:
\[\beta = \gamma + \beta_{i+k} \ \Rightarrow \ \gamma < \beta < \beta_{i+k}  \ \text{ou} \  \beta_{i+k} < \beta < \gamma \]
Par l'hypothèse de récurrence, on a $\beta_{i} < \gamma   < \beta_{i+p} $ et, dans les deux cas, on en déduit que $\beta_{i} < \beta < \beta_{i+p}$.
\end{description}
\end{preuve}
\begin{preuve2}
Soient $i$ et $j$ deux entiers tels que $1 \leq i < j \leq N$. Comme les monômes ordonnés en les $X_{\beta}, \beta \in \Phi^+$, forment une base de $\mathcal{U}_q^+(\mathfrak{g})$ (Théorème \ref{PBW}), on a :
\[X_{\beta_j}X_{\beta_i} = \sum_{\gamma_{1} < ... < \gamma_{p}, \ p \geq 1, \  k_i \in \mathbb{N}}  c_{\mathbf{\overline{k},\overline{\gamma}}} X_{\gamma_{1}}^{k_{1}}...X_{\gamma_{p}}^{k_{p}} \ \ (\star) \]
Rappelons que $\mathcal{U}_q^+(\mathfrak{g})$ est graduée par $\mathbb{Z}\Phi$ et que $wt(X_{\beta}) = \beta$ pour tout $\beta \in \Phi^+$. On a donc $c_{\mathbf{\overline{k},\overline{\gamma}}} \neq 0 \Rightarrow wt(X_{\gamma_{1}}^{k_{1}}...X_{\gamma_{p}}^{k_{p}}) := k_{1} \times \gamma_{1}+...+k_{p} \times \gamma_{p} = \beta_i+\beta_j$. Il reste donc à démontrer que :
\[\begin{array}{p{4cm}llrp{2cm}r}
&c_{\mathbf{\overline{k},\overline{\gamma}}} = 0 & \text{si} \ \gamma_{1} < \beta_{i}; &&& (1)\\
&c_{\mathbf{\overline{k},\overline{\gamma}}} = 0 & \text{si} \ \beta_{j}  < \gamma_{p}; &&& (2)\\
&c_{\mathbf{\overline{k},\overline{\gamma}}} = 0 & \text{si} \ \gamma_{1} = \beta_{i} \ \text{et} \  \gamma_{p} < \beta_{j}; &&& (3)\\
&c_{\mathbf{\overline{k},\overline{\gamma}}} = 0 & \text{si} \ \gamma_{1} > \beta_{i} \  \text{et} \  \gamma_{p} = \beta_{j}; &&& (4)\\
&c_{\mathbf{\overline{k},\overline{\gamma}}} = q^{-(\beta_{i},\beta_{j})} & \text{si} \ p = 2, \gamma_{1} = \beta_{i} \ \text{et} \  \gamma_{2} = \beta_{j};&&& (5)
\end{array}\]
Dans la suite, pour $k \in \llbracket 1, N \rrbracket$, on note $w_{k} := s_{i_{1}}...s_{i_{k}}$. Soit $t$ le plus petit entier tel qu'il existe un monôme du second membre de l'égalité $(\star)$ ( coefficient non nul) commençant par $X_{\gamma_{1}}$ avec $\gamma_1 = \beta_t$. De même, soit $u$ le plus grand entier tel qu'il existe un monôme du second membre de l'égalité $(\star)$ finissant par $X_{\gamma_{p}}$ avec $\gamma_p = \beta_u$. Comme dans \cite{J}, on note $\mathcal{U}_q(\mathfrak{b}^{+})$ (resp $\mathcal{U}_q(\mathfrak{b}^-)$) la sous algèbre de  $\mathcal{U}_q(\mathfrak{g})$ engendrée par les éléments $E_{\alpha}$ (resp. $F_{\alpha}$) et $K_{\alpha}^{\pm}$ ($\alpha \in \Pi$).
\begin{itemizedot}
\item Si $t < i$, en appliquant $T_{w_{t}}^{-1}$ aux deux membres de $(\star)$, on obtient :
\[T_{w_{t}}^{-1}(X_{\beta_j}X_{\beta_i}) = \sum_{\beta_{t} = \gamma_1 < \gamma_{2} < ... < \gamma_{p}}  c_{\mathbf{\overline{k},\overline{\gamma}}} T_{w_{t}}^{-1}(X_{\beta_{t}}^{k_{1}}...X_{\gamma_{p}}^{k_{p}}) + \sum_{\beta_t < \gamma_{1} < ... < \gamma_{p}}  c_{\mathbf{\overline{k},\overline{\gamma}}} T_{w_{t}}^{-1}(X_{\gamma_{1}}^{k_{1}}...X_{\gamma_{p}}^{k_{p}}) \ (\star \star) \] 
Observons que si $l > t$, on a $T_{w_{t}}^{-1}(X_{\beta_l}) = T_{w_{t}}^{-1}T_{w_{l-1}}(E_{\alpha_{i_{l}}})= T_{\alpha_{i_{t+1}}}...T_{\alpha_{i_{l-1}}}(E_{\alpha_{i_{l}}})$. Comme l'écriture $s_{i_{t+1}}...s_{i_{l}}$ est réduite, on a $s_{i_{t+1}}...s_{i_{l-1}}(\alpha_{i_{l}}) \in \Phi^{+}$. On en déduit, par la proposition \ref{auto lusztig2}, que $T_{w_{t}}^{-1}(X_{\beta_l})  \in \mathcal{U}_q^+(\mathfrak{g})$. Puisque $i > t$, il en résulte que le premier membre ainsi que la deuxième somme T du second membre de $(\star \star)$ sont dans $\mathcal{U}_q^+(\mathfrak{g})$. Comme ci-dessus, on a $T_{w_{t}}^{-1}(X_{\beta_{t}}) = T_{\alpha_{i_{t}}}^{-1}(E_{\alpha_{i_{t}}}) = \tau \circ T_{\alpha_{i_{t}}} \circ \tau (E_{\alpha_{i_{t}}})$  (\cite[section 8.18]{J}) $=   -K_{\alpha_{i_{t}}}^{-1}F_{\alpha_{i_{t}}} = -K_{\alpha_{i_{t}}}^{-1}Y_{\alpha_{i_{t}}}$. La première somme $S$ du second membre de $(\star \star)$ est donc une somme finie de produits $ Y_{\alpha_{i_{t}}}^{k_{1}}K_{\alpha_{i_{t}}}^{-k_{1}}S_{k_{1}}$ avec $S_{k_{1}} \in  \mathcal{U}_q^+(\mathfrak{g})$ et $k_{1} \in \mathbb{N}^{\star}$. D'autre part, on a $S =  T_{w_{t}}^{-1}(X_{\beta_j}X_{\beta_i}) - T \in \mathcal{U}_q^+(\mathfrak{g})$. Par le théorème \ref{PBW}, on a donc S=0, ce qui contredit la définition de $t$. De l, on a $t \geq i$, ce qui démontre (1).
\par
\item Si j < u, on applique alors $T_{w_{j}}^{-1}$ aux deux membres de $(\star)$ :
\[T_{w_{j}}^{-1}(X_{\beta_j}X_{\beta_i}) = \sum_{\gamma_{1} < \gamma_{2} < ... < \gamma_{p} \leq \beta_{j}}  c_{\mathbf{\overline{k},\overline{\gamma}}} T_{w_{j}}^{-1}(X_{\gamma_{1}}^{k_{1}}...X_{\gamma_{p}}^{k_{p}}) + \sum_{\gamma_{1} < ... < \gamma_{r-1} <\beta_{j} < \gamma_{r} < ... < \gamma_{p}}  c_{\mathbf{\overline{k},\overline{\gamma}}} T_{w_{j}}^{-1}(X_{\gamma_{1}}^{k_{1}}...X_{\gamma_{p}}^{k_{p}}) \]
On déduit du lemme \ref{Uqb-} que le premier membre et la première somme du second membre de l'égalité précédente sont dans $\mathcal{U}_q(\mathfrak{b}^-)$. La seconde somme de cette même égalité n'est pas un élément de $\mathcal{U}_q(\mathfrak{b}^-)$ car la proposition \ref{auto lusztig2} implique (par un raisonnement analogue à celui du premier point) que $T_{w_{j}}^{-1}(X_{\gamma_{r}}^{k_{r}}...X_{\gamma_{p}}^{k_{p}}) \in \mathcal{U}_q^+(\mathfrak{g})$ pour $\beta_{j} <  \gamma_{r} ... < \gamma_{p}$. Cette troisième somme est donc nulle par le théorème \ref{PBW}. En composant par  $T_{w_{j}}$, on en déduit que :
\[ \sum_{\gamma_{1} < ... < \gamma_{r-1} <\beta_{j} < \gamma_{r} < ... < \gamma_{p}}  c_{\mathbf{\overline{k},\overline{\gamma}}} X_{\gamma_{1}}^{k_{1}}...X_{\gamma_{p}}^{k_{p}} = 0.\]
Du théorème \ref{PBW} on déduit à nouveau que chacun des monômes de cette somme est nul, ce qui contredit la définition de u et démontre le point (2).
\item Pour le point (3), on raisonne par l'absurde. Si un tel $c_{\mathbf{\overline{k},\overline{\gamma}}}$ est non nul, un monôme de la forme $X_{\beta_{i}}^{k_{1}}...X_{\gamma_{p}}^{k_{p}}$ avec $\beta_i < \gamma_{p} < \beta_{j}$ apparat dans la relation. Par suite $(k_{1}-1) \times \beta_{i} + ... + k_{p} \times \beta_{p} = \beta_{j}$, ce qui contredit le lemme \ref{convexe}.
\item Le point (4) se démontre comme le point (3). 
\item A ce stade, on a donc démontre qu'il existe $a \in \mathbb{K}$ tel que 
\[X_{\beta_j}X_{\beta_i} - aX_{\beta_i}X_{\beta_j} = \sum_{\beta_{i} < \gamma_{1} < ... < \gamma_{p} < \beta_{j}}  c_{\mathbf{\overline{k},\overline{\gamma}}} X_{\gamma_{1}}^{k_{1}}...X_{\gamma_{p}}^{k_{p}} \]
Il reste à calculer la valeur de $a$. On compose chaque membre par $T_{w_{i}}^{-1}$ :\\
\begin{equation}\label{eg7}
T_{w_{i}}^{-1}(X_{\beta_j})T_{w_{i}}^{-1}(X_{\beta_i}) - aT_{w_{i}}^{-1}(X_{\beta_i})T_{w_{i}}^{-1}(X_{\beta_j}) = \sum_{\beta_{i} < \gamma_{1} < ... < \gamma_{p} \leq \beta_{j}}  c_{\mathbf{\overline{k},\overline{\gamma}}} T_{w_{i}}^{-1}(X_{\gamma_{1}}^{k_{1}})...T_{w_{i}}^{-1}(X_{\gamma_{p}}^{k_{p}})
\end{equation}
et on note $X:=T_{w_{i}}^{-1}(X_{\beta_j})$. C'est un élément de $\mathcal{U}_q^{+}(\mathfrak{g})$ par la proposition \ref{auto lusztig2}. De même, le second membre $X'$ est un élément de $\mathcal{U}_q^+(\mathfrak{g})$ et puisque $T_{w_{i}}^{-1}(X_{\beta_i}) := T_{\alpha_{i}}^{-1}(E_{\alpha_{i}}) = -K_{\alpha_{i}}^{-1}F_{\alpha_{i}}$, on a :
\[-XK_{\alpha_{i}}^{-1}F_{\alpha_{i}}+aK_{\alpha_{i}}^{-1}F_{\alpha_{i}}X = X'\]
Par le lemme \ref{Uqb-}, $X$ est homogène de poids $w_i^{-1}(\beta_j)$, de sorte que 
\[-XK_{\alpha_{i}}^{-1}F_{\alpha_{i}} = -q^{(\alpha_{i}, w_i^{-1}(\beta_j))}K_{\alpha_{i}}^{-1}XF_{\alpha_{i}}\]
Compte tenu des relations de commutation entre les $F_{\alpha}$ et les $E_{\alpha}$ (cf définition \ref{Uqg}), on a : 
\[XF_{\alpha_{i}} = F_{\alpha_{i}}X + K_{1}X_{1} + X_{2}, \ \text{avec} \ K_{1} \in \mathcal{U}_q^0(\mathfrak{g}), X_{1} \ \text{et} \ X_{2} \in \mathcal{U}_q^+(\mathfrak{g})\]
Comme $F_{\alpha_{i}} = Y_{\alpha_{i}}$,  l'égalité (\ref{eg7}) s'écrit :
\[-q^{(\alpha_{i}, w_i^{-1}(\beta_j))}K_{\alpha_{i}}^{-1}Y_{\alpha_{i}}X + -q^{(\alpha_{i}, w_i^{-1}(\beta_j))}K_{\alpha_{i}}^{-1}K_{1}X_{1} + -q^{(\alpha_{i}, w_i^{-1}(\beta_j))}K_{\alpha_{i}}^{-1}X_{2} +aK_{\alpha_{i}}^{-1}Y_{\alpha_{i}}X = X',\]
de sorte que 
\[-q^{(\alpha_{i}, w_i^{-1}(\beta_j))}K_{\alpha_{i}}^{-1}Y_{\alpha_{i}}X+aK_{\alpha_{i}}^{-1}Y_{\alpha_{i}}X = q^{(\alpha_{i}, w_i^{-1}(\beta_j))}K_{\alpha_{i}}^{-1}K_{1}X_{1} + q^{(\alpha_{i}, w_i^{-1}(\beta_j))}K_{\alpha_{i}}^{-1}X_{2} + X' \]
Comme $X, X', X_{1}$ et $X_{2}$ sont dans $\mathcal{U}_q^{+}(\mathfrak{g})$, il résulte du théorème \ref{PBW} que les deux membres de cette égalité sont nuls, de sorte que $a = q^{(\alpha_{i}, w_i^{-1}(\beta_j))} = q^{-( w_{i}^{ -1}(\beta_{i}) , w_{i}^{ -1}(\beta_{j}))} = q^{-(\beta_{i},\beta_{j})}$.
Ceci démontre le point (5).
\end{itemizedot}
\end{preuve2}
\subsection{Relations de commutation entre les $X_{\gamma}$ dans les plans admissibles}
Le but de cette partie est de démontrer que les $X_{\gamma}$ vérifient des relations de commutation analogues aux relations de commutation entre les $E_{\gamma}$ rappelées dans le paragraphe 3.2.
\subsubsection{Construction d'un troisième système de générateurs}
\begin{conv}
\begin{itemizedot}
\item Soit i un entier de $\llbracket 1, n \rrbracket$. $\tau$ désignant l'automorphisme de   $\mathcal{U}_q(\mathfrak{g})$ défini au lemme \ref{tau}, notons $T'_{\alpha_{i}} = \tau \circ T_{\alpha_{i}} \circ \tau$. Ceci est un  automorphisme de $\mathcal{U}_q(\mathfrak{g})$ qui vérifie immédiatement les relations suivantes :
\[T'_{\alpha_i}E_{\alpha_i} = -K^{-1}_{\alpha_i}F_{\alpha_i}, \   T'_{\alpha_i}F_{\alpha_i} = -E_{\alpha_i} K_{\alpha_i}\ T'_{\alpha_i}K_{\alpha_j} = K_{\alpha_j}K_{\alpha_i}^{-a_{ij}} \ (j \in \llbracket 1,n \rrbracket) \]
et pour $j \neq i$ :
\[  T'_{\alpha_i}E_{\alpha_j} = \sum_{r + s = - a_{ij}}(-1)^r q^{d_i s} E_{\alpha_i}^{(s)}E_{\alpha_j}E_{\alpha_i}^{(r)} \]
\[ T'_{\alpha_i}F_{\alpha_j} = \sum_{r + s = - a_{ij}}(-1)^r q^{-d_i s} F_{\alpha_i}^{(r)}F_{\alpha_j}F_{\alpha_i}^{(s)}  \]
\item De même, si $w_{p}  \in W$ a pour décomposition réduite $w_{p} = s_{i_{1}}...s_{i_{p}}$, on pose $T'_{w_{p}} :=  \tau \circ T_{w_{p}} \circ \tau$ et on observe que l'on a encore  $T'_{w_{p}} = T'_{\alpha_{i_{1}}}...T'_{\alpha_{i_{p}}} $.
\item Si $\beta \in \Phi^+$, on pose $w_{\beta} := s_{i_{1}}...s_{i_{\beta}-1}$ et on définit $X'_{\beta} := T'_{w_{\beta}}(E_{\alpha_{i_{\beta}}})$ et  $Y'_{\beta} := T_{w_{\beta}}(F_{\alpha_{i_{\beta}}})$. On a immédiatement $X'_{\alpha} = E_{\alpha}$ et par analogie $Y'_{\alpha} = F_{\alpha}$  pour $\alpha \in \Pi$.
\end{itemizedot}
\end{conv}
Le théorème de Levendorski{\v\i} et Soibelman prend alors la forme suivante :
\begin{prop} \label{LS2}
Avec les notations précédentes, si $i$ et $j$ sont deux entiers tels que $1 \leq i < j \leq N$
\begin{equation*}
	X'_{\beta_i}X'_{\beta_j} - q^{-(\beta_i , \beta_j)}X'_{\beta_j}X'_{\beta_i} = \sum_{\beta_{i} < \gamma_{1} < ... < \gamma_{p} < \beta_{j}, \ p \geq 1, \ k_i \in \mathbb{N}}  c_{\mathbf{\overline{k},\overline{\gamma}}} X_{\gamma_{1}}^{\prime k_{1}}...X_{\gamma_{p}}^{\prime k_{p}}
\end{equation*}
 avec  $c_{\mathbf{\overline{k},\overline{\gamma}}} \in \mathbb{K}$ et $c_{\mathbf{\overline{k},\overline{\gamma}}} \neq 0 \Rightarrow wt(X_{\gamma_{1}}^{k_{1}}...X_{\gamma_{p}}^{k_{p}}) := k_{1} \times \gamma_{1}+...+k_{p} \times \gamma_{p} = \beta_i+\beta_j$
\end{prop}
\begin{preuve}
Commen\c cons par exprimer $X_{\beta}$ en fonction de $X'_{\beta}$ ($\beta \in \Phi^+$):
\[X_{\beta} = T_{w'_{\beta}}(E_{\alpha_{i_{\beta}}}) = \tau \circ T'_{w'_{\beta}} \circ \tau (E_{\alpha_{i_{\beta}}}) = \tau (X'_{\beta}).\]
Ainsi, compte tenu du théorème de Levendorski{\v\i} et Soibelman (théorème \ref{LS}), on a :
\[\begin{array}{ll}
X'_{\beta_i}X'_{\beta_j} - q^{-(\beta_i , \beta_j)}X'_{\beta_j}X'_{\beta_i} &= \tau (X_{\beta_j}X_{\beta_i} - q^{-(\beta_i , \beta_j)}X_{\beta_i}X_{\beta_j})\\
&= -q^{-(\beta_i , \beta_j)} \tau (X_{\beta_i}X_{\beta_j} - q^{(\beta_i , \beta_j)}X_{\beta_j}X_{\beta_i})\\
&\displaystyle = -q^{-(\beta_i , \beta_j)} \tau (\sum_{\beta_{i} < \gamma_{1} < ... < \gamma_{p} < \beta_{j}, \ p \geq 1, \ k_i \in \mathbb{N}}  c_{\mathbf{\overline{k},\overline{\gamma}}} X_{\gamma_{1}}^{k_{1}}...X_{\gamma_{p}}^{k_{p}})\\
&\displaystyle = \sum_{\beta_{i} < \gamma_{1} < ... < \gamma_{p} < \beta_{j}, \ p \geq 1, \ k_i \in \mathbb{N}}  c'_{\mathbf{\overline{k},\overline{\gamma}}} \tau(X_{\gamma_{p}})^{k_{p}}...\tau(X_{\gamma_{1}})^{k_{1}}\\
&\displaystyle = \sum_{\beta_{i} < \gamma_{1} < ... < \gamma_{p} < \beta_{j}, \ p \geq 1, \ k_i \in \mathbb{N}}  c'_{\mathbf{\overline{k},\overline{\gamma}}} X_{\gamma_{p}}^{\prime k_{p}}...X_{\gamma_{1}}^{\prime k_{1}}
\end{array}\]
Pour ré-ordonner, les monômes du second membre on utilise alors le lemme ci-dessous. 
\end{preuve}
\begin{lemme}
On suppose que,  pour tout $ 1<i<j<N$, on ait une relation de la forme $X'_{\beta_i}X'_{\beta_j} - q^{-(\beta_i , \beta_j)}X'_{\beta_j}X'_{\beta_i} = P_{ij} $, où $P_{ij} \in \mathbb{K}\langle X'_{\beta_{i+1}},...,X'_{\beta_{j-1}} \rangle$.  Alors pour tout $1<k<l<N$, tout polynôme M en les $X'_{\beta_{k}},...,X'_{\beta_{l}}$ peut s'écrire comme une combinaison linéaire de monômes ordonnés en les $X'_{\beta_{k}},...,X'_{\beta_{l}}$.
\end{lemme}
\begin{preuve}
Il suffit de démontrer le résultat lorsque $M$ est un monôme de $ \mathbb{K}\langle X'_{\beta_{k}},...,X'_{\beta_{l}} \rangle$. On raisonne par récurrence sur $m := l-k$, le nombre de variables apparaissant dans l'écriture de $M$.
\begin{itemize}
\item $m = 1$ : Le résultat est évident.\\

\item m-1 $\Rightarrow$ m : On raisonne par récurrence sur le nombre $t$ d'occurrences de $X'_{\beta_{k}}$ dans $M$.
\begin{itemize}
\item Si $t=1$ : On écrit $M = M_{1} X'_{\beta_{k}} M_{2}$ avec $M_{1}$ et $M_{2}$ deux monômes en $X'_{\beta_{k+1}},...,X'_{\beta_{l}}$ . On raisonne alors par récurrence sur $p = d M_{1}$.
\begin{itemize}
\item p = 0 : $M = X'_{\beta_{k}}M_{2}$ avec $M_{2} \in \mathbb{K} \langle X'_{\beta_{k+1}},...,X'_{\beta_{l}} \rangle$. En appliquant l'hypothèse de récurrence sur le nombre de variables  $M_{2}$, on obtient le résultat voulu.\\
\item p-1 $\Rightarrow$ p : On a alors $M_{1} = M'_{1} X'_{\beta_{s}}$ avec $s \in \llbracket k+1,l \rrbracket$ et $M'_{1} \in \mathbb{K}\langle X'_{\beta_{k+1}},...,X'_{\beta_{l}} \rangle$, de degré p-1 . Ainsi, en utilisant la relation de commutation entre $X'_{\beta_{k}}$ et $X'_{\beta_{s}}$ donne par hypothèse on a :
\[M = M'_{1}X'_{\beta_{s}}X'_{\beta_{k}}M_{2} = q^{(\beta_{k},\beta_{s})}M'_{1}X'_{\beta_{k}}X'_{\beta_{s}}M_{2} + M'_{1}P_{ks}M_{2}, \  P_{ks} \in \mathbb{K} \langle X'_{\beta_{k+1}},...,X'_{\beta_{s-1}} \rangle.\]
Par l'hypothèse de récurrence sur le degré de $M_{1}$, le monôme  $M'_{1}X'_{\beta_{k}}X'_{\beta_{s}}M_{2}$ est une combinaison linaire de monômes ordonnés en $X'_{\beta_{k}},...,X'_{\beta_{l}}$.  Par l'hypothèse de récurrence sur le nombre de variable, chaque monôme de $M'_{1}P_{ks}M_{2}$ est une combinaison linaire de monômes ordonnés en $X'_{\beta_{k+1}},...,X'_{\beta_{l}}$. Duo le résultat.
\end{itemize}
\item t-1 $\Rightarrow$ t: Si la première occurrence de $X'_{\beta_{k}}$ est en première position dans l'écriture de \\
$M = X'_{\beta_{k}}M_{1}X'_{\beta_{k}}...M_{t}X'_{\beta_{k}}M_{t+1}$, le monôme $M_{1}X'_{\beta_{k}}...M_{t}X'_{\beta_{k}}M_{t+1}$ est une combinaison linaire de monômes ordonnés en $X'_{\beta_{k}},...,X'_{\beta_{l}}$ d'où le résultat.\\
Sinon, on a $M = M_{1}X'_{\beta_{k}}...M_{n}X'_{\beta_{k}}M_{n+1}$ et, par l'hypothèse de récurrente applique  $M_{1}X'_{\beta_{k}}...M_{n}$ on a $M = \sum_{i} X_{\beta_{k}}^{\prime a_{i}}M'_{i}X'_{\beta_{k}}M_{n+1}$, avec $a_{i} \in \mathbb{N}$ et $M'_{i} $ monômes (ordonnés) en $X'_{\beta_{k+1}},...,X'_{\beta_{l}}$.\\
Chaque monôme $M'_{i}X'_{\beta_{k}}M_{n+1}$ présente une seule occurrence de $X'_{\beta_{k}}$. C'est donc (voir le cas $t =1$) une combinaison linaire de monômes ordonnés en  $X'_{\beta_{k+1}},...,X'_{\beta_{l}}$. D'où le résultat.
\end{itemize}
\end{itemize}
\end{preuve}
\subsubsection{Relations entre les $E_{\beta}$ et les $X'_{\beta}$.}
Comme dans les paragraphes précédents, $\Phi^+$ est muni d'un ordre de Lusztig associé  une décomposition réduite $w_{0} = s_{i_{1}}...s_{i_{N}}$. Dans ce cas, on peut préciser comme suit la version du théorème de Levendorski{\v\i} et Soibelman (proposition  \ref{LS2}).
\begin{theo} \label{1.4}
Si $i$ et $j$ sont deux entiers tels que $1 \leq i < j \leq N$, on a :
\begin{equation*}
	X'_{\beta_i}X'_{\beta_j} - q^{-(\beta_i , \beta_j)}X'_{\beta_j}X'_{\beta_i} = \sum_{\beta_{i} < \gamma_{1} < ... < \gamma_{p} < \beta_{j}}  C_{\mathbf{\overline{k},\overline{\gamma}}} X_{\gamma_{1}}^{\prime k_{1}}...X_{\gamma_{p}}^{\prime k_{p}}.
\end{equation*}
Les monômes du second membres apparaissant avec un coefficient $C_{\mathbf{\overline{k},\overline{\gamma}}}$ non nul vérifient :
\begin{itemizedot}
\item $wt(X_{\gamma_{1}}^{\prime k_{1}}...X_{\gamma_{p}}^{\prime k_{p}}) = \beta_i+\beta_j$;
\item $\gamma_{1}$ n'est pas dans la même boîte que $\beta_{i}$;
\item $\gamma_{p}$ n'est pas dans la même boîte que $\beta_{j}$.
\end{itemizedot}
\end{theo}
La preuve de ce théorème repose essentiellement sur le résultat suivant :
\begin{lemme}\label{1.2}
Soit $B=\{\beta_p, ..., \beta_{p+l}\}$ une boîte et soit $\alpha_{i_{p}}, ..., \alpha_{i_{p+l}}$ les racines simples correspondantes. Alors $\forall k \in \llbracket 0,l \rrbracket$, 
\[ T'_{\alpha_{i_{p}}} ... T'_{\alpha_{i_{p+k-1}}}(E_{\alpha_{i_{p+k}}}) =  E_{\alpha_{i_{p+k}}} = T'_{\alpha_{i_{p}}} ...T'_{i_{\alpha_{p+k-1}}}T'_{i_{\alpha_{p+k+1}}}...T'_{\alpha_{i_{p+l}}}(E_{\alpha_{i_{p+k}}}) \]
\end{lemme}
\begin{preuve}
On sait que si $\alpha_{1}$ et $\alpha_{2}$ sont deux racines simples orthogonales, on a $T_{\alpha_{1}}(E_{\alpha_{2}}) = E_{\alpha_{2}} = \tau(E_{\alpha_{2}})$ de sorte que $T'_{\alpha_{1}}(E_{\alpha_{2}}) = E_{\alpha_{2}}$. Comme $\alpha_{i_{p}}, ..., \alpha_{i_{p+l}}$ sont 2 à 2 orthogonales (la proposition \ref{boîte}), on a immédiatement les formules ci-dessus.
\end{preuve}

\begin{preuve3}
Le premier point résulte de la proposition \ref{LS2}.
Si dans la décomposition de $w_{0}$, on change l'ordre des réflexions  associes  des racines simples correspondant à une même boîte B, on trouve une nouvelle décomposition réduite de $w_{0}$. Les racines positives de B sont permutes de la même manière et les autres sont inchangés. Par le lemme ci-dessus, les $X'_{\beta}$, $\beta \in B$, sont aussi permutes de la même manière mais non modifies, et les $X'_{\gamma}$, $\gamma \notin B$, sont inchangés. Sans perte de généralité, on peut donc supposer que $\beta_i$ est maximale dans sa boîte et que $\beta_{j}$ est minimale dans la sienne. Il en résulte que si $\beta_i < \gamma_{1} < ... < \gamma_{p} < \beta_{j}$ alors $\gamma_{1}$ n'est pas dans la même boîte que $\beta_{i}$ et  $\gamma_{p}$ n'est pas dans le même boîte que $\beta_{j}$.
\end{preuve3}
\begin{rmq}
La preuve du théorème précédent est aussi valable pour les éléments $X_\beta, \beta \in \Phi_+$ de sorte qu'on appliquera aussi le théorème \ref{1.4} à ces éléments.
\end{rmq}
On peut à présent démontrer : 
\begin{theo}\label{theo3.6}
\[\forall \beta \in \Phi^+, \exists \lambda_\beta \in \mathbb{K}\setminus\{0\} \ \ \text{tel que} \ \ X'_\beta = \lambda_\beta E_\beta\]
\end{theo}
\begin{preuve}
Soit $\beta$ et $\beta'$ deux racines positives telles que $\beta > \beta'$. On note (par analogie avec $[E_{\beta}, E_{\beta'}]_q$) $[X'_{\beta}, X'_{\beta'}]_q = X'_{\beta}X'_{\beta'} - q^{(\beta,\beta')} X'_{\beta'}X'_{\beta}$.\\
Traitons d'abord le cas $\Phi = G_2$ et reprenons les conventions de la proposition \ref{propG2}. On sait (3.4.1. conventions) que, puisque $\beta_1$ et $\beta_6$ sont simples, on a $X'_{\beta_1} = E_{\beta_1}$ et $X'_{\beta_6} = E_{\beta_6}$.\\
On a donc 
\[[X'_{\beta_6},X'_{\beta_1}]_q = [E_{\beta_6},E_{\beta_1}]_q = \lambda E_{\beta_5} \ \text{avec} \ \lambda \in \mathbb{K}^\star. \]
Par le théorème \ref{1.4},  on a aussi $[X'_{\beta_6},X'_{\beta_1}]_q= \mu X'_{\beta_5}$ avec $\mu \in \mathbb{K}$ et, par suite, $X'_{\beta_5} = \lambda_{\beta_5} E_{\beta_5}$ avec $\lambda_{\beta_5} \in \mathbb{K}^\star$.\\
De là, $[X'_{\beta_5},X'_{\beta_1}]_q = \lambda_{\beta_5} [E_{\beta_5},E_{\beta_1}]_q = \nu E_{\beta_3}$ avec $\nu \in \mathbb{K}^\star$. On en déduit comme ci-dessus que $X'_{\beta_3} = \lambda_{\beta_3} E_{\beta_3}$ avec $\lambda_{\beta_3} \in \mathbb{K}^\star$.\\
De la même manière, en considérant $[X'_{\beta_3},X'_{\beta_1}]_q = \lambda_{\beta_3} [E_{\beta_3},E_{\beta_1}]_q$, on montre que $X'_{\beta_2} = \lambda_{\beta_2} E_{\beta_2}$ avec $\lambda_{\beta_2} \in \mathbb{K}^\star$.\\
Enfin, on a $[X'_{\beta_5},X'_{\beta_3}]_q = \lambda_{\beta_5}\lambda_{\beta_3} [E_{\beta_5},E_{\beta_3}]_q = \nu E_{\beta_4}$ avec $\nu \in \mathbb{K}^\star$, d'où on déduit que $X'_{\beta_4} = \lambda_{\beta_4} E_{\beta_4}$ avec $\lambda_{\beta_4} \in \mathbb{K}^\star$.\\
\par
Supposons maintenant $\Phi \neq G_2$, considérons une colonne $C_{t}$ $(t \in \llbracket 1, n \rrbracket)$ et démontrons le théorème pour toutes les racines de $C_{t}$.\\
\par
On étudie d'abord le cas des racines non exceptionnelles.\\
Soit donc $\beta \in C_{t}$, $\beta$ non exceptionnelle. Raisonnons par récurrence sur $h(\beta)$.
\begin{description}
\item[Initialisation :] Si $h(\beta) = 1$, alors $\beta = \alpha_{t}$ et d'après les conventions de  $X'_{\alpha_{t}} = E_{\alpha_{t}}$.
\item[Hérédité :] Supposons $h(\beta) > 1$ et le résultat vrai pour chaque $\delta \in C_{t}$, $\delta$ non exceptionnelle telle que $h(\delta) < h(\beta)$.
Par la proposition \ref{prop2.10}, il existe une racine simple $\alpha_i$ ($i<t$) telle que $\beta - \alpha_i = \gamma \in C_t$. De plus, $\gamma$ n'est pas exceptionnelle car, sinon $\beta = \gamma + \alpha_{i}$ le serait. $P := < \alpha_i, \beta>$  est alors un plan admissible de \hyperref[type 2]{type (2.1), (2.2) ou (2.3)} et on a $[E_\gamma,E_{\alpha_i}]_q = c E_\beta$ $(c \in \mathbb{K} \setminus\{ 0\})$ (cf. paragraphe 3.2 page \pageref{type 2}).\\
Comme $h(\gamma) = h(\beta)-1 < h(\beta)$,  on a $X'_\gamma = \lambda_{\gamma} E_{\gamma}$  $(\lambda_{\gamma } \in \mathbb{K}\setminus\{0\}$), et comme $E_{\alpha_i} = X'_{\alpha_i}$, on a :
\[[X'_\gamma,X'_{\alpha_i}]_q = \lambda_{\gamma} [E_\gamma,E_{\alpha_i}]_q = \lambda_\gamma c E_\beta . \]
Par le théorème \ref{1.4},  $E_\beta$ est une combinaison linaire de monômes $X'_{\delta_1}...X'_{\delta_s}$ avec \\
$\alpha_i < \delta_{1} \leq ... \leq \delta_{s} < \gamma$, $\delta_{s}$ n'appartient pas à la boîte de $\gamma$, $\delta_{1}$ n'appartient pas à la boîte de $\alpha_{i}$ et 
\[ \delta_1+...+\delta_s = \alpha_i + \gamma  = \beta  \ \ \ (\star) .\]
Pour chacun de ces monômes, $\delta_s \in C_{t}$ et $\delta_{s}$ n'est pas exceptionnelle (puisque $\beta \in C_{t}$ et $\beta$ n'est pas exceptionnelle). Comme $\delta_s < \gamma $ et $\delta_{s}$ n'appartient pas à la boîte de $\gamma$, on a  $h(\delta_s) > h(\gamma) \Rightarrow h(\delta_s) \geq h(\beta) \Rightarrow  (s = 1$ et $\delta_1 = \beta) \Rightarrow E_\beta =  aX'_\beta $ avec $ a \in \mathbb{K}\setminus\{0\}$, d'où le résultat.
\end{description}
\begin{tabular}{@{}lc@{}}
\begin{minipage}{.8\linewidth}
Supposons à présent que $\beta$ est la racine exceptionnelle de $C_{t}$ (cf figure ci-contre). Soit $\gamma$ l'élément de $C_{t}$ qui précède $\beta$ dans l'ordre de Lusztig et soit $\delta = s_{D}(\gamma)$, de sorte que $\delta + \gamma = \beta$. Par le proposition \ref{h'}, on a $h'(\beta) = m+\frac{1}{2}$ avec $m \in \mathbb{N}^{\star}$ et $h'(C_{t})= \llbracket 1, 2m \rrbracket$. Si B est la boîte de $C_{t}$ qui précède $\beta$, on a donc $h'(B) = h(B) = t+1$. Comme la boîte de $\beta$ est réduite  un élément, on a $\gamma \in B$, donc $h(\gamma) = m+1$ et, par suite, $h(\delta) = m$.  $P =Vect(\gamma , \delta)$ est alors un plan admissible de type (1.1) ou (1.2), et $[E_\delta,E_{\gamma}]_q = c E_\beta (c \neq 0)$ (cf. paragraphe 3.2).\\
Comme $\gamma$ et $\delta$ sont non exceptionnelles, on sait déjà (voir ci dessus) que $X'_\gamma =  \lambda_\gamma E_\gamma$ et $X'_{\delta} = \lambda_\delta E_{\delta}$ $(\lambda_{\gamma} \neq 0, \lambda_{\delta} \neq 0)$. Ainsi, on a :
\[[X'_\delta,X'_{\gamma}]_q =\lambda_\delta \lambda_\gamma   [E_\delta,E_{\gamma}]_q = \lambda_\delta \lambda_\gamma  c E_\beta \ (\lambda_{\gamma} \neq 0, \lambda_{\delta} \neq 0).\] 

\end{minipage} &
 \begin{tabular}{c|c|}
\cline{2-2}
  & $\beta_k$ \\ 
 \cline{2-2}
     $C_t$ :  & \vdots \\ 
\cline{2-2}
     &$\gamma$ \\ 
 \cline{2-2}
    & $\beta$ \\ 
\cline{2-2}
     & $\delta $\\ 
   \cline{2-2}
     & \vdots \\ 
\cline{2-2}
     & $\beta_N$ \\ 
\cline{2-2}
  \end{tabular}\\
\end{tabular}\\
\par
Comme ci-dessus,  $E_\beta$ est une combinaison linaire de monômes $X'_{\delta_1}...X'_{\delta_s}$ avec \\
$\gamma < \delta_{1} \leq ... \leq \delta_{s} < \delta$, $\delta_{s}$ n'appartient pas à la boîte de $\delta$ et $\delta_{1}$ n'appartient pas à la boîte de $\gamma$. Comme $\beta$  est la seule racine de $C_{t}$ vérifiant  $\gamma <\beta < \delta$, $\beta$ n'appartient pas à la boîte de $\delta$ et $\beta$ n'appartient pas à la boîte de $\gamma$, on en déduit que $s=1$ et $\delta_{1} = \beta$. Donc $E_\beta = a X'_\beta$ avec $a \in \mathbb{K}\setminus\{0\}$.
\end{preuve}
Des théorèmes \ref{1.4} et \ref{theo3.6}, on en déduit :
\begin{cor} \label{LSE}
Si $i$ et $j$ sont deux entiers tels que $1 \leq i < j \leq N$, on a :
\begin{equation*}
	E_{\beta_i}E_{\beta_j} - q^{-(\beta_i , \beta_j)}E_{\beta_j}E_{\beta_i} = \sum_{\beta_{i} < \gamma_{1} < ... < \gamma_{p} < \beta_{j}, \ p \geq 1, \ k_i \in \mathbb{N}}  C'_{\mathbf{\overline{k},\overline{\gamma}}} E_{\gamma_{1}}^{k_{1}}...E_{\gamma_{p}}^{k_{p}}.
\end{equation*}
Les monômes du second membres apparaissant avec un coefficient $C'_{\mathbf{\overline{k},\overline{\gamma}}}$ non nul vérifient :
\begin{itemizedot}
\item $wt(X_{\gamma_{1}}^{\prime k_{1}}...X_{\gamma_{p}}^{\prime k_{p}}) = \beta_i+\beta_j$;
\item $\gamma_{1}$ n'est pas dans la même boîte que $\beta_{i}$;
\item $\gamma_{p}$ n'est pas dans la même boîte que $\beta_{j}$.
\end{itemizedot}\end{cor}

\subsubsection{Lien avec la construction de Jantzen}
\begin{prop}\label{relcomX}
Soient $\beta_{1} < \beta_{2}$ deux racines positives.
\begin{enumerate}
\item Si $E_{\beta_{1}}E_{\beta_{2}}-q^{-(\beta_{1},\beta_{2})}E_{\beta_{2}}E_{\beta_{1}} = k E^{m}_{\gamma}$ ($k \neq 0, m \geq 1$ et  $\gamma \in \Phi^+$), alors\\
$X_{\beta_{1}}X_{\beta_{2}}-q^{+(\beta_{1},\beta_{2})}X_{\beta_{2}}X_{\beta_{1}} = k' X^{m}_{\gamma}$ ($k' \neq 0$).
\item Si $E_{\beta_{1}}E_{\beta_{2}}-q^{-(\beta_{1},\beta_{2})}E_{\beta_{2}}E_{\beta_{1}} = k E_{\gamma}E_{\delta}$ ($k \neq 0, \gamma, \delta \in \Phi^+, \gamma$ et $\delta$ appartenant  une même boîte), alors $X_{\beta_{1}}X_{\beta_{2}}-q^{+(\beta_{1},\beta_{2})}X_{\beta_{2}}X_{\beta_{1}} = k' X_{\gamma}X_{\delta}$ ($k' \neq 0$).
\end{enumerate}
\end{prop}
\begin{preuve}
Soit $\beta \in \Phi^{+}$. Rappelons (Section 3.4.1) que $X_{\beta} := T_{w'_{\beta}}(E_{\alpha_{i_{\beta}}})$, $X'_{\beta} := T'_{w'_{\beta}} (E_{\alpha_{i_{\beta}}})$, et que $T_{w'_{\beta}} = \tau \circ T'_{w'_{\beta}} \circ \tau $. On a donc $X_{\beta} = \tau \circ T'_{w'_{\beta}} \circ \tau (E_{\alpha_{i_{\beta}}}) = \tau (X'_{\beta})$. Rappelons également (Théorème \ref{theo3.6}) que  $X'_{\beta} = \lambda_{\beta} E_{\beta}$ avec $\lambda_{\beta} \in \mathbb{K}^{\star}$.\\
Soient $\beta_{1} < \beta_{2}$ deux racines positives.
\begin{enumerate}
\item Si $E_{\beta_{1}}E_{\beta_{2}}-q^{-(\beta_{1},\beta_{2})}E_{\beta_{2}}E_{\beta_{1}} = k E^m_{\gamma}$ ($k \neq 0, \gamma \in \Phi^+$), alors :
\[\begin{array}{rl}
X_{\beta_{1}}X_{\beta_{2}}-q^{+(\beta_{1},\beta_{2})}X_{\beta_{2}}X_{\beta_{1}} &= \tau(X'_{\beta_{1}})\tau(X'_{\beta_{2}})-q^{+(\beta_{1},\beta_{2})}\tau(X_{\beta_{2}})\tau(X_{\beta_{1}})\\
\tiny\\
 &= \tau(X'_{\beta_{2}}X'_{\beta_{1}}-q^{+(\beta_{1},\beta_{2})}X'_{\beta_{1}}X'_{\beta_{2}})\\
 \tiny \\
 &=-q^{+(\beta_{1},\beta_{2})} \tau(X'_{\beta_{1}}X'_{\beta_{2}}-q^{-(\beta_{1},\beta_{2})}X'_{\beta_{2}}X'_{\beta_{1}})\\
 \tiny \\
 &=-q^{+(\beta_{1},\beta_{2})}\lambda_{\beta_{1}}\lambda_{\beta_{2}} \tau(E_{\beta_{1}}E_{\beta_{2}}-q^{-(\beta_{1},\beta_{2})}E_{\beta_{2}}E_{\beta_{1}})\\
 \tiny \\
 &=-q^{+(\beta_{1},\beta_{2})}\lambda_{\beta_{1}}\lambda_{\beta_{2}} \tau(k E^m_{\gamma})\\
 \tiny \\
 &=\displaystyle \frac{-q^{+(\beta_{1},\beta_{2})}\lambda_{\beta_{1}}\lambda_{\beta_{2}} k}{\lambda_{\gamma}}\tau((X'_{\gamma})^m)\\
 \tiny \\
 &= k' X^m_{\gamma} \ \text{avec} \ k' \in \mathbb{K}^{\star}
 \end{array}\]
 \item Si $E_{\beta_{1}}E_{\beta_{2}}-q^{-(\beta_{1},\beta_{2})}E_{\beta_{2}}E_{\beta_{1}} = k E_{\gamma}E_{\delta}$ ($k \neq 0, \gamma, \delta \in \Phi^+, \gamma$ et $\delta$ appartenant  une même boîte) alors, en effectuant le même calcul qu'en 1., on trouve :
\[X_{\beta_{1}}X_{\beta_{2}}-q^{(\beta_{1},\beta_{2})}X_{\beta_{2}}X_{\beta_{1}} = k' \tau (X'_{\gamma}X'_{\delta}) = k' X_{\delta}X_{\gamma} \ (k' \neq 0)\]
Comme $\gamma$ et $\delta$ sont dans la môme boîte, on sait (proposition \ref{boîte}) que $(\delta, \gamma) = 0$, de sorte que, par le théorème \ref{LS}, $X_{\gamma}X_{\delta} = X_{\gamma}X_{\delta} $. Ce qui termine la démonstration.
\end{enumerate}

\end{preuve}

\section{Effacement des dérivation dans $\mathcal{U}_q^+(\mathfrak{g})$}
\subsection{Notations et rappels}
Dans cette partie, on note $A := \mathcal{U}_q^+(\mathfrak{g})$, $X_i := X_{\beta_{i}}$ pour $1 \leq i \leq N$, et $\lambda_{i,j}:= q^{-(\beta_{j},\beta_{i})}$ pour $1 \leq i,j \leq N$. On a vu (théorème \ref{LS}) que, si $1 \leq i < j \leq N$, on a :
\begin{equation}
X_{j}X_{i} - \lambda_{j,i}X_{i}X_{j} = P_{j,i}
\end{equation}
avec
\begin{equation}\label{Pji}
P_{j,i} = \sum_{\overline{k} = (k_{i+1},...,k_{j-1})} c_{\bar{k}} X_{i+1}^{k_{i+1}}...X_{j-1}^{k_{j-1}} \ (c_{\bar{k}} \in \mathbb{K}).
\end{equation}
De plus, compte tenu de la $\Phi$-graduation de  $\mathcal{U}_q^+(\mathfrak{g})$, on a 
\begin{equation}\label{impli}
c_{\bar{k}} \neq 0 \Rightarrow \lambda_{l,i+1}^{k_{i+1}} ... \lambda_{l,j-1}^{k_{j-1}} = \lambda_{l,j}\lambda_{l,i} \ \text{pour} \ 1 \leq l \leq N
\end{equation}
Ainsi, $A$ vérifie les égalités (1) et l'hypothèse 6.1.1 de \cite[Section 6 page 503-504]{MR1967309}.\\
\par
Compte tenu du théorème \ref{PBW}, les monômes ordonnés en les $X_{i}$ forment une base de $A$ de sorte l'on a (\cite[Proposition 6.1.1]{MR1967309}): 
\begin{prop}\label{ore}
\begin{enumerate}
\item $A$ est une extension de Ore itérée que l'on peut écrire :
\[A = \mathbb{K}[X_{1}][X_2;\sigma_2,\delta_2]...[X_N;\sigma_N,\delta_N]\]
\item Si $1 \leq m \leq N$, il existe un (unique) automorphisme $h_{m}$ de l'algèbre $A$ qui vérifie $h_{m}(X_{i}) = \lambda_{m,i} X_{i}$ pour  $1 \leq i \leq N$.
\end{enumerate}
\end{prop}
De plus, (\cite[Proposition 6.1.2]{MR1967309})
\begin{prop}\label{conv}
\begin{enumerate}
\item $A$ vérifie les conventions de la Section 3.1 de \cite{MR1967309}, à savoir :
\begin{itemizedot}
\item Pour tout $j \in \llbracket 2, N \rrbracket$, $\sigma_j$ est un automorphisme $\mathbb{K}$-linéaire et $\delta_j$ est $\sigma_{j}$-dérivation (à gauche) $\mathbb{K}$-linéaire et localement nilpotente. 
\item Pour tout $j \in \llbracket 2, N \rrbracket$, on a $\sigma_{j} \circ \delta_{j} =q_{j} \delta_{j} \circ \sigma_{j}$ avec $q_{j} = \lambda_{j,j} = q^{-||\beta_{j}||^2}$, et pour tout $i \in \llbracket 1, j-1 \rrbracket$, $\sigma_{j}(X_{i}) = \lambda_{j,i}X_{i}$.
\item Aucun des $q_{j} \ (2 \leq j \leq N )$ n'est une racine de l'unité.
\end{itemizedot}
\item $A$ vérifie l'hypothèse 4.1.2 de \cite{MR1967309}, à savoir :\\
Le sous-groupe $H$ du groupe des automorphismes de A engendré par les $h_{l}$ vérifie 
\begin{itemizedot}
\item Pour tout h dans H, les indéterminées $X_{1}, ..., X_{N}$ sont des h-vecteurs propres.
\item L'ensemble $\{\lambda \in K^{*} | (\exists h \in H)(h(X_{1}) = \lambda X_{1}\}$ est infini.
\item Si $m \in \llbracket 2, N \rrbracket$, il existe $h_{m} \in H$ tel que $h_{m}(X_{i}) = \lambda_{m,i}X_{i}$ si $1 \leq i < m$ et $h_{m}(X_{m}) = q_{m}X_{m}$. 
\end{itemizedot}
\end{enumerate}
\end{prop}
Ceci nous permet d'utiliser la théorie de l'effacement des dérivations (\cite{MR1967309}) que nous décrivons dans le paragraphe suivant.
\subsection{L'algorithme d'effacement}
Compte tenu des propositions \ref{ore} et \ref{conv}, $A$ est une algèbre intègre est n\oe thérienne. Notons $F$ son corps des fractions. On définit récursivement les familles $X^{(l)} = (X_{i}^{(l)})_{1 \leq i \leq N}$ d'éléments de $F^{\star} := F \setminus \{0\}$, et les algèbres $A^{(l)} := \mathbb{K}<X_{1}^{(l)},...,X_{N}^{(l)}>$ lorsque $l$ décroît de $N + 1$ à 2 comme dans \cite[Section 3.2]{MR1967309}. On a donc, pour chaque $l \in \llbracket 2, N+1 \rrbracket$ :

\begin{lemme}\label{4.1.1}
Si $1 \leq i < j \leq N$, on a :
\begin{equation}\label{relationl}
X_{j}^{(l)}X^{(l)}_{i} - \lambda_{j,i}X^{(l)}_{i}X^{(l)}_{j} = P^{(l)}_{j,i}
\end{equation}
 avec 
 \begin{equation}\label{Pjil}
P^{(l)}_{j,i} = \left\{ 
\begin{array}{lc}
0 & \text{si} \ j \geq l \\
\displaystyle \sum_{\overline{k} = (k_{i+1},...,k_{j-1})} c_{\bar{k}} (X_{i+1}^{(l)})^{k_{i+1}}...(X_{j-1}^{(l)})^{k_{j-1}} & \text{si} \ j < l,
\end{array}
\right.
\end{equation}
où les $c_{\bar{k}}$ sont \textbf{les mêmes} que ceux de la formule (\ref{Pji}), de sorte que l'on dispose encore de l'implication (\ref{impli}).
\end{lemme}
\begin{preuve}
Cela résulte immédiatement de \cite[Théorème 3.2.1]{MR1967309}.
\end{preuve}
\begin{lemme}\label{4.1.2}
Les monômes ordonnés en $X_{1}^{(l)},...,X_{N}^{(l)}$ forment une base de $A^{(l)}$ comme $\mathbb{K}$ espace vectoriel.
\end{lemme}
\begin{preuve}
Cela résulte également de \cite[Théorème 3.2.1]{MR1967309}.
\end{preuve}
Des lemmes \ref{4.1.1} et \ref{4.1.2} ci-dessus et de \cite[Proposition 6.1.1]{MR1967309}, on déduit :
\begin{lemme}\label{4.1.3}
\begin{enumerate}
\item $A^{(l)}$ est une extension de Ore itérée que l'on peut écrire :
\[A^{(l)} = \mathbb{K}[X^{(l)}_{1}][X_2^{(l)};\sigma^{(l)}_2,\delta^{(l)}_2]...[X^{(l)}_N;\sigma^{(l)}_N,\delta^{(l)}_N]\]
où les $\sigma_{j}^{(l)}$ sont des automorphismes $\mathbb{K}$-linéaires et les $\delta_{j}^{(l)}$ sont des $\sigma_{j}^{(l)}$-dérivations (à gauche) $\mathbb{K}$-linéaires tels que, pour $1 \leq i < j \leq N$, $\sigma_{j}^{(l)}(X_{i}^{(l)}) = \lambda_{j,i} X_{i}^{(l)}$ et $\delta_{j}^{(l)}(X_{i}^{(l)}) = P_{j,i}^{(l)}$.
\item  $A^{(l)}$ est la $\mathbb{K}$ algèbre engendrée par les indéterminées $X_{1}^{(l)},...,X_{N}^{(l)}$ soumises aux relations (\ref{relationl}).
\end{enumerate}
\end{lemme}
Rappelons que les automorphismes $h_{m}$ ($1 \leq m \leq N$ ) de l'algèbre $A$ définis dans la proposition \ref{ore} se prolongent (de manière unique) en des automorphismes, encore notés $h_{m}$, du corps $F$.  
\begin{lemme}\label{4.1.4}
Si $1 \leq m,i \leq N$, on a  $h_{m}(X_{i}^{(l)}) = \lambda_{m,i} X_{i}^{(l)}$ de sorte que $h_{m}$ induit (par restriction) un automorphisme de l'algèbre $A^{(l)}$, noté $h_{m}^{(l)}$.
\end{lemme}
\begin{preuve}
Comme $A$ vérifie l'hypothèse 4.1.2. de \cite{MR1967309} (Proposition \ref{conv}) et comme $h_{m}(X_{i}) =  \lambda_{m,i} X_{i}$ , le lemme résulte de \cite[Lemme 4.2.1]{MR1967309}.
\end{preuve}
\begin{convsanss}
On notera $H^{(l)}$ le sous-groupe  engendré par les $h_{m}^{(l)}$ ($1 \leq m \leq N$) du groupe des automorphismes de $A^{(l)}$.
\end{convsanss}
Par \cite[Proposition 6.1.2]{MR1967309}, on a :
\begin{lemme}\label{4.1.5}
L'extension de Ore itérée $A^{(l)} = \mathbb{K}[X^{(l)}_{1}][X_2^{(l)};\sigma^{(l)}_2,\delta^{(l)}_2]...[X^{(l)}_N;\sigma^{(l)}_N,\delta^{(l)}_N]$ vérifie les conventions de la section 3.1 de \cite{MR1967309} avec, comme ci-dessus, $\lambda_{i,j} = q^{-(\beta_{i},\beta_{j})}$ et $q_{i} = \lambda_{i,i} = q^{-||\beta_{i}||^{2}}$ pour $ 1\leq i,j \leq N$. Elle vérifie également l'hypothèse 4.1.2 de \cite{MR1967309} en rempla\c cant $H$ par $H^{(l)}$.
\end{lemme}
\begin{cor}\label{Hprime}
Si $J$ est un idéal $H^{(l)}$-premier de $A^{(l)}$ au sens de \cite[II.1.9]{MR1898492}, alors J est complètement premier.
\end{cor}
\begin{preuve}
On a :
\begin{itemizedot}
\item $A^{(l)} = \mathbb{K}[X^{(l)}_{1}][X_2^{(l)};\sigma^{(l)}_2,\delta^{(l)}_2]...[X^{(l)}_N;\sigma^{(l)}_N,\delta^{(l)}_N]$ est une extension de Ore itérée (lemme \ref{4.1.3}).
\item $X^{(l)}_{1},X_2^{(l)},...,X^{(l)}_N$ sont des $H^{(l)}$-vecteurs propres (lemme \ref{4.1.4}).
\item Si $1 \leq i < j \leq N$, on a $h_{i}^{(l)}(X_{j}^{(l)}) = \lambda_{j,i} X_{i}^{(l)} = \sigma_{j}^{(l)}(X_{i}^{(l)})$ et $h_{j}^{(l)}(X_{j}^{(l)}) = q_{j} X_{j}^{(l)}$ avec $q_{j} = \lambda_{j,j} \in \mathbb{K}^{\star}$ non racine de l'unité (lemmes \ref{4.1.3} et \ref{4.1.4}). 
\end{itemizedot}
On en déduit, par \cite[Theorem II.5.12]{MR1898492}, que $J$ est complètement premier.
\end{preuve}
Il résulte de la construction de l'algorithme d'effacement ( \cite[Section 3.2]{MR1967309}) que l'on a :
\begin{lemme}\label{4.1.6}
\begin{enumerate}
\item $X_{i}^{(N+1)} = X_{i}$ pour tout $i \in \llbracket 1, N \rrbracket$
\item Si $2 \leq l \leq N$ et si $i \in \llbracket 1, N \rrbracket$, on a 
\begin{equation}\label{Xil}
X_{i}^{(l)} = 
\left\{
\begin{array}{lc}
 X_{i}^{(l+1)} & \text{si} \ i \geq l    \\
\displaystyle \sum_{n=0}^{+ \infty} \left[ \frac{(1-q_{l})^{-n}}{[n]!_{q_{l}}} \left(\delta_{l}^{(l+1)}\right)^{n}\circ \left(\sigma_{l}^{(l+1)}\right)^{-n}\left(X_{i}^{(l+1)} \right) \right]\left(X_{l}^{(l+1)} \right)^{-n} &   \text{si} \ i < l 
\end{array}
\right.
\end{equation}

\end{enumerate}

\end{lemme}
\begin{lemme}\label{4.1.7}
Soit $J$ un idéal (bilatère) $H^{(l)}$-invariant de $A^{(l)}$. Considérons un entier $j \in \llbracket 2, N \rrbracket$ et notons\\ $B = \mathbb{K}[X^{(l)}_{1}][X_2^{(l)};\sigma^{(l)}_2,\delta^{(l)}_2]...[X^{(l)}_{j-1};\sigma^{(l)}_{j-1},\delta^{(l)}_{j-1}]$ la sous-algèbre de $A^{(l)}$ engendrée par $X^{(l)}_{1},...,X^{(l)}_{j-1}$.  Alors $\sigma_{j}^{(l)}(B \cap J) = B \cap J$ et $\delta_{j}^{(l)}(B \cap J) \subset B \cap J$. 
\end{lemme}
\begin{preuve}
Par les lemmes \ref{4.1.3} et \ref{4.1.4}, on a pour $1 \leq i < j$, 
\begin{equation}\label{HS}
\sigma_{j}^{(l)}(X_{i}^{(l)}) = \lambda_{j,i}X_{i}^{(l)} = h_{j}^{(l)}(X_{i}^{(l)})
\end{equation}
Il en résulte que, pour tout $b \in B$, on a $\sigma_{j}^{(l)}(b) = h_{j}^{(l)}(b)$. Comme $J$ est $H^{(l)}$-invariant, et comme $B$ est  $\sigma_{j}^{(l)}$-invariant, on en déduit que, pour tout $b \in B \cap J$, on a $\sigma_{j}^{(l)}(b) \in B \cap J$. Ainsi, $\sigma_{j}^{(l)}(B \cap J) \subset B \cap J$. De l'égalité $\ref{HS}$, on déduit :
\begin{equation}
\left( \sigma_j^{(l)} \right)^{-1}( X_{i}^{(l)} ) = \lambda_{j,i}^{-1} X_{i}^{(l)}  = \left( h_{i}^{(l)} \right)^{-1} ( X_{i}^{(l)} )
\end{equation}
Comme ci-dessus, il en résulte que $\left( \sigma_{j}^{(l)} \right)^{-1}(B \cap J) \subset B \cap J$, de sorte que \\
$\sigma_{j}^{(l)}(B \cap J) = B \cap J$. \\
Si $b \in B \cap J$, on a alors $\delta_{j}^{(l)}(b) = X_{j}^{(l)}b - \sigma_j^{(l)}(b)X_{j}^{(l)} \in B \cap J$.
\end{preuve}
Si $l \in \llbracket 2, N \rrbracket$, il résulte de (\ref{Xil}) que $X_{l}^{(l)} = X_{l}^{(l+1)}$. Ceci est donc un élément non nul commun aux algèbres $A^{(l)}$ et $A^{(l+1)}$ (on rappelle que tous les $X_{i}^{(l)}$ sont non nuls). Ainsi, l'ensemble $S_{l} := \{(X_{l}^{(l)})^{n} | n \in \mathbb{N} \}$ est un système multiplicatif d'éléments réguliers de $A^{(l)}$ et aussi de $A^{(l+1)}$.
De \cite[Théorème 3.2.1]{MR1967309}, on déduit alors :
\begin{lemme}\label{4.1.8}
Soit $l \in \llbracket 2, N \rrbracket$. $S_{l}$ vérifie la condition de Ore (des deux côtés) dans $A^{(l)}$ et aussi dans $A^{(l+1)}$. De plus, on a :
\[A^{(l)}S_{l}^{-1} = A^{(l+1)}S_{l}^{-1}\]
\end{lemme}

\subsection{Le spectre premier et les diagrammes}
Rappelons que l'on a posé, par convention, $X_{i} = X_{\beta_{i}}$ pour $1 \leq i \leq N$. Notons $\overline{A} := A^{(2)} = \mathbb{K}<T_{\beta_{1}},...,T_{\beta_{N}}>$ avec $T_{\beta_{i}} = X_{i}^{(2)}$ pour tout $i$. Par les lemmes \ref{4.1.1} et \ref{4.1.3}, $\overline{A}$ est   l'espace affine quantique engendré par les indéterminés $T_{\beta_{i}}$ ($1 \leq i \leq N$) soumises aux relations $T_{\beta_{j}}T_{\beta_i} = \lambda_{j,i} T_{\beta_i}T_{\beta_j}$ pour $1 \leq i < j \leq N$. \\
Considérons un entier $l \in \llbracket 2, N \rrbracket$ et un idéal premier $P \in \text{Spec}(A^{(l+1)})$. \\
\par
\begin{itemizedot}
\item Supposons $X_{l}^{(l+1)} \notin P$.\\
Par \cite[Lemme 4.2.2 et 4.3.1]{MR1967309},
$S_{l} \cap P = \emptyset $ et $Q = A^{(l)} \cap P S_{l}^{-1} \in \text{Spec}(A^{(l)})$.
\item Supposons $X_{l}^{(l+1)} \in P$.\\
Par \cite[Lemme 4.3.2]{MR1967309}, il existe un (unique) homomorphisme d'algèbres surjectif 
\[g:A^{(l)} \rightarrow \frac{A^{(l+1)}}{(X_{l}^{(l+1)})}\] qui vérifie, pour tout i, $g(X_{i}^{(l)} )= \overline{X_{i}^{(l+1)}} (:= X_{i}^{(l+1)} + (X_{l}^{(l+1)}))$, de sorte que $Q = g^{-1}(\frac{P}{(X_{l}^{(l+1)})}) \in \text{Spec}(A^{(l)})$.
\end{itemizedot}
\par
On définit ainsi une application $\phi_{l} : \text{Spec}(A^{(l+1)}) \rightarrow \text{Spec}(A^{(l)}) : P \mapsto Q$ puis, par composition, une application $\phi = \phi_{2} \circ ... \circ \phi_{N} $ : Spec$(A) \rightarrow  \text{Spec}(\bar{A})$. Par \cite[Proposition 4.3.1]{MR1967309}, on a :
\begin{lemme}
Chaque $\phi_{l}$ ($2 \leq l \leq N$) est injective, de sorte que $\phi$ est injective.
\end{lemme}
\begin{dfn}\label{diagcauchon}
\begin{enumerate}
\item On appelle diagramme toute partie $\Delta $ de l'ensemble $\Phi^+$ des racines positives, et on note :
 $$\text{Spec}_{\Delta}(\overline{A}) := \{Q \in \text{Spec}(\overline{A}) \ | \ Q \cap \{T_{\beta_1},...,T_{\beta_N}\} = \{T_{\beta} \ | \ \beta \in \Delta\}\}.$$
\item Le diagramme $\Delta$ est un diagramme de Cauchon s'il existe $P \in \text{Spec}(A)$ tel que $\phi(P) \in \text{Spec}_{\Delta}(\overline{A})$, c'est à dire, si $\phi(P) \cap  \{T_{\beta_1},...,T_{\beta_N}\} =  \{T_{\beta} \ |  \ \beta \in \Delta\}$, et on note :
$$\text{Spec}_{\Delta}(A) = \{P \in \text{Spec}(A) \ | \ \phi(P) \in \text{Spec}_{\Delta}(\overline{A})\}.$$
\end{enumerate}
\end{dfn}
On a alors, par \cite[Théorèmes 5.1.1, 5.5.1 et 5.5.2]{MR1967309} : 
\begin{prop}\label{Hstrat}
\begin{enumerate}
\item Si $\Delta$ est un diagramme de Cauchon, on a  $\phi (\text{Spec}_{\Delta}(A)) = \text{Spec}_{\Delta}(\overline{A})$ et $\phi$ induit un homéomorphisme bi-croissant de $\text{Spec}_{\Delta}(A)$ sur $\text{Spec}_{\Delta}(\overline{A})$.
\item La famille $\text{Spec}_{\Delta}(A)$ (avec $\Delta$ diagramme de Cauchon) coïncide avec la H-stratification de Goodearl-Letzter de Spec$(A)$ (\cite{MR1898492}).
\end{enumerate}
\end{prop}
Dans la suite, on va chercher à décrire plus précisément les diagrammes de Cauchon, pour cela on aura besoin du critère de la proposition suivante :
\begin{prop} \label{critereimage}
Soit $P\it^{(m)}$ un idéal premier $H$-invariant de $A^{(m)}$.\\
$P\it^{(m)}$ $ $ $\in$ $ $ $Im(\phi_{m})$ $ $ si et seulement si l'une des conditions suivantes est satisfaite 
\begin{enumerate}
\item $ $ $X_{m}^{(m)}$ $ $ $\notin$ $ $ $P\it^{(m)}$. 
\item $ $ $X_{m}^{(m)}$ $ $ $\in$ $ $ $P\it^{(m)}$ $ $ et $ $ $\Theta^{(m)}(\delta_{m}^{(m+1)}(X_{i}^{(m+1)})) $  $\in$ $ $ $P\it^{(m)}$ $ $ pour $ $ $1$  $\leq$ $i$ $\leq$ $m-1$.\\
(où $ $ $\delta_{m}^{(m+1)}$($X_{i}^{(m+1)}$) =  $P_{m,i}^{(m+1)}$ ($X_{i+1}^{(m+1)}, \mbox{ } ... \mbox{ } , \mbox{ } X_{m-1}^{(m+1)}$) (lemme \ref{4.1.1}) et $ $ où $ $ $\Theta^{(m)}: \mbox{ } k<X_{1}^{(m+1)}, \mbox{ } ... \mbox{ } , \mbox{ } X_{m-1}^{(m+1)}> \mbox{ } \rightarrow \mbox{ } k<X_{1}^{(m)}, \mbox{ } ... \mbox{ } , \mbox{ } X_{m-1}^{(m)}>$ $ $ est l'homomorphisme qui transforme chaque $ $ $X_{l}^{(m+1)}$ $ $ en $ $ $X_{l}^{(m)}$)
\end{enumerate}
\end{prop}
\begin{preuve}
Supposons que $ $ $P\it^{(m)}$ $ $ $\in$ $ $ $Im(\phi_{m})$, $ $ de sorte que $ $ $P\it^{(m)}$ $ $ = $ $ $\phi_{m}$($P\it^{(m+1)}$) $ $ avec  $ $ $P\it^{(m+1)}$ $ $ $\in$ $ $ $Spec$($A^{(m+1)}$), $ $ et on suppose que la condition 1. n'est pas vérifiée. Ceci implique que $ $ $P\it^{(m)}$ $ $ = $ $ $ker(g)$ $ $ où $ $ $g: \mbox{ } A^{(m)} \mbox{ } \rightarrow \mbox{ } A^{(m+1)} / P\it^{(m+1)}$ $ $ est l'homomorphisme qui transforme chaque $ $ $X_{i}^{(m)}$ $ $ en $ $ $x_{i}^{(m+1)}$ $ $ = $ $ $X_{i}^{(m+1)}$ $ $ + $ $ $P\it^{(m+1)}$.
On considère $ $ $1$  $\leq$ $i$ $\leq$ $m-1$. \\ Rappelons que $ $ $\delta_{m}^{(m+1)}$($X_{i}^{(m+1)}$) $ $ = $ $ $P_{m,i}^{(m+1)}$($X_{i+1}^{(m+1)}, \mbox{ } ... \mbox{ } , \mbox{ } X_{m-1}^{(m+1)}$) et que $ $ $\Theta^{(m)}: \mbox{ } k<X_{1}^{(m+1)}, \mbox{ } ... \mbox{ } , \mbox{ } X_{m-1}^{(m+1)}> \mbox{ } \rightarrow \mbox{ } k<X_{1}^{(m)}, \mbox{ } ... \mbox{ } , \mbox{ } X_{m-1}^{(m)}>$ $ $ est l'homomorphisme qui transforme chaque $ $ $X_{l}^{(m+1)}$ $ $ en $ $ $X_{l}^{(m)}$. $ $ Comme $ $ $X_{m}^{(m)}$ $ $ $\in$ $ $ $P\it^{(m)}$, $ $ on a $ $ $X_{m}^{(m+1)}$ $ $ $\in$ $ $ $P\it^{(m+1)}$ $ $ (\cite{MR1967309}, proposition 4.3.1.) et donc, $ $ $\delta_{m}^{(m+1)}$($X_{i}^{(m+1)}$) $ $ $\in$ $ $ $P\it^{(m+1)}$. $ $ A présent, on a \\ $g$($\Theta^{(m)}$($\delta_{m}^{(m+1)}$($X_{i}^{(m+1)}$))) $ $ = $ $ $g$($\Theta^{(m)}$($P_{m,i}^{(m+1)}$($X_{i+1}^{(m+1)}, \mbox{ } ... \mbox{ } , \mbox{ } X_{m-1}^{(m+1)}$))) $ $ = $ $ $g$($P_{m,i}^{(m+1)}$($X_{i+1}^{(m)}, \mbox{ } ... \mbox{ } , \mbox{ } X_{m-1}^{(m)}$)) $ $ = $ $ $P_{m,i}^{(m+1)}$($x_{i+1}^{(m+1)}, \mbox{ } ... \mbox{ } , \mbox{ } x_{m-1}^{(m+1)})$ $ $ = $ $ $P_{m,i}^{(m+1)}$($X_{i+1}^{(m+1)}, \mbox{ } ... \mbox{ } , \mbox{ } X_{m-1}^{(m+1)})$ $ $ + $ $ $P\it^{(m+1)}$ $ $ = $ $ $0$. \\ On en déduit que $ $ $\Theta^{(m)}$($\delta_{m}^{(m+1)}$($X_{i}^{(m+1)}$))) $ $ = $ $ $ker(g)$ $ $ = $ $ $P\it^{(m)}$. \\ $ $ \\
Si la condition 1. est vérifiée, alors $ $ $P\it^{(m)}$ $ $ $\in$ $ $ $Im(\phi_{m})$ $ $ par (\cite{MR1967309}, lemma 4.3.1.). \\ $ $ \\
Supposons que la condition 2. est vérifiée. Alors si $ $ $1$  $\leq$ $i$ $\leq$ $m-1$, $ $ on a, comme précédemment, $ $  $P_{m,i}^{(m+1)}$($X_{i+1}^{(m)}, \mbox{ } ... \mbox{ } , \mbox{ } X_{m-1}^{(m)}$) $ $ = $ $ $\Theta^{(m)}$($\delta_{m}^{(m+1)}$($X_{i}^{(m+1)}$)) $ $ $\in$ $ $ $P\it^{(m)}$. $ $ Donc, dans $ $ $Q^{(m)}$ $ $ = $ $ $A^{(m)}/P\it^{(m)}$, $ $  on a  $ $  $P_{m,i}^{(m+1)}$($x_{i+1}^{(m)}, \mbox{ } ... \mbox{ } , \mbox{ } x_{m-1}^{(m)}$) $ $ = $ $ $0$. \\ Comme $ $ $P_{m,i}^{(m)}$ $ $ = $ $ $0$ $ $ (à cette étape, $X_m^{(m)}$ est normalisant), on peux écrire $ $ $ $ $x_{m}^{(m)}x_{i}^{(m)}$ $ $ $-$ $ $ $\lambda_{m,i}x_{i}^{(m)}x_{m}^{(m)}$ $ $ = $ $ $P_{m,i}^{(m)}$($x_{i+1}^{(m)}, \mbox{ } ... \mbox{ } , \mbox{ } x_{m-1}^{(m)}$) $ $ = $ $ $0$ $ $ = $ $ $P_{m,i}^{(m+1)}$($x_{i+1}^{(m)}, \mbox{ } ... \mbox{ } , \mbox{ } x_{m-1}^{(m)}$). \\  Si  $ $ $1$  $\leq$ $i$ $\leq$ $j-1$ $ $ avec $j$ $ $ $\neq$ $ $ $m$, $ $ on a (d'après le lemme \ref{4.1.1}) \\ $ $ $x_{j}^{(m)}x_{i}^{(m)}$ $ $ $-$ $ $ $\lambda_{j,i}x_{i}^{(m)}x_{j}^{(m)}$ $ $ = $ $ $P_{j,i}^{(m)}$($x_{i+1}^{(m)}, \mbox{ } ... \mbox{ } , \mbox{ } x_{j-1}^{(m)}$) $ $ = $ $ $P_{j,i}^{(m+1)}$($x_{i+1}^{(m)}, \mbox{ } ... \mbox{ } , \mbox{ } x_{j-1}^{(m)}$). \\ Donc, par la propriété universelle des algèbres définies par générateurs et relations, il existe un (unique) homomorphisme $ $ $\epsilon: \mbox{ } A^{(m+1)} \mbox{ } \rightarrow \mbox{ } Q^{(m)}$ $ $ qui transforme chaque $ $ $X_{l}^{(m+1)}$ $ $ en $ $ $x_{l}^{(m)}$. $ $ Cet homomorphisme est surjectif, et son noyau $ $ $ker(\epsilon)$ $ $ = $ $ $P\it^{(m+1)}$ $ $ est un ideal premier de $ $ $A^{(m+1)}$. $ $ On observe que, comme $ $ $X_{m}^{(m)}$ $ $ $\in$ $ $ $P\it^{(m)}$, $ $ on a $ $ $X_{m}^{(m+1)}$ $ $ $\in$ $ $ $P\it^{(m+1)}$, $ $  et que $ $ $\epsilon$ $ $ induit un isomorphisme
$$ \overline{\epsilon}: \mbox{ } A^{(m+1)} / P\it^{(m+1)} \mbox{ } \rightarrow \mbox{ } Q^{(m)} \mbox{ } = \mbox{ } A^{(m)} / P\it^{(m)}$$ qui transforme chaque  $ $ $x_{l}^{(m+1)}$ $ $ en $ $ $x_{l}^{(m)}$. $ $ Rappelons que $ $ $f_{m}: \mbox{ } A^{(m)} \mbox{ } \rightarrow \mbox{ } A^{(m)} / P\it^{(m)}$  désigne la projection canonique. \\ Ainsi, $g \mbox{ } = \mbox{ } (\overline{\epsilon})^{-1} \circ f_{m}: \mbox{ } A^{(m)} \mbox{ } \rightarrow \mbox{ } A^{(m+1)} / P\it^{(m+1)}$ $ $ est l'homomorphisme qui transforme chaque $ $ $X_{l}^{(m)}$ $ $ en $ $ $x_{l}^{(m+1)}$. $ $ Comme $ $ $ker(g)$ $ $ = $ $ $ker(f_{m})$ $ $ = $ $ $P\it^{(m)}$, $ $ on en déduit que $ $ $P\it^{(m)}$ $ $ = $ $ $\phi_{m}$($P\it^{(m+1)}$).
\end{preuve}
\section{Diagrammes de Cauchon dans $\mathcal{U}_q^+(\cat{g})$}
Dans \cite{MR1967310}, G. Cauchon utilise un moyen combinatoire pour décrire les diagrammes "admissibles" (appelés ici "diagrammes de Cauchon") pour $O\sb q(M\sb n(k))$. A l'aide de la théorie des plans admissibles de Lusztig (voir section 3.2) et des résultats démontrés dans la section 3.3, nous décrivons ces diagrames pour $\mathcal{U}_q^+(\cat{g})$ (avec $\cat{g}$ algèbre de Lie simple complexe de dimension finie). Le but de cette partie est de démontrer le résultat suivant :\\
Un diagramme $\Delta$ vérifie toutes les contraintes provenant des plans admissibles (notion à définir) si et seulement si $\Delta$ est un diagramme de Cauchon (au sens de la définition \ref{diagcauchon}) .\\
\subsection{Contraintes dans un diagramme}

\begin{lemme}\label{lemme4.1}
Soit $j \in \llbracket 1, N \rrbracket$. Soient $l \in \llbracket 2, N \rrbracket$, $P^{(l+1)}$ un idéal premier de $A^{(l+1)}$ et $P^{(l)} = \varphi_l (P^{(l+1)})$.
\begin{enumerate}
\item Si $X_{j}^{(l+1)} \in P^{(l+1)}$ alors $X_{j}^{(l)} \in P^{(l)}$.
\item Si $X_{j}^{(l+1)} = X_{j}^{(l)}$ (donc, en particulier, si $j \geq l$), on a $X_{j}^{(l+1)} \in P^{(l+1)} \Leftrightarrow X_{j}^{(l)} \in P^{(l)}$.
\end{enumerate}
\end{lemme}
\begin{preuve}
Le 2. se démontre comme \cite[lemme 4.3.4]{MR1967309}. Montrons donc le 1. lorsque $j < l$.
\begin{description}
\item[$1^{er}$ cas :] Si le pivot $\varpi:=X_{l}^{(l+1)} \in P^{(l+1)}$, rappelons (cf. Section 4.2) qu'il existe un  homomorphisme d'algèbres surjectif 
\[g:A^{(l)} \rightarrow \frac{A^{(l+1)}}{(X_{l}^{(l+1)})}\] qui vérifie $g(X_{i}^{(l)} )= \overline{X_{i}^{(l+1)}} (:= X_{i}^{(l+1)} + (X_{l}^{(l+1)}))$ pour tout $i \in \llbracket 1, N \rrbracket$. Comme $X_{j}^{(l+1)} \in P^{(l+1)}$, on a $g(X_{j}^{(l)}) \in \frac{P^{(l+1)}}{(X_{l}^{(l+1)})}$ de sorte que $X_{j}^{(l)} \in g^{-1}(\frac{P^{(l+1)}}{(X_{l}^{(l+1)})}):= P^{(l)}$.

\item[$2^{nd}$ cas :]  Supposons que le pivot $\varpi:=X_{l}^{(l+1)} \notin P^{(l+1)}$ et notons $S_l := \{\varpi^n | n \in \mathbb{N}\}$. Rappelons (cf Section 4.2) que l'on a alors $P^{(l)} = A^{(l)} \cap (P^{(l+1)}S_{l}^{-1})$.\\
Notons  $\displaystyle J := \bigcap_{h \in H^{(l+1)}} h( P^{(l+1)})$ et observons que J est un idéal bilatère $H^{(l+1)}$-invariant par construction. Comme $A^{(l+1)}$ vérifie l'hypothèse 4.1.2 de  \cite{MR1967309} (Lemme \ref{4.1.5}), $X_{j}^{(l+1)}$  est un $H^{(l+1)}$-vecteur propre. De là, puisque $X_{j}^{(l+1)}$ appartient à $P^{(l+1)}$, il appartient aussi à J.\\
Par le lemme \ref{4.1.7}, on en déduit que $\left(\delta_{l}^{(l+1)}\right)^{n}\circ \left(\sigma_{l}^{(l+1)}\right)^{-n}\left(X_{j}^{(l+1)} \right) \in J \subset P^{(l+1)}$ pour tout $n \in \mathbb{N}$. Il en résulte que :
\[X_{j}^{(l)} = \sum_{n=0}^{+ \infty} \left[ \frac{(1-q_{l})^{-n}}{[n]!_{q_{l}}}\left(\delta_{l}^{(l+1)}\right)^{n}\circ \left(\sigma_{l}^{(l+1)}\right)^{-n}\left(X_{j}^{(l+1)} \right) \right] \left(X_{l}^{(l+1)} \right)^{-n} \in  P^{(l+1)}S_{l}^{-1}.\] 
Ainsi, $X_{j}^{(l)} \in A^{(l)} \cap (P^{(l+1)}S_{l}^{-1}) = P^{(l)}$. 
\end{description}
\end{preuve}
\begin{lemme}\label{5.2}
Soient $l \in \llbracket 2, N \rrbracket$, $P^{(l+1)}$ un idéal premier de $A^{(l+1)}$. Considérons un entier $j$ avec $2 \leq j < l$ et notons $P^{(j)} = \varphi_j \circ ... \circ \varphi_l (P^{(l+1)})$.
\begin{enumerate}
\item Supposons que $\beta_j$ est dans la boîte de $\beta_l$ ou dans la boîte précédente. Alors
\begin{itemizedot}
\item $X_j^{(j+1)} = X_j^{(j+2)} = ... = X_j^{(l+1)}$,
\item $X_j^{(j+1)} \in P^{(j+1)} \Rightarrow X_j^{(j+2)} \in P^{(j+2)} \Rightarrow ... \Rightarrow X_j^{(l+1)} \in P^{(l+1)}$.
\end{itemizedot}   
\item Supposons que les boites $B$ et $B'$ de $\beta_j$ et $\beta_l$ (respectivement) soient séparées par une boite $B''$ réduite à un élément $\beta_e$ tel que $X_e^{(e+1)} \in P^{(e+1)}$. Alors $X_j^{(j+1)} \in P^{(j+1)} \Rightarrow X_j^{(l+1)} \in P^{(l+1)}$.
\end{enumerate}
\end{lemme}
\begin{preuve}
\begin{enumerate}
\item Soit $k \in  \llbracket j+1, l\rrbracket$ de sorte que $\beta_k$ est, soit dans la boîte de $\beta_j$, soit dans celle de $\beta_l$. Comme ces boîtes sont consécutives, on a $X_kX_j = q^{-<\beta_k,\beta_j>}X_jX_k$, de sorte que par le lemme \ref{4.1.1}, $X_k^{(k+1)}X_j^{(k+1)} = q^{-<\beta_k,\beta_j>}X_j^{(k+1)}X_k^{(k+1)}$. On a donc $\delta_k^{(k+1)}(X_j^{(k+1)}) = 0$ et, par \cite[Section 3.2]{MR1967309}, on a :
\[X_j^{(k)} = \sum_{s=0}^{+ \infty} \lambda_s \left( \delta_k^{(k+1)}\right)^s \circ \left( \sigma_k^{(k+1)} \right)^{-s}(X_j^{(k+1)})\left( X_k^{(k+1)}\right)^{-s} = \sum_{s=0}^{+ \infty} \lambda'_s \left( \delta_k^{(k+1)}\right)^s(X_j^{(k+1)})\left( X_k^{(k+1)}\right)^{-s} =  X_j^{(k+1)}\]
($\lambda_s, \lambda'_s \in \mathbb{K}$).\\
Ceci montre le premier point. Le second point résulte du lemme \ref{lemme4.1}.
\item $B$ et $B''$ étant consécutives, il résulte du point 1. que $X_{j}^{(j+1)} = ... = X_{j}^{(e+1)}$ et que $X_j^{(j+1)} \in P^{(j+1)} \Rightarrow ... \Rightarrow X_j^{(e+1)} \in P^{(e+1)}$. Pour terminer montrons, par récurrence sur k, que :
\[X_j^{(k)} \in P^{(k)} \Rightarrow X_j^{(k+1)} \in P^{(k+1)} \ \text{pour} \ e+1 \leq k \leq l\] 
On écrit comme dans le 1. :
\[X_j^{(k)} = X_j^{(k+1)} + \sum_{s=1}^{+ \infty} \lambda_s \left( \delta_
k^{(k+1)}\right)^s \circ \left( \sigma_k^{(k+1)} \right)^{-s}(X_j^{(k+1)})\left( X_k^{(k+1)}\right)^{-s} \ \ (\lambda_s \in \mathbb{K})\]
\begin{itemizedot}
\item Si $\delta_k^{(k+1)}(X_j^{(k+1)}) = 0$, alors on a $X_j^{(k)} = X_j^{(k+1)}$ et on conclut par le lemme \ref{lemme4.1}.
\item Sinon, on a $\delta_k^{(k+1)}(X_j^{(k+1)}) = \lambda \left(X_e^{(k+1)}\right)^m$ $(m \in \mathbb{N}^\star, \lambda \in \mathbb{K}^\star$) par le lemme \ref{4.1.1} et, puisque $B'$ et $B''$ sont adjacentes,   
\[ \begin{array}{cl}
\delta_k^{(k+1)}\left( X_e^{(k+1)} \right) = 0 & \Rightarrow \left( \delta_k^{(k+1)}\right)^s(X_j^{(k+1)}) =  \lambda \left( \delta_k^{(k+1)}\right)^{s-1}\left( \left(X_e^{(k+1)}\right)^m \right) =0 \ \text{pour} \ s>1\\
 & \Rightarrow X_j^{(k)} = X_j^{(k+1)} + \lambda' \left(X_e^{(k+1)}\right)^m\left( X_k^{(k+1)}\right)^{-1} \ \text{avec} \ \lambda' \in \mathbb{K}^\star.
 \end{array}\]
\begin{itemizedot}
\item Si $X_k^{k+1} \in P^{(k+1)}$, on considère l'homomorphisme $g : A^{(k)} \rightarrow \frac{A^{(k+1)}}{(X_k^{(k+1)})}$ qui vérifie $g(X_i^{(k)}) = \overline{X_i^{(k+1)}}$ pour $i \in \llbracket 1, N \rrbracket$ (cf Section 4.2). Par définition de $\phi_k$ (\cite[Notation 4.3.1.]{MR1967309}), on a $P^{(k)} = g^{-1}\left( \frac{P^{(k+1)}}{(X_k^{(k+1)})} \right)$. Donc  
\[X_j^{(k)} \in P^{(k)} \Rightarrow g(X_j^{(k)}) = \overline{X_j^{(k+1)}} \in \frac{P^{(k+1)}}{(X_k^{(k+1)})} \Rightarrow  X_j^{(k+1)} \in P^{(k+1)}\]
\item  Par le 1., on a $X_e^{(e+1)}  = ... = X_e^{(k)}= X_e^{(k+1)}$ et $X_e^{(e+1)} \in P^{(e+1)} \Rightarrow ... \Rightarrow X_e^{(k)} \in P^{(k)} \Rightarrow X_e^{(k+1)} \in P^{(k+1)}$. Posons, comme dans \cite[Théorème 3.2.1]{MR1967309}, $S_k := \{\left(X_k^{(k+1)}\right)^n | n \in \mathbb{N} \}$ de sorte que $P^{(k+1)} = A^{(k+1)}\cap (P^{(k)}S_k^{-1}) $ par définition de $\varphi_k$ \cite[Notation 4.3.1.]{MR1967309}. On a alors :\\
\[X_j^{(k+1)} = X_j^{(k)} - \lambda' \left(X_e^{(k+1)}\right)^m\left( X_k^{(k+1)}\right)^{-1} = X_j^{(k)} - \lambda' \left(X_e^{(k)}\right)^m\left( X_k^{(k+1)}\right)^{-1} \in P^{(k)}S_k^{-1}.\]
Comme $X_j^{(k+1)}$ est aussi dans $A^{(k+1)}$, on a bien $X_j^{(k+1)} \in P^{(k+1)}$.
\end{itemizedot} 
\end{itemizedot}
\end{enumerate}
\end{preuve}
Pour trouver la forme des diagrammes de Cauchon, on va se servir de \cite[proposition 5.2.1]{MR1967309}. En tenant compte des notations utilisées ici, cette proposition s'écrit : 
\begin{prop}\label{prop4.2}
Soit $\Delta$ un diagramme de Cauchon et soit $P \in \text{Spec}(A)$. Pour que $P$ appartienne à Spec$_{\Delta}(A)$ il faut et il suffit qu'il vérifie le critère ci-dessous :
\[ (\forall \ l \in \llbracket 1, N \rrbracket ) \ \ (X_{l}^{(l+1)} \in P^{(l+1)} \Leftrightarrow \beta_{l} \in \Delta)  \]
\end{prop}
Ceci permet de démontrer la proposition suivante :
\begin{prop}\label{prop4.4.3}
Soient $\Delta$ un diagramme de Cauchon et $\beta_{l} \in \Delta$ ($1 \leq l \leq N$). Supposons qu'il existe un entier $k \in \llbracket 1, l-1\rrbracket$ vérifiant $X_{\beta_{l}}X_{\beta_{k}}-q^{-(\beta_{l},\beta_{k})}X_{\beta_{k}}X_{\beta_{l}} = c X_{\beta_{i_{1}}}...X_{\beta_{i_{s}}}$ avec $c \in \mathbb{K}^{\star}, s \geq 1$ et $k < i_{1} \leq ... \leq i_{s} < l$. Alors l'un des $\beta_{i_{r}} \ (1 \leq r \leq s)$ appartient à $\Delta$.
\end{prop}
\begin{preuve}
Soit $P \in \text{Spec}_{\Delta}(A)$. Par le lemme  \ref{4.1.1}, on a:
 $$X_{l}^{(l+1)}X^{(l+1)}_{k}-q^{-(\beta_{l},\beta_{k})}X^{(l+1)}_{k}X^{(l+1)}_{l} = c X^{(l+1)}_{\beta_{i_{1}}}...X^{(l+1)}_{\beta_{i_{s}}} := M.$$
Par la proposition \ref{prop4.2}, on a $X_{l}^{(l+1)} \in P^{(l+1)}$ de sorte que $M \in P^{(l+1)}$. Comme $P^{(l+1)}$ est un idéal premier de $A^{(l+1)}$, on sait (\cite[II.6.9]{MR1898492}) que $P^{(l+1)}$ est complètement premier de sorte qu'il existe $r \in \llbracket 1, s \rrbracket$ tel que $X_{\beta_{i_{r}}}^{(l+1)} \in P^{(l+1)}$. Par le lemme \ref{lemme4.1}, on en déduit $X_{i_{r}}^{(i_{r}+1)} \in P^{(i_{r}+1)}$ et, par la proposition \ref{prop4.2}, $\beta_{i_{r}} \in \Delta$.
\end{preuve}
\begin{conv}
On dit qu'un diagramme $\Delta$ vérifie la contrainte
\begin{enumerate}
\item $\beta_{j_{0}} \rightarrow \beta_{j_{1}}$  si $\beta_{j_{0}} \in \Delta \Rightarrow \beta_{j_{1}} \in \Delta$.
\item \begin{tabular}{cc}
\scalebox{.8}{\begin{psmatrix}[rowsep=.2cm, colsep=1cm]
& $\beta_{j_{1}}$\\
$\beta_{j_{0}}$ & \vdots \\
& $\beta_{j_{s}}$
\end{psmatrix}
\psset{nodesep=2pt, linestyle=dashed, arrows=->}
\ncline{2,1}{1,2}
\ncline{2,1}{2,2}
\ncline{2,1}{3,2}}
&
si $\beta_{j_{0}} \in \Delta \Rightarrow ( \beta_{j_{1}} \in \Delta$ ou ... ou $\beta_{j_{1}} \in \Delta )$.
\end{tabular}
\end{enumerate}
\end{conv}
La proposition \ref{prop4.4.3} peut donc se réécrire :
\begin{prop} \label{prop4.4}
Soient $\Delta$ un diagramme de Cauchon et $\beta_{l} \in \Delta$ ($1 \leq l \leq N$). Supposons qu'il existe un entier $k \in \llbracket 1, l-1\rrbracket$ vérifiant $X_{\beta_{l}}X_{\beta_{k}}-q^{-(\beta_{l},\beta_{k})}X_{\beta_{k}}X_{\beta_{l}} = c X_{\beta_{i_{1}}}^{m_{1}}...X_{\beta_{i_{s}}}^{m_{s}}$ avec $c \in \mathbb{K}^{\star}, s \geq 1, k < i_{1} < ... < i_{s} < l$ et  $m_{1},...,m_{s} \in \mathbb{N}^{\star}$. Alors 
\begin{enumerate}
\item Si $s = 1$, la flèche pleine $\beta_{l} \rightarrow \beta_{i_{1}}$ est une contrainte 
\item Si $s \geq 2$, le système 
\begin{center}
\scalebox{.8}{\begin{psmatrix}[rowsep=.2cm, colsep=1cm]
& $\beta_{i_{1}}$\\
$\beta_{l}$ & \vdots \\
& $\beta_{i_{s}}$
\end{psmatrix}
\psset{nodesep=2pt, linestyle=dashed, arrows=->}
\ncline{2,1}{1,2}
\ncline{2,1}{2,2}
\ncline{2,1}{3,2}}
\end{center}
est une contrainte.
\end{enumerate}
\end{prop}
Dans les trois propriété qui suivent, on considère $\Delta$ un diagramme de Cauchon.
\begin{prop}\label{contrainteelementaire}
Soient $1 \leq l \leq n$ et $\beta \in C_{l}$. S'il existe $i \in \llbracket 1, l-1 \rrbracket$ tel que $\beta + \alpha_{i} = m \beta'$ avec $m \in \mathbb{N}^{\star}$ et $\beta' \in \Phi^{+}$, alors $\beta \rightarrow \beta'$ est une contrainte. 
\end{prop}
\begin{preuve}
On sait (proposition \ref{propG2} lorsque $\Phi = G_2$, corollaire \ref{corplanadm} lorsque $\Phi \neq G_2$) que l'on dispose alors d'une formule de redressement de la forme $E_{\beta}E_{\alpha_{i}} - q^{(\beta,\alpha_{i})}E_{\alpha_{i}}E_{\beta} = k E^{m}_{\beta'}$ avec $k \neq 0$ (où les $E_{\gamma}$ sont les générateurs de Lusztig de $\mathcal{U}_{q}^{+}(\mathfrak{g}))$.\\
Il résulte alors de la proposition \ref{relcomX} que $X_{\beta}X_{\alpha_{i}} - q^{-(\beta,\alpha_{i})}X_{\alpha_{i}}X_{\beta} = k' X^{m}_{\beta'}$ avec $k' \neq 0$ puis, de la proposition \ref{prop4.4}, que $\beta \rightarrow \beta'$ est une contrainte.
\end{preuve}
\begin{prop}\label{contrainteelementaireavantex}
Soit $C_{l}$ ($1 \leq l \leq n$) une colonne exceptionnelle. Si $\beta \in C_{l}$ est dans la boîte qui suit celle de la racine exceptionnelle $\beta_{ex}$, alors $\beta \rightarrow \beta_{ex}$ est une contrainte.  
\end{prop}
\begin{preuve}
Supposons d'abord que $\Phi$ est de type $G_2$. Avec les notations de la proposition \ref{propG2}, on a donc $l=2$, $\beta_{ex} = \beta_4$, $\beta = \beta_5$ et on dispose d'une formule de redressement de la forme $E_{\beta_5}E_{\beta_3} - q^{(\beta_3,\beta_5)}E_{\beta_3}E_{\beta_5} = k E_{\beta_{4}}$ avec $k \in \mathbb{K}^\star$. On en déduit, par la proposition \ref{prop4.4} que $\beta = \beta_5 \rightarrow \beta_{ex}= \beta_4$ est une contrainte.\\
\par 
Supposons maintenant $\Phi \neq G_2$. On sait (proposition \ref{h'}) que $h'(\beta_{ex}) = t + \frac{1}{2}$ ($t \in \mathbb{N}^{\star}$), de sorte que $h'(\beta) = h(\beta) = t$. On sait également (proposition \ref{prop2.8}) que si $D = $Vect$(\beta_{ex})$, on a $\beta' = s_{D}(\beta) = \beta_{ex} - \beta \in C_{l}$, de sorte que $h'(\beta') = h(\beta') = h(\beta_{ex}) - h(\beta) =  t+1$. Il en résulte (définition \ref{planadm} et remarque \ref{systracine}) que $P = $Vect$(\beta, \beta')$ est un plan admissible de type (1.1) ou (1.2). Par la proposition \ref{rellusztig} on dispose alors d'une formule de redressement de la forme $E_{\beta}E_{\beta'} - q^{(\beta,\beta')}E_{\beta'}E_{\beta} = k E_{\beta_{ex}}$ avec $k \in \mathbb{K}^\star$. On en déduit, par la proposition \ref{prop4.4}, que $\beta \rightarrow \beta_{ex}$ est une contrainte.  
\end{preuve}
\begin{prop}\label{contrainteelementaireapresex}
Soit $C_{l}$ ($1 \leq l \leq n$) une colonne exceptionnelle et $\beta_{ex}$ sa racine exceptionnelle. Supposons qu'il existe $i \in \llbracket 1, l \llbracket$ tel que $\beta_{ex} + \alpha_{i} = \beta'_{i_{1}} + \beta'_{i_{2}}$ avec $\beta'_{i_{1}} \neq \beta'_{i_{2}}$ dans la boîte qui précède celle de $\beta_{ex}$. Alors le système \begin{center}
\scalebox{.8}{\begin{psmatrix}[rowsep=.2cm, colsep=1cm]
& $\beta'_{i_{1}}$\\
$\beta_{ex}$ & \\
& $\beta'_{i_{2}}$
\end{psmatrix}
\psset{nodesep=2pt, linestyle=dashed, arrows=->}
\ncline{2,1}{1,2}
\ncline{2,1}{3,2}}
\end{center}
est une contrainte.
\end{prop}
\begin{preuve}
Comme, par hypothèse, $\beta'_{i_{1}} \neq \beta'_{i_{2}}$ sont dans la boîte qui précède celle de $\beta_{ex}$, le système de racines n'est pas de type $G_2$ (cf proposition \ref{propG2}).\\
Comme dans la preuve proposition \ref{contrainteelementaire}, il suffit de montrer que :
$$[E_{\beta_{ex}},E_{\alpha_{i}}]_{q}:=E_{\beta_{ex}}E_{\alpha_{i}} - q^{(\beta_{ex},\alpha_{i})}E_{\alpha_{i}}E_{\beta_{ex}} = \lambda E_{\beta'_{i_{1}}}E_{\beta'_{i_{2}}} \ \text{avec} \ \lambda \in \mathbb{K}^{\star}$$ 
Rappelons (proposition \ref{prop2.8}) que $\beta_{ex} \ \bot \ \alpha_{i}$, de sorte que  :
$$(\alpha_{i} , \beta'_{i_{1}} + \beta'_{i_{2}} ) = (\alpha_{i}, \beta_{ex} + \alpha_{i}) = ||\alpha_{i}||^2 \Rightarrow (\alpha_{i} , \beta'_{i_{1}} ) > 0 \ \text{ou} \ (\alpha_{i},\beta'_{i_{2}} ) > 0.$$
On peut donc supposer, sans perte de généralité, que $(\alpha_{i},\beta'_{i_{2}} ) > 0$, de sorte que (corollaire \ref{corplanadm}) $[E_{\beta'_{i_{2}}},E_{\alpha_{i}}]_{q} = 0$.\\
\par
Comme dans la preuve de la proposition précédente, on a :
\begin{itemizedot}
\item $h'(\beta_{ex}) = t + \frac{1}{2}$ ($t \in \mathbb{N}^{\star}$) et $h'(\beta'_{i_{1}}) = h'(\beta'_{i_{2}}) = t + 1$,
\item  $\beta_{i_{1}} = s_{D} (\beta'_{i_{1}})$ et $\beta_{i_{2}} = s_{D} (\beta'_{i_{2}})$ sont alors dans $C_{l}$ et vérifient $h'(\beta_{i_{1}}) = h'(\beta_{i_{2}}) = t$,
\item $E_{\beta_{i_{2}}}E_{\beta'_{i_{2}}} - q^{(\beta_{i_{2}},\beta'_{i_{2}})}E_{\beta'_{i_{2}}}E_{\beta_{i_{2}}} = k E_{\beta_{ex}}$ avec $k \neq 0$. \hspace{2cm} $(\star)$ 
\end{itemizedot}
$\mbox{ }$\\
Par définition de $\beta_{i_{2}}$, on a $\beta_{ex} = \beta_{i_{2}} + \beta'_{i_{2}}$, de sorte que \\
$$\beta'_{i_{1}} + \beta'_{i_{2}} =\beta_{ex} + \alpha_{i}  = \beta_{i_{2}} + \beta'_{i_{2}} + \alpha_{i} \Rightarrow \beta'_{i_{1}} = \beta_{i_{2}} + \alpha_{i}.$$
Par le corollaire \ref{corplanadm}, on a donc $\ [E_{\beta_{i_{2}}},E_{\alpha_{i}}]_{q}:= h E_{\beta'_{i_{1}}} \ (h \neq 0).$
\par
On sait que $\mathcal{U}_q^+(\mathfrak{g})$ est $\mathbb{Z}\Phi$-graduée. Il existe donc un (unique) automorphisme $\sigma$ de $\mathcal{U}_q^+(\mathfrak{g})$ 
qui vérifie pour tout $Q \in \mathcal{U}_q^+(\mathfrak{g})$, homogène de degré $\beta$, $\sigma(Q) = q^{(\beta,\alpha_{i})} Q$.\\
Notons $\delta$ la $\sigma$-dérivation intérieure à droite associée à $E_{\alpha_{i}}$, de sorte que \\
 $$\delta(Q) = QE_{\alpha_{i}} - E_{\alpha_{i}}\sigma(Q) \ \ \ \ \ (\forall Q \in \mathcal{U}_q^+(\mathfrak{g}))$$
Si $\beta \in C_{l}$, on a $\delta(E_{\beta}) = E_{\beta}E_{\alpha_{i}} - q^{(\beta,\alpha_{i})}E_{\alpha_{i}}E_{\beta} = [E_{\beta},E_{\alpha_{i}}]_{q}$ et, de là :
$\delta(E_{\beta'_{i_{2}}}) = 0$ et  $\delta(E_{\beta_{i_{2}}}) = h E_{\beta'_{i_{1}}}$.
On déduit alors de $(\star)$ que 
\[\begin{array}{rl}
k[E_{\beta_{ex}},E_{\alpha_{i}}]_{q} & = k \delta(E_{\beta_{ex}}) = \delta(E_{\beta_{i_{2}}}E_{\beta'_{i_{2}}})-q^{(\beta_{i_{2}},\beta'_{i_{2}})}\delta(E_{\beta'_{i_{2}}}E_{\beta_{i_{2}}}) \\
& = E_{\beta_{i_{2}}}\delta(E_{\beta'_{i_{2}}}) + \delta(E_{\beta_{i_{2}}})\sigma(E_{\beta'_{i_{2}}}) - q^{(\beta_{i_{2}},\beta'_{i_{2}})}(E_{\beta'_{i_{2}}}\delta(E_{\beta_{i_{2}}}) + \delta(E_{\beta'_{i_{2}}})\sigma(E_{\beta_{i_{2}}}))\\
& = h [q^{(\beta'_{i_{2}},\alpha_{i})}E_{\beta'_{i_{1}}}E_{\beta'_{i_{2}}} - q^{(\beta_{i_{2}},\beta'_{i_{2}})}E_{\beta'_{i_{2}}}E_{\beta'_{i_{1}}}]
\end{array}
\]
Comme $\beta'_{i_{2}}$ et $\beta'_{i_{1}}$ sont dans une même boîte, on sait (corollaire \ref{LSE}) que $E_{\beta'_{i_{1}}}E_{\beta'_{i_{2}}} = E_{\beta'_{i_{2}}}E_{\beta'_{i_{1}}}$, de sorte que :
\[k[E_{\beta_{ex}},E_{\alpha_{i}}]_{q} = h (q^{(\beta'_{i_{2}},\alpha_{i})} - q^{(\beta_{i_{2}},\beta'_{i_{2}})}) E_{\beta'_{i_{1}}}E_{\beta'_{i_{2}}}\]
Puisque $\beta_{i_{2}} + \beta'_{i_{2}} = \beta_{ex}$,  $P =$ Vect($\beta_{i_{2}},\beta'_{i_{2}}$) est un plan admissible de type (1.1) ou (1.2) (remarque \ref{systracine}) avec $\{\beta_{i_{2}}, \beta'_{i_{2}}\} = \{\beta, \beta' \}$, de sorte que $(\beta_{i_{2}},\beta'_{i_{2}}) \leq 0$.  Comme on a supposé $(\alpha_{i},\beta'_{i_{2}} ) > 0$, on a $(\beta_{i_{2}},\beta'_{i_{2}}) \neq (\alpha_{i},\beta'_{i_{2}} )$ de sorte que le coefficient ci-dessus de $E_{\beta'_{i_{1}}}E_{\beta'_{i_{2}}}$ est non nul. De là :
$$[E_{\beta_{ex}},E_{\alpha_{i}}]_{q} = \lambda E_{\beta'_{i_{1}}}E_{\beta'_{i_{2}}} \ \text{avec} \ \lambda \neq 0.$$ 
\end{preuve}

\subsection{Contraintes provenant d'un plan admissible}
On définit ici la notion de contrainte provenant d'un plan admissible P et on va vérifier que les diagrammes de Cauchon vérifient toute les contraintes provenant des plans admissibles. Rappelons d'abord les notations utilisées.
\begin{nota}
\begin{itemizedot}
\item On note $C_1, ..., C_n$ les sous ensembles de $\Phi^+$ représentant les colonnes.
\item Dans la suite, on considère un diagramme $\Delta$, c'est à dire une partie de $\Phi^+$
\item On considère un entier $j \in \llbracket 1, n \rrbracket$, on note   $\Delta_{j} := \Delta \cap C_j = \{ \beta_{u_{1}}, ..., \beta_{u_{l}} \} \subset C_{j} = \{\beta_{k},...,\beta_{r}\}$. On note ensuite, lorsqu'elles existent, $\beta_{ex}$ la racine exceptionnelle et $B_{ex} := \{\beta_{ex}\}$ la boîte qui la contient. S'il en est ainsi, on note $B_1$ la boîte de $C_{j}$ qui précède $B_{ex}$ et, $B'_1$ celle qui suit $B_{ex}$ ; de sorte que $s_{D}(B_1) = B_1'$.
\end{itemizedot}
\end{nota}
Dans les propositions \ref{contrainteelementaire}, \ref{contrainteelementaireavantex} et  \ref{contrainteelementaireapresex} , on a établi l'existence de contraintes en utilisant des plans admissibles. On va formaliser ce fait dans la définition de "contrainte provenant d'un plan admissible" suivante :
\begin{dfn}\label{contrainteplanadm}
Soit $\beta$ une racine de la colonne $C_j$ telle que $h'(\beta) = l$ et P un plan admissible (cf G. Lusztig \cite{Lu2} et \cite{Lu3}) de l'un des types suivants avec $\Phi_P^+ = \Phi^+ \cap P$ :
\begin{enumerate}
\item $\Phi_P^+ = \{\beta, \beta + \alpha_i, \alpha_i\}$ avec $i < j$. \hyperref[systracine]{type (2.1)}
\item $\Phi_P^+ = \{\beta, \beta + \alpha_i, \beta + 2\alpha_i, \alpha_i\}$ avec $i < j$. \hyperref[systracine]{type (2.3)}
\item $\Phi_P^+ = \{\beta, \beta + \beta' = \beta_{ex}, \beta'\}$, $\beta' \in C_j$ et $h(\beta') = h(\beta) +1$. \hyperref[systracine]{type (1.1)}
\item $\Phi_P^+ = \{\alpha_i , \alpha_i + \beta, \alpha_i + 2\beta = \beta_{ex}, \beta\}$ avec $i < j$, $h'(\alpha_i + 2\beta) = \frac{2l +1}{2}$ et $h(\beta) = l$. \hyperref[systracine]{type (1.2)} ou \hyperref[systracine]{type (2.2)} 
\item $\Phi_P^+ = \{\beta = \beta_{ex} , \alpha_i\}$ avec $i < j$, $\alpha_i \ \bot \ \beta$ et il existe $\beta_1$ et $\beta_2$ dans $C_j$ telles que $\beta + \alpha_i = \beta_1 + \beta_2$. \hyperref[systracine]{type (2.4)} \\
\item $\Phi_P^+ = \Phi^+ = \{\beta_1, ..., \beta_6\}$ est la partie positive d'un système de racines de type $G_2$ (notations de la proposition \ref{propG2}).
\end{enumerate}
Dans chacun des cas ci-dessus, on appelle \textit{contrainte(s) provenant de P}, la ou les contraintes suivantes :
\begin{enumerate}
\item $\beta \rightarrow \beta + \alpha_i$.
\item $\beta \rightarrow \beta + \alpha_i$ et $\beta + \alpha_i \rightarrow \beta + 2\alpha_i$.
\item $\beta \rightarrow \beta + \beta'$.
\item $\beta \rightarrow \alpha_i + \beta $, $\beta \rightarrow \alpha_i + 2\beta$ et $\alpha_i + 2\beta \rightarrow \alpha_i +\beta $.
\item \scalebox{.8}{\begin{psmatrix}[rowsep=.2cm, colsep=1cm]
& $\beta_{1}$\\
$\beta_{ex}$ & \\
& $\beta_{2}$
\end{psmatrix}
\psset{nodesep=2pt, linestyle=dashed, arrows=->}
\ncline{2,1}{1,2}
\ncline{2,1}{3,2}}

\item $\beta_6 \rightarrow \beta_5, \beta_5 \rightarrow \beta_4, \beta_5 \rightarrow \beta_3, \beta_4 \rightarrow \beta_3, \beta_3 \rightarrow \beta_2$.
\end{enumerate}
\end{dfn}
\begin{lemme}\label{contraintedeP}
Soit $\beta \in C_j$.
\begin{enumerate}
\item Si $\beta$ appartient à la boîte qui suit $\{\beta_{ex}\}$, alors $\beta \rightarrow \beta_{ex}$ est une contrainte provenant d'un plan admissible.
\item S'il existe $i < j$ tel que $\gamma = \beta + \alpha_i \in \Phi^+$ alors $\beta \rightarrow \gamma$ est une contrainte provenant d'un plan admissible.
\end{enumerate}
\end{lemme}
\begin{preuve}
Le résultat est vrai dans le cas où $\Phi$ est de type $G_2$. Dans la suite de la démonstration on considère donc un système de racines qui n'est pas de type $G_2$.
\begin{enumerate}
\item Soit $P = <\beta, \beta_{ex}>$. C'est un plan admissible de type 3 ou 4 dans la définition précèdente et dans chaque cas de figure, $\beta \rightarrow \beta_{ex}$ est une contrainte provenant de P.
\item Soit $P = <\beta, \alpha_i>$. C'est un plan admissible de type 1,2 ou 4 dans la définition précédente et dans tous les cas $\beta \rightarrow \gamma$ est une contrainte provenant de P.
\end{enumerate}
\end{preuve}
\begin{prop}
Soit $\Delta$ un diagramme de Cauchon. $\Delta$ vérifie toutes les contraintes (au sens de la section 4.1.1) provenant des plans admissibles.
\end{prop}
\begin{preuve}
Soit $\beta \in \Delta$ et soit P un plan admissible contenant $\beta$. Notons, comme dans la définition \ref{contrainteplanadm}, $\Phi_P^+ = \Phi^+ \cap P$.
\begin{enumerate}
\item Si $\Phi_P^+ = \{\beta, \beta + \alpha_i, \alpha_i\}$ avec $i < j$, alors, d'après la proposition \ref{contrainteelementaire}, $\Delta$ vérifie la contrainte $\beta \rightarrow \beta + \alpha_i$.
\item Si $\Phi_P^+ = \{\beta, \beta + \alpha_i, \beta + 2\alpha_i, \alpha_i\}$ avec $i < j$,  alors, d'après la proposition \ref{contrainteelementaire} appliquée à $\beta$ et $\beta + \alpha_i$, $\Delta$ vérifie les contraintes $\beta \rightarrow \beta + \alpha_i$ et $\beta + \alpha_i \rightarrow \beta + 2\alpha_i$.
\item Si $\Phi_P^+ = \{\beta, \beta + \beta' = \beta_{ex}, \beta'\}$ avec $i < j $, $\beta' \in C_j$ et $h(\beta') = h(\beta) +1$, alors, d'après la proposition \ref{contrainteelementaireavantex}, $\Delta$ vérifie la contrainte $\beta \rightarrow \beta + \beta'$.
\item Si $\Phi_P^+ = \{\alpha_i , \alpha_i + \beta, \alpha_i + 2\beta = \beta_{ex}, \beta\}$ avec $i < j$ et $h'(\alpha_i + 2\beta) = \frac{2l +1}{2}$, alors, d'après les propositions \ref{contrainteelementaire}, \ref{contrainteelementaireavantex} et \ref{contrainteelementaireapresex}, $\Delta$ vérifie respectivement les contraintes $\beta \rightarrow \alpha_i + \beta $, $\beta \rightarrow \alpha_i + 2\beta$ et $\alpha_i + 2\beta \rightarrow \alpha_i +\beta $.
\item Si $\Phi_P^+ = \{\beta = \beta_{ex} , \alpha_i\}$ avec $i < j$, $\alpha_i \ \bot \ \beta$ et s'il existe $\beta_1$ et $\beta_2$ dans $C_j$ tels que $\beta + \alpha_i = \beta_1 + \beta_2$, alors, d'après la proposition \ref{contrainteelementaireapresex}, $\Delta$ vérifie la contrainte \scalebox{.8}{\begin{psmatrix}[rowsep=.2cm, colsep=1cm]
& $\beta_{1}$\\
$\beta$ & \\
& $\beta_{2}$
\end{psmatrix}
\psset{nodesep=2pt, linestyle=dashed, arrows=->}
\ncline{2,1}{1,2}
\ncline{2,1}{3,2}}.  \\ 
\item Si $\Phi_P^+ = \Phi^+$ est de type $G_2$, la proposition \ref{contrainteelementaire} implique que $\Delta$ vérifie les contraintes $\beta_6 \rightarrow \beta_5, \beta_5 \rightarrow \beta_3, \beta_4 \rightarrow \beta_3, \beta_3 \rightarrow \beta_2$. La proposition \ref{contrainteelementaireavantex} implique quant à elle que $\Delta$ vérifie la contrainte $\beta_5 \rightarrow \beta_4$.
\end{enumerate}
\end{preuve}

\subsection{La réciproque}
Le but de cette partie est de démontrer la réciproque du résultat démontré dans la partie précédente, c'est à dire :
\begin{theo}\label{reciproque contrainte}
Si $\Delta$ est un diagramme qui vérifie toutes les contraintes provenant des plans admissibles, alors $\Delta$ est un diagramme de Cauchon.
\end{theo}
Soit $\beta \in \Phi^+$, une racine positive de la colonne $C_j$ avec $1 < j < n$. On note $B_0$ la boîte de $\beta$, $B_1$ la boîte qui précède $B_0$ (si elle existe) et $B_2$ celle qui précède $B_1$ (si elle existe) dans la colonne $C_j$.\\
On pose $\Phi_\beta^+ = \{ \alpha_i \ | \ i < j\} \cup \{\gamma < \beta \ | \ \gamma \in B_0 \} \cup B_1 \cup (B_2$ si $B_1 = \{\beta_{ex}\})$. Si $\gamma \in \Phi^+$, il existe $k \in \llbracket 1, N \rrbracket$ tel que $\gamma = \beta_k$ et on rappelle (cf section 4.1) que $X_\gamma = X_k$.\\
On pose enfin $D_\beta = \mathbb{K}<X_\gamma \ | \ \gamma < \beta >$.
\begin{lemme} \label{Dbeta}
\[D_\beta = \mathbb{K}<X_\gamma \ | \ \gamma \in \Phi_\beta^+>\]
\end{lemme}
\begin{preuve}
On pose $D'_\beta = \mathbb{K}<X_\gamma \ | \ \gamma \in \Phi_\beta^+> \ \subset \ D_\beta $. Commen\c cons par montrer que pour i < j, $\{X_{\gamma}, \gamma \in C_{i}\} \subset D'_\beta$. Si $\Phi$ est de type $G_2$, $\{X_{\gamma}, \gamma \in C_{i}\}$ est soit vide, soit réduit à $X_{\alpha_1} \in D'_\beta$. Pour $\Phi \neq G_2$, on démontre ce résultat par récurrence sur $h(\gamma)$. 
\begin{description}
\item[Si $h(\gamma) = 1$ :] alors $\gamma = \alpha_{i}$ et $X_{\gamma} \in  D'_\beta$ par définition de $ \Phi_\beta^+$.
\item[Si $h(\gamma) > 1$ et $\gamma$ non exceptionnelle :] Par la proposition \ref{prop2.10}, il existe $l < i$ tel que $\gamma' = \gamma - \alpha_l \in \Phi^+$, de sorte que, par le corollaire \ref{corplanadm} et la proposition \ref{relcomX}, on a $X_\gamma \in \mathbb{K}<X_{\gamma'}, X_{\alpha_l}> \subset D'_\beta$ (hypothèse de récurrence).
\item[Si $h(\gamma) > 1$ et $\gamma$ exceptionnelle :] On sait (proposition \ref{prop2.8}) que dans ce cas, il existe deux racines non exceptionnelles de $C_{i}$, $\eta_1$ et $\eta_2$,  telles que $\eta_1+ \eta_2 = \gamma$ et $h(\eta_2) = h(\eta_1) + 1$. On en déduit par le corollaire \ref{corplanadm} et la proposition \ref{relcomX} que $X_\gamma \in \mathbb{K}<X_{\eta_{1}}, X_{\eta_{2}}> \subset D'_\beta$ ($X_{\eta_{1}}$ et $X_{\eta_{2}}$ sont dans $D'_\beta$ car $\eta_1$ et $\eta_2$ sont non exceptionnelle).
\end{description}
Reste à montrer que $\{X_{\gamma} | \gamma \in C_{j}, \gamma < \beta \} \subset D'_\beta$.\\
Si $h(\gamma) = h(\beta)$ avec $\gamma < \beta$, alors $\gamma \in \Phi_\beta^+$. Donc $X_{\gamma} \in D'_\beta$.\\
On raisonne à nouveau par récurrence (sur $h(B)$) pour montrer que pour toute boîte non exceptionnelle $B$ de $C_{j}$ telle que $B < B_0$ (i.e. toute racine $\beta$ de $B$ est inférieure strictement à toutes les racines de $B_0$), on a $\{X_{\gamma} | \gamma \in B\} \subset D'_{\beta}$.
\begin{itemizedot}
\item \textbf{Supposons $B_1$ non exceptionnelle.}
\begin{description}
\item[Initialisation :] Le résultat est vrai pour la boîte $B_{1}$ puisque $B_1 \subset \Phi_\beta^+$.
\item[Hérédité :] Soient B une boîte non exceptionnelle de $C_j$ telle que $h(B) > h(B_{1})$ et $\gamma \in B$. Par la proposition \ref{prop2.10}, il existe $\alpha_{l} \in \Pi \ (l < j)$ tel que $\gamma - \alpha_{l} \in \Phi^+$. $\gamma' := \gamma - \alpha_{l}$ est dans une boîte $B'$ non exceptionnelle de $C_j$ telle que $h(B) = h(B') + 1 > h(B') \geq h(B_{1}) > h(B_0)$ et on a $X_{\gamma'} \in D'_{\beta}$ par l'hypothèse de récurrence. \\
Si $\Phi \neq G_2$, on déduit du corollaire \ref{corplanadm} et de la proposition \ref{relcomX} que $[X_{\gamma'},X_{\alpha_l}]_q = k X_\gamma$ avec $k \in \mathbb{K}^\star$. Comme $X_{\alpha_{l}} \in D'_{\beta}$, on en déduit que $X_{\gamma} \in D'_{\beta}$.\\
Si $\Phi = G_2$, on déduit des propositions \ref{propG2} et \ref{relcomX} que $[X_{\gamma'},X_{\alpha_l}]_q = k X_\gamma$ avec $k \in \mathbb{K}^\star$. Comme $X_{\alpha_{l}} \in D'_{\beta}$, on en déduit que $X_{\gamma} \in D'_{\beta}$.\\
\end{description}
\item \textbf{Supposons $B_1$ exceptionnelle.}
\begin{description}
\item[Initialisation :] Le résultat est vrai pour la boîte $B_{2}$ puisque, dans ce cas, $B_2 \subset \Phi_\beta^+$.
\item[Hérédité :] On fait le même raisonnement que ci-dessus en rempla\c cant $B_{1}$ par $B_{2}$.
\end{description}
\end{itemizedot}
Reste à montrer que si $B = \{\beta_{ex}\}$ est une boîte exceptionnelle de $C_{j}$ telle  $B < B_{0}$, on a $X_{\beta_{ex}} \in D'_{\beta}$.\\
Si $B = B_1$, on a $B \subset \Phi_\beta^+$, d'où le résultat.\\
Supposons donc $B < B_1$. Comme ci-dessus, on a $\beta_{ex} = \eta_{1} + \eta_{2}$ avec $\eta_{1}$ et $\eta_{2}$ deux racines non exceptionnelle de $C_j$ telles que $h(\eta_{2}) = h(\eta_1) +1$. Les boîtes de $\eta_1$ et $\eta_2$ sont non exceptionnelles, de part et d'autre de $B$, donc inférieures ou égales à $B_1$, donc inférieures strictements à $B_0$. Le résultat étant vrai pour les boîtes non exceptionnelles, $X_{\eta_{1}} \in D'_{\beta}$ et $X_{\eta_{2}} \in D'_{\beta}$.\\
Si $\Phi \neq G_2$, on déduit (comme ci-dessus) du corollaire \ref{corplanadm} et de la proposition \ref{relcomX} que $X_{\beta_{ex}} \in D'_{\beta}$.\\
Si $\Phi = G_2$, on déduit (comme ci-dessus) des propositions \ref{propG2} et \ref{relcomX} que $X_{\beta_{ex}} \in D'_{\beta}$.\\

On a donc bien $D_\beta = D'_\beta$.
\end{preuve}
On rappelle que $A = \mathcal{U}_q^+(\mathfrak{g}) = \mathbb{K}<X_{\beta_i} | i \in \llbracket 1 , N \rrbracket > := \mathbb{K}<X_i | i \in \llbracket 1 , N \rrbracket >$. Soient $\beta_{r}$ et $\beta_{r+1}$ ($1 \leq r \leq N-1$) deux racines consécutives de $\Phi^+$ ($\beta_{r} < \beta_{r+1}$). On rappelle que $A^{(r+1)} = \mathbb{K}<X_{i}^{(r+1)}>$ et $A^{(r)} = \mathbb{K}<X_{i}^{(r)}>$ (1 < r < N) désignent les algèbres déduites de $A$ par l'algorithme d'effacement des dérivations de la section 4.
\begin{lemme} \label{relations}
Soit  $\beta_r \in \Phi^+$, une racine positive de la colonne $C_j$ et $D_{\beta_r}^{(r+1)} = \mathbb{K}<X^{(r+1)}_\gamma \ | \ \gamma < \beta_r >$.
Alors \[D^{(r+1)}_{\beta_r} = \mathbb{K}<X^{(r+1)}_\gamma \ | \ \gamma \in \Phi_{\beta_r}^+>\]
\end{lemme}
\begin{preuve}
Par le lemme \ref{4.1.1}, les formules de redressement entre les $X^{(r+1)}_\gamma$ lorsque $\gamma \leq \beta_r$ sont les mêmes que les formules de redressement entre les $X_\gamma$ avec $\gamma \leq \beta_r$. Il suffit donc de réécrire la démonstration du lemme \ref{Dbeta}, en remplaçant chaque $X_\gamma$ par $X^{(r+1)}_\gamma$.
\end{preuve}
\textit{Notons, comme dans le paragraphe 4, $\varphi : $ Spec $ A \hookrightarrow$ Spec $(\overline{A})  (\overline{A} = A^{(2)})$ l'injection canonique, c'est à dire la composée des injections canoniques $\varphi_{r} : $ Spec $(A^{(r+1)}) \hookrightarrow$ Spec $(A^{(r)})$ pour $r \in \llbracket 2, N \rrbracket$, et rappelons qu'un sous-ensemble $\Delta$ de $\Phi^+$ est un diagramme de Cauchon si et seulement si $(\exists P \in \text{Spec}(A))(\varphi (P) = <T_\gamma | \gamma \in \Delta>$).}\\ 
\begin{preuve4}
Soit $\Delta \subset \Phi^+$ un diagramme vérifiant les contraintes provenant des plans admissibles. Posons $Q:= <T_\gamma | \gamma \in \Delta>$. Par \cite[Section 5.5.]{MR1967309}, c'est un idéal $H^{(2)}$-premier, donc complètement premier, de $A^{(2)} = \overline{A}$ et si $\beta \in \Phi^+ \setminus \Delta$, $T_\beta$ est régulier modulo Q. Par suite, $Q \cap \Phi^+ = \{ T_\gamma | \gamma \in \Delta \}$.\\
Montrons par récurrence, que pour chaque $r \in \llbracket 2, N+1 \rrbracket $, il existe $P^{(r)} \in $ Spec $(A^{(r)})$ tel que $Q = \varphi_{2} \circ ... \circ \varphi_{r-1} (P^{(r)}).$\\
\textbf{Initialisation :} Si $r= 2$, dans ce cas, on a $\varphi_{2} \circ ... \circ \varphi_{r-1} = Id_{\text{Spec}(\overline{A})}$ et $P^{(2)} = Q$ convient.\\ 
\textbf{Considérons un entier $r \in \llbracket 2, N \rrbracket $, supposons qu'il existe $P^{(r)} \in \text{Spec} (A^{(r)})$ tel que $\varphi_{2} \circ ... \circ \varphi_{r-1}(P^{(r)}) = Q$ et montrons qu'il existe $P^{(r+1)} \in \text{Spec} (A^{(r+1)})$ tel que $\varphi_{r} (P^{(r+1)}) = P^{(r)}$. Il en résultera que $\varphi_{2} \circ ... \circ \varphi_{r}(P^{(r+1)}) = Q$.}\\
$ $\\
\begin{itemizedot}
\item Si $X_{r}^{(r)} \notin P^{(r)}$, cela résulte du premier point de la proposition \ref{critereimage}.
\item Supposons donc que $X_{r}^{(r)} \in P^{(r)}$. D'après le second point de la proposition \ref{critereimage}, tout revient à montrer que 
$$\Theta^{(r)}\left(\delta_{r}^{(r+1)}(X_{i}^{(r+1)})\right) \in P\it^{(r)} \ \ \text{pour} \ 1  \leq i \leq  r-1.$$
\end{itemizedot}
\begin{obs}
Il suffit de démontrer que $\Theta^{(r)}\left(\delta_{r}^{(r+1)}(X_{i}^{(r+1)})\right) $  $\in$ $ $ $P\it^{(r)}$ $ $ pour chaque $i \in \llbracket 1, r-1 \rrbracket$ tel que $\beta_i \in \Phi_{\beta_r}^+$.
\end{obs}
\begin{preuve5}
Soit $i \in \llbracket 1, r-1 \rrbracket$. D'après le lemme \ref{relations}, on a :
\[ X_i^{(r+1)} = \sum_{j_1,...,j_s \in \Gamma} m_{j_1, ..., j_s} X_{j_1}^{(r+1)}...X_{j_s}^{(r+1)} \ \ \text{où} \ \ \Gamma := \{ j \in \llbracket 1, r-1 \rrbracket | \beta_j \in \Phi_{\beta_r}^+ \} \ \text{et} \ m_{j_1,...,j_s} \in \mathbb{K}.\]
On a alors :
\[\begin{array}{rl}
\delta_{r}^{(r+1)}(X_{i}^{(r+1)}) 	&= \sum  m_{j_1, ..., j_s} \delta_{r}^{(r+1)}(X_{j_1}^{(r+1)}...X_{j_s}^{(r+1)})\\
									&= \sum m_{j_1, ..., j_s} \left[\delta_{r}^{(r+1)}(X_{j_1}^{(r+1)})X_{j_2}^{(r+1)}...X_{j_s}^{(r+1)} + \sigma_{r}^{(r+1)}(X_{j_1}^{(r+1)})\delta_{r}^{(r+1)}(X_{j_2}^{(r+1)})...X_{j_s}^{(r+1)} + \right. \\
									&\left. ... + \sigma_{r}^{(r+1)}(X_{j_1}^{(r+1)}...X_{j_{s-1}}^{(r+1)})\delta_{r}^{(r+1)}(X_{j_s}^{(r+1)}) \right]\\
									&= \sum m_{j_1, ..., j_s} \left[  \delta_{r}^{(r+1)}(X_{j_1}^{(r+1)})X_{j_2}^{(r+1)}...X_{j_s}^{(r+1)} + \lambda_{r,j_1}X_{j_1}^{(r+1)}\delta_{r}^{(r+1)}(X_{j_2}^{(r+1)})...X_{j_s}^{(r+1)} + \right. \\
									& \left. ... + \lambda_{r,j_1}...\lambda_{r,j_{s-1}}X_{j_1}^{(r+1)}...X_{j_{s-1}}^{(r+1)}\delta_{r}^{(r+1)}(X_{j_s}^{(r+1)})\right].
\end{array}\]
Par suite 
\[\begin{array}{rl}
\Theta^{(r)}\left(\delta_{r}^{(r+1)}(X_{i}^{(r+1)})\right) =
									& \sum m_{j_1, ..., j_s} \left[  \Theta^{(r)}\left(\delta_{r}^{(r+1)}(X_{j_1}^{(r+1)})\right)X_{j_2}^{(r)}...X_{j_s}^{(r)} \right.\\
									& + \lambda_{r,j_1}X_{j_1}^{(r)}\Theta^{(r)}\left(\delta_{r}^{(r+1)}(X_{j_2}^{(r+1)})\right)...X_{j_s}^{(r)} +  ... \\
									& \left. + \lambda_{r,j_1}...\lambda_{r,j_{s-1}}X_{j_1}^{(r)}...X_{j_{s-1}}^{(r)}\Theta^{(r)}\left(\delta_{r}^{(r+1)}(X_{j_s}^{(r+1)})\right)\right].
\end{array}\]
Comme chaque $\Theta^{(r)}\left(\delta_{r}^{(r+1)}(X_{j_l}^{(r+1)})\right) \in P^{(r)}$ par hypothèse, on a bien $\Theta^{(r)}\left(\delta_{r}^{(r+1)}(X_{i}^{(r+1)})\right) \in P^{(r)}$.
\end{preuve5}
\begin{tabular}{lc}
\begin{minipage}{.8\linewidth}
\bf{Retour à la démonstration du théorème \ref{reciproque contrainte} :} \rm \\
Pour chaque $s \in \llbracket 2, r-1 \rrbracket$, on notera $P^{(s)} = \varphi_s \circ ... \varphi_{r-1}(P^{(r)})$.
\begin{obs}
$\beta_{r} \in \Delta$. 
\end{obs}
En effet, comme $X_{r}^{(r)} \in P^{(r)}$, on déduit du lemme \ref{lemme4.1} que l'on a, successivement  $X_{r}^{(r-1)} \in P^{(r-1)}$, ..., $X_{r}^{(2)}\in P^{(2)} = Q$. Comme $Q \cap \{T_{\beta} | \beta \in \Phi^+\} = \{T_\gamma | \gamma \in \Delta \}$, on a $T_{\beta_r} (= X_{r}^{(2)}) \in Q \Rightarrow \beta_{r} \in \Delta$.\\

On rappelle que, si $\beta_r \in C_j$, $\Phi_{\beta_r}^+ = \{ \alpha_i \ | \ i < j\} \cup \{\gamma < \beta_r \ | \ \gamma \in B_0\} \cup B_1 \cup (B_2$ si $B_1 = \{\beta_{ex}\})$ ($B_0$ est la boîte qui contient $\beta_r$, $B_1$ est la boîte qui précède $B_0$ dans $C_j$ quand elle existe et $B_2$ la boîte qui précède $B_1$ dans $C_j$ lorsqu'elle existe).
\end{minipage} 
&
\begin{minipage}{.2\linewidth}
\hspace{.5cm} $C_j$\\
\begin{tabular}{l|c|c}
\cline{2-2}
&&\\
\cline{2-2}
&& \\
&$\vdots$ & \\
\cline{2-2}\cline{2-2}\cline{2-2}\cline{2-2}\cline{2-2}\cline{2-2}
 && \\
 \cline{2-2}
& & \\
 \cline{2-2}
 && \multirow{-3}{2cm}{ $\left. \begin{minipage}{0.1cm} $ $\\ $ $ \\ $ $ \\ \end{minipage} \right\} B_2$} \\
 \cline{2-2}\cline{2-2}\cline{2-2}\cline{2-2}\cline{2-2}\cline{2-2}\cline{2-2}
& & \\
 \cline{2-2}
 && \\
 \cline{2-2}
 && \multirow{-3}{2cm}{ $\left. \begin{minipage}{0.1cm} $ $\\ $ $ \\ $ $ \\ \end{minipage} \right\} B_1$} \\
 \cline{2-2}\cline{2-2}\cline{2-2}\cline{2-2}\cline{2-2}\cline{2-2}\cline{2-2}
 && \\
 \cline{2-2}
 &$\beta_{r}$& \\
 \cline{2-2}
 && \multirow{-3}{2cm}{ $\left. \begin{minipage}{0.1cm} $ $\\ $ $ \\ $ $ \\ \end{minipage} \right\} B_0$} \\
 \cline{2-2}\cline{2-2}\cline{2-2}\cline{2-2}\cline{2-2}\cline{2-2}\cline{2-2}
 &$\vdots $ \\
 \cline{2-2}
& $\alpha_j$ & \\
 \cline{2-2}
\end{tabular}
\end{minipage}
\end{tabular}\\
Soit $i \in \llbracket 1, r-1 \rrbracket$ tel que $\beta_i \in \Phi_{\beta_r}^+$.\\
\begin{itemizedot}
\item Si $\beta_i \in B_0 \cup B_1$, il résulte du théorème \ref{1.4} (et de la remarque qui suit) que $\delta_{r}^{(r+1)}(X_{i}^{(r+1)}) = 0$. Par suite $\Theta^{(r)}\left(\delta_{r}^{(r+1)}(X_{i}^{(r+1)})\right) = 0 \in P^{(r)}$.
\item Supposons que $B_1 = \{\beta_{ex} \}$ avec $\beta_{ex}= \beta_e$ ($e < r$), et que $\beta_i \in B_2$.\\ 
D'après le théorème \ref{1.4}, $\delta_{r}^{(r+1)}(X_{i}^{(r+1)}) = P_{r,i}^{(r+1)}$ est homogène de poids $\beta_r + \beta_i$ et les variables $X_l^{(r+1)}$ apparaissant dans $P_{r,i}^{(r+1)}$ sont telles que $\beta_l \in B_1 = \{\beta_{e}\}$. Par suite $P_{r,i}^{(r+1)}$ est nul ou de la forme $\lambda X_{e}^m$ avec $\lambda \in \mathbb{K}^{\star}$ et $m\beta_{ex} = \beta_r + \beta_i$, de sorte que (en identifiant les coefficients de $\alpha_j$) on a $m = 1$. \\
Si $P_{r,i}^{(r+1)} = 0$, on a $\Theta^{(r)}\left(\delta_{r}^{(r+1)}(X_{i}^{(r+1)})\right) = 0 \in P^{(r)}$.\\
Sinon, supposons $P_{r,i}^{(r+1)} = \lambda X_{e}^m$. Comme $\Delta$ vérifie les contraintes provenant des plans admissibles, il résulte du lemme \ref{contraintedeP} que $\Delta$ vérifie la contrainte $\beta_r \rightarrow \beta_{ex}$ et, puisque $\beta_r \in \Delta$, on a $\beta_{ex} \in \Delta$.\\
On a donc $X_e^{(2)} \in Q = P^{(2)}$ et par le lemme \ref{lemme4.1}, $X_e^{(e+1)} \in P^{(e+1)}$. Comme $\beta_{e}$ et $\beta_r$ sont dans des boîtes adjacentes par construction, on déduit du lemme \ref{5.2} que 
\[X_e^{(e+1)} \in P^{(e+1)} \Rightarrow X_e^{(r)} \in P^{(r)}.\]
Par suite, $\Theta^{(r)}\left(\delta_{r}^{(r+1)}(X_{i}^{(r+1)})\right) = \Theta^{(r)}(\lambda  X_e^{(r+1)}) = \lambda X_e^{(r)} \in P^{(r)}$.
\item Considèrons à présent le cas $\beta_i = \alpha_k$ avec $k<j$.\\ 
Si $\delta_{r}^{(r+1)}(X_{i}^{(r+1)}) = 0$, on a $\Theta^{(r)}\left(\delta_{r}^{(r+1)}(X_{i}^{(r+1)})\right) = 0 \in P^{(r)}$.\\
Supposons $\delta_{r}^{(r+1)}(X_{i}^{(r+1)}) \neq 0$. D'après le théorème \ref{1.4}, on a
\[\delta_{r}^{(r+1)}(X_{i}^{(r+1)}) = \sum_{i < j_1 \leq ... \leq j_s < r} c_{j_1,...,j_s} X_{j_1}^{(r+1)}...X_{j_s}^{(r+1)} \ \ (c_{j_1,...,j_s} \in \mathbb{K})\]
\[c_{j_1,...,j_s} \in \mathbb{K}^* \Rightarrow \left( \beta_{j_1} + ... + \beta_{j_s} = \beta_r + \alpha_k \ \text{et} \ \beta_{j_1},...,\beta_{j_s} \notin B_0 \right) \]
On a alors
\[\Theta^{(r)}\left(\delta_{r}^{(r+1)}(X_{i}^{(r+1)})\right) = \sum_{i < j_1 \leq ... \leq j_s < r} c_{j_1,...,j_s} X_{j_1}^{(r)}...X_{j_s}^{(r)}\]
et tout revient à démontrer que, si $c_{j_1,...,j_s} \in \mathbb{K}^*$, on a $X_{j_1}^{(r)}...X_{j_s}^{(r)} \in P^{(r)}$.\\
\par 
Supposons donc $c_{j_1,...,j_s} \neq 0 $. En considérant le coefficient de $\alpha_j$ dans l'égalité 
\begin{equation}\label{equaracine}
\beta_{j_1} + ... + \beta_{j_s} = \beta_r + \alpha_k
\end{equation}
on voit que $\beta_{j_s} \in C_j$. Comme $\beta_{j_s} \notin B_0$ et $j_s < r$, la boîte $B_1$ existe.\\
Distinguons plusieurs cas.
\begin{itemizedot}
\item $B_0$ et $B_1$ sont ordinaires. Comme $j_s < r$ et $\beta_{j_s} \notin B_0$, on a $h(\beta_r) < h(\beta_{j_s})$. Par (\ref{equaracine}), on a aussi $h(\beta_{j_s}) \leq h(\beta_r + \alpha_k) = h(\beta_r) + 1$. Il en résulte que $s=1$ et que $\beta_{j_s} \in B_1$. On dispose alors (lemme \ref{contraintedeP}) de la contrainte $\beta_r \rightarrow \beta_{j_s}$. Comme $\beta_r \in \Delta$, on a $\beta_{j_s} \in \Delta$ et, comme ci-dessus, $X_{j_s}^{(j_s+1)} \in P^{(j_s+1)}$ puis, de là, $X_{j_s}^{(r)} \in P^{(r)}$. Le monôme considéré est donc bien dans $P^{(r)}$.\\
\item $B_0$ est ordinaire et $B_1$ est exceptionnelle de sorte que $B_2$ existe. Comme dans le cas précédent, on vérifie que $s =1$ et que $\beta_{j_s} \in B_2$. On dispose donc encore (lemme \ref{contraintedeP}) de la contrainte $\beta_r \rightarrow \beta_{j_s}$. On dispose aussi (lemme \ref{contraintedeP}) de la contrainte $\beta_r \rightarrow \beta_e$. Comme $\beta_r \in \Delta$, on a $\beta_e, \beta_{j_s} \in \Delta$, de sorte que $X_e^{(e+1)} \in P^{(e+1)}$ et $X_{j_s}^{(j_s+1)} \in P^{(j_s+1)}$. Par le second point du lemme \ref{5.2}, on en déduit que $X_{j_s}^{(r)} \in P^{(r)}$. Le monôme considéré est donc bien dans $P^{(r)}$.\\
\item $B_0$ est exceptionnelle. Puisque $\beta_{j_s} \notin B_0$, $\beta_{j_s}$ est ordinaire dans $C_j$. Par l'égalité (\ref{equaracine}), on a $s \geq 2$ et $\beta_{j_{s-1}}$ ordinaire dans $C_j$. Si on note $h(\beta_r) = 2l +1 $ ($l \geq 1$), on sait que $h(\beta_{j_{s-1}}) \geq l+1, h(\beta_{j_{s}}) \geq l+1$ et $h(\beta_r + \alpha_k) = 2l+2$. On en déduit que $s = 2$ et que $\beta_{j_{s-1}}, \beta_{j_s} \in B_1$. L'égalité (\ref{equaracine}) s'écrit alors $\beta_r + \alpha_k = \beta_{j_{s-1}} + \beta_{j_s}$.
\begin{itemizedot}
\item Supposons $\beta_{j_{s-1}} \neq \beta_{j_s}$, de sorte que $\beta_{j_{s-1}}$ et $\beta_{j_s}$ étant dans la même boîte $B_1$, elles sont orthogonales. Il en résulte que $\Phi$ n'est pas de type $G_2$ (dans le cas $G_2$, les boîtes sont réduites à un élément). Notons $P = <\beta_{j_s},\beta_{j_{s-1}}>$ le plan engendré par $\beta_{j_{s}}$, $\beta_{j_{s-1}}$, et supposons $\Phi_P^+ \neq \{\beta_{j_{s-1}},\beta_{j_s}\}$. Alors, puisque $\Phi_P \neq G_2$, $\Phi_P$ est de type $A_2$ ou $B_2$. Comme $\beta_{j_{s-1}}$ et $\beta_{j_s}$ sont orthogonales, $\Phi_P$ est de type $B_2$ et il existe $\beta \in \Phi^+$ telle que $\beta_r + \alpha_k = \beta_{j_{s-1}} + \beta_{j_s} = m \beta$ avec $m = 1$ ou $2$. \\
Si $m = 1$, alors $\beta$ et $\beta_r$ sont deux racines exceptionnelles distinctes de $C_j$, ce qui est impossible.\\
Supposons $m = 2$. On a alors $\beta_r + \alpha_k = \beta_{j_{s-1}} + \beta_{j_s} = 2 \beta$. On en déduit que $h(\beta) = l + 1$, de sorte que $\beta$ est aussi un élément de $B_1$ différent de $\beta_{j_{s-1}}$ et $\beta_{j_s} $. Il en résulte que $\beta, \beta_{j_{s-1}}, \beta_{j_s}$ sont deux à deux orthogonales, ce qui est contradictoire avec l'égalité $\beta_{j_{s-1}} + \beta_{j_s} = 2 \beta$.\\
On a donc $\Phi_P^+  =  \{\beta_{j_{s-1}},\beta_{j_s}\}$ et il résulte de la définition \ref{contrainteplanadm} que l'on dispose de la contrainte \scalebox{.8}{\begin{psmatrix}[rowsep=.2cm, colsep=1cm]
& $\beta_{j_{s-1}}$\\
$\beta_{r}$ & \\
& $\beta_{j_s}$
\end{psmatrix}
\psset{nodesep=2pt, linestyle=dashed, arrows=->}
\ncline{2,1}{1,2}
\ncline{2,1}{3,2}}. L'une des deux racines $\beta_{j_s}$ ou $\beta_{j_{s-1}}$ est alors dans $\Delta$. Si, par exemple, $\beta_{j_s} \in \Delta$, on a successivement, comme dans le premier cas, $X_{j_s}^{(j_s+1)} \in P^{(j_s+1)}$ et $X_{j_s}^{(r)} \in P^{(r)}$. Le monôme considéré est donc bien dans $P^{(r)}$.
\item Si $\beta_{j_{s-1}} = \beta_{j_s}$ alors l'égalité (\ref{equaracine}) s'écrit alors $\beta_r + \alpha_k = 2 \beta_{j_s}$. Si on pose $\beta = s_D(\beta_{j_s}) = \beta_r - \beta_{j_s} \in \Phi^+$ et si on retranche $\beta_{j_s}$ à chaque membre de l'égalité précèdente, on obtient $\beta + \alpha_k = \beta_{j_s}$. Notons $P = <\beta_r,\beta_{j_s}>$ le plan engendré par $\beta_r$ et $\beta_{j_s}$. \\
Supposons $\Phi$ de type $G_2$. Alors on a $\beta_r = \beta_4$, $\alpha_k = \beta_1$ et $\beta_{j_s} = \beta_3$. Par la définition \ref{contrainteplanadm}, on dispose donc de la contrainte $\beta_r \rightarrow \beta_{j_s}$.\\
Supposons que $\Phi$ n'est pas de type $G_2$. Il résulte de l'égalité $\beta_r + \alpha_k = 2 \beta_{j_s}$ que $\Phi_P$ est de type $B_2$, de sorte que $\Phi_P^+ = \{\alpha_k, \alpha_k + \beta = \beta_{j_s}, \alpha_k + 2\beta = \beta_r, \beta \}$ avec $h(\beta) = h(\beta_r) - h(\beta_{j_s}) = 2l+1 - (l+1) = l$. P est donc un plan admissible de type 4 au sens de la définition  \ref{contrainteplanadm}. On dispose donc à nouveau de la contrainte $\beta_r \rightarrow \beta_{j_s}$. \\
Ainsi, dans tous les cas, on a $\beta_{j_s} \in \Delta$. On a donc successivement, comme dans le premier cas, $X_{j_s}^{(j_s+1)} \in P^{(j_s+1)}$ et $X_{j_s}^{(r)} \in P^{(r)}$. Le monôme considéré est donc bien dans $P^{(r)}$, ce qui termine la démonstration.
\end{itemizedot}
  
\end{itemizedot}
\end{itemizedot}
\end{preuve4}

\section{Forme des diagrammes pour une décomposition de $w_{0}$ particulière.}
Dans ce paragraphe, on donne une description explicite des diagrammes de Cauchon pour une décomposition particulière de $w_{0}$ dans chaque type d'algèbre de Lie simple de dimension finie. 
\subsection{Le cas des familles infinies}
\subsubsection{Le cas $A_n$, $n \geq 1$}

\begin{conv}
On numérote les racines simples de manière que le diagramme de Dynkin soit : 
\begin{center}
	$ \alpha_1 - \alpha_2 - \cdots - \alpha_{n-1} - \alpha_n $.
\end{center}
On sait (voir par exemple \cite[Section 5]{MR1628449}) que l'écriture $ w_0 = s_{\alpha_1} \circ (s_{\alpha_2}\circ s_{\alpha_1})   \cdots \circ (s_{\alpha_n} \circ s_{\alpha_{n-1}} \circ \cdots \circ s_{\alpha_1} )$ est une décomposition réduite qui induit l'ordre ci-dessous sur les racines positives.

\[\begin{array}{c|c|c|c|c|}
& C_{1} & C_{2} & & C_{n}\\
\cline{2-5}
 &   \beta_1 = \alpha_1 &  
 \beta_2 = \alpha_1 + \alpha_2  & \cdots & 
\beta_{N-n+1} = \alpha_1 + \cdots + \alpha_{n-1} + \alpha_n  \\
  \cline{2-5}
\multicolumn{2}{c|}{} & \multicolumn{1}{|c|}{
 \beta_3 = \alpha_2} & \multicolumn{1}{|c|}{\cdots} & \multicolumn{1}{|c|}{\vdots} \\
 \cline{3-5}
\multicolumn{3}{c|}{} & \multicolumn{1}{|c|}{\ddots} & \multicolumn{1}{|c|}{\vdots} \\
\cline{4-5}
\multicolumn{4}{c|}{} & \multicolumn{1}{|c|}{
\beta_N = \alpha_n}  \\
\cline{5-5}
\end{array}\]
\end{conv}
\par
On constate que ceci est un ordre de Lusztig et qu'aucune des colonnes $C_{1}, ..., C_{n}$ n'est exceptionnelle.\\
\begin{lemme}\label{contraintesan}
Si $1 \leq l \leq n$ et si $\beta_{j}, \beta_{j+1}$ sont deux racines consécutives de $C_{l}$ alors $\beta_{j+1} \rightarrow \beta_{j} $ est une contrainte.
\end{lemme}
\begin{preuve}
Si $\beta_{j}$ et $\beta_{j+1}$ sont deux racines consécutives de $C_{l}$ alors il existe une racine simple $\alpha_i$ ($1 \leq i < l \leq n$) telle que $\beta_{j} = \beta_{j+1} + \alpha_{i}$. D'après la proposition \ref{contrainteelementaire}, $\beta \rightarrow \beta + \alpha_{i}$ est une contrainte.

\end{preuve}
On en déduit :
\begin{cor}\label{contraintecol}
Si $C_{l} = \{\beta_{s}, \beta_{s+1},...,\beta_{r} = \alpha_{l}\}$ est la colonne $l$ avec $1 \leq l \leq n$, les contraintes provenant des plans admissibles sont :
\[ \beta_{r} \rightarrow \beta_{r-1} \rightarrow ... \rightarrow \beta_{s+1} \rightarrow \beta_{s}\]
\end{cor}
\begin{conv}
\begin{itemize}
\item Si  $C_{l} = \{\beta_{s}, \beta_{s+1},...,\beta_{r} = \alpha_{l}\}$ est la colonne $l$ avec $1 \leq l \leq n$, les colonnes tronquées contenues dans $C_l$ sont les parties de la forme $\{\beta_{s}, \beta_{s+1}, ..., \beta_{t}\}$, $t \in \llbracket s, r \rrbracket$.
\item On note $\mathcal{D}$ l'ensemble de tous les diagrammes qui vérifient les contraintes du corollaire \ref{contraintecol} pour $1 \leq l \leq n$.
\end{itemize}
\end{conv}
On remarque que $\mathcal{D}$ est  l'ensemble de tous les diagrammes $\Delta$ qui sont des réunions de colonnes tronquées (voir figure ci-dessous).

\renewcommand{\PuzzleBlackBox}{\rule{.75\PuzzleUnitlength}%
{.75\PuzzleUnitlength}} 
\PuzzleUnsolved

\begin{center}
\begin{tabular}{cp{2cm}c}
\begin{Puzzle}{5}{5}%
|* |1 |*|* | * |.
|{} |1|*|*| * |.
|{}|{} |* | 1|* |.
|{} |{}|{} |1| * |.
|{} |{} |{}|{}|1 |.
\end{Puzzle}  
& &
\begin{Puzzle}{5}{5}%
|* |1 |*|* | * |.
|{} |1|* |*| * |.
|{}|{} |* | 1|* |.
|{} |{}|{} |*| * |.
|{} |{} |{}|{}|1 |.
\end{Puzzle}\\
$\in \mathcal{D} $& &
$\notin \mathcal{D}$
\end{tabular}
\end{center}
\begin{rmq}
L'ensemble des diagrammes de Cauchon noté $\mathcal{D}$  a même cardinal que le groupe de Weyl $W$.
\end{rmq}
\begin{preuve}
Comme $\mathcal{D}$ est l'ensemble des diagrammes $\Delta$ qui sont des réunions de colonnes tronquées, on a $|\mathcal{D}| = (n+1)! = |W|$. 
\end{preuve}
\subsubsection{Le cas $B_n$, $n \geq 2$}
\begin{conv}
On numérote les racines simples de manière que le diagramme de Dynkin soit :
\[\alpha_1 \Leftarrow \alpha_2 - \cdots - \alpha_{n-1} - \alpha_n \]
On sait  (voir par exemple \cite[Section 6]{MR1628449}) que l'écriture 
\[w_0 = s_{\alpha_1} \circ (s_{\alpha_2}\circ s_{\alpha_1}\circ s_{\alpha_2})   \cdots \circ (s_{\alpha_n} \circ s_{\alpha_{n-1}} \circ \cdots \circ s_{\alpha_2} \circ s_{\alpha_1}\circ s_{\alpha_2} \circ \cdots \circ s_{\alpha_n})\]
 est une décomposition réduite qui induit l'ordre ci-dessous sur les racines positives.
\end{conv}
\begin{tabular}{c|c|c|p{3cm}|c|}
\cline{5-5}
\multicolumn{4}{c|}{} & \multicolumn{1}{|c|}{$\beta_{(n-1)^2+1} = 2\alpha_1 + \cdots  + 2\alpha_{n-1} + \alpha_n$}  \\
\cline{4-5}
\multicolumn{3}{c|}{} & \multicolumn{1}{|c|}{} & \multicolumn{1}{|c|}{$\vdots$} \\
 \cline{3-5}
\multicolumn{2}{c|}{} & \multicolumn{1}{|c|}{$
 \beta_2 =2 \alpha_1 + \alpha_2$} & \multicolumn{1}{|c|}{} & \multicolumn{1}{|c|}{$\beta_{N-n} = 2\alpha_1 + \alpha_{2} + \cdots + \alpha_{n-1} + \alpha_n $} \\ 
\cline{2-5}
 &  $ \beta_1 = \alpha_1 $& $ \beta_3 = \alpha_1 + \alpha_2 $ &  & $
\beta_{N-n+1} = \alpha_1 + \cdots + \alpha_{n-1} + \alpha_n $ \\
  \cline{2-5}
\multicolumn{2}{c|}{} & \multicolumn{1}{|c|}{$
 \beta_4 = \alpha_2$} & \multicolumn{1}{|c|}{} & \multicolumn{1}{|c|}{$\beta_{N-n+2} = \alpha_2 + \cdots + \alpha_{n-1} + \alpha_n $} \\
 \cline{3-5}
\multicolumn{3}{c|}{} & \multicolumn{1}{|c|}{} & \multicolumn{1}{|c|}{$\vdots$} \\
\cline{4-5}
\multicolumn{4}{c|}{} & \multicolumn{1}{|c|}{$
\beta_N = \alpha_n$}  \\
\cline{5-5}
\end{tabular}\\
\par
On constate que ceci est un ordre de Lusztig et qu'aucune des colonnes n'est exceptionnelle.
\begin{lemme}\label{contraintesbn}
Si $1 \leq l \leq n$ et si $\beta_{j}, \beta_{j+1}$ sont deux racines consécutives de $C_{l}$ alors $\beta_{j+1} \rightarrow \beta_{j} $ est une contrainte.
\end{lemme}
\begin{preuve}
Si $\beta_{j}$ et $\beta_{j+1}$ sont deux racines consécutives de $C_{l}$ alors il existe une racine simple $\alpha_i$ ($1 \leq i < l \leq n$) telle que $\beta_{j} = \beta_{j+1} + \alpha_{i}$. D'après la proposition \ref{contrainteelementaire}, $\beta \rightarrow \beta + \alpha_{i}$ est une contrainte.
\end{preuve}
On en déduit :
\begin{cor}\label{contraintecolbn}
Si $C_{l} = \{\beta_{s}, \beta_{s+1},...,\beta_{r} = \alpha_{l}\}$ est la colonne $l$ avec $1 \leq l \leq n$,les contraintes provenant des plans admissibles sont :
\[ \beta_{r} \rightarrow \beta_{r-1} \rightarrow ... \rightarrow \beta_{s+1} \rightarrow \beta_{s}\]
\end{cor}
\begin{conv}
\begin{itemize}
\item Si  $C_{l} = \{\beta_{s}, \beta_{s+1},...,\beta_{r} = \alpha_{l}\}$ est la colonne $l$ avec $1 \leq l \leq n$, les colonnes tronquées contenues dans $C_l$ sont les parties de la forme $\{\beta_{s}, \beta_{s+1}, ..., \beta_{t}\}$, $t \in \llbracket s, r \rrbracket$.
\item On note $\mathcal{D}$ l'ensemble de tous les diagrammes qui vérifient les contraintes du corollaire \ref{contraintecolbn} pour $1 \leq l \leq n$.
\end{itemize}
\end{conv}
On remarque que $\mathcal{D}$ est  l'ensemble de tous les diagrammes $\Delta$ qui sont des réunions de colonnes tronquées (voir figure ci-dessous).
\renewcommand{\PuzzleBlackBox}{\rule{.75\PuzzleUnitlength}%
{.75\PuzzleUnitlength}} 
\PuzzleUnsolved
\begin{center}
\begin{tabular}{cccc}
\begin{Puzzle}{7}{4}
|{}|{}|{}| *|.
|{}|{}| *| *| .
|{}| 1|* | 1| .
|*| 1| *|1 | . 
|{}| 1|* | 1| .
|{}|{}| 1| 1|.
|{}|{}|{}| 1| .
\end{Puzzle}  

& &

\begin{Puzzle}{7}{4}
|{}|{}|{}| 1| .
|{}|{}| *| 1| .
|{}| 1|* | 1|.
|*| 1| *|1 | . 
|{}| 1|* | *| .
|{}|{}| 1| *|.
|{}|{}|{}| 1|.
\end{Puzzle}  \\
&&&\\
&&&\\
&&&\\
&&&\\
&&&\\
 $\in \mathcal{D}$ & & $ \notin \mathcal{D}$
\end{tabular}
\end{center}
\begin{rmq}
L'ensemble des diagrammes de Cauchon noté $\mathcal{D}$  a même cardinal que le groupe de Weyl $W$.
\end{rmq}
\begin{preuve}
Comme $\mathcal{D}$ est l'ensemble des diagrammes $\Delta$ qui sont des réunions de colonnes tronquées, on a $|\mathcal{D}| = 2^{n+1}(n+1)! = |W|$.
\end{preuve}
\subsubsection{$C_n$, $n \geq 3$}
\begin{conv}
On numérote les racines simples de manière que le diagramme de Dynkin soit :
\[\alpha_1 \Rightarrow \alpha_2 - \cdots - \alpha_{n-1} - \alpha_n  \]
On sait  (voir par exemple \cite[Section 6]{MR1628449}) que l'écriture 
\[w_0 = s_{\alpha_1} \circ (s_{\alpha_2}\circ s_{\alpha_1}\circ s_{\alpha_2})   \cdots \circ (s_{\alpha_n} \circ s_{\alpha_{n-1}} \circ \cdots \circ s_{\alpha_2} \circ s_{\alpha_1}\circ s_{\alpha_2} \circ \cdots \circ s_{\alpha_n})\]
 est une décomposition réduite qui induit l'ordre ci-dessous sur les racines positives.
\end{conv}
\begin{tabular}{c|c|c|p{3cm}|c|}
\cline{5-5}
\multicolumn{4}{c|}{} & \multicolumn{1}{|c|}{$
\beta_{(n-1)^2+1} = \alpha_1 + 2\alpha_2 \cdots  + 2\alpha_{n-1} + \alpha_n$}  \\
\cline{4-5}
\multicolumn{3}{c|}{} & \multicolumn{1}{|c|}{} & \multicolumn{1}{|c|}{$\vdots$} \\
 \cline{3-5}
\multicolumn{2}{c|}{} & \multicolumn{1}{|c|}{$ \beta_2 = \alpha_1 + \alpha_2 $} & \multicolumn{1}{|c|}{} & \multicolumn{1}{|c|}{$\beta_{N-n} = \alpha_1 + \alpha_{2} + \cdots + \alpha_{n-1} + \alpha_n $} \\ 
\cline{2-5}
 &  $ \beta_1 = \alpha_1$& $ \beta_3 = \alpha_1 + 2\alpha_2 $ &  & $
\beta_{N-n+1} = \alpha_1 +2\alpha_2 \cdots + 2\alpha_{n-1} + 2\alpha_n$ \\
  \cline{2-5}
\multicolumn{2}{c|}{} & \multicolumn{1}{|c|}{$
 \beta_4 = \alpha_2 $} & \multicolumn{1}{|c|}{} & \multicolumn{1}{|c|}{$\beta_{N-n+2} = \alpha_2 + \cdots + \alpha_{n-1} + \alpha_n$} \\
 \cline{3-5}
\multicolumn{3}{c|}{} & \multicolumn{1}{|c|}{} & \multicolumn{1}{|c|}{$\vdots$} \\
\cline{4-5}
\multicolumn{4}{c|}{} & \multicolumn{1}{|c|}{$\beta_N = \alpha_n$}  \\
\cline{5-5}
\end{tabular}\\

\par
On constate que ceci est un ordre de Lusztig et que toutes les colonnes exceptée la première sont exceptionnelles.
\begin{lemme}\label{contraintescn}
Si $1 \leq l \leq n$ et si $\beta_{j}, \beta_{j+1}$ sont deux racines consécutives de $C_{l}$ alors $\beta_{j+1} \rightarrow \beta_{j} $ est une contrainte.
\end{lemme}
\begin{preuve}
Chaque colonne comporte un nombre impair d'éléments, notons $C_{l} = \{\beta_{u_{1}},...,\beta_{u_{2k+1}}\}$, $\beta_{u_{k+1}}$ est une racine exceptionnelle et toutes les boîtes comportent une unique racine.\\ 
\begin{itemizedot}
\item Soient $j \in \llbracket 1, k-1 \rrbracket \cup \llbracket k+1, 2k+1 \rrbracket$, $\beta_{u_{j}}$ et $\beta_{u_{j+1}}$ deux racines consécutives de $C_{l}$. alors il existe une racine simple $\alpha_{i_{j}}$ ($1 \leq i_{j} < l \leq n$) telle que $\beta_{u_{j}} = \beta_{u_{j+1}} + \alpha_{i_{j}}$. D'après la proposition \ref{contrainteelementaire}, $\beta_{u_{j+1}} \rightarrow \beta_{u_{j}}$ est une contrainte.
\item $\beta_{u_{k+2}}$ est dans la boîte qui suit la racine exceptionnelle, d'après la proposition \ref{contrainteelementaireavantex}, $\beta_{u_{k+2}} \rightarrow \beta_{u_{k+1}}$ est une contrainte.
\item Comme $\beta_{u_{k+1}} + \beta_{1} = 2 \beta_{u_{k-1}}$ alors d'après la proposition \ref{contrainteelementaireapresex}, $\beta_{u_{k+1}} \rightarrow \beta_{u_{k}}$ est une contrainte.
\end{itemizedot}
\end{preuve}
On en déduit :
\begin{cor}\label{contraintecolcn}
Si $C_{l} = \{\beta_{s}, \beta_{s+1},...,\beta_{r} = \alpha_{l}\}$ est la colonne $l$ avec $1 \leq l \leq n$, les contraintes provenant des plans admissibles sont :
\[ \beta_{r} \rightarrow \beta_{r-1} \rightarrow ... \rightarrow \beta_{s+1} \rightarrow \beta_{s}\]
\end{cor}
\begin{conv}
\begin{itemize}
\item Si  $C_{l} = \{\beta_{s}, \beta_{s+1},...,\beta_{r} = \alpha_{l}\}$ est la colonne $l$ avec $1 \leq l \leq n$, les colonnes tronquées contenues dans $C_l$ sont les parties de la forme $\{\beta_{s}, \beta_{s+1}, ..., \beta_{t}\}$, $t \in \llbracket s, r \rrbracket$.
\item On note $\mathcal{D}$ l'ensemble de tous les diagrammes qui vérifient les contraintes du corollaire \ref{contraintecolcn} pour $1 \leq l \leq n$.
\end{itemize}
\end{conv}
On remarque que $\mathcal{D}$ est  l'ensemble de tous les diagrammes $\Delta$ qui sont des réunions de colonnes tronquées (voir figure ci-dessous).
\renewcommand{\PuzzleBlackBox}{\rule{.75\PuzzleUnitlength}%
{.75\PuzzleUnitlength}} 
\PuzzleUnsolved
\begin{center}
\begin{tabular}{cccc}
\begin{Puzzle}{7}{4}
|{}|{}|{}| *|.
|{}|{}| *| *| .
|{}| 1|* | 1| .
|*| 1| *|1 | . 
|{}| 1|* | 1| .
|{}|{}| 1| 1|.
|{}|{}|{}| 1| .
\end{Puzzle}  

& &

\begin{Puzzle}{7}{4}
|{}|{}|{}| 1| .
|{}|{}| *| 1| .
|{}| 1|* | 1|.
|*| 1| *|1 | . 
|{}| 1|* | *| .
|{}|{}| 1| *|.
|{}|{}|{}| 1|.
\end{Puzzle}  \\
&&&\\
&&&\\
&&&\\
&&&\\
&&&\\
 $\in \mathcal{D}$ & & $ \notin \mathcal{D}$
\end{tabular}
\end{center}
\begin{rmq}
L'ensemble des diagrammes de Cauchon noté $\mathcal{D}$  a même cardinal que le groupe de Weyl $W$.
\end{rmq}
\begin{preuve}
Comme $\mathcal{D}$ est l'ensemble des diagrammes $\Delta$ qui sont des réunions de colonnes tronquées, on a $|\mathcal{D}| = 2^{n+1}(n+1)! = |W|$.
\end{preuve}
\subsubsection{$D_n$, $n \geq 4$}
\begin{conv}
On numérote les racines simples de manière que le diagramme de Dynkin soit :
\[\begin{array}{ccccccccccc}
\alpha_1 && & &   & &&&  &     &   \\ 
& \diagdown & & &   & &&&  &     & \\
&& \alpha_3 &\text{---~}& \alpha_{4} &\text{---~}& \cdots &\text{---~}& \alpha_{n-1}&\text{---~} & \alpha_n \\
& \diagup & & &   & &&&  &     & \\
\alpha_2&& & &  &  &&&&  &
\end{array}\]
On sait  (voir par exemple \cite[Section 6]{MR1628449}) que l'écriture :
 \[w_0 = s_{\alpha_1} \circ s_{\alpha_2} \circ (s_{\alpha_3}\circ s_{\alpha_1} \circ s_{\alpha_2} \circ s_{\alpha_3})   \cdots \circ (s_{\alpha_n} \circ s_{\alpha_{n-1}} \circ \cdots \circ s_{\alpha_3} \circ s_{\alpha_1} \circ s_{\alpha_2}\circ s_{\alpha_3} \circ \cdots \circ s_{\alpha_n}) \]
 est une décomposition réduite qui induit l'ordre ci-dessous sur les racines positives.\\
\begin{tabular}{c|c|c|c|c|}
\cline{5-5}
\multicolumn{4}{c|}{} & \multicolumn{1}{|c|}{$\beta_{N-2n+1} = \alpha_1 + \alpha_2 + 2 \alpha_3 \cdots  + 2\alpha_{n-1} + \alpha_n  $}  \\
\cline{4-5}
\multicolumn{3}{c|}{} & \multicolumn{1}{|c|}{} & \multicolumn{1}{|c|}{$\vdots$} \\
 \cline{3-5}
\multicolumn{2}{c|}{} & \multicolumn{1}{|c|}{$
 \beta_3 = \alpha_1 + \alpha_2 + \alpha_3   $} & \multicolumn{1}{|c|}{} & \multicolumn{1}{|c|}{$
\beta_{N-n-1} = \alpha_1 + \alpha_2 + \alpha_3 \cdots + \alpha_{n-1} + \alpha_n$} \\ 
\cline{2-5}
 &  $
 \beta_1 = \alpha_1 $& $  \beta_4 = \alpha_2 + \alpha_3 $ &  & $
\beta_{N-n} = \alpha_1 \ (\textrm{ou} \ \alpha_2) \ + \alpha_3 \cdots + \alpha_{n-1} + \alpha_n$\\

\cline{2-5}
 &   $\beta_2 = \alpha_2$& $ \beta_5 = \alpha_1 + \alpha_3 $ &  & $
\beta_{N-n+1} = \alpha_2 \ (\textrm{ou} \ \alpha_1) \ + \alpha_3 \cdots + \alpha_{n-1} + \alpha_n$ \\
  \cline{2-5}
\multicolumn{2}{c|}{} & \multicolumn{1}{|c|}{$
 \beta_6 = \alpha_3 $} & \multicolumn{1}{|c|}{} & \multicolumn{1}{|c|}{$
\beta_{N-n+2} = \alpha_3 \cdots + \alpha_{n-1} + \alpha_n$} \\
 \cline{3-5}
\multicolumn{3}{c|}{} & \multicolumn{1}{|c|}{} & \multicolumn{1}{|c|}{$\vdots$} \\
\cline{4-5}
\multicolumn{4}{c|}{} & \multicolumn{1}{|c|}{$
\beta_N = \alpha_n$}  \\
\cline{5-5}
\end{tabular}\\
On a réunit ici les colonnes 1 et 2 pour une meilleur lisibilité de l'ordre sur les racines.
\end{conv}
On constate que ceci est un ordre de Lusztig et qu'aucune des colonnes n'est exceptionnelle.
\begin{lemme}\label{contraintesdn}
Soit une colonne $C_l = \{\beta_{u_{1}},...,\beta_{u_{s}},\beta_{u_{s+1}}, ... ,\beta_{u_{2s}} \}$ (avec $s = l-1$, $l \geq 3$), on a les contraintes :\\
\begin{center}
\begin{psmatrix}[rowsep=.3cm, colsep=.6cm]
 				&				&		&			&$\beta_{u_{s}}$ \\ 
$\beta_{u_{2s}}$ &$ \beta_{u_{2s-1}} $& ... &$\beta_{u_{s+2}}$& 			&	$\beta_{u_{s-1}}$ &  ...& $\beta_{u_{2}}$  &$\beta_{u_{1}}$ \\
				& 						 &		&             &  $\beta_{u_{s+1}}$ 
\end{psmatrix}
\psset{arrows=->, nodesep=2pt}
\ncline{2,1}{2,2}
\ncline{2,2}{2,3}
\ncline{2,3}{2,4}
\ncline{2,4}{1,5}
\ncline{2,4}{3,5}
\ncline{1,5}{2,6}
\ncline{3,5}{2,6}
\ncline{2,6}{2,7}
\ncline{2,7}{2,8}
\ncline{2,8}{2,9}
\end{center}
\end{lemme}
\begin{preuve}
Observons tout d'abord que seule la boîte $B = \{\beta_{u_{s}}, \beta_{u_{s+1}}\}$ comporte plus d'un élément.
\begin{itemizedot}
\item Pour $j \in \llbracket 1, s-2 \rrbracket$, on a $\beta_{u_{j+1}} + \alpha_{s-j} = \beta_{u_{j}}$ et $\beta_{u_{2s-j+1}} + \alpha_{s-j+1} = \beta_{u_{2s-j}}$. Ainsi, par la proposition \ref{contrainteelementaire}, on a les contraintes $\beta_{u_{j+1}} \rightarrow \beta_{u_{j}}$ et $\beta_{u_{2s-j+1}}\rightarrow \beta_{u_{2s-j}}$
\item On a $\beta_{u_{s+2}} + \alpha_{1} \in B$ et $\beta_{u_{s+2}} + \alpha_{2} \in B$ de sorte que par la proposition \ref{contrainteelementaire}, on a les contraintes  $\beta_{u_{s+2}} \rightarrow \beta_{u_{s+1}}$ et $\beta_{u_{s+2}}\rightarrow \beta_{u_{s}}$.
\item On a $\beta_{u_{s-1}} - \alpha_{1} \in B$ et $\beta_{u_{s-1}} - \alpha_{2} \in B$ de sorte que par la proposition \ref{contrainteelementaire}, on a les contraintes  $\beta_{u_{s+1}} \rightarrow \beta_{u_{s-1}}$ et $\beta_{u_{s}}\rightarrow \beta_{u_{s-1}}$.
\end{itemizedot}
\end{preuve}
\begin{conv}
\begin{itemize}
\item Si  $C_{l} = \{\beta_{s}, \beta_{s+1},...,\beta_{r} = \alpha_{l}\}$ est la colonne $l$ avec $1 \leq l \leq n$, les colonnes tronquées privées de l'élément $\beta_{m}$ ($s \leq m \leq r$) contenues dans $C_l$ sont les parties de la forme $\{\beta_{s}, \beta_{s+1}, ...,\beta_{m-1},\beta_{m+1}, \beta_{t}\}$, $t \in \llbracket s, r \rrbracket$.
\item On note $\mathcal{D}$ l'ensemble de tous les diagrammes qui vérifient les contraintes du lemme \ref{contraintesdn} pour $1 \leq l \leq n$.
\end{itemize}
\end{conv}
On en déduit que $\mathcal{D}$ est l'ensemble des réunions de colonnes tronquées avec la possibilité d'avoir des colonnes tronquées avec un "trou" sur la ligne correspondant à $\alpha_1$ lorsque la colonne tronquée ne dépasse pas la ligne correspondant à $\alpha_2$.  \\
Ainsi :
\begin{center}
\begin{tabular}{cp{2cm}c}
\begin{Puzzle}{8}{4}
|{}|{}| *| .
|{}| *|* | .
|1| *| 1|. 
|*| *| *|. 
|{}| *|1 |.
|{}|{}| 1|.
\end{Puzzle}  
& &
\begin{Puzzle}{8}{4}
|{}|{}| *| .
|{}| *|* |.
|1| *| 1|.
|*| 1| *|.  
|{}| *|* | .
|{}|{}| *| .
\end{Puzzle}  \\
&&\\
&&\\
&& \\
 $\in \mathcal{D}$  & &$ \notin \mathcal{D}$
\end{tabular}\\
\end{center}
\begin{rmq}
L'ensemble des diagrammes de Cauchon noté $\mathcal{D}$  a même cardinal que le groupe de Weyl $W$.
\end{rmq}
\begin{preuve}
Comme $\mathcal{D}$ est l'ensemble des diagrammes $\Delta$ qui sont des réunions de colonnes tronquées qui possèdent éventuellement un trou, on a 
\[|\mathcal{D}| = 4 \times 6 \times 8 \times ... \times 2n = 2^{n-1}n! = |W|.\]
\end{preuve}
\subsection{Les cas exceptionnels}
\subsubsection{$G_2$}
\begin{conv}
On numérote les racines simples de manière que le diagramme de Dynkin soit :
\begin{center}
	$ \alpha_1 \Lleftarrow \alpha_2 $
\end{center}
On sait que l'écriture  $w_0 = s_{\alpha_1} \circ s_{\alpha_2} \circ s_{\alpha_1} \circ s_{\alpha_2} \circ s_{\alpha_1} \circ s_{\alpha_2} $ est une décomposition réduite de $w_{0}$ qui induit l'ordre suivant sur les racines positives.
\begin{center}
\begin{tabular}{c|c|c|}
\cline{3-3}
\multicolumn{2}{c|}{} & \multicolumn{1}{|c|}{$\beta_2 = 3\alpha_1 + \alpha_2$}\\
\cline{3-3}
\multicolumn{2}{c|}{} & \multicolumn{1}{|c|}{$\beta_3 = 2\alpha_1 + \alpha_2$} \\
\cline{2-3}
 &  $ \beta_1 = \alpha_1$ & $ \beta_4 = 3\alpha_1 + 2\alpha_2$ \\
 \cline{2-3}
\multicolumn{2}{c|}{} & \multicolumn{1}{|c|}{$\beta_5 = \alpha_1 + \alpha_2$}\\
 \cline{3-3}
\multicolumn{2}{c|}{} & \multicolumn{1}{|c|}{$\beta_6 = \alpha_2$}\\
\cline{3-3}
\end{tabular}
\end{center}
\end{conv}
\begin{lemme}
On a les contraintes :
\scalebox{.8}{\begin{psmatrix}[rowsep=.2cm, colsep=1cm]
$\beta_{6}$ & $\beta_{5}$& $\beta_{4}$& $\beta_{3}$& $\beta_{2}$
\end{psmatrix}
\psset{nodesep=2pt, arrows=->}
\ncline{1,1}{1,2}
\ncline{1,2}{1,3}
\ncline{1,3}{1,4}
\ncline{1,4}{1,5}}
\end{lemme}
\begin{preuve}\label{contraintecolg2}
Pour trouver ces contraintes, on applique les propositions \ref{contrainteelementaire}, \ref{contrainteelementaireavantex} et \ref{contrainteelementaireapresex} car on a les égalités suivantes ($\beta_{4}$ est une racine exceptionnelle):
\[\beta_{6} + \alpha_{1} = \beta_{5},\  h'(\beta_{5}) + 1 = \beta_{4},\  \beta_{4} + \alpha_{1} = 2 \beta_{3},\hspace{.2cm} \beta_{3} + \alpha_{1} = \beta_{2}.\]
\end{preuve}
\begin{convsanss}
On note $\mathcal{D}$ l'ensemble de tous les diagrammes qui vérifient les contraintes du lemme \ref{contraintecolg2} .
\end{convsanss}
\begin{rmq}
L'ensemble des diagrammes de Cauchon noté $\mathcal{D}$ a même cardinal que le groupe de Weyl $W$.
\end{rmq}
\begin{preuve}
Comme $\mathcal{D}$ est l'ensemble des diagrammes $\Delta$ qui sont des réunions de colonnes tronquées 
\[|\mathcal{D}| = 2 \times 6 = 12 = |W|.\]
\end{preuve}
\subsubsection{$F_4$}
\begin{conv}
On numérote les racines simple de sorte que le diagramme de Dynkin soit : 
	\[\alpha_1 \ \text{---~} \alpha_2 \Rightarrow \alpha_3 \ \text{---~} \alpha_4\]
On choisit la décomposition réduite de l'élément de plus grande longueur du groupe de Weyl suivante : 
 $$w_0 = s_4s_3s_4s_2s_3s_4s_2s_3s_2s_1s_2s_3s_4s_2s_3s_1s_2 s_1 s_3 s_4 s_2s_3s_2s_1 $$
Cette décomposition induit l'ordre convexe suivant sur les racines positives :\\
\begin{tabular}{cc}
\begin{tabular}{c|c|c|c|c|}
\cline{5-5}
\multicolumn{4}{c|}{} & \multicolumn{1}{|c|}{$
\beta_{10} (1, 3, 4, 2)$}\\
\cline{5-5}
\multicolumn{4}{c|}{} & \multicolumn{1}{|c|}{$
\beta_{11} (1, 2, 4, 2) $}\\
\cline{5-5}
\multicolumn{4}{c|}{} & \multicolumn{1}{|c|}{$
\beta_{12} (1, 2, 3, 2)$}\\
\cline{5-5}
\multicolumn{4}{c|}{} & \multicolumn{1}{|c|}{$
\beta_{13} (1, 2, 3, 1)$}\\
\cline{4-5}
\multicolumn{3}{c|}{} & \multicolumn{1}{|c|}{$\beta_4(0, 1, 2, 2 )$}&\multicolumn{1}{|c|}{$
\beta_{14} (1, 2, 2, 2)$}\\
\cline{4-5}
\multicolumn{3}{c|}{} & \multicolumn{1}{|c|}{$\beta_5( 0, 1, 2, 1)$} & \multicolumn{1}{|c|}{$\beta_{15} (1, 2, 2, 1)$} \\
 \cline{3-5}
\multicolumn{2}{c|}{} & \multicolumn{1}{|c|}{$\beta_2 (0, 0, 1, 1)$} & \multicolumn{1}{|c|}{$\beta_6 (0, 1, 1, 1)$} & \multicolumn{1}{|c|}{$\beta_{16} (1, 1, 2, 2)$} \\ 
\cline{2-5}
 &  $\beta_1 ( 0, 0, 0, 1) $& $ \beta_3 (0, 0, 1, 0)$  & $\beta_7 (0, 1, 2, 0) $& $
\beta_{17} (2, 3, 4, 2)$ \\
  \cline{2-5}
\multicolumn{3}{c|}{} & \multicolumn{1}{|c|}{$\beta_8 (0, 1, 1, 0)$} & \multicolumn{1}{|c|}{$\beta_{18} (1, 2, 2, 0)$} \\
 \cline{4-5}
\multicolumn{3}{c|}{} & \multicolumn{1}{|c|}{$\beta_9 (0, 1, 0, 0)$} & \multicolumn{1}{|c|}{$\beta_{19} (1, 1, 2, 1)$} \\
\cline{4-5}
\multicolumn{4}{c|}{} & \multicolumn{1}{|c|}{$\beta_{20} (1, 1, 1, 1)$}  \\
\cline{5-5}
\multicolumn{4}{c|}{} & \multicolumn{1}{|c|}{$\beta_{21} (1, 1, 2, 0)$}  \\
\cline{5-5}
\multicolumn{4}{c|}{} & \multicolumn{1}{|c|}{$\beta_{22} (1, 1, 1, 0)$}  \\
\cline{5-5}
\multicolumn{4}{c|}{} & \multicolumn{1}{|c|}{$\beta_{23} (1, 1, 0, 0)$}  \\
\cline{5-5}
\multicolumn{4}{c|}{} & \multicolumn{1}{|c|}{$\beta_{24} (1, 0, 0, 0)$}  \\
\cline{5-5}
\end{tabular}
&
\begin{tabular}{c|c|c|c|c|}
\cline{5-5}
\multicolumn{4}{c|}{} & \multicolumn{1}{|c|}{10}\\
\cline{5-5}
\multicolumn{4}{c|}{} & \multicolumn{1}{|c|}{9}\\
\cline{5-5}
\multicolumn{4}{c|}{} & \multicolumn{1}{|c|}{8}\\
\cline{5-5}
\multicolumn{4}{c|}{} & \multicolumn{1}{|c|}{7}\\
\cline{4-5}
\multicolumn{3}{c|}{} & \multicolumn{1}{|c|}{5}&\multicolumn{1}{|c|}{7}\\
\cline{4-5}
\multicolumn{3}{c|}{} & \multicolumn{1}{|c|}{4} & \multicolumn{1}{|c|}{6} \\
 \cline{3-5}
\multicolumn{2}{c|}{} & \multicolumn{1}{|c|}{2} & \multicolumn{1}{|c|}{3} & \multicolumn{1}{|c|}{6} \\ 
\cline{2-5}
 &  1&  1  & 3 & 11/2 \\
  \cline{2-5}
\multicolumn{3}{c|}{} & \multicolumn{1}{|c|}{2} & \multicolumn{1}{|c|}{5} \\
 \cline{4-5}
\multicolumn{3}{c|}{} & \multicolumn{1}{|c|}{1} & \multicolumn{1}{|c|}{5} \\
\cline{4-5}
\multicolumn{4}{c|}{} & \multicolumn{1}{|c|}{4}  \\
\cline{5-5}
\multicolumn{4}{c|}{} & \multicolumn{1}{|c|}{4}  \\
\cline{5-5}
\multicolumn{4}{c|}{} & \multicolumn{1}{|c|}{3}  \\
\cline{5-5}
\multicolumn{4}{c|}{} & \multicolumn{1}{|c|}{2}  \\
\cline{5-5}
\multicolumn{4}{c|}{} & \multicolumn{1}{|c|}{1}  \\
\cline{5-5}
\end{tabular}
\end{tabular}\\
\end{conv}
On vérifie que chaque colonne est ordinaire ou exceptionnelle puis on refait le même tableau en calculant $h'(\beta_i)$ pour tout i pour vérifier que l'ordre choisi correspond bien à celui de Lusztig. On connaît déjà la forme des diagrammes sur les 2 premières colonnes. Grâce aux relations de commutations et à la proposition \ref{contrainteelementaire}, on a les contraintes suivantes pour l'avant dernière colonne :
\begin{center}
\begin{psmatrix}[rowsep=.3cm, colsep=.5cm]
&  & 6 \\
 9 & 8 & & 5& 4 \\ 
 & & 7  \\  
\end{psmatrix}
\psset{arrows=->, nodesep=2pt}
\ncline{2,1}{2,2}
\ncline{2,2}{1,3}
\ncline{2,2}{3,3}
\ncline{1,3}{2,4}
\ncline{3,3}{2,4}
\ncline{2,4}{2,5}
\end{center}
Il y a donc 8 possibilités de remplir cette colonne en suivant les contraintes.\\

Au moyen des propositions \ref{contrainteelementaire}, \ref{contrainteelementaireavantex} et \ref{contrainteelementaireapresex}, on obtient les contraintes suivantes pour la dernière colonne :\\
\begin{tabular}{lc}
\begin{tabular}{c|c}
Diagrammes commençant par la case : & Nombre de diagrammes\\
\hline
Pas de case & 1\\
10					& 1	\\
11					& 1\\
12					& 1\\
13					& 1\\
14				& 2\\
15					& 1\\
16					& 3\\
17					& 2\\
18					& 2\\
19				& 2\\
20					& 2\\
21					& 2\\
22					& 1\\
23					& 1\\
24					& 1\\
\hline
 & \\
TOTAL 			& 24 
\end{tabular}
&
\begin{minipage}{3cm}
\begin{psmatrix}[rowsep=.3cm, colsep=.2cm]
 & 24 \\ 
 & 23 \\ 
 & 22 \\ 
20 &  & 21 \\ 
19 & &18 \\ 
&   17 \\ 
 16& & 15 \\ 
14 & & 13 \\ 
   & 12 \\ 
 & 11 &\\
 & 10&  
\end{psmatrix}
\psset{arrows=->, nodesep=2pt}
\ncline{10,2}{11,2}
\ncline{9,2}{10,2}
\ncline{8,3}{9,2}
\ncline{8,1}{9,2}
\ncline{7,3}{8,3}
\ncline{7,3}{8,1}
\ncline{7,1}{8,1}
\ncline{6,2}{7,3}
\ncline{5,1}{7,1}
\ncline{5,3}{6,2}
\ncline{5,1}{6,2}
\ncline{4,3}{5,3}
\ncline{4,1}{5,1}
\ncline{4,3}{5,1}
\ncline{3,2}{4,1}
\ncline{3,2}{4,3}
\ncline{2,2}{3,2}
\ncline{1,2}{2,2}
\end{minipage}
\end{tabular}\\
\par
On dénombre donc : $2 \times 3 \times 8 \times 24 = 2^7\times 3^2$ diagrammes soit le cardinal du groupe de Weyl.
\subsubsection{$E_6$}
\begin{conv}
On numérote les racines simple de sorte que le diagramme de Dynkin soit : 
\[\begin{array}{ccccccccc}
 					& 					&   				 & 					 & \alpha_2	&&  &     &   \\ 
 					& 					&   				 & 					 & |					&&  &     & \\
 \alpha_1 &\text{---~}& \alpha_{3} &\text{---~}& \alpha_4 &\text{---~}& \alpha_{5}&\text{---~} & \alpha_6 \\
\end{array}\]
Pour décrire $w_0$, on va remarquer que les racines $\alpha_1$ à $\alpha_5$ forment un système de racines de type $D_5$. Si on appelle $\tau$, l'élément de plus grande longueur utilisé précédemment pour le type $D_5$ alors la décomposition 
 $$w_0 = \tau s_{6} s_5 s_4 s_2 s_3  s_1  s_4  s_3 s_5  s_4  s_6 s_2 s_5  s_4 s_3 s_1$$
induit l'ordre convexe suivant sur les racines positives : \\
 \small
\begin{minipage}{20cm}
\begin{tabular}{p{0cm}@{}|c|c|c|c|c|}
\cline{6-6}
\multicolumn{5}{c|}{} & \multicolumn{1}{|c|}{ $ 
\beta_{21} = (1,2,2,3,2,1)$}\\
\cline{6-6}
\multicolumn{5}{c|}{} & \multicolumn{1}{|c|}{$ 
\beta_{22} = (1,1,2,3,2,1)$}\\
\cline{6-6}
\multicolumn{5}{c|}{} & \multicolumn{1}{|c|}{$
\beta_{23} = (1,1,2,2,2,1)$}\\
\cline{6-6}
\multicolumn{5}{c|}{} & \multicolumn{1}{|c|}{$
\beta_{24} = (1,1,2,2,1,1)$}\\
\cline{5-6}
\multicolumn{4}{c|}{} & \multicolumn{1}{|c|}{\tiny $\beta_{13}=\alpha_1+\alpha_2 + 2\alpha_3 + 2\alpha_4 + \alpha_5$} & \multicolumn{1}{|c|}{$\beta_{25} = (1,1,1,2,2,1)$} \\
 \cline{4-6}
\multicolumn{3}{c|}{} & \multicolumn{1}{|c|}{\tiny $\beta_7 =  \alpha_2 + \alpha_3 + 2\alpha_4 + \alpha_5$} & \multicolumn{1}{|c|}{\tiny $\beta_{14} = \alpha_1 + \alpha_2 + \alpha_3 + 2\alpha_4 + \alpha_5$} & \multicolumn{1}{|c|}{$\beta_{26} = (0,1,1,2,2,1)$} \\ 
\cline{3-6} 
\multicolumn{2}{c|}{} & \multicolumn{1}{|c|}{\tiny $\beta_3 =  \alpha_2 + \alpha_4+\alpha_5 $} & \multicolumn{1}{|c|}{\tiny $\beta_{8} =  \alpha_2 + \alpha_3 + \alpha_4+ \alpha_5 $} & \multicolumn{1}{|c|}{\tiny $\beta_{15} = \alpha_1 + \alpha_2 + \alpha_3 + \alpha_4 + \alpha_5$}& \multicolumn{1}{|c|}{$\beta_{27} = (1,1,1,2,1,1)$} \\ 
\cline{2-6}
 & \tiny $\beta_1 = \alpha_2 $& \tiny $ \beta_4 = \alpha_4 + \alpha_5$  & \tiny $\beta_9 = \alpha_2 + \alpha_3 + \alpha_4 $& \tiny $  \beta_{16} = \alpha_1 + \alpha_{3} + \alpha_4 + \alpha_5$ & $\beta_{28} =(0,1,1,2,1,1)$ \\
  \cline{2-6}
& \tiny $\beta_2 = \alpha_5 $& \tiny $ \beta_5 = \alpha_2 + \alpha_4$  & \tiny $\beta_{10} =  \alpha_3 + \alpha_4 +\alpha_5 $&\tiny  $ \beta_{17} = \alpha_1+ \alpha_2 + \alpha_{3} + \alpha_4 $ & $\beta_{29} = (1,1,1,1,1,1)$ \\
  \cline{2-6}
\multicolumn{2}{c|}{} & \multicolumn{1}{|c|}{\tiny $\beta_6 = \alpha_4$} & \multicolumn{1}{|c|}{\tiny $\beta_{11} = \alpha_3 + \alpha_4$} & \multicolumn{1}{|c|}{\tiny $\beta_{18} =\alpha_1 + \alpha_3 + \alpha_4  $} & \multicolumn{1}{|c|}{$\beta_{30} =  (0,1,1,1,1,1)$} \\
 \cline{3-6}
\multicolumn{3}{c|}{} & \multicolumn{1}{|c|}{\tiny $\beta_{12} = \alpha_3$} & \multicolumn{1}{|c|}{\tiny $\beta_{19} = \alpha_1 + \alpha_3$} & \multicolumn{1}{|c|}{ $\beta_{31} =  (1,0,1,1,1,1)$} \\
\cline{4-6}
\multicolumn{4}{c|}{} & \multicolumn{1}{|c|}{\tiny $\beta_{20} = \alpha_1$} & \multicolumn{1}{|c|}{$\beta_{32} =(0,1,0,1,1,1)$}  \\
\cline{5-6}
\multicolumn{5}{c|}{} & \multicolumn{1}{|c|}{$\beta_{33} = (0,0,1,1,1,1) $}  \\
\cline{6-6}
\multicolumn{5}{c|}{} & \multicolumn{1}{|c|}{$\beta_{34} =  (0,0,0,1,1,1)$}  \\
\cline{6-6}
\multicolumn{5}{c|}{} & \multicolumn{1}{|c|}{$\beta_{35} = (0,0,0,0,1,1)$}  \\
\cline{6-6}
\multicolumn{5}{c|}{} & \multicolumn{1}{|c|}{$\beta_{36} = (0,0,0,0,0,1)$}  \\
\cline{6-6}
\end{tabular}\\
\end{minipage}
\end{conv}
On obtient les contraintes ci-dessous au moyen de la proposition \ref{contrainteelementaire}, on dénombre alors les diagrammes qui vérifient ces contraintes :\\
\begin{tabular}{lccr}
\begin{minipage}{3cm}
\begin{psmatrix}[rowsep=.5cm, colsep=.2cm]
 & 36 \\ 
 & 35 \\ 
 & 34 &  \\ 
33& &32\\
31 &  & 30 \\ 
29 &  & 28 \\ 
27 &  & 26 \\ 
24 &  & 25 \\
&23 \\
 & 22 \\ 
 & 21 
\end{psmatrix}
\psset{arrows=->, nodesep=2pt}
\ncline{1,2}{2,2}
\ncline{2,2}{3,2}
\ncline{3,2}{4,1}
\ncline{3,2}{4,3}
\ncline{4,1}{5,1}
\ncline{4,1}{5,3}
\ncline{4,3}{5,3}
\ncline{5,1}{6,1}
\ncline{5,3}{6,1}
\ncline{5,3}{6,3}
\ncline{6,1}{7,1}
\ncline{6,3}{7,1}
\ncline{6,3}{7,3}
\ncline{7,1}{8,3}
\ncline{7,1}{8,1}
\ncline{7,3}{8,3}
\ncline{8,3}{9,2}
\ncline{8,1}{9,2}
\ncline{9,2}{10,2}
\ncline{10,2}{11,2}
\end{minipage}
& &&
\begin{minipage}{6cm}
\begin{tabular}{c|c}
Diagrammes commençant par la case : & Nombre de diagrammes\\
\hline
Pas de case & 1\\
21					& 1	\\
22					& 1\\
23					& 1\\
24					& 1\\
25					& 2\\
26					& 2\\
27					& 2\\
28					& 1\\
29					& 3\\
30					& 1\\
31					& 4\\
32					& 2\\
33					& 2\\
34					& 1\\
35					& 1\\
36					& 1\\
\hline
 & \\
TOTAL 			& 27 = $3^3$
\end{tabular} \\
\end{minipage}
\end{tabular} \\
\par
Pour les 5 première colonnes, on retrouve un système de type $D_5$, soit $4 \times 6 \times 8 \times 10 = 2^7 \times 3 \times 5$ diagrammes.\\
En ajoutant la dernière colonne, on obtient $2^7 \times 3^4 \times 5$ diagrammes soit le cardinal du groupe de Weyl associé. On vérifie donc encore que le cardinal de $\mathcal{D}$ (ensemble des diagrammes de Cauchon) est le même que celui de $W$.
\subsubsection{$E_7$}
\begin{conv}
On numérote les racines simple de sorte que le diagramme de Dynkin soit : 
\[\begin{array}{ccccccccccc}
 					& 					&   				 & 					 & \alpha_2	&&  &     &&&  \\ 
 					& 					&   				 & 					 & |					&&  &     & &&\\
 \alpha_1 &\text{---~}& \alpha_{3} &\text{---~}& \alpha_4 &\text{---~}& \alpha_{5}&\text{---~} & \alpha_6 &\text{---~}& \alpha_7 \\
\end{array}\]
Comme les racines $\alpha_1$ à $\alpha_6$ forment un système de racines de type $E_6$, on appelle $\sigma$, l'élément de plus grande longueur utilisé précédemment pour le type $E_6$. La décomposition 
\[ w_0 = \sigma s_{7} s_6 s_5 s_4 s_2 s_3  s_1  s_4  s_3 s_5  s_4  s_6 s_2 s_5  s_7 s_4 s_6 s_3 s_5 s_1 s_4 s_2 s_3 s_4 s_5 s_6 s_7\]
 induit alors l'ordre convexe suivant sur les racines positives. Pour pouvoir stocker toutes les informations dans un tableau, on va remplacer les vecteurs par leurs coordonnées sur la base $\alpha_1, \cdots , \alpha_7$ en ne mettant pas les coordonnées inutiles. \\
\begin{tabular}{p{0cm}@{}|c|c|c|c|c|c|}
\multicolumn{6}{c|}{} & \multicolumn{1}{|c|}{\tiny $ 
\mathbf{\alpha_1,\alpha_{2},\alpha_3,\alpha_4,\alpha_5,\alpha_6, \alpha_7}$}\\
\cline{7-7}
\multicolumn{6}{c|}{} & \multicolumn{1}{|c|}{\tiny $\beta_{37} (2, 2, 3, 4, 3, 2, 1)$}\\
\cline{7-7}
\multicolumn{6}{c|}{} & \multicolumn{1}{|c|}{\tiny $\beta_{38} (1, 2, 3, 4, 3, 2, 1)$}\\
\cline{7-7}
\multicolumn{6}{c|}{} & \multicolumn{1}{|c|}{\tiny $\beta_{39} (1, 2, 2, 4, 3, 2, 1)$}\\
\cline{7-7}
\multicolumn{6}{c|}{} & \multicolumn{1}{|c|}{\tiny $\beta_{40} (1, 2, 2, 3, 3, 2, 1)$}\\
\cline{7-7}
\multicolumn{6}{c|}{} & \multicolumn{1}{|c|}{\tiny $\beta_{41} ( 1, 1, 2, 3, 3, 2, 1 )$}\\
\cline{7-7}
\multicolumn{5}{c|}{} & \multicolumn{1}{|c|}{\tiny $ 
\mathbf{\alpha_1,\alpha_{2},\alpha_3,\alpha_4,\alpha_5,\alpha_6}$}& \multicolumn{1}{|c|}{\tiny $\beta_{42} (1, 2, 2, 3, 2, 2, 1)$}\\
\cline{6-7}
\multicolumn{5}{c|}{} & \multicolumn{1}{|c|}{\tiny $ 
\beta_{21} (1, 2, 2, 3, 2, 1)$}& \multicolumn{1}{|c|}{\tiny $\beta_{43} (1, 2, 2, 3, 2, 1, 1)$}\\
\cline{6-7}
\multicolumn{5}{c|}{} & \multicolumn{1}{|c|}{\tiny$ 
\beta_{22} (1, 1, 2, 3, 2, 1) $}& \multicolumn{1}{|c|}{\tiny $\beta_{44} (1, 1, 2, 3, 2, 2, 1)$} \\
\cline{6-7}
\multicolumn{5}{c|}{} & \multicolumn{1}{|c|}{\tiny$
\beta_{23} (1, 1, 2, 2, 2, 1)$}& \multicolumn{1}{|c|}{\tiny $\beta_{45} (1, 1, 2, 3, 2, 1, 1)$}\\
\cline{6-7}
\multicolumn{4}{c|}{} & \multicolumn{1}{|c|}{\tiny $ 
\mathbf{\alpha_1,\alpha_{2},\alpha_3,\alpha_4,\alpha_5,\alpha_6}$}& \multicolumn{1}{|c|}{\tiny$
\beta_{24}(1, 1, 2, 2, 1, 1)$}& \multicolumn{1}{|c|}{\tiny $\beta_{46} (1, 1, 2, 2, 2, 2, 1)$}\\
\cline{5-7}
\multicolumn{3}{c|}{} & \multicolumn{1}{|c|}{\tiny $ 
\mathbf{\alpha_2,\alpha_{3},\alpha_4,\alpha_5}$}& \multicolumn{1}{|c|}{\tiny $\beta_{13}(1,1,2,2,1)$} & \multicolumn{1}{|c|}{\tiny$\beta_{25} (1, 1, 1, 2, 2, 1)$}& \multicolumn{1}{|c|}{\tiny $\beta_{47} (1, 1, 2, 2, 2, 1, 1)$} \\
 \cline{4-7}
\multicolumn{2}{c|}{}& \multicolumn{1}{|c|}{\tiny $ 
\mathbf{\alpha_2,\alpha_{4},\alpha_5}$} & \multicolumn{1}{|c|}{\tiny $\beta_7 (1,1,2,1)$} & \multicolumn{1}{|c|}{\tiny $\beta_{14} (1,1,1,2,1)$} & \multicolumn{1}{|c|}{\tiny$\beta_{26} (0, 1, 1, 2, 2, 1)$}& \multicolumn{1}{|c|}{\tiny $\beta_{48}(1, 1, 1, 2, 2, 2, 1) $} \\ 
\cline{3-7} 
\multicolumn{1}{c|}{}& \multicolumn{1}{|c|}{\tiny $ 
\mathbf{\alpha_2,\alpha_{5}}$} & \multicolumn{1}{|c|}{\tiny $\beta_3 (1,1,1) $} & \multicolumn{1}{|c|}{\tiny $\beta_{8} (1,1,1,1) $} & \multicolumn{1}{|c|}{\tiny $\beta_{15} (1,1,1,1,1)$}& \multicolumn{1}{|c|}{\tiny$\beta_{27} (1, 1, 1, 2, 1, 1)$}& \multicolumn{1}{|c|}{\tiny $\beta_{49} (1, 1, 2, 2, 1, 1, 1)$} \\ 
\cline{2-7}
 & \tiny $\beta_1 (1,0) $& \tiny $ \beta_4 (0,1,1)$  & \tiny $\beta_9 (1,1,1,0) $& \tiny $  \beta_{16} (1,0,1,1,1)$ & \tiny $\beta_{28} (0, 1, 1, 2, 1, 1)$ & \tiny $\beta_{50} (1, 1, 1, 2, 2, 1, 1)$ \\
  \cline{2-7}
& \tiny $\beta_2 (0,1) $& \tiny $ \beta_5 (1,1,0)$  & \tiny $\beta_{10} (0,1,1,1) $&\tiny  $ \beta_{17} (1,1,1,1,0)$ & \tiny $\beta_{29} (1,1,1,1,1,1)$& \tiny $\beta_{51} (0, 1, 1, 2, 2, 2, 1)$ \\
  \cline{2-7}
\multicolumn{2}{c|}{} & \multicolumn{1}{|c|}{\tiny $\beta_6 (0,1,0)$} & \multicolumn{1}{|c|}{\tiny $\beta_{11} (0,1,1,0)$} & \multicolumn{1}{|c|}{\tiny $\beta_{18} (1,0,1,1,0)$} & \multicolumn{1}{|c|}{\tiny $\beta_{30} (0, 1, 1, 1, 1, 1)$}& \multicolumn{1}{|c|}{\tiny $\beta_{52} (1, 1, 1, 2, 1, 1, 1)$} \\
 \cline{3-7}
\multicolumn{3}{c|}{} & \multicolumn{1}{|c|}{\tiny $\beta_{12} (0,1,0,0)$} & \multicolumn{1}{|c|}{\tiny $\beta_{19}(1,0,1,0,0)$} & \multicolumn{1}{|c|}{\tiny $\beta_{31} (1, 0, 1, 1, 1, 1)$}& \multicolumn{1}{|c|}{\tiny $\beta_{53} (0, 1, 1, 2, 2, 1, 1)$} \\
\cline{4-7}
\multicolumn{4}{c|}{} & \multicolumn{1}{|c|}{\tiny $\beta_{20}(1,0,0,0,0)$} & \multicolumn{1}{|c|}{\tiny$\beta_{32}(0, 1, 0, 1, 1, 1) $}& \multicolumn{1}{|c|}{\tiny $\beta_{54}(1, 1, 1, 1, 1, 1, 1) $}  \\
\cline{5-7}
\multicolumn{5}{c|}{} & \multicolumn{1}{|c|}{\tiny$\beta_{33} (0, 0, 1, 1, 1, 1)$}& \multicolumn{1}{|c|}{\tiny $\beta_{55}(0, 1, 1, 2, 1, 1, 1) $}  \\
\cline{6-7}
\multicolumn{5}{c|}{} & \multicolumn{1}{|c|}{\tiny$\beta_{34} (0, 0, 0, 1, 1, 1)$}& \multicolumn{1}{|c|}{\tiny $\beta_{56}(1, 0, 1, 1, 1, 1, 1) $}  \\
\cline{6-7}
\multicolumn{5}{c|}{} & \multicolumn{1}{|c|}{\tiny$\beta_{35} ( 0, 0, 0, 0, 1, 1)$}& \multicolumn{1}{|c|}{\tiny $\beta_{57} (0, 1, 1, 1, 1, 1, 1)$}  \\
\cline{6-7}
\multicolumn{5}{c|}{} & \multicolumn{1}{|c|}{\tiny$\beta_{36} (0,0,0,0,0,1)$}& \multicolumn{1}{|c|}{\tiny $\beta_{58} (0, 1, 0, 1, 1, 1, 1)$}  \\
\cline{6-7}
\multicolumn{6}{c|}{} & \multicolumn{1}{|c|}{\tiny $\beta_{59} (0, 0, 1, 1, 1, 1, 1)$}\\
\cline{7-7}
\multicolumn{6}{c|}{} & \multicolumn{1}{|c|}{\tiny $\beta_{60} (0, 0, 0, 1, 1, 1, 1)$}\\
\cline{7-7}
\multicolumn{6}{c|}{} & \multicolumn{1}{|c|}{\tiny $\beta_{61} (0, 0, 0, 0, 1, 1, 1)$}\\
\cline{7-7}
\multicolumn{6}{c|}{} & \multicolumn{1}{|c|}{\tiny $\beta_{62} (0,0,0,0,0,1,1)$}\\
\cline{7-7}
\multicolumn{6}{c|}{} & \multicolumn{1}{|c|}{\tiny $\beta_{63} (0,0,0,0,0,0,1)$}\\
\cline{7-7}
\end{tabular}\\
\end{conv}
On connaît déjà la forme des diagrammes sur les 6 premières colonnes. On utilise la même méthode que pour $E_6$ (proposition \ref{contrainteelementaire}) pour trouver les contraintes dans la dernière colonne de $E_7$. On obtient : \\
\begin{tabular}{cc}
\begin{minipage}{3cm}
\begin{psmatrix}[rowsep=.3cm, colsep=.2cm]
 &  & 63 \\ 
 &  & 62 \\ 
 &  & 61 \\ 
 &  & 60 \\ 
 & 59 &  & 58 \\ 
&56 &  & 57 \\ 
 & 54 &  & 55 \\ 
 & 52 &  & 53 \\ 
  49 &  & 50 &  & 51 \\ 
  & 47 &  & 48 \\ 
 & 45 &  & 46 \\ 
&43 &  & 44 \\ 
 & 42 &  & 41 \\ 
 &  & 40 \\ 
 &  & 39 \\ 
 &  & 38 \\
 & & 37
\end{psmatrix}
\psset{arrows=->, nodesep=2pt}
\ncline{1,3}{2,3}
\ncline{2,3}{3,3}
\ncline{3,3}{4,3}
\ncline{4,3}{5,2}
\ncline{4,3}{5,4}
\ncline{5,2}{6,4}
\ncline{5,2}{6,2}
\ncline{5,4}{6,4}
\ncline{6,4}{7,4}
\ncline{6,4}{7,2}
\ncline{6,2}{7,2}
\ncline{7,4}{8,4}
\ncline{7,4}{8,2}
\ncline{7,2}{8,2}
\ncline{8,4}{9,5}
\ncline{8,4}{9,3}
\ncline{8,2}{9,3}
\ncline{8,2}{9,1}
\ncline{9,5}{10,4}
\ncline{9,3}{10,4}
\ncline{9,3}{10,2}
\ncline{9,1}{10,2}
\ncline{10,4}{11,4}
\ncline{10,2}{11,4}
\ncline{10,2}{11,2}
\ncline{11,4}{12,4}
\ncline{11,2}{12,4}
\ncline{11,2}{12,2}
\ncline{12,4}{13,2}
\ncline{12,4}{13,4}
\ncline{12,2}{13,2}
\ncline{13,2}{14,3}
\ncline{13,4}{14,3}
\ncline{14,3}{15,3}
\ncline{15,3}{16,3}
\ncline{16,3}{17,3}
\end{minipage}
&
\begin{tabular}{c|c}
Diagrammes commençant par la case : & Nombre de diagrammes\\
\hline
Pas de case & 1\\
37					& 1	\\
38					& 1\\
39					& 1\\
40					& 1\\
41					& 1\\
42					& 2\\
43					& 2\\
44					& 2\\
45					& 1\\
46					& 3\\
47					& 1\\
48					& 4\\
49					& 2\\
50					& 2\\
51					& 6\\
52					& 2\\
53					& 4\\
54					&	3\\
55					& 2\\
56					& 4\\
57					&	2\\
58					& 2\\
59					& 2\\
60					& 1\\
61					& 1\\
62					& 1\\
63					&	1\\
\hline
 & \\
TOTAL 			& 56 = $2^3 \times 7$
\end{tabular} \\
\end{tabular}\\
\par
Les diagrammes de Cauchon pour les 6 première colonnes proviennent du cas $E_{6}$, soit $2^7 \times 3^4 \times 5$ diagrammes. En ajoutant la dernière colonne, on a $2^{10} \times 3^4 \times 5 \times 7$ diagrammes soit le cardinal du groupe de Weyl associé. On a encore $| \mathcal{D}| = |W|$.

\subsubsection{$E_8$}
\begin{conv}
On numérote les racines simple de sorte que le diagramme de Dynkin soit : 
\[\begin{array}{ccccccccccccc}
 					& 					&   				 & 					 & \alpha_2	&&  &     &&& &&  \\ 
 					& 					&   				 & 					 & |					&&  &     &&& &&\\
 \alpha_1 &\text{---~}& \alpha_{3} &\text{---~}& \alpha_4 &\text{---~}& \alpha_{5}&\text{---~} & \alpha_6 &\text{---~}& \alpha_7& \text{---~}& \alpha_8\\
\end{array}\]
Comme les racines $\alpha_1$ à $\alpha_7$ forment un système de racines de type $E_7$, on note $\sigma_{7}$, l'élément de plus grande longueur utilisé précédemment pour le type $E_7$. On vérifie que la décomposition réduite suivante 
\par
$w_0 = \sigma_{7} s_8 s_{7} s_6 s_5 s_4 s_2 s_3  s_1  s_4  s_3 s_5  s_4  s_6 s_2 s_5  s_7 s_4 s_6
 s_8 s_3 s_5 s_7 s_1 s_4 s_6 s_3 s_2 s_5 s_4 s_{5} s_2 s_3 s_6 s_1 s_4  s_7  s_3  s_5 s_8s_4  s_6 s_2 s_5  s_7s_4 $\\$s_6 s_3 s_5 s_1 s_4 s_2 s_3 s_4 s_5 s_6 s_7 s_8$.\\
induit l'ordre convexe suivant sur les racines positives (on ne fait apparaître ici que la dernière colonne). Dans la deuxième colonne on calcul $h'(\beta_i)$ : \\
\[
\begin{array}{ccc}
\begin{array}{c|c}
\beta_i & h'(\beta_i) \\
\hline
\beta_{64}( 2, 3, 4, 6, 5, 4, 3, 1 ) &	28 \\
\hline
\beta_{65}( 2, 3, 4, 6, 5, 4, 2, 1 ) &	27 \\
\hline
\beta_{66}( 2, 3, 4, 6, 5, 3, 2, 1 ) &	26 \\
\hline
\beta_{67}( 2, 3, 4, 6, 4, 3, 2, 1 ) &	25 \\
\hline
\beta_{68}( 2, 3, 4, 5, 4, 3, 2, 1 ) &	24 \\
\hline
\beta_{69}( 2, 2, 4, 5, 4, 3, 2, 1 ) &	23 \\
\beta_{70}( 2, 3, 3, 5, 4, 3, 2, 1 ) &	23 \\
\hline
\beta_{71}( 1, 3, 3, 5, 4, 3, 2, 1 ) & 	22 \\
\beta_{72}( 2, 2, 3, 5, 4, 3, 2, 1 )	& 	22 \\
\hline
\beta_{73}(1, 2, 3, 5, 4, 3, 2, 1 ) &	21 \\
\beta_{74}( 2, 2, 3, 4, 4, 3, 2, 1 ) 	& 	21 \\
\hline
\beta_{75}( 1, 2, 3, 4, 4, 3, 2, 1)&	20 \\
\beta_{76}(2, 2, 3, 4, 3, 3, 2, 1 )& 	20 \\
\hline
\beta_{77}( 1, 2, 2, 4, 4, 3, 2, 1 )&	19 \\
\beta_{78}( 1, 2, 3, 4, 3, 3, 2, 1 )& 	19 \\
\beta_{79}( 2, 2, 3, 4, 3, 2, 2, 1 )&	19 \\
\hline
\beta_{80}(1, 2, 2, 4, 3, 3, 2, 1 )&	18  \\
\beta_{81}( 1, 2, 3, 4, 3, 2, 2, 1 )&	18 \\
\beta_{82} (2, 2, 3, 4, 3, 2, 1, 1 )&18\\
\hline
\beta_{83}( 1, 2, 2, 3, 3, 3, 2, 1 )& 	17 \\
\beta_{84}( 1, 2, 2, 4, 3, 2, 2, 1 )&	17 \\
\beta_{85}( 1, 2, 3, 4, 3, 2, 1, 1 )& 	17 \\
\hline
\beta_{86}( 1, 1, 2, 3, 3, 3, 2, 1 )& 	16 \\
\beta_{87}( 1, 2, 2, 3, 3, 2, 2, 1 )&	16 \\
\beta_{88}( 1, 2, 2, 4, 3, 2, 1, 1 )&	16 \\
\hline
\beta_{89}( 1, 1, 2, 3, 3, 2, 2, 1 )&	15 \\
\beta_{90}( 1, 2, 2, 3, 2, 2, 2, 1 )& 	15 \\
\beta_{91}( 1, 2, 2, 3, 3, 2, 1, 1 )&	15 \\
\hline
\beta_{92}( 2, 3, 4, 6, 5, 4, 3, 2 )&29/2 \\
\hline
\end{array}
&&
\begin{array}{c|c}
\beta_i & h'(\beta_i) \\
\hline
\beta_{93}( 1, 1, 2, 3, 2, 2, 2, 1 )&	14 \\
\beta_{94}( 1, 1, 2, 3, 3, 2, 1, 1)&	14\\
\beta_{95}( 1, 2, 2, 3, 2, 2, 1, 1 )&	14\\
\hline 
\beta_{96}( 1, 1, 2, 2, 2, 2, 2, 1 )&	13\\
\beta_{97}( 1, 2, 2, 3, 2, 1, 1, 1 )&	13\\
\beta_{98}( 1, 1, 2, 3, 2, 2, 1, 1 )&	13 \\
\hline 
\beta_{99}( 1, 1, 1, 2, 2, 2, 2, 1 )& 	12\\
\beta_{100}	(1, 1, 2, 3, 2, 1, 1, 1)&	12\\
\beta_{101}	( 1, 1, 2, 2, 2, 2, 1, 1 )&	12\\
\hline
\beta_{102}( 0, 1, 1, 2, 2, 2, 2, 1 )&	11\\
\beta_{103}( 1, 1, 2, 2, 2, 1, 1, 1 )&	11\\
\beta_{104}	( 1, 1, 1, 2, 2, 2, 1, 1 )&	11\\
\hline 
\beta_{105}( 1, 1, 2, 2, 1, 1, 1, 1 )&	10 \\
\beta_{106}	(1, 1, 1, 2, 2, 1, 1, 1 )&10 \\
\beta_{107}	( 0, 1, 1, 2, 2, 2, 1, 1 )&	10 \\
\hline
\beta_{108}( 1, 1, 1, 2, 1, 1, 1, 1 )&	9 \\
\beta_{109}	(0, 1, 1, 2, 2, 1, 1, 1 )&	9 \\
\hline
\beta_{110}	(1, 1, 1, 1, 1, 1, 1, 1 )& 	8 \\
\beta_{111}	( 0, 1, 1, 2, 1, 1, 1, 1)&	8 \\
\hline
\beta_{112}	( 1, 0, 1, 1, 1, 1, 1, 1 )&	7  \\
\beta_{113}	( 0, 1, 1, 1, 1, 1, 1, 1 )&	7 \\
\hline
\beta_{114}	( 0, 1, 0, 1, 1, 1, 1, 1 )&	6 \\ 
\beta_{115}	( 0, 0, 1, 1, 1, 1, 1, 1 )&	6 \\
\hline
\beta_{116}	( 0, 0, 0, 1, 1, 1, 1, 1 )&	5 \\
\hline 
\beta_{117}	( 0, 0, 0, 0, 1, 1, 1, 1)&	4 \\
\hline
\beta_{118}	( 0, 0, 0, 0, 0, 1, 1, 1 )& 	3\\
\hline
\beta_{119}	( 0, 0, 0, 0, 0, 0, 1, 1 )&	2 \\
\hline
\beta_{120}	( 0, 0, 0, 0, 0, 0, 0, 1 )&	1 \\
\end{array}
\end{array}\]
\end{conv}
On connaît déjà la forme des diagrammes sur les 7 premières colonnes. 
Au moyen des propositions \ref{contrainteelementaire}, \ref{contrainteelementaireavantex} et \ref{contrainteelementaireapresex}, on obtient les contraintes pour la dernière colonne.\\
En particulier, on va avoir des contraintes du type $\mathbf{(i \Rightarrow j \ \text{ou} \ k)}$ :
\[(92 \Rightarrow 91 \ \text{ou} \ 90) \ \text{et} \ (92 \Rightarrow 90 \ \text{ou} \ 89) \ \text{et} \ (92 \Rightarrow 91 \ \text{ou} \ 89).\]
On note ces contraintes par des pointillés dans le graphe ci-dessous : \\
\begin{center}
\begin{psmatrix}[rowsep=.25cm, colsep=.2cm]
 &  && 120 \\ 
 &  & &119 \\ 
 &  & &118 \\ 
 &  & &117 \\ 
 &  & &116 \\ 
 & &115 &  & 114 \\ 
& &112 & &113 \\ 
 & &110  & &111 \\ 
 &  & 108 &  & 109 \\
 & 105 &  & 106 &  & 107 \\ 
 & 103 &  & 104 &  & 102 \\ 
 & 100 &  & 101 &  & 99 \\ 
 & 97 &  & 98 &  & 96 \\ 
 & 95 &  & 94 &  & 93 \\ 
 &  &  & 92 \\ 
 & 91 &  & 90 &  & 89 \\ 
 & 88 &  & 87 &  & 86 \\ 
 & 85 &  & 84 &  & 83 \\ 
 & 82 &  & 81 &  & 80 \\
 & 79 &  & 78 &  & 77 \\ 
 && 76 &  & 75 \\ 
 && 74 &  & 73 \\ 
 && 72 &  & 71 \\ 
 && 69 &  & 70 \\ 
 &&  & 68 \\ 
 &&  & 67 \\ 
 &&  & 66 \\ 
 &&  & 65 \\ 
 &&  & 64 
\end{psmatrix}
\psset{arrows=->, nodesep=2pt}
\ncline{1,4}{2,4}
\ncline{2,4}{3,4}
\ncline{3,4}{4,4}
\ncline{4,4}{5,4}
\ncline{5,4}{6,3}
\ncline{5,4}{6,5}
\ncline{6,3}{7,5}
\ncline{6,3}{7,3}
\ncline{6,5}{7,5}
\ncline{7,5}{8,5}
\ncline{7,5}{8,3}
\ncline{7,3}{8,3}
\ncline{8,5}{9,5}
\ncline{8,5}{9,3}
\ncline{8,3}{9,3}
\ncline{9,5}{10,6}
\ncline{9,5}{10,4}
\ncline{9,3}{10,4}
\ncline{9,3}{10,2}
\ncline{10,6}{11,4}
\ncline{10,6}{11,6}
\ncline{10,4}{11,4}
\ncline{10,4}{11,2}
\ncline{10,2}{11,2}
\ncline{11,4}{12,4}
\ncline{11,4}{12,6}
\ncline{11,2}{12,4}
\ncline{11,2}{12,2}
\ncline{11,6}{12,6}
\ncline{12,4}{13,4}
\ncline{12,4}{13,6}
\ncline{12,2}{13,4}
\ncline{12,2}{13,2}
\ncline{12,6}{13,6}
\ncline{13,4}{14,2}
\ncline{13,4}{14,4}
\ncline{13,4}{14,6}
\ncline{13,2}{14,2}
\ncline{13,6}{14,6}
\ncline{14,2}{15,4}
\ncline{14,2}{16,2}
\ncline{14,2}{16,4}
\ncline{14,4}{15,4}
\ncline{14,4}{16,2}
\ncline{14,4}{16,6}
\ncline{14,6}{15,4}
\ncline{14,6}{16,4}
\ncline{14,6}{16,6}
\ncline{16,2}{17,2}
\ncline{16,2}{17,4}
\ncline{16,4}{17,4}
\ncline{16,6}{17,4}
\ncline{16,6}{17,6}
\ncline{17,2}{18,2}
\ncline{17,2}{18,4}
\ncline{17,4}{18,4}
\ncline{17,4}{18,6}
\ncline{17,6}{18,6}
\ncline{18,2}{19,2}
\ncline{18,2}{19,4}
\ncline{18,4}{19,4}
\ncline{18,4}{19,6}
\ncline{18,6}{19,6}
\ncline{19,2}{20,2}
\ncline{19,4}{20,2}
\ncline{19,4}{20,4}
\ncline{19,6}{20,4}
\ncline{19,6}{20,6}
\ncline{20,2}{21,3}
\ncline{20,4}{21,3}
\ncline{20,4}{21,5}
\ncline{20,6}{21,5}
\ncline{21,3}{22,3}
\ncline{21,5}{22,3}
\ncline{21,5}{22,5}
\ncline{22,3}{23,3}
\ncline{22,5}{23,3}
\ncline{22,5}{23,5}
\ncline{23,3}{24,5}
\ncline{23,3}{24,3}
\ncline{23,5}{24,5}
\ncline{24,5}{25,4}
\ncline{24,3}{25,4}
\ncline{25,4}{26,4}
\ncline{26,4}{27,4}
\ncline{27,4}{28,4}
\ncline{28,4}{29,4}
\psset{nodesep=1pt, linestyle=dashed}
\ncline{15,4}{16,2}
\ncline{15,4}{16,4}
\ncline{15,4}{16,6}
\end{center}
On utilise alors les contraintes ci-dessus pour trouver le nombre de diagrammes possibles.\\
\begin{tabular}{cc}
\begin{tabular}{c|c}
Diagrammes commençant  & Nombre de\\
par la case :					&  diagrammes \\
\hline
Pas de case & 1\\
64					& 1	\\
65					& 1\\
66					& 1\\
67					& 1\\
68					& 1\\
69					& 1\\
70					& 2\\
71					& 2\\
72					& 2\\
73					& 1\\
74					& 3\\
75					& 1\\
76					& 4\\
77					& 2\\
78					& 2\\
79					& 7\\
80					& 2\\
81					&	3\\
82					& 11\\
83					& 5\\
84					&	4\\
85					& 6\\
86					& 9\\
87					& 6\\
88					& 5\\
89					& 4\\
90					&	12\\
91					& 6\\
92					& 8
\end{tabular}
&
\begin{tabular}{c|c}
Diagrammes commençant  & Nombre de\\
par la case :					&  diagrammes \\
\hline
93					& 5\\
94					& 3\\
95					& 6\\
96					& 8\\
97					& 8\\
98					&4 \\
99					& 12\\
100					& 3\\
101					& 6\\
102					& 16\\
103					& 3\\
104					& 8\\
105					& 5\\
106					& 4\\
107					& 7\\
108					& 3\\
109					& 3\\
110					& 4\\
111					& 2\\
112					& 5\\
113					& 2\\
114					& 2\\
115					& 2\\
116					& 1 \\
117					& 1\\
118					& 1\\
119					& 1\\
120					& 1\\
\hline
 & \\
TOTAL 			& 240 = $2^4 \times 3 \times 5$
\end{tabular} 
\end{tabular} \\
\par
$ $ \\
$ $ \\
Pour les 7 première colonnes, on retrouve un système de type $E_7$, soit $2^{10} \times 3^4 \times 5$ diagrammes.\\
En ajoutant la dernière colonne, on dénombre $2^{10} \times 3^4 \times 5 \times 7$ diagrammes, soit le cardinal du groupe de Weyl associé. On a, comme dans les autres cas, $|\mathcal{D}|=|W|$.

\nocite{*}
\pagebreak
\tableofcontents
\pagebreak
\bibliographystyle{alpha}
\bibliography{biblio.bib}